%% file: main.tex
\documentclass[11pt, oneside]{Thesis} 

\graphicspath{{Figures/}{../Figures/}}
\usepackage{enumitem}
\usepackage{amsmath}
\usepackage{mathtools}
\usepackage{mathrsfs}
\usepackage{amsthm}
\usepackage[usenames, dvipsnames]{color}
\usepackage[T1]{fontenc}
\usepackage[utf8]{inputenc}
\usepackage[square, numbers, comma, sort&compress]{natbib} 
\usepackage[nottoc,notlof,notlot]{tocbibind} 
\renewcommand\bibname{References}
\usepackage{epigraph} 
\usepackage{makeidx} 
\makeindex 
\usepackage{subfiles}
\usepackage{graphicx}
\usepackage{hyperref}
\hypersetup{urlcolor=blue, colorlinks=true, plainpages=true, hypertexnames=true,} 
\title{\ttitle} 

\newcommand{\R}{\mathbb{R}}
\newcommand{\B}{\mathbb{B}}
\newcommand{\N}{\mathbb{N}}
\newcommand{\norm}[1]{\left \| #1  \right \|}
\newcommand{\paren}[1]{\left (  #1 \right ) }
\newcommand{\smap}[3]{#1 : \mathbb{R}^{#2} \longrightarrow \mathbb{R}^{#3} }
\newcommand{\mmap}[3]{#1 : \mathbb{R}^{#2} \rightrightarrows \mathbb{R}^{#3} }
\newcommand{\gph}[1]{ \mathrm{gph} \, #1 }
\newcommand{\rfp}[2]{ (\bar{#1}, \bar{#2}) }
\newcommand{\rfpp}[2]{ (\bar{#1}\, | \, \bar{#2}) }
\newcommand{\NUV}{neighborhoods $U$ of $\bar{x} $ and $V $ of $\bar{y}$ }
\newcommand{\mfa}{ \mathrm{~~~~for~all~} }
\newcommand{\at}[2]{at $ \bar{#1} $ for $ \bar{#2} $}
\newcommand{\inp}[2]{\langle\,#1 \, , #2 \, \rangle}

\theoremstyle{plain}
\newtheorem{thm}{\bf Theorem}[section]
\newtheorem{lem}[thm]{\bf  Lemma}
\newtheorem{prop}[thm]{\bf Proposition}
\newtheorem{cor}[thm]{\bf Corollary}

\theoremstyle{definition}
\newtheorem{note}[thm]{\bf  Note}
\newtheorem{defn}[thm]{\bf Definition}
\newtheorem{eg}[thm]{\bf Example}

\newtheorem{rem}[thm]{\bf Remark}

\renewenvironment{proof}{\emph{\textbf{Proof.}} \normalfont}{\hfill $\Box$}

\allowdisplaybreaks

\begin{document}
\frontmatter 

\setstretch{1.3} 

\fancyhead{} 
\rhead{\thepage} 
\lhead{} 

\pagestyle{fancy} 

\newcommand{\HRule}{\rule{\linewidth}{0.5mm}} 

\hypersetup{pdftitle={\ttitle}}
\hypersetup{pdfsubject=\subjectname}
\hypersetup{pdfauthor=\authornames}
\hypersetup{pdfkeywords=\keywordnames}


\begin{titlepage}
\begin{center}
\textsc{\LARGE \univname}\\[1cm] 
\begin{figure}[ht]
	\centering
		\includegraphics[width=0.28\textwidth]{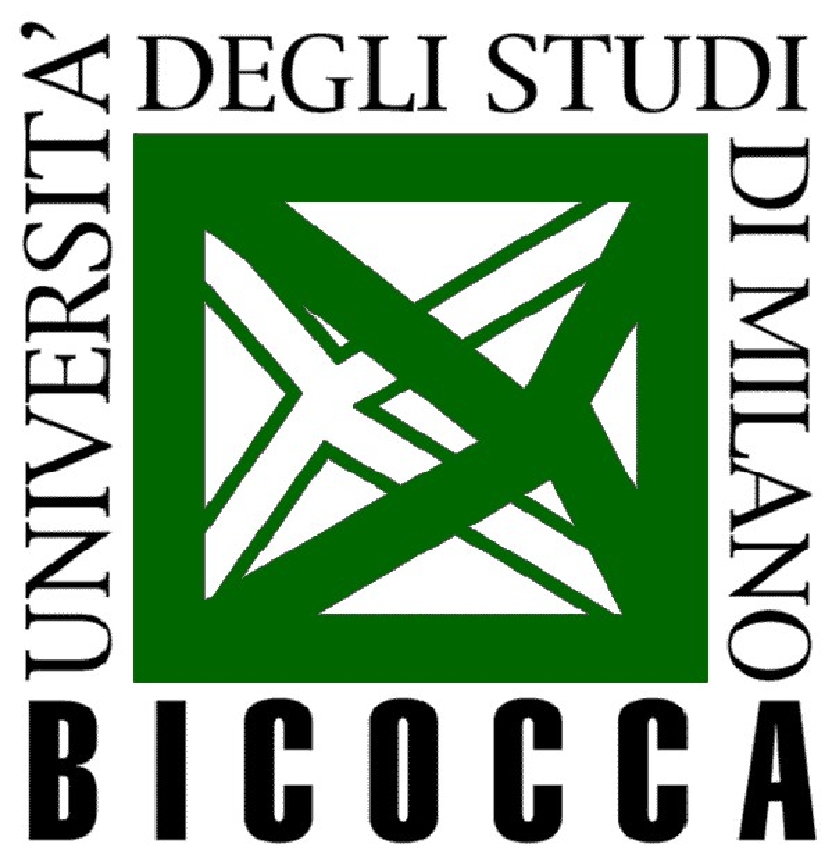}
\end{figure}
\vspace*{1cm}
\textsc{\Large Doctoral Thesis}\\[0.5cm] 

\HRule \\[0.4cm] 
{\Large \bfseries \ttitle}\\[0.4cm] 
\HRule \\[1.5cm] 
 
\begin{minipage}{0.4\textwidth}
\begin{flushleft} \large
\emph{Author:}\\
\authornames
\end{flushleft}
\end{minipage}
\begin{minipage}{0.4\textwidth}
\begin{flushright} \large
\emph{Supervisor:} \\
\supname
\end{flushright}
\end{minipage}\\[2.5cm]
 
\large \textit{A thesis submitted in fulfilment of the requirements\\ for the degree of \degreename}\\[0.3cm] 
\textit{in the}\\[0.4cm]
\deptname\\ [0.3cm]
\MakeUppercase{\large \univname} \\[2cm]

{\large March 2017}\\[1cm] 

\vfill
\end{center}

\end{titlepage}


\Declaration{

\addtocontents{toc}{\vspace{1em}} 

I, \authornames, declare that this thesis titled, \lq\lq Metrically Regular Generalized Equations: A Case Study in Electronic Circuits\rq\rq \,and the work presented in it are my own. I confirm that:

\begin{itemize} 
\item[\tiny{$\blacksquare$}] This work was done wholly or mainly while in candidature for a research degree at this University.
\item[\tiny{$\blacksquare$}] Where any part of this thesis has previously been submitted for a degree or any other qualification at this University or any other institution, this has been clearly stated.
\item[\tiny{$\blacksquare$}] Where I have consulted the published work of others, this is always clearly attributed.
\item[\tiny{$\blacksquare$}] Where I have quoted from the work of others, the source is always given. With the exception of such quotations, this thesis is entirely my own work.
\item[\tiny{$\blacksquare$}] I have acknowledged all main sources of help.
\item[\tiny{$\blacksquare$}] Where the thesis is based on work done by myself jointly with others, I have made clear exactly what was done by others and what I have contributed myself.
\end{itemize}
 
Signed:\\
\rule[2.5em]{25em}{0.5pt} 
~Date: 10/March/2017 \\
\rule[2.5em]{25em}{0.5pt} }
\clearpage 






\addtotoc{Abstract} 

\abstract{\addtocontents{toc}{\vspace{0.5em}} 

\vspace*{2cm}
In this thesis, set-valued maps are considered to model the $ i- v $ characteristics of semiconductors like diode, and transistor. Using the circuit theory laws, a generalized equation is obtained. The main concern of the thesis is to investigate how perturbing the input signal will affect the output variables. The problem is studied in two cases: the static case, where the input signal is a DC source; and the dynamic case, where there exists an AC source in the circuit. \\
In the static case, 
the problem can be reduced to the existence or absence of local stability properties of the solution map, like the Aubin property, isolated calmness, and calmness, or metric regularity for the inverse map. Some tools from variational analysis are used to provide necessary and/or sufficient conditions that guarantee these properties. \\
In the dynamic case, those pointwise results are used to obtain descriptions for regularity properties of the solution trajectories in function spaces.
}
\clearpage 


\setstretch{1.3} 

\vspace*{3.5cm}
\acknowledgements{\addtocontents{toc}{\vspace{0.5em}} 
\lhead{\emph{Acknowledgements}}
\vspace*{2.5cm}

Imagine each human as a sphere in $\R^n$ (with $ n \geq 4$), whose mass is concentrated in a human body shape at the center of this sphere and his/her effect radiates centrifugally. Thus, whatever a human being does, would affect all the other humans, more or less. Of course, such an effect may decrease when you get farther from the center, but it never reaches zero. \\
You may define an \lq\lq effective radius\rq\rq, inside which, people are physically close enough to see each other, to start a conversation, to communicate and share experiences and emotions, ..., to live with each other. Such an effective radius, if it exists, is different from one person to another. Also there could be some disagreements in defining this radius. You may see someone, everyday for a while, but never have the courage to say something to him/her. On the other hand, it is possible to read a note and feel a great and in depth sympathy and closeness to the person who has written it. \\
Such a model, if one accepts, can explain many things in human relations and particularly, the difficulty I may encounter to write an acknowledgement. There has been many people in my life, who had helped me to be where I am. Consciously or unconsciously, their manners and deeds have affected me, even though I may not know them in person. To be honest, I should thank all of humans. \\
After this general introduction, I may name a few people who had stronger influence on me, at least in the last three years, which would not mean being rude or ungrateful to those I may not name here. \\
Above all, I want to express my sincere appreciation to my advisers, Prof.ssa Rita Pini, and Prof. Amos Uderzo. English was not the mother tongue of me nor theirs, but their kindness allows me to feel powerful enough to talk, their patience provides me a safety to ask my questions, their expertise permits me to run faster in my path, and their openness to new experiences enables me to bridge a gap between my bachelor and master studies. I learned a lot from their mathematical knowledge and manner. \\
I would like to thank Prof. Radek Cibulka for his kindness and trust to give me drafts of his unpublished works; his patience to answer my endless questions and doubts; his hospitality to let me visit him in Pilsen; and his ideas to enrich my work. The spark of most of the new results in this thesis has been produced during my stay in Pilsen. Tomas Roubal also helped me to have a more pleasant and fruitful visiting period, thanks to his valuable friendship, which I hope lasts for a life-time.\\
During my stay in Italy, I had very good moments with my friends that I may not forget. Thanks to Simone, Davide, Chiara, Elena, Jessica, Daniela, Federico, Morteza, Alberto, Mani, Masoud, Maryam, Reza, ... and all the others.\\
Last but not least, I shall thank my family: my parents and my sister for their continuous love and support; my brothers who empowered and encouraged me for higher education. 
} 
\clearpage 


\pagestyle{fancy} 

\lhead{\emph{Contents}} 
\tableofcontents 

\lhead{\emph{List of Figures}} 
\listoffigures 

\mainmatter 

\pagestyle{fancy} 

\subfile{Chapters/Introduction}

\subfile{Chapters/Chapter01}

\subfile{Chapters/Chapter02}

\subfile{Chapters/Chapter03-a}

\subfile{Chapters/Chapter03-b}

\subfile{Chapters/Chapter04}

\newpage
\label{Bibliography}

\lhead{\emph{References}} 

\bibliographystyle{amsplain-nodash}
\nocite{*}
\bibliography{mainbibadapted} 


\clearpage 

\lhead{\emph{Notation}} 

\listofnomenclature{ll} 
{
$\R$ & the real numbers \hspace*{0.2cm} \\
$\R_{+}$ & the non-negative real numbers   \\
$\R^m$ & the $m$-dimensional Euclidean space \\ 
$\N$ & the natural numbers: $1, 2, 3, . . .$ \\[0.6em]
$ \B_{r} (\bar{x}) $ & closed ball of radius $ r > 0 $ centred at $\bar{x}$ \\
$ \B $ & closed unit ball, $ \B_{1} (0) $\\
$ \mathrm{int\,} C $ & interior of the set $C$ \\
$ \mathrm{cl\,} C $ & closure of the set $C$ \\
$ x_k \longrightarrow x $ & the sequence $ (x_k) $ is convergent to $x$ \\
$ t_k \downarrow 0 $ & a sequence of positive numbers $ t_k $ tending to $ 0 $\\[0.6em]
$ d(x, C) $ & distance from $x$ to the set $C$ \\
$e(C, D) $ & excess of the set $C$ beyond the set $D$ \\
$ | x | $ & absolute value of $ x \in \R $ \\
$\norm{x}$ & norm of $x$ \\
$ \inp{x}{y} $ & canonical inner product, bilinear form \\
$ | H |^{+} $ & outer norm \\
$| H |^{-}$ & inner norm \\[0.6em]
$ T(\bar{x};\,C) $ & Bouligand-Severi tangent cone (contingent cone) to the set $C$ at $\bar{x}$ \\
$\widetilde{T}(\bar{x};\, C)$ & Bouligand paratingent cone to the set $C$ at $\bar{x}$ \\
$\widehat{N}(\bar{x}; \,C) $ &  Fr\'echet normal cone (regular normal cone) to the set $C$ at $\bar{x}$ \\
$N(\bar{x}; \,C)$ & Mordukhovich normal cone (limiting normal cone)  to the set $C$ at $\bar{x}$ \\[0.6em]
$ F : X \rightrightarrows Y $ & set-valued mapping from $X$ into the subsets of $Y$ \\
$ f : X \to Y $ & function $f$ from $X$ into $Y$ \\
$ A^T $ & transpose of the matrix $A$ \\
$ \mathrm{ker\,} A$ & kernel of the the linear operator $A$ \\
$ \mathrm{det\, } A$ & determinant of the matrix $A$\\
$ \gph{F} $ & graph of the mapping $F$ \\
$\mathrm{dom\,} F $ & domain of the mapping $F$ \\
$\mathrm{rge\,} F $ & range of the mapping $F$ \\[0.6em]
$\partial_F f(\bar{x}) $ & Fr\'{e}chet subdifferential (F-subdifferential) of the function $f$ at $\bar{x}$ \\
$ \partial f(\bar{x}) $ & Mordukhovich (limiting) subdifferential of the function $f$ at $\bar{x}$ \\
$ \partial_{>} f(\bar{x}) $ & outer subdifferential of the function $f$ at $ \bar{x} $ \\
$ \partial_B h (\bar{u}) $ & Bouligand's limiting Jacobian of the function $ h $ at $\bar{u}$ \\
$ \partial h (\bar{u}) $ & Clarke's generalized Jacobian of $h$ at $\bar{u}$ \\[0.6em]
$ \nabla f (\bar{x}) $ & Jacobian matrix of the function $f$ at $\bar{x}$ \\
$D f(\bar{x}) $ & derivative of the function $f$ at $\bar{x}$ \\
$ DF\rfpp{x}{y} $ & graphical derivative of the mapping $F$ \at{x}{y} \\ 
$ D^* F\rfpp{x}{y} $ & coderivative of the mapping $F$ \at{x}{y} \\ 
$ \widetilde{D} F \rfpp{x}{y} $ & strict graphical derivative of the mapping $F$ \at{x}{y} \\[0.6em]
$ \mathrm{clm\,}( f; \bar{x} ) $ & calmness modulus of the function $f$ at $\bar{x}$ \\
$ \mathrm{clm}_x ( f ; \rfp{p}{x} ) $ & partial calmness modulus of $f$ with respect to $x$ at $\rfp{p}{x}$ \\
$ \widehat{\mathrm{clm}}_x \, (f; (\bar{p}, \bar{x})) $ & uniform partial calmness modulus of $f$ with respect to $x$ \\
$ \mathrm{clm\,}(S; \bar{y}|\bar{x})$ & calmness modulus of $S$ at $\bar{y}$ for $\bar{x}$  \\
$ \mathrm{lip\,}(f; \bar{x}) $ & Lipschitz modulus of the function $f$ at $\bar{x}$ \\
$ \widehat{\mathrm{lip}}_x \, (f; (\bar{p}, \bar{x})) $ & uniform partial Lipschitz modulus of $f$ with respect to $x$ \\
$ \mathrm{lip\,}(S; \bar{y}| \bar{x}) $ & Lipschitz modulus of $ S $ at $ \bar{y}$ for $\bar{x}$ \\
$ \mathrm{reg\,}(F; \bar{x}|\bar{y})$ & regularity modulus of $F$ at $\bar{x}$ for $\bar{y}$ \\
$ \mathrm{subreg\,} (F; \bar{x} | \bar{y}) $ & modulus of metric sub-regularity of $F$ at $\bar{x}$ for $\bar{y}$
}

\newpage
\clearpage
\lhead{\emph{Index}}
\phantomsection 
\printindex

\end{document}

%% file: Chapters/Introduction.tex
\chapter*{Introduction} 
\addcontentsline{toc}{chapter}{Introduction} \markboth{INTRODUCTION}{}
\label{Introduction} 

\lhead{ \emph{Introduction}} 
\setlength{\epigraphwidth}{.49\textwidth}
\epigraph{\emph{The definition of a good mathematical problem is the mathematics it generates rather than the problem itself.}}{\textit{Andrew J. Wiles}}

It would be useful to compare learning mathematics with practising arts like painting and drawing. Imagine you want to draw a landscape. As the first steps in mastering as a painter are to watch carefully, and try to make a simplified sketch of what the painter is seeing, then going slowly toward the details to create more accurate copy of the landscape (in other words, closer to what is \lq\lq real\rq\rq ) by improving his/her knowledge of the painting instruments, mastering hand skills, and an endless process of trial and error to become better and better, so should be the approach to use mathematics in studying the nature.\\
Modelling a natural phenomenon is a way to study it in a more abstract way than doing experiments on it (if possible), and a goal of modelling is to \lq\lq predict\rq\rq \,the behaviour of the phenomenon under study with respect to changes of different variables in the model. 
The more appropriate a model would be,  a \emph{better} prediction of the phenomenon it would provide, and by \lq\lq better\rq\rq \,one can think of different desired properties: a wider range of change for the input variables, an easier way to find the equilibrium point, a more detailed view of small changes, or ...\\
Our study starts with modelling electrical circuits. As it would be explained with more details in Chapter \ref{Chapter2}, every component in the circuit could be understood as a relation between the current passing through it and the drop of potential difference (voltage, for short) over it. For many components such a relation is a function, but some electronic components like diodes, and transistors, need a set-valued map to describe their $ i - v $ characteristic. Understanding the behaviour of a circuit depends on knowing the current of each branch and voltage of each component in the circuit, and to do so, we need to use the $ i - v $ characteristics of components, together with Kirchhoff’s current and voltage laws to form a system of $n$ equations, $n$ variables, and then, solve it. When there exists a set-valued map in the model, we will obtain a \emph{generalized equation} like:
\begin{equation} \label{GE-int}
0 \in f(x) + F(x),
\end{equation}
where $\smap{f}{n}{n}$, and $\mmap{F}{n}{n}$ could be determined specifically for each circuit.\\
Generalized equations are interesting structures in their own. First, observe that $ F( \cdot ) \equiv \{ 0 \} $, reduces \eqref{GE-int} to an ordinary equation $ 0 = f(x) $, hence the name. Second, consider a closed and convex set $ \Omega $ in $\R^n$, and let $F$ be the normal cone (in the sense of convex analysis) to $\Omega$ at a point $ x \in \R^n$, that is,
\begin{equation}
F(x) = N(x; \Omega) =
\left \{ \begin{matrix}
\big\{ v ~|~ \inp{v}{c - x} \leq 0 \mfa c \in \Omega \big\} & & x \in \Omega,\\
\emptyset ~~~~~~~~~~~~~~~~~~~~~~~~~~~~~~~~~~~~~~~~~~~  & & x \not \in \Omega.
\end{matrix}\right.
\end{equation}
Geometrically, this is the cone of all outward normals to $\Omega$ at $x$. Note that, if \eqref{GE-int} holds with this particular $F$, then the sum on the right hand side is non-empty (it contains $0$), so $ N(x; \Omega) \neq \emptyset $, which means $ x \in \Omega $. Also, $ - f(x) $ must belong to $N(x; \Omega)$, so for each $ c \in \Omega$,
\begin{equation*} 
\inp{ -f(x) }{ c - x } \, \leq \, 0.
\end{equation*}
Thus, we can see that \eqref{GE-int} holds if and only if $x$ satisfies the so called \emph{variational inequality}:
\begin{equation} \label{VI-int}
x \in \Omega, \mathrm{~and~}~ \inp{ f(x) }{ c - x } \, \geq \, 0 \mathrm{~~for~each~} c \in \Omega, 
\end{equation}
and this, geometrically, means that $f(x)$ is an inward normal to $\Omega$ at $x$. Robinson in \citep{robinson1979generalized, robinson1982generalized, robinson1983generalized, robinson} has studied this particular type of generalized equations in details and found the setting of generalized equations as an appropriate way to express and analyse problems in complementarity, mathematical programming, and variational inequalities.\\
It is worthwhile mentioning that, although one could write \eqref{GE-int} as $ 0 \in G(x) $, by defining a new set-valued map $ \mmap{G}{n}{n} $ as $ G(x) = f(x) + F(x) $, 
keeping the single-valued and set-valued parts separated is more useful in practice, as it turns out that $f$ is often a fairly smooth function, while 
$F$ involves \lq\lq corners\rq\rq.\\
In the study of electrical circuits, power supplies (that is, both current and voltage sources) play an important role. Not only their failure in providing the minimum voltage level for other components to work would be a problem, but also small changes in the provided voltage level will affect the whole circuit and the goal it has been designed for. This small changes around a desired value could happen mainly because of failure in precise measurements, ageing process, and thermal effects (which are explained in Chapter \ref{Chapter2}, see also \citep{johns2008analog, Sedra}). \\
Thus, based on the type of voltage/current sources in the circuits, we consider two different cases:
\begin{enumerate}[topsep=-1ex, itemsep=0ex, partopsep=1ex, parsep=1ex, leftmargin = 4ex]
\item[1.] \textbf{static case:}\index{problem ! static case}\\ 
This is the situation when the signal sources in the circuit are DC (that is, its value is not changing with respect to time). 
For practical reasons, we would prefer to rewrite \eqref{GE-int} as $ p \in f(x) + F(x) $, where $ p \in \R^n $ is a fixed vector representing the voltage or current sources in the circuit. Then, we define the corresponding \emph{solution mapping} as follows:
\begin{equation} \label{solution mapping-int}
p \longmapsto S(p) := \left \{ x \in \R^n ~|~ p \in f(x) \, + F(x) \right \}. 
\end{equation}
In this framework, the desire to have small deviations of $ x $ with respect to perturbations of $p$ around a presumed point $ \rfp{p}{x} \in \gph{S} $ could be investigated as the local stability properties (like the Aubin property, isolated calmness, and calmness) of $S$ at $\bar{p}$ for $\bar{x}$. Or equivalently, we can ask for metric regularity features of $ \Phi = f + F $ at $\bar{x}$ for $\bar{p}$. The details of this analogy is given in Chapter \ref{Chapter2}, while the definition of these local properties and their equivalences (one for $S$, the other for $\Phi = S^{-1}$) is provided in Chapter \ref{Chapter1}. We would provide handy theorems to check these properties in Chapter \ref{Chapter3}.\\
Let us mention that by the term \emph{metrically regular generalized equation}, we refer to a generalized equation \ref{GE-int} where the right hand side is a metrically regular map. 
In general, the stability properties under parameter perturbations is a very important topic in engineering, not only for the determination of the behaviour of a system with respect to perturbations but also for the construction of algorithms for the numerical simulation of the problem. 
\item[2.] \textbf{dynamic case:} \index{problem ! dynamic case}\\
When an AC signal source (that is, its value is a function of time) is in the circuit, problem could be more complicated. First of all, all the other variables of the model would become a function of time, too. Second, it is not appropriate any more to formulate the solution mapping as $S(p)$. One can consider a parametric generalized equation like:
\begin{equation} \label{PGE-int}
0 \in h(t, x) + F(x),
\end{equation}
where $ h : \R \times \R^n \longrightarrow \R^n $ now depends on a scalar parameter $ t \in [0, 1]$\footnote{
In fact, $t$ can belong to any finite interval like $ [0, T]$ for a $T>0$. The starting point $ t = 0$ is considered as the moment that the circuit starts working, in other words, when the circuit is connected to the signal sources and is turned on with a key. We keep the time interval as $ [0, 1] $ in the entire thesis for simplicity.\\ 
}, 
and define the solution mapping corresponding to \eqref{PGE-int} as
\begin{equation}
S : t \mapsto S(t) = \{ x \in \R^n ~|~ h(t, x) + F(x) \ni 0 \},
\end{equation}
where $ h(t, x) = - p(t) + f(x) $ corresponding to the previous notation. The third difficulty rises here; the study of the effects of perturbations of $p $ is not equivalent any more to searching the local stability properties of $S$. 
The \emph{strongly regular point} criteria of Robinson \citep{robinson} comes into play now, which guarantees a good behaviour of the problem solutions (cf. Theorem \ref{theorem 2B.7} and the explanations therein). \\
One can encounter this framework with an overlook to the previous case, and consider this problem as an iteration process, that is, for any $t \in [0, 1]$ we have a static case problem. This approach is well known and well studied in the literature, both as a pointwise study (see for example 
\citep{mordukhovich1994stability, kyparisis1987sensitivity, kummer1984generalized, shapiro1994sensitivity, Bianchi2013279, Artacho20101149, uderzo2009some}), 
or as a numerical method and for designing algorithms (see for example 
\citep{dontchev2010newton, ferreira2016kantorovich, ferreira2016unifying, adly2016newton, adly2015newton} and references therein). \\
However, we have an intuition (coming from experimental observations and simulation results) that the variable $x$ could be seen as a function (of time) in this case, and there is some relation between the signal source function (we call it input signal) and this function (we call it output signal). Thus, instead of looking at the sets $S(t)$, we focus on \emph{solution trajectories}, functions like $ x : [0, 1] \to \R^n $ such that 
\begin{equation} \label{solution trajectory for time-varying case-int}
x(t) \in S(t),  \mfa t \in [0, 1],
\end{equation}
that is, $x(\cdot) $ is a selection for $S$ over $ [0,1] $; and search for their regularity properties in function spaces. This study would be the main concern of Chapter \ref{Chapter4} which is mostly our contribution to the subject.\\
There, we provide a smoothness relation between the input signal and the output signal, and study the perturbation effect of the input signal on solution trajectories. 
\end{enumerate}
All the statements, definitions, and examples from books or papers has been cited, though some explanations and changes has been added to adapt them with our setting. Most of the proofs have been rewritten, partly or entirely, to make it more easier to follow\footnote{
This is a risky note, as it uncovers the level of my mathematical knowledge by showing what was assumed not clear or difficult to understand. I am eager to take the risk and also cheerful to learn more.\\
}.
Figures, especially in Chapter \ref{Chapter2}, are obtained from internet, and we did not consider them something crucial to cite. However, when they contain an original idea, they have been explained in the context and truly cited. \emph{Wikipedia} is always a good source to find appropriate photos for electronic components.\\
The thesis is organized in four chapters. We will provide a brief review of their contents here.\\
In Chapter \ref{Chapter1}, the main mathematical preliminaries of the thesis are presented. Most of the definitions and notations that would be used throughout this thesis are introduced in this chapter. After mentioning basic concepts of set-valued maps in Section \ref{Multifunctions: First Properties}, we begin to present some local properties for a general set-valued map in the next two sections. We gather these properties under the general name of \emph{local stability} properties in Section \ref{Local Stability Properties} and \emph{regularity} in Section \ref{section1.3}.\\
We tried to familiarize the reader with these properties by providing several examples to distinguish the delicate differences in the definitions and by describing the relation between these notions (see, for example, Theorems \ref{MR-AP}, \ref{LO-MR relation}, and \ref{Characterization by Calmness of the Inverse}). \\
Section \ref{A Review on Perturbation Results} will be an omen of what would be the main question of this thesis: we have a set-valued map $F$ with a particular stability-like\index{stability-like}\footnote{
We use this phrase in our general descriptions in order to refer to any of the local stability properties of Section \ref{Local Stability Properties} (that is, Aubin property, calmness, linear openness, and isolated calmness) or any of the regularity terms of Section \ref{section1.3} (that is, metric regularity, strong metric regularity, metric sub-regularity, and strong metric sub-regularity). \\
The reason behind this choice is that, first of all, these properties provide a local and not general description of the map under study. Secondly, the term \lq\lq stability\rq\rq \,has already been used in analysis, dynamical systems, and many other areas of mathematics with certain definitions, which is different from what we want to study here.\\ 
}
 property. We consider a perturbing function $ f $ affecting $F$ in a summation form (the reason for using this model will be discussed in Chapter \ref{Chapter2}). Then, we wonder under which conditions the new set-valued map $ G = f + F $, would still have the stability-like property?\\
We can also model a closely related situation: we start from $ f + F $, and then the perturbed model would be of the form $ g + F $. Having the set-valued map fixed, one might think of a possible relation between $ g $ and $ f $. Then, we get the idea to approximate $ f $ in the model which may not be smooth enough, with a smooth function. This is the reason we introduce \emph{approximations} and \emph{estimators} in that section.\\
Chapter \ref{Chapter2} starts by a short review on the theory of electrical circuits, and a brief introduction to some electronic components, which is a wide and deep topic. It would be hard to talk about so many things in few pages, and it might cause more confusion than clarification; but since the aim of this chapter is only to provide the framework of how to model properly a real world situation into a generalized equations setting, we tried to avoid any further details. \\
Some engineering terms may not be explained completely, importance and necessity of studying some circuits are not highlighted well enough, and not all interesting examples are included in the text, but this is the price one should pay eagerly to write a thesis in pure mathematics about a concrete problem.\\
Although we have been inspired by the works already done in this direction (see, for example, \citep{adly, adly2, adly3, Roubal:Thesis:2015}), we found some ambiguity and misinterpretations in the modelling of some circuits in the literature, and thus, tried to write a self-reliant chapter about it. Most of the materials in this chapter are new or, at least, reorganized in a new form to serve the purpose of this thesis\footnote{
This sentence should not be understood in the way that we have discovered the Kirchhoff’s laws, or invented the diodes, or things like this which could be found in electronic books, but simply means that we noticed the $ i -v $ relation of diodes and other semiconductors as set-valued maps and used circuit theory to form a proper generalized equation, not an approximated equation, nor a variational inequality. 
It is the idea of viewing the old problem from a new stand point. \\
}.\\
Section \ref{Formulating the problem} has several important aims. Firstly, it provides the general form of generalized equations that would be considered in the rest of the thesis as an outcome of modelling process. Secondly, it provides physical explanation for the importance of studying the small perturbations of source signals. And in third place, it gives a meaning to the stability-like properties introduced in Chapter \ref{Chapter1}: in the static case, the question about sensitivity of the circuit to small perturbations of the power source, is translated to the question whether the solution mapping of the obtained generalized equation has some stability-like properties or not. For the dynamic case, which is the situation where sources are changing with time known as AC voltages or currents, we only present the framework of the generalized equation. How they could be related to stability-like concepts would be the subject of Chapter \ref{Chapter4}. \\
This chapter ends with a review on an alternative formulation method in Subsection \ref{Alternative Formulation, Pros and Cons}, and examples of circuits we would like to discuss stability of the solution mappings related to them. There, we tried to explain how this formulation could open a door to variational inequality study, and why this model could not be used for all circuits (see for instance, the circuit in Example \ref{Driving an LED with AC, on the importance of understanding the behaviour of the circuit}).\\
Chapter \ref{Chapter3} starts with introducing some tools from variational analysis. Starting from a geometrical interpretation of derivatives, the concept of graphical derivative (with three different approaches) for set-valued maps is introduced in Section \ref{Variational Geometry}, and some calculus rules for these derivatives is provided in subsection \ref{Calculus Rules}. Although there are many results in this area (see for example 
\citep{mordukhovich1994generalized, MORDUKHOVICH19973059, MORDUKHOVICH19951401, durea2012openness, bianchi2015linear, voisei2008sum}), we only present the most needed ones and adapt them to our setting. \\
Then we provide two norm-like maps that assign an extended real number (i.e. $ \R \cup \{- \infty, + \infty \} $) to a \emph{positively homogeneous} set-valued map (cf. Definitions \ref{Positively Homogeneous Mappings}, and \ref{Outer and Inner Norms}) in Subsection \ref{Inner and Outer Norms}. These definitions allow us to characterize the stability-like
properties of a set-valued map in terms of the \emph{inner} and \emph{outer norm} of their graphical derivatives (see for example Proposition \ref{Norm Characterizations}, and Theorems \ref{GDCMR}, \ref{CDCMR theorem}, and \ref{SGDCSMR}).\\
Subsection \ref{Subdifferentials} will introduce \emph{subdifferentials} and \emph{generalized Jacobians} for functions that may not be differentiable. It might be better to introduce subdifferentials before graphical derivatives. Firstly, because they extend the idea of differentiability to a non-smooth function in a graphical way. Secondly, a subdifferential is usually a set-valued map. So, one might be more ready for defining the concept of graphical derivatives for set-valued maps afterwards, but we preferred to keep them where they are as we only use them in Section \ref{Results about Calmness} and Subsection \ref{Strong Metric (Sub-) Regularity with a Nonsmooth Single-valued Part} to show a possible way to continue the stability study when dealing with a generalized equation with a non-smooth single-valued part. \\
The remaining sections of this chapter have a common structure. Each section will discuss one of the stability-like properties for the specific generalized equation $ f(\cdot) + B F(C \cdot) \ni p $, with the following assumptions:
\begin{enumerate}[topsep=-1ex, itemsep=-1ex, partopsep=1ex, parsep=1ex, leftmargin = 9ex]
\item[(A1)] $B$ is injective;
\item[(A2)] $f$ is continuously differentiable in $\R^n$; 
\item[(A3)] $F$ has closed graph;
\item[(A4)] $C$ is surjective; and 
\item[(A5)] there are $ \mmap{F_i}{}{},~ i \in \{ 1,... \, , m \} $ such that $ F(x) = \prod\limits_{i=1}^{m} F_i(x_i) $ whenever $ x = (x_1, ..., x_m)^T \in \R^m $.
\end{enumerate}
We try to use the norm characterization obtained for an arbitrary set-valued mapping in Subsection \ref{Inner and Outer Norms} for each property, and calculus rules of Subsection \ref{Calculus Rules} to go step by step toward a pointwise easy-to-check criteria for the local stability of the solution mapping.\\
The main theorem of each section uses only assumptions (A1) - (A3) (see, for example, Theorem \ref{theorem 3.1}, and \ref{Theorem 4.1}). Results using (A4) and/or (A5) are expressed as corollaries (see, for instance, Corollary \ref{AP2}, and \ref{Corollary 4.1}). Although in many circuits these assumptions may hold, there are specific situations where these conditions may not be satisfied.\\
Among the first three assumptions, (A3) is not such a strong requirement and holds for the $i - v$ characteristics of semiconductors like diodes in our study. In Subsection \ref{Isolated Calmness Without Injectivity Assumption} we try to consider the case where $B$ is not injective, and use the following condition instead of (A1) to provide some similar statements for isolated calmness in this case (cf. Theorem \ref{IC2}).
\begin{itemize}[nolistsep]
\item[]$(\widetilde{\mathrm{A1}})$ Suppose that there is $ \bar{v} \in F(C\bar{z}) $ such that
\begin{equation*}
\bar{p} = f(\bar{z}) + B \, \bar{v} ~\mathrm{~~and~~}~ \overline{ \bigcup_{t > 0} \dfrac{ \mathrm{rge~} F_C - \bar{v} } {t} } \,  \bigcap \mathrm{~ker~} B = \{ 0_{\R^m} \}.
\end{equation*}
\end{itemize}
The process is almost the same as before. Example \ref{IC without A1} shows such a situation in a circuit.\\
In Subsection \ref{Strong Metric (Sub-) Regularity with a Nonsmooth Single-valued Part} we focus on functions $f$ which are not smooth enough to satisfy (A2). A possible approach would be considering the generalized Jacobians and replace (A2) with
\begin{itemize}[nolistsep]
\item[]$(\widetilde{\mathrm{A2}})$ $f$ is locally Lipschitz continuous on $\R^n$.
\end{itemize}
In Section \ref{Results in terms of Metric Regularity}, using the relations between the regularity terms and local stability properties of the inverse map expressed in Section \ref{section1.3}, we restate the results of previous sections in terms of metric regularities in Theorems \ref{SMR-Theorem 2}, and \ref{Summary of MR Criteria for GE in Case of a Nonsmooth Single-valued Part}. Not only gathering all separated statements together would ease the future referring, but also expressing results in terms of metric regularity would be more useful in view of the method we choose for our study in Chapter \ref{Chapter4}.\\
In Chapter \ref{Chapter4} we try to answer the question of perturbation effect when the vector $p$ becomes a function of time, which corresponds to the electronic circuits working with AC sources. In Section \ref{Obtaining a Proper Model} we obtain a proper model for the electrical circuit in this case, and explain the shift of our goal from studying the solution mapping to solution trajectories. In Subsection \ref{Review on an existence theorem}, we express an important theorem (cf. Theorem \ref{theorem 6G.1}) to build our structure. Theorem \ref{theorem 6G.1} and its backbone theorem (Theorem \ref{theorem 5G.3}), were first stated in \citep{dontchev2013}, and then appear in \citep{implicit} with a small modification. Being assured that there does exist a solution trajectory, we search for the smoothness and other properties of the trajectories. \\
Section \ref{Results} starts with a proposition that has a simple proof, but contains a novel idea which considerably eases the study of the problem in the time varying case (cf. Proposition \ref{different auxiliary maps}). In Subsection \ref{Continuity of Solution Trajectories}, we provide results that highlight smoothness dependence of trajectories on the input signal (cf. Proposition \ref{claim1}, and Corollary \ref{claim1-corollary}). \\
In Subsection \ref{Uniform Strong Metric Regularity-subsection}, we deviate from the study about solution trajectories for a while to present a uniform strong metric regularity result along a trajectory, with assumptions weaker than Theorem \ref{theorem 6G.1}. The reason we present Theorems \ref{claim2}, and \ref{claim2-second version} in this subsection is that they are necessary for obtaining the results about perturbation effect in Subsection \ref{Perturbations of the Input Signal}. \\
The final result of this chapter would be Theorem \ref{Existence of a Continuous Trajectory for the Perturbed Problem} that guarantees the existence of a solution trajectory for the perturbed problem which is continuous, and whose distance from the solution trajectory of the (non-perturbed) generalized equation is controlled by the distance between the input signal and its perturbed function.

%% file: Chapters/Chapter01.tex
\chapter{Mathematical Preliminaries} 

\label{Chapter1} 

\lhead{Chapter 1. \emph{Mathematical Preliminaries}} 

\setlength{\epigraphwidth}{.55\textwidth}
\epigraph{\emph{Life is the only art that we are required to practice without preparation, and without being allowed the preliminary trials, the failures and botches, that are essential for training.}}{\textit{Lewis Mumford}}

The main mathematical preliminaries of this thesis are presented in this chapter. 
We start with mentioning basic concepts of set-valued maps in Section \ref{Multifunctions: First Properties}, and then present some local properties for a general set-valued map under the general name of \emph{local stability} properties in Section \ref{Local Stability Properties}. \\
In Section \ref{section1.3} we will introduce the regularity terms, and will clarify how these concepts are related to local stability properties.
The last section of this chapter, Section \ref{A Review on Perturbation Results}, will be devoted to the study of perturbing a set-valued map $F$, having a particular stability-like property, with a function $ f $. We wonder under which conditions the new set-valued map $ G = f + F $, would still have the stability-like property.

\section{Multifunctions: First Properties} \label{Multifunctions: First Properties}
In this section, we review the most needed definitions, terminology, and notations, that we will use throughout the thesis.\\ 
In order to avoid confusion, after this section we would refer to multifunctions as set-valued maps and use \lq\lq function\rq\rq \, only for single-valued maps. The different notation will also help us to separate the maps more easily.

\defn \textbf{(Set-valued Maps)} \label{definition of Multifunctions} \citep[p. 63]{Schirotzek} 
Let $X$ and $Y$ be vector spaces. A mapping $ F : X \longrightarrow 2^{Y} $, which associates to $ x \in X $ a (possibly empty) subset 
$ F(x) $ of $Y$, is called a \textit{multifunction}\footnote{
\textbf{Historical note}. 
The usage of terms \emph{multifunction}, and \emph{multivalued function} as a sort of allowing \lq\lq function\rq\rq\,in mathematics to mean also \lq\lq multivalued function\rq\rq, had some ups and downs in history and finally dropped out of usage at some point in the first half of the twentieth century. \\
It seems that the theory of multivalued functions (with this name) was fairly systematically developed for the first time in 1959 in the book \lq\lq \emph{Topological Spaces: including a treatment of multi-valued functions, vector spaces, and convexity}\rq\rq\,by the French mathematician Claude Berge (cf. \citep{berge1963}).
}
 or \textit{set-valued mapping}, and is denoted by $ F : X \rightrightarrows Y $. \\
The \textit{graph}, \textit{domain}, and \textit{range} of $ F $ are defined, respectively, by
\begin{center}
gph $ F := \{ (x,y) \in X \times Y ~|~ x \in X ~,~ y \in F (x) \}$,\\
dom $ F := \{ x \in X  ~|~  F (x) \neq \emptyset \}$, ~~~~
\hspace*{2.06cm}
\\rge $ F := \{ y \in Y ~|~ y \in F(x) ~~\mathrm{for~some}~ x \in X \}$.~~~
\end{center}
If $A$ is a subset of $X$, we write
$F(A) := \bigcup_{x \in A} F(x) $.

\note
Although the above and many other definitions in this thesis could be expressed in general vector spaces (see, for example, \citep{Mordukhovich}), we keep our attention only on 
$ \R^n $ for some positive natural number $n$.

\note
The multifunction $ F $ is said to be \emph{closed} or \emph{convex} if $ \gph{F} $  is closed or convex, respectively. 
We call $ F $ \emph{closed-valued} or \emph{bounded-valued} if $ F(x) $ is, respectively, a closed, or a bounded subset of $ Y $ for any $ x \in X $. \\
Notice that a closed multifunction is closed-valued but the converse is not true, for example consider the map $\mmap{F}{}{}$ with $F(0) = [-1, 0] $ and $F(x) = [0, 1]$ 
for $ x \neq 0$.

\begin{rem}
When $F$ assigns more than one element to $x$ we say it is \textit{multi-valued} at $x$, and when it assigns no element at all, it is \textit{empty-valued} at $x$. When it assigns exactly one element $y$ to $x$, it is\textit{ single-valued} at $x$, in which case we allow ourselves to write $ F(x) = y $ instead of $ F(x) = \{ y \}$.\\
Thus, a mapping $\smap{f}{n}{m}$ can be identified with the (single-valued) multifunction $\mmap{F}{n}{m} $ defined by $ F(x) = \{ f(x) \}, ~ x \in \R^n $. 
Hence the name \emph{multifunction} is explained. 
Moreover, concepts defined below for set-valued mappings will be also applicable to a function $\smap{f}{n}{m}$ according to this identification.
\end{rem}

Though, at the beginning the concept of set-valued maps may seem to be superficial and an unnecessary complication; the following example shows that it is deeply hidden in the heart of analysis, as in many cases after formulating the problem in terms of a relation (or function) between the input and output data, we are interested in the properties of the inverse map. 

\eg \textbf{(Inverse Functions)}\index{function ! inverse}
Consider a (single valued) function $ \smap{f}{n}{m} $ as a relation that describes the behaviour of system, a relation between inputs $ x \in \R^n $ and outputs $ y \in \R^m$.\\ In many applications we are interested in finding the accurate input which gives us a desired output. Mathematically speaking, we are interested in the \emph{solution mapping} 
\begin{center}
$ S(y) := \{ x \in \R^n ~|~ y = f(x) \} $.
\end{center}
This map is generally a set-valued map and if we have the further assumption on $f$ to be injective, it would be single valued and known as the inverse function
$ S = f^{-1}$.\\
In fact, the inverse of a set-valued map $\mmap{F}{n}{m}$ at a point $\bar{y} \in \R^m$ is defined as 
$$ F^{-1} (\bar{y}) := \{ x \in \R^n ~|~ \bar{y} \in F(x) \}. $$
In this manner, a function $f$ always has an inverse $ f^{-1} $ as a set-valued mapping. The question of an \emph{inverse function} comes down then, to passing to some \lq\lq piece\rq\rq \,of the graph of $ f^{-1} $. Soon, we will define selection and localization for a set-valued map (ref. Definitions \ref{Selection} and \ref{Single-valued Localization}) to describe this \lq\lq piecing\rq\rq \,process. 

\textbf{Terminology and Notation} \citep[p. 7]{implicit}
\begin{itemize}
\item[\textbf{(a)}]
In working with $ \R^n $ we will denote by $ \| x \| $ the Euclidean norm associated with the canonical inner product
\begin{equation*}
<x,x'> = \sum_{j=1}^{n} x_j x'_j ~~\mathrm{ for }~~ x= (x_1 , ... , x_n) ~~ \mathrm {and}~~ x'= (x'_1 , ... , x'_n),
\end{equation*}
namely
\begin{equation*}
\| x \| = \sqrt{ <x,x> } = 
\begin{bmatrix}
 \sum_{j=1}^{n} x_j ^2
\end{bmatrix} ^{1/2}.
\end{equation*}
\item[\textbf{(b)}]
The closed ball around $ \bar{x} $  with radius $ r $ is $ \mathbb{B}_r(\bar{x}) = \{ x ~|~~ \| x - \bar{x} \|  \leq r \} $. We denote the closed unit ball $ \B_1(0) $ by $ \B $. 
\item[\textbf{(c)}]
A \emph{neighborhood} of $\bar{x} $ is any set $U$ such that $ \mathbb{B} _r (\bar{x}) \subset U $ for some $ r > 0 $. (Often the neighborhoods can conveniently be taken to be closed balls themselves.) 
\item[\textbf{(d)}]
The \emph{interior} and \emph{closure} of a set $ C \subset \R^n $ will be denoted by int $ C $ and cl $ C $. Thus, $ \mathrm{int\,} \B _r (x) $ will indicate the open ball around 
$ \bar{x} $  with radius $ r $. 
\item[\textbf{(e)}]
The \emph{distance} from a point $ x \in \R^d $ to a set $B$ is denoted by $ d(x,B) $ and defined as
\begin{equation*}
d (x, B) := \inf_{y \in B} d(x, y),
\end{equation*}
with the convention that $ d(x,\emptyset) = + \infty $. \\
The \emph{Hausdorff distance}\index{Haudorff distance} of two non-empty sets $ A $ and $ B $ is then defined as
\begin{equation*}
d_H \, (A, B) : = \max \{ \, e( A, B) , \: e( B,A) \}
\end{equation*}
where $ e(A, B)  $ is the \emph{excess}\index{excess} of $ A $ beyond $ B $ defined as 
\begin{equation}\label{excess}
e \, (A, B) : = \sup_{ x \in A} d(x, B) = \sup_{ x \in A} \, \inf_{y \in B} d(x, y),
\end{equation}
with the following convention that
\begin{equation*}
e \, (\emptyset, B) = \left \{ \begin{matrix}
0 & & \mathrm{when~} B \neq \emptyset, \\
\infty & & \mathrm{otherwise}.
\end{matrix} \right.
\end{equation*}
\item[\textbf{(f)}] 
Given two set-valued mappings $ F : X \rightrightarrows  Y $
and $ G : Y \rightrightarrows Z $, the \emph{restriction}\index{restriction} $ G\mid_F : Y  \rightrightarrows  Z $ of $G$ on $F$ is defined by
\begin{equation}\label{restriction}
G|_F (y) := 
\left\{ \begin{matrix}
G(y) && \mathrm{~if~} y \in F(x),\\
\emptyset ~~ && \mathrm{otherwise}.
\end{matrix}\right.
\end{equation}
\end{itemize}

We shortly review three basic properties of functions, mostly for the sake of notation and ease of future reference. The set-valued counterparts of calmness and Lipschitz continuity would be introduced in next section, while the derivative definitions for set-valued maps are postponed to Chapter \ref{Chapter3}.

\defn \textbf{(Calmness)} \citep[p. 25]{implicit} \label{calmness for single-valued} \index{calmness ! for functions}
A function $ \smap{f}{n}{m} $ is said to be \emph{calm} at $\bar{x}$ relative to a set $D$ in $\R^n$ if  $\bar{x} \in D \cap \mathrm{dom }\, f $ and there exists a constant 
$ \kappa \geq 0 $ such that
\begin{equation} \label{calmness for single-valued maps}
\norm{f(x) - f(\bar{x}) } \, \leq \, \kappa \norm{ x - \bar{x} } \mfa x \in D \cap \mathrm{dom}\, f.
\end{equation}
The calmness property (\ref{calmness for single-valued maps}) can alternatively be expressed in the form of the inclusion
$$  f(x) \in f(\bar{x}) + \kappa \norm{x - \bar{x}} \B  \mfa x \in D \cap \mathrm{dom }\, f. $$
That expression connects with the generalization of the definition of calmness to set-valued mappings, which we will discuss at length in next section.\\
The \emph{calmness modulus} of $f$ at $\bar{x}$, denoted by $ \mathrm{clm }( f; \bar{x} ) $, is the infimum of the set of values $ \kappa \geq 0 $ for which there exists a neighborhood $D$ of $\bar{x} $ such that (\ref{calmness for single-valued maps}) holds.

\note
According to this definition, if $\bar{x}$ is an isolated point, we have $ \mathrm{clm }( f; \bar{x} ) = 0 $. As long as $\bar{x}$ is not an isolated point of $\mathrm{dom }\, f$, the calmness modulus satisfies
\begin{equation*}
\mathrm{clm }( f; \bar{x} ) = \limsup_{ \genfrac{}{}{0pt}{}{x \, \in \, D \, \cap \, \mathrm{dom} \, f, ~ x \to \bar{x}}{x \neq \bar{x}} } \dfrac{\norm{f(x) - f(\bar{x}) } }{\norm{ x - \bar{x} }}.
\end{equation*}
When $f$ is not calm at $\bar{x}$, from the definition we get $  \mathrm{clm }( f; \bar{x} ) = \infty $. In this way,
$$ f \mathrm{~is~calm~at~} \bar{x} \Longleftrightarrow  \mathrm{clm }( f; \bar{x} ) < \infty .$$

\defn \textbf{(Lipschitz Continuous Functions)} \index{function ! Lipschitz continuous} \hfill \\ \citep[p. 29]{implicit}
A function $ \smap{f}{n}{m} $ is said to be \emph{Lipschitz continuous} relative to a set $D$, or on a set $D$, if $ D \subset \mathrm{dom }\, f$ and there exists a constant $ \kappa \geq 0 $ (Lipschitz constant) such that
\begin{equation} \label{Lipschitz continuity for single-valued maps}
\norm{f(x') - f(x) } \, \leq \, \kappa \norm{ x' - x } \mfa x', x \in D \cap \mathrm{dom }\, f.
\end{equation}
It is said to be Lipschitz continuous \textit{around} $\bar{x}$ when this inequality holds for some neighborhood $D$ of $\bar{x}$.\\
The \emph{Lipschitz modulus} of $f$ at $\bar{x}$, denoted by $ \mathrm{lip }(f; \bar{x}) $, is the infimum of the set of values of $\kappa$ for which there exists a neighborhood $D$ of $\bar{x}$ such that (\ref{Lipschitz continuity for single-valued maps}) holds. Equivalently,
\begin{equation*}
\mathrm{lip }( f; \bar{x} ) = \limsup_{ \genfrac{}{}{0pt}{}{ x, x' \to \bar{x} }{x \neq x'} } \dfrac{\norm{f(x') - f(x) } }{\norm{ x' - x }}. 
\end{equation*}

\note
Note that, by this definition, for the Lipschitz modulus we have $ \mathrm{lip }( f; \bar{x} ) = \infty $ precisely in the case where, for every $ \kappa > 0 $ and every neighborhood $D$ of $\bar{x}$, there are points $x', x \in D $ violating (\ref{Lipschitz continuity for single-valued maps}). Thus,
$$ f \mathrm{~is~Lipschitz~continuous~around~} \bar{x} \Longleftrightarrow  \mathrm{lip }( f; \bar{x} ) < \infty .$$
For an open set $C$, a function $f$ is \emph{locally Lipschitz continuous} on $C$ exactly when $ \mathrm{lip }( f; \bar{x} ) < \infty$ for every $ x \in C $. \\
Every continuously differentiable function on an open set $C$ is locally Lipschitz continuous on $C$.\\

\note \label{notation}
A function $ \smap{f}{n}{m} $ is differentiable\index{function ! differentiable} at a point $\bar{x}$, when $ \bar{x} \in \mathrm{~int~dom~} f$ and there is a linear mapping $ \smap{A}{n}{m} $ with the property that for every $ \epsilon > 0 $ there exists $ \delta > 0 $ with
$$ \norm{ f(\bar{x}+h) -  f(\bar{x}) - A h } \, \leq \, \epsilon \norm{h} ~~\mathrm{for~every~} h \in \R^n \mathrm{~with~} \norm{h} < \delta. $$
If such a mapping $A$ exists at all, it is unique; it is denoted by $Df(\bar{x})$ and is called the \emph{derivative} of $f$ at $\bar{x}$.
The $ m \times n $ matrix that represents the derivative $Df(\bar{x})$ is called the \emph{Jacobian}\index{function ! Jacobian} of $f$ at $\bar{x}$ and is denoted by $ \nabla f(\bar{x}) $.\\
In distinguishing between $Df(\bar{x})$ as a linear mapping and $ \nabla f(\bar{x}) $ as its matrix, one can guard better against ambiguities which may arise in some situations. Also it would provide a better form for representing some results coming latter (see for instance, Example \ref{coderivative of strict differentiable maps} and Theorem \ref{sum rule for coderivatve}). \\

\defn \textbf{(Strict Differentiability)} \label{strict differentiability}\index{function ! strict differentiable} \hfill \\ \citep[p. 34]{implicit}
A function $ \smap{f}{n}{m} $ is said to be \emph{strictly differentiable} at a point $\bar{x}$ if there is a linear mapping $\smap{A}{n}{m} $ such that
\begin{equation}
\mathrm{lip}(e; \bar{x}) = 0 \mathrm{~~~for~~} e(x) : = f(x) - [f(\bar{x})+A(x - \bar{x})].
\end{equation}

In particular, in this case we have that clm$(e; \bar{x}) = 0$ and hence $f$ is differentiable at $\bar{x}$ with $ A = \nabla f(\bar{x}) $, but strictness imposes a requirement on the difference 
$$ e(x) - e(x') = f(x) - [ f(x')+ \nabla f(\bar{x})(x - x')] $$
also when $ x' \not = \bar{x} $. Specifically, it demands for each $ \epsilon > 0 $, the existence of a neighborhood $U$ of $\bar{x}$ such that
\begin{equation} \label{strict diff. relation}
\norm{ \, f(x) - [f(x') + \nabla f(\bar{x})(x - x') ] \, } \, \leq \, \epsilon \, \norm{ x - x'} \mathrm{~~~~for~every~~} x,x' \in U. 
\end{equation}

\begin{eg}
For the function $\smap{f}{}{}$ defined as below, simple calculations show that $f$ is differentiable at $ \bar{x} = 0 $, but there is no $ \epsilon > 0 $ such that \eqref{strict diff. relation} holds. 
\begin{equation*}
f(x) := 
\left \{ \begin{matrix}
x^2 \, \sin( \dfrac{1}{x}) & &  x \neq 0,\\
0~~~~~~~~~ & &   x = 0,
\end{matrix} \right.
\end{equation*}
To observe that, consider the sequences $ x_n := \dfrac{1}{(n + \frac{1}{2}) \pi} $, and $ x'_n := \dfrac{1}{(n + \frac{3}{2}) \pi} $. Therefore, $f$ is not strictly differentiable at $0$.
\end{eg}

\rem \label{strict differentiability properties}
The following statements which could be easily obtained from the above definition, are useful in the sequel, especially when dealing with coderivatives.
\begin{itemize}[noitemsep,topsep=0pt]
\item Every function $f$ that is continuously differentiable in a neighborhood of $\bar{x}$ is strictly differentiable at $\bar{x}$ (cf. \citep[p. 35]{implicit}).
\item Every mapping $f$ strictly differentiable at $\bar{x}$ is Lipschitz continuous around $\bar{x}$, or \emph{locally Lipschitzian} around this point  (\citep[p. 19]{Mordukhovich}) , that is, there is a
neighborhood $U$ of $\bar{x}$ and a constant $ l \geq 0 $ such that 
\begin{equation}
\norm{ f (x) - f (u) } \, \leq \, l \norm{ x - u } \mfa x, u \in U.
\end{equation}
\end{itemize}

\defn \textbf{(Selection)} \citep[p. 64]{Schirotzek} \label{Selection}\index{selection}
Given a set-valued map $ \mmap{F}{n}{m} $, a function $ \smap{f}{n}{m} $ is said to be a \emph{selection} of $ F $ if $ f(x) \in F(x) $ for each 
$ x \in \mathrm{dom} \, F  $.

\begin{eg}
Consider the set-valued mapping $\mmap{F}{}{}$ defined for every $ x \in \R_{+}$ with $F(x) := \{ y \in \R ~|~ y \geq x^3 \} $; whose graph is the \textit{epigraph} of the function $ \smap{f}{}{} $, $f(x) = x^3 $ for $ x \geq 0 $, that is, $ \gph{F} = \mathrm{epi}\,f$.
Then, the function $s$ defined as 
\begin{equation*}
s(x) := 
\left \{ \begin{matrix}
x^2 & & 0 \, \leq \, x \, \leq \, 1,\\
x^4 & &  ~~~~~~ x > 1,
\end{matrix} \right.
\end{equation*}
is a \emph{continuous} selection of $F$.
\end{eg}

\defn \textbf{(Graphical Localization)} \label{Graphical-Localization}\index{localization ! graphical localization} \hfill \\ \citep[p. 4]{implicit} 
Given $ \mmap{F}{n}{m} $ and a pair $ ( \bar{x}, \bar{y}) \in \gph{F} $, a \emph{graphical localization} of $F$ at $\bar{x}$ for $ \bar{y} $ is a set-valued mapping $ \tilde{F} $ such that 
\begin{center}
$ \gph{\tilde{F} } = (U \times V) \cap \gph{F} $ for some neighborhoods $U$ of $\bar{x}$ and $V$ of  $ \bar{y} $, 
\end{center}
so that 
\begin{center}
$ \tilde{F} : x \mapsto 
\left\{\begin{matrix}
F(x) \cap V & \mathrm{when~} x \in U \\ 
\emptyset & \mathrm{otherwise.}
\end{matrix}\right. $
\end{center}

\hspace*{0.5cm}
\note
\textbf{(a)} The inverse of $ \tilde{F} $ has the form
\begin{center}
$ \tilde{F}^{-1} (y) =
\left\{\begin{matrix}
 F^{-1}(y) \cap U & \mathrm{ when~} y \in V \\ 
\emptyset & \mathrm{otherwise,}
\end{matrix}\right. $
\end{center}
and therefore, is a graphical localization of the set-valued mapping $F^{-1} $ at $ \bar{y} $ for $ \bar{x} $. \hfill \\
\textbf{(b)} The neighborhoods $ U $ and $ V $ can often be taken, conveniently, as closed balls $ \mathbb{B}_{a} (\bar{x}) $ and $ \mathbb{B}_{b} (\bar{y}) $. \hfill \\
\textbf{(c)} The domain of a graphical localization $ \tilde{F} $ of $F$ with respect to $U$ and $V$ may differ from $ U \cap \mathrm{dom } F $, and in general depends on the choice of $V$ (like the situation described in  the following example). \hfill \\
To avoid this, specially when one deals with the inverse mapping, some authors will consider a slightly different definition which would be \emph{graphical localization \textbf{around}} $\bar{x}$ for $\bar{y}$ that takes into account only neighborhoods $U \subset \mathrm{ dom } \tilde{F}$ in the above definition (see, for example, \cite{cibulkaRoubal}). \hfill \\ 
This condition guarantees that $\bar{x}$ is an interior point of dom $\tilde{F}$. 

\begin{figure}[ht]
  \begin{minipage}{.49\textwidth}
    \begin{eqnarray*}
		F (x) :=
		\left \{\begin{array}{ll}		
		$\big[$f_1(x), \, f_2(x)$\big]$ & x \neq \bar{x},\\[1em]
		$\big[$f_1(x), \, f_2(x)$\big]$ ~ \cup  ~ \{ \bar{y} \} & x = \bar{x}.
		\end{array} \right.
	\end{eqnarray*} 
	\hspace*{1cm}
  \end{minipage}
  \begin{minipage}{.50\textwidth}
    \centering
		\includegraphics[width=0.95\textwidth]{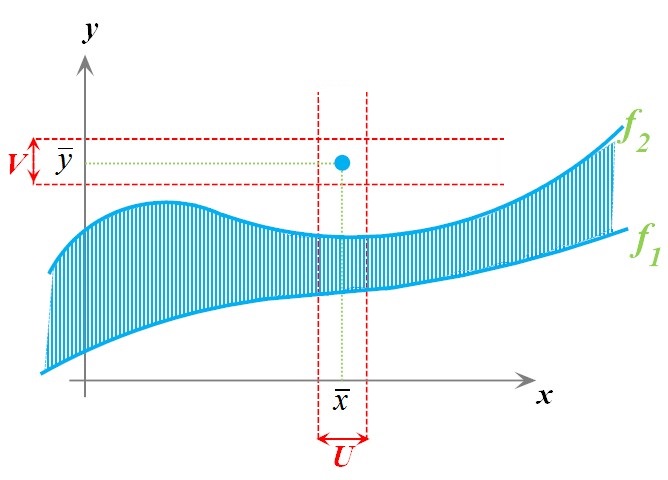}
  \end{minipage}
\caption{The role of neighborhood $V$ in graphical localization}  
\label{fig: graphical localization}
\end{figure}

\eg
Consider the set valued mapping $ \mmap{F}{}{} $ defined in Figure (\ref{fig: graphical localization}), in terms of smooth functions $f_1$, and $f_2$.

One can easily check that based on the chosen $V$, for any $ x \in U \setminus  \left \{ \bar{x} \right \} $, $ \tilde{F}(x) = \emptyset $, and $ \left \{ \bar{x} \right \} = \mathrm{dom } \tilde{F} \not = U \cap \mathrm{dom } F$.

\defn \textbf{(Single-valued Localization)} \citep[p. 4]{implicit} \label{Single-valued Localization}\index{localization ! single-valued localization}
By a \textit{single-valued localization} of $ F $ at $ \bar{x} $  for $ \bar{y} $ will be meant a graphical localization that is a function, its domain not necessarily being a neighborhood of $ \bar{x} $. \\
The case where the domain is indeed a neighborhood of $ \bar{x} $ will be indicated by referring to as single-valued localization of $ F $ \textbf{\textit{around}} $ \bar{x} $ for $ \bar{y} $ instead of just \textbf{\textit{at}} $ \bar{x} $ for $ \bar{y} $.

\eg \textbf{(Classical Inverse Function Theorem)} \citep[p. 11]{implicit} \label{classical inverse function theorem}
Considering a function $\smap{f}{n}{n}$, one might be interested in the possibility of having a single-valued inverse map. The problem is well studied in the literature under the title of ``Inverse Mapping Theorems\rq\rq.  \\
For future reference, we would like to mention an old classical theorem of this type.\\
Let $\smap{f}{n}{n}$ be continuously differentiable in a neighborhood of a point $\bar{x}$ and let $ \bar{y} := f (\bar{x}) $. If $ \nabla f (\bar{x}) $ is non-singular, then $f^{-1}$ has a single-valued localization $s$ around $\bar{y}$ for $\bar{x}$. Moreover, the function $s$ is continuously differentiable in a neighborhood $V$ of $\bar{y}$, and its Jacobian satisfies
$$ \nabla s(y) = \nabla f (s(y)) ^{-1} \mathrm{~~~for~every~~}  y \in V. $$

\eg \textbf{(Implicit Functions)} \citep[p. 4]{implicit}\index{function ! implicit}
In passing from inverse functions to implicit functions, we need to pass from an equation $ f(x) = y $ to one of the form 
\begin{equation} \label{Dini eq. 1}
g(p,x) = 0 ~~ \mathrm{ for~a~function } ~~ g : \R^d \times \R^n \longrightarrow \R^m,
\end{equation}
in which $ p $ acts as a parameter. \\
The question is no longer about inverting $ f $, but the framework of set-valuedness is valuable nonetheless because it allows us to immediately introduce the \emph{solution mapping} 
\begin{equation} \label{Solution mapping1}
\mmap{S}{d}{n}~~ \mathrm{ with }~~ S(p) =\{ x~ | ~ g(p,x) = 0 \}.
\end{equation}
We can then look at pairs $ (\bar{p},\bar{x}) $ in $ \gph{S} $ and ask whether $ S $ has a single-valued localization $ s $ around $ \bar{p} $ for $ \bar{x} $. Such a localization is exactly what constitutes an implicit function coming out of the equation\footnote{
\textbf{Historical note}. 
A very early result in this direction was introduced and proved by Ulisse Dini (1845–1918) in his lecture notes of 1877-78, which is now known as the \textit{classical implicit function theorem} or \textit{Dini’s theorem}; though the set-valued solution mapping $ S $ in (\ref{Solution mapping1}) never enters the picture directly. We express the theorem with our notations as below: \\

\textbf{Dini classical implicit function theorem}.
Let the function  $ f: \R^d \times \R^n \longrightarrow \R^n $ in (\ref{Dini eq. 1}) be continuously differentiable in a neighborhood of $ (\bar{p},\bar{x}) $ and such that $ f (\bar{p},\bar{x}) = 0 $, and let the partial Jacobian of $ f $ with respect to $ x $ at $ (\bar{p},\bar{x}) $, namely $ \nabla_{x} f (\bar{p},\bar{x}) $, be non-singular. \\
Then the solution mapping $ S $ defined in (\ref{Solution mapping1}) has a single-valued localization $ s $ around $ \bar{p} $  for $ \bar{x} $ which is continuously differentiable in a neighborhood $ Q $ of $ \bar{p} $ with Jacobian satisfying 
\begin{equation*}
\nabla s(p) = - \nabla_{x} f( p, s(p) ) ^{-1} \nabla_{p} f( p, s(p) ) ~~\mathrm{ for~every } ~p \in Q.
\end{equation*}
}.

\defn  \textbf{(Monotone Mappings) } \label{Monotone Mappings} \index{monotone map} \hfill \\ \citep[p. 195]{implicit} \citep[p. 104]{Aubin}
A mapping $ \mmap{F}{n}{n} $ is said to be \emph{monotone} if
\begin{equation} \label{monotonicity}
 \inp{y' - y}{x' - x} \geq 0 ~\mathrm{~~~~whenever~~~~} (x',y'), (x,y) \in \gph{F}.
\end{equation}
It is called \emph{maximal monotone}\index{monotone ! maximal monotone} when no more points can be added to $\gph{F}$ without running into a violation of (\ref{monotonicity}).\\
In other words, a monotone set-valued map $F$ is \emph{maximal} if there is no other monotone set-valued map whose graph strictly contains the graph of $F$.\\
$F$ is called \emph{locally monotone}\index{monotone ! locally monotone} \at{x}{y} if $ \rfp{x}{y} \in \gph{F} $ and for some neighborhood $W$ of $\rfp{x}{y}$, one has
\begin{equation} \label{locally monotone criteria}
 \inp{y' - y}{x' - x} \geq 0 ~\mathrm{~~~~whenever~~~~} (x',y'), (x,y) \in \gph{F}  \cap W. 
\end{equation}

\begin{eg}
Consider the set-valued mappings $ \mmap{F_1, F_2}{}{}$ defined as
\begin{equation*}
F_1(x) := 
\left \{ \begin{matrix}
x^3 -1 & & x < 0,\\
\{-1, +1 \} & & x = 0,\\
x^3 + 1 & & x > 0.
\end{matrix} \right.
~~~~~~~\mathrm{and}~~~~~~~
F_2 (x) := 
\left \{ \begin{matrix}
x^3 -1 & & x < 0,\\
[-1, +1 ] & & x = 0,\\
x^3 + 1 & & x > 0.
\end{matrix} \right.
\end{equation*}
It is easy to check that both mappings are monotone. Since $ \gph{F_1} \subset \gph{F_2} $, one concludes that $F_1$ is not a maximal monotone map.
\end{eg}

\rem \label{maximal monotone maps properties}
The following statements, which could be easily obtained from the above definition, are used in the sequel.
\begin{itemize}[noitemsep,topsep=0pt]
\item If $A$ and $B$ are monotone maps and $\lambda,  \mu >0$ are scalars, then $ \lambda A + \mu B $ is also monotone (cf. \citep[p. 105]{Aubin}).
\item Since monotonicity is a property bearing on the graph of $F$, a set-valued map $F$ is monotone (or maximal monotone) if and only if its inverse $F^{-1}$ is monotone (or maximal monotone) (cf. \citep[p. 105]{Aubin}).
\item A necessary and sufficient condition for a set-valued map $F$ to be maximal monotone is that the property
 $$ \inp{u - v}{x - y} \, \geq \, 0, \mfa (y,v) \in \gph{F} $$
is equivalent to $ u \in F(x) $ (cf. \citep[p. 107]{Aubin}).
\end{itemize}

\begin{eg} \textbf{(Single-valued Localization and Selection for Locally Monotone Maps)} \citep{cibulkaRoubal} 
This example will discuss an interesting property of locally monotone maps. The claim is the following: \\
A set-valued mapping $ \mmap{S}{l}{l} $, which is locally monotone at $ \rfp{y}{x} \in \gph{S} $, has a single-valued Lipschitz continuous localization around $\bar{y}$ for $\bar{x}$ if and only if it has a Lipschitz continuous selection around $\bar{y}$ for $\bar{x}$. \\
Find $W$ such that (\ref{locally monotone criteria}) holds. Let $s$ be a local selection for $S$, defined on $ \mathrm{int\,} \B_r (\bar{y}) $ for some $ r > 0 $; and Lipschitz continuous there such that $ \mathrm{int\,} \B_r (\bar{y}) \times \mathrm{int~} \B_{\kappa r} (\bar{x}) \subset W $, where $ \kappa > 0 $ is the corresponding Lipschitz constant.\\
Fix any $ y \in \mathrm{int\,} \B_r (\bar{y}) $. By definition, $ s(\bar{y}) = \bar{x} $, and 
$$ \norm{ s(y) - s( \bar{y}) } \, \leq \, \kappa \, \norm{ y - \bar{y}} \, \leq \, \kappa r, $$
implies that $ s(y) \in  \mathrm{int\,} \B_{\kappa r} (\bar{x})$. Therefore, the point $ s(y)$ lies in $ S(y) \cap  \mathrm{int\,} \B_{\kappa r} (\bar{x}) $. 
It suffices to show that the latter set is singleton. Suppose that this is not the case. Find $ x \in \R^l $ such that 
$$ x \in S(y) \cap  \mathrm{int\,} \B_{\kappa r} (\bar{x}) \mathrm{~~~with~~~} x \neq s(y). $$
Let $ b := \norm{ x - s(y) } > 0 $, and $ c := \frac{ x - s(y) }{b}$. Thus,
\begin{equation} \label{eg-formula01}
\begin{split}
\inp{x}{c} & =  \inp{x}{\dfrac{ x - s(y) }{b}} = \dfrac{1}{b} \inp{x}{x - s(y)} = \dfrac{1}{b} \inp{x - s(y) + s(y)}{x - s(y) } \\
& =  \dfrac{1}{b} \inp{x - s(y)}{x - s(y)} + \dfrac{1}{b} \inp{ s(y)}{x - s(y) } \\
& = \, b + \inp{s(y)}{c}
\end{split}
\end{equation}
Find $ \tau > 0 $ such that $ \kappa \tau < b $ and $ y + \tau c \in \mathrm{int\,} \B_r (\bar{y}) $. Since $ \norm{c} = 1 $, the Cauchy-Schwartz inequality and the Lipschitz continuity of $s$ imply that
\begin{equation} \label{eg-formula02}
\inp{ s(y + \tau c) - s(y) }{c} \, \leq \, \norm{ s(y + \tau c) - s(y) } \, \norm{c} \, \leq \, \kappa \, \norm{ (y + \tau c) - y } = \kappa \tau .
\end{equation}
Since $ \big( y + \tau c, s(y + \tau c) \big) $ and $ (y, x) $ are in $ \gph{S} \cap W $, local monotonicity of $S$, Condition (\ref{locally monotone criteria}), reveals that 
$ 0 \, \leq \, \inp{s(y + \tau c) - x}{y + \tau c - y} \, = \, \tau \inp{s(y + \tau c) - x}{c} $. Thus,
\begin{equation} \label{eg-formula03}
\inp{s(y + \tau c)}{c} \, \geq \, \inp{x}{c}.
\end{equation}
Now, one can use these equations to obtain 
\begin{equation*}
b + \inp{s(y)}{c} \overset{\mathrm{\eqref{eg-formula01}}}{=} \inp{x}{c}  \overset{\mathrm{\eqref{eg-formula03}}}{ \, \leq \,} \inp{s(y + \tau c)}{c} \overset{\mathrm{\eqref{eg-formula02}}}{ \, \leq \,} \inp{s(y)}{c} + \kappa \tau \, < \, \inp{s(y)}{c} + b.
\end{equation*}
We arrived at a contradiction, therefore $ S(y) \cap  \mathrm{int\,} \B_{\kappa r} (\bar{x}) = \{ s(y) \} $ for each $ y \in \mathrm{int\,} \B_{r} (\bar{y}) $. The opposite implication is trivial.
\end{eg}

\section{Local Stability Properties}\label{Local Stability Properties}
This section is devoted to introducing the main local properties we would consider in this thesis. These properties which mostly has a well known counterpart for single-valued maps, let us stud the local behaviour of a set-valued map. \\
Later in Chapter \ref{Chapter2}, we will see how slight modifications in formulating the problem will result to each of these properties. 

\defn \textbf{(Locally Closed Set) } \cite[p. 172]{implicit} \index{locally closed set}
A set $C$ is said to be \emph{locally closed} at $x \in C$ if there exists a neighborhood $U$ of $x$ such that the intersection $ C \cap U $ is closed.\\
It could be equivalently defined as the existence of a scalar $ r > 0$ such that the set $ C \cap \B_{r} (x) $ is closed.

\defn \textbf{(Aubin Property)} \label{AP} \index{Aubin property} \hfill \\ \citep[p. 172]{implicit}
A mapping $ \mmap{S}{m}{n} $ is said to have the \emph{Aubin property\footnote{
\textbf{Historical note}. 
Actually, the property that was introduced by Aubin in \citep{Aubin} under the name of \lq\lq pseudo-Lipschitz\rq\rq, with the formulation \eqref{Aubin2}. 
Under closedness of the values of the mapping both formulations and their constants agree, see \citep{artacho2007regular}; but without this assumption, the constant $\kappa$ in (\ref{Aubin2}) might be slightly larger than $\kappa$ in (\ref{Aubin1}). For instance, consider the following mapping $ \mmap{S}{}{} $ defined as
\begin{equation*}
S(y) := 
\left \{ \begin{matrix}
(- \infty, y] & & \mathrm{~for~} y \mathrm{~rational},~~\\
(- \infty, y) & & \mathrm{~for~} y \mathrm{~irrational}.
\end{matrix} \right.
\end{equation*}
Here (\ref{Aubin1}) holds around $(0, 0)$ for $ \kappa = 1 $ while (\ref{Aubin2}) is only valid for $ \kappa > 1 $.\\
}} 
at $ \bar{y} \in \R^m $ for $ \bar{x} \in \R^n $ if $ \bar{x} \in S(\bar{y}) $, the graph of $ S $ is locally closed at $ (\bar{y}, \bar{x}) $, and there is a constant $ \kappa \geq 0 $ together with neighborhoods $ U $ of $ \bar{x} $ and $ V $ of $ \bar{y} $ such that\\
\begin{equation}\label{Aubin1}
e(S(y') \cap U, S(y) ) \leq \kappa ~\| y'- y \| \mathrm{~~~for~all~~} y', y \in V,
\end{equation}
or equivalently, there exist $ \kappa , U $ and $ V $, as described, such that 
\begin{equation}\label{Aubin2}
S(y') \cap U \subset S(y) + \kappa ~\| y' - y \| \mathbb{B}  \mathrm{~~~for~all~~} y',y \in V.
\end{equation}
The infimum of $ \kappa $ over all such combinations of $ \kappa, U,$ and $V, $ is called the \emph{Lipschitz modulus} of $ S $ at $ \bar{y}$ for $\bar{x}$ and is denoted by lip $ (S; \bar{y}| \bar{x}) $.

\note
\textbf{(a)}  The absence of this property is signalled by lip $ (S; \bar{y}| \bar{x})  = \infty$.\\
\textbf{(b)} When $ S $ is single-valued on a neighborhood of $ \bar{y} $, the Lipschitz modulus lip $ (S; \bar{y}| ~S(\bar{y}) ) $ equals the usual Lipschitz modulus lip $ (S; \bar{y}) $ for functions.\\
\textbf{(c)} It is not claimed that (\ref{Aubin1}) and (\ref{Aubin2}) are themselves equivalent, although this is true when $S(y)$ is closed for every $y \in V$. Nonetheless, the infimum furnishing lip $ (S; \bar{y}| \bar{x}) $ is the same whichever formulation is adopted.\\

\rem
A mapping $\mmap{S}{m}{n}$ is said to be \emph{Lipschitz continuous} relative to a (non-empty) set $D$ in $\R^m$ if $ D \subset \mathrm{dom}\, S $, $S$ is closed-valued on $D$, and there exists $ \kappa \geq 0 $ (Lipschitz constant) such that
$$S(y') \subset S(y) + \kappa \, \| y' - y \| \, \B \mfa y',y \in D,$$
or equivalently, there exists $ \kappa \geq 0 $ such that
$$ d(x,S(y)) \leq \kappa \, d(y,S^{-1}(x) \cap D) \mfa x \in \R^n \mathrm{~and~} y \in D. $$
The similarity between \emph{Aubin property} and this definition is the reason why some authors refer to having the \emph{Aubin property} as being ``\emph{Lipschitz-like}'' or ``\emph{pseudo-Lipschitz}'' \footnote{ 
This definition is in fact an equivalent way to define Lipschitz continuity for set-valued maps, based on the Proposition 3C.1 \citep[ p.161]{implicit}. This form suits better for the purpose of comparison.}.\\

\eg \citep[p. 172]{implicit}
Consider the set-valued mapping $ \mmap{S}{}{} $ defined as 
\begin{eqnarray}
S(y) =  \left\{\begin{matrix}
\{0, 1 + \sqrt{y} \}  & \mathrm{~~for~} y \geq 0, \\ 
0~~~~~~~~~~~~ & \mathrm{~~for~} y < 0.
\end{matrix}\right.
\end{eqnarray}
At $0$, the value $S(0)$ consists of two points, $0$ and $1$. This mapping has the Aubin property at $0$ for $0$ but not at $0$ for $1$. 
To see the latter, consider sequences $y_n = \frac{1}{n} $, and $ y'_n = \frac{2}{n}$ for $ n \in \N $. For any neighborhood $V$ of $\bar{y} = 0 $, there is $ N \in \N $ such that for $ n > N$, points $y_n$, and $y'_n$ are inside $V$. Then, no matter how small would be the neighborhood $U$ of $\bar{x} = 1$, $S(y'_n) \cap U$ would include $ 1+ \sqrt{y'_n}$ and thus, $ e \big( S(y'_n) \cap U, S(y_n) \big) = \dfrac{y'_n - y_n}{\sqrt{y'_n} + \sqrt{y_n}}$. 
To have the Aubin property, $\kappa$ must satisfy the following inequality
$$ \dfrac{\sqrt{n}}{1 + \sqrt{2}} \, \leq \, \kappa, $$
which is absurd. This example shows that the Aubin property is tied to a particular point in the graph of the mapping.

The Aubin property could alternatively be defined with one variable ``free,” as shown in the next proposition.

\begin{prop} [\textbf{Alternative Description of Aubin Property}] \label{alternative Aubin} \hfill \\ \citep[p. 176]{implicit}
A mapping $ \mmap{S}{m}{n} $ has the Aubin property at $\bar{y}$ for $\bar{x}$ with constant $ \kappa > 0 $ if and only if $ \bar{x} \in S(\bar{y}) $, $ \gph{S} $ is locally closed at $\rfp{y}{x} $, and there exist neighborhoods $U$ of $\bar{x}$ and $V $ of $ \bar{y}$ such that 
\begin{equation}\label{propformula}
e \, (S(y') \cap U, S(y) ) ~ \leq ~ \kappa ~ \| y' - y \| \mathrm{~~~~~for~all~} y' \in  \R^m \mathrm{~and~} y \in V.
\end{equation}
\end{prop}
\begin{proof}
Clearly, (\ref{propformula}) implies (\ref{Aubin1}). Assume (\ref{Aubin1}) with corresponding $ U $ and $V$ and choose positive $a$ and $b$ such that $\B_a (\bar{x}) \subset U $ and $ \B_b (\bar{y}) \subset V $. Let $ 0 < a' < a $ and $ 0 < b'< b $ be such that 
\begin{equation}\label{propformula2}
2 \kappa  b' \, + \, a' \leq ~ \kappa  b.
\end{equation}
For any $ y \in \B_{b'} (\bar{y}) $ we have from (\ref{Aubin1}) that
\begin{equation*}
d (\bar{x}, S(y) ) ~\leq ~ e \, (S(\bar{y}) \cap U, S(y) )  ~\leq ~ \kappa ~ \| y - \bar{y} \| ~\leq ~ \kappa b' ,
\end{equation*}  
hence
\begin{equation} \label{propformula3}
e \, (\B_{a'}(\bar{x}), S(y) ) ~ \leq ~ \kappa b'+a' . 
\end{equation}
Take any $ y' \in \R^m $. If $ y' \in \B_b (\bar{y}) $ the inequality in (\ref{propformula}) comes from (\ref{Aubin1}) and there is nothing more to prove. 
Assume $ \|y' - \bar{y} \| > b $. Then $ \| y -y' \| > b - b' $ and from (\ref{propformula2}),  
\begin{equation*}
\kappa b' + a' ~ \leq ~ \kappa (b - b') ~ \leq ~ \kappa \|y - y' \| .
\end{equation*}
Using this in (\ref{propformula3}) we obtain
\begin{equation}
e \, (\B_{a'}(\bar{x}), S(y) ) ~ \leq ~ \kappa ~ \| y' - y \| 
\end{equation} 
and since $ S(y') \cap \B_{a'} (\bar{x}) $ is obviously a subset of $ \B_{a'}(\bar{x}) $, we come again to (\ref{propformula}).
\end{proof}

The following proposition will highlight an interesting fact of a map with Aubin property, which we will use later in Chapter \ref{Chapter3}.

\begin{prop} [\textbf{Local Non-emptiness}] \label{AP properties1}\hfill \\ \citep[p. 173]{implicit}
If $ \mmap{S}{m}{n} $ has the Aubin property at $\bar{y}$ for $\bar{x}$, then for every neighborhood $U$ of $\bar{x}$ there exists a neighborhood $V$ of $\bar{y}$ such that $ S(y) \cap U \not = \emptyset \mathrm{~~for all~} y \in V $.
\end{prop}
\begin{proof}
The inclusion (\ref{Aubin2}) for $ y'= \bar{y} $ yields
\begin{equation*}
 \bar{x} \in S (y) + \, \kappa ~ \| y - \bar{y} \| \, \B \mathrm{~~~for~every~} y \in V, 
\end{equation*}
which is the same as
\begin{equation*}
 (\bar{x} + ~ \kappa \,  \| y - \bar{y} \| \, \B) \cap S(y) \not = \emptyset \mathrm{~~~for~every~} y \in V.
\end{equation*}
That is, $ S(y) $ intersects every neighborhood of $ \bar{x} $ when $y$ is sufficiently close to $ \bar{y}$.
\end{proof}

The next property would be the \lq\lq one-variable\rq\rq \,version of the Aubin property. We have already introduced a calm function (cf. Definition \ref{calmness for single-valued}); it is now time to define the set-valued counterpart.

\defn \textbf{(Calmness)} \label{C definition} \index{calmness ! for set-valued maps} \hfill \\ \citep[p. 197]{implicit}
A mapping $ \mmap{S}{m}{n} $ is said to be \emph{calm} at $ \bar{y} $ for $ \bar{x} $ if $ \rfp{y}{x} \in \gph{S} $, and there is a constant $ \kappa \geq 0 $ along with neighborhoods $ U $ of $ \bar{x} $ and $ V $ of $ \bar{y} $ such that
\begin{equation}\label{calmness1}
e(S(y) \cap U, S(\bar{y})) \leq \kappa ~ \| y - \bar{y} \| \mathrm{~~~~for~all~} y \in V.
\end{equation}
Equivalently, the property in (\ref{calmness1}) can be also written as
\begin{equation}\label{calmness2}
S(y) \cap U \subset S(\bar{y}) +  \kappa ~ \| y - \bar{y} \| \, \B \mathrm{~~~~for~all~} y \in V
\end{equation}
although perhaps with larger constant $\kappa $.\\
The infimum of $\kappa $ over all such combinations of $\kappa $, $U$ and $V$ is called the \textit{calmness modulus} of $S$ at $\bar{y}$ for $\bar{x}$ and is denoted by clm $(S; \bar{y}|\bar{x})$.\\

\note
\textbf{(a)} The absence of this property is signalled by clm $(S; \bar{y}|\bar{x}) = \infty $.\\
\textbf{(b)} As in the case of the Lipschitz modulus, it is not claimed that (\ref{calmness1}) and (\ref{calmness2}) are themselves equivalent; anyway, the infimum furnishing clm $(S; \bar{y}|\bar{x}) $ is the same with respect to either of them.

\defn \textbf{(Isolated Calmness)} \label{IC definition}\index{calmness ! isolated calmness} \hfill \\ \citep[p. 201]{implicit}
A mapping $ \mmap{S}{m}{n} $ is said to have the \textit{isolated calmness} property if it is calm at $ \bar{y} $ for $\bar{x}$ and, in addition, $S$ has a graphical localization at $ \bar{y} $ for $\bar{x}$ that is single-valued at $\bar{y}$ itself (with value $\bar{x}$). Specifically, this refers to the existence of a constant $ \kappa \geq 0 $ and neighborhoods $ U$ of $\bar{x}$ and $V$ of $\bar{y}$ such that
\begin{equation}
\| x - \bar{x} \| ~\leq \kappa~ \| y - \bar{y} \| \mathrm{~~~~when~} x \in S(y) \cap U \mathrm{~~and~~}  y \in V\footnote{
Isolated calmness was formally introduced by A. L. Dontchev in \citep{dontchev1995characterizations} under the name \lq\lq local upper Lipschitz continuity at a point\rq\rq.
}.
\end{equation}

\note
\textbf{(a)} Observe that in this definition $S(\bar{y}) \cap U $ is a singleton, namely the point $\bar{x}$, so $\bar{x}$ is an isolated point in $S(\bar{y})$, hence the terminology. \\
\textbf{(b)} Isolated calmness can equivalently be defined as the existence of a (possibly slightly larger) constant $\kappa$ and neighborhoods $U$ of $\bar{x}$ and $V$ of $\bar{y}$ such that
\begin{equation}
S(y) \cap U \subset  ~ \bar{x} + \kappa ~  \| y- \bar{y} \| \, \B \mathrm{~~~when~} y \in V.
\end{equation}
 
The last property of this subsection would be Linear openness, for which, we remind the definition of an open map.
\defn \textbf{(Openness)} \citep[p. 180]{implicit}
A mapping $ \mmap{F}{n}{m} $ is said to be \emph{open} at $\bar{x}$ for $\bar{y}$ if $ \bar{y} \in F(\bar{x}) $, and for every neighborhood $U$ of $\bar{x}$, $ F(U) $ is a neighborhood of $ \bar{y} $.

\defn \textbf{(Linear Openness)} \label{linear openness-def}\index{open at linear rate}\index{linear openness} \hfill \\ \citep[p. 180]{implicit} 
A mapping $ \mmap{F}{n}{m} $ is said to be \emph{linearly open} at $\bar{x}$ for $\bar{y}$  when $ \bar{y} \in F(\bar{x}) $, the graph of $F$ is locally closed at $\rfp{x}{y}$, and there is a constant $ \kappa > 0 $ together with neighborhoods $U$ of $ \bar{x} $ and $V$ of $\bar{y} $ such that
\begin{equation}\label{linear openness}
F \, (x + \, \kappa r \, \mathrm{int\,} \B )  \supset \, [ F(x) + \, r \, \mathrm{int\,} \B ]  \cap V \mathrm{~~~~~for~all~} x \in U \mathrm{~~and~all~} r > 0.
\end{equation}

\note
Linear openness is in fact, a particular case of openness and follows from (\ref{linear openness}) for $ x =\bar{x} $, and the following example shows the reverse is not true in general. Linear openness postulates openness \textit{around} the reference point with balls having proportional radii. \\
For example, consider the function $\smap{f}{}{}$ defined as $ f(x) = x^3 $ for every $ x \in \R $. Obviously, $f$ is open around the origin, but in order to have the linear openness at $ \bar{x} = 0$ for $ \bar{y} =0 $, we must find 
$ \kappa > 0 $ together with \NUV such that \eqref{linear openness} holds. \\
In particular, let $ x = \bar{x} = 0 $. Thus, 
$ f(x) + \, r \, \mathrm{int\,} \B = \mathrm{int\,} \B_r (0) $, and we can also compute $ f (x + \, \kappa r \, \mathrm{int\,} \B ) = \mathrm{int\,} \B_{\kappa^3 \, r^3} (0) $. Hence, we must have $ \kappa^3 > \dfrac{1}{r^2}$, which implies $ \kappa \longrightarrow \infty $ as $ r \longrightarrow 0 $. Therefore, it does not have the linear openness property at the origin. 

\section{Regularity Terms and Their Relation with Local Stability Properties} \label{section1.3}

In this section, we will introduce a regularity concept in terms of the distance of sets, hence a \lq\lq metric\rq\rq \,regularity. The first benefit of these metric regularity terms (the four different versions) is to gather closely related properties which describe the local behaviour of a set-valued map under a single roof, that is metric regularity. \\
The small differences between these definitions and their relation with local stability properties defined in the previous section is studied throughout this section. 
It is worth mentioning that Ioffe in \citep{ioffe2016metric, ioffe2016metric2} has provided a comprehensive note on metric regularity, with many old and new results.

\defn \textbf{(Metric Regularity)}\label{def. mr}\index{metric regularity ! metric regularity} \hfill \\ \citep[p. 176]{implicit}
A mapping $ \mmap{F}{n}{m} $ is said to be \emph{metrically regular}\footnote{
\textbf{Historical note}. 
Although the term \lq\lq regularity\rq\rq \,has been already used in the literature to describe a similar property (see, for example \citep{ioffe1979regular, robinson}), it was firstly Borwein \citep{borwein1986stability} who called this property \emph{metric regularity}.
Let $f$ be a $ C ^1 $ function from a Banach space $X$ (or from an open set in $X$) into another Banach space $Y$. It was usually called \emph{regular} at $z$ if $f' (z) $ maps $X$ onto $Y$.\\
} 
at $\bar{x}$ for $\bar{y}$ when $ \bar{y} \in F(\bar{x})$, the graph of $F$ is locally closed at $ \rfp{x}{y} $, and there is a constant $ \kappa \geq 0 $ together with neighborhoods $U$ of $\bar{x}$ and $V$ of $\bar{y} $ such that
\begin{equation}\label{metric regularity}
d \big(x, F^{-1} (y) \big) \leq \kappa~ d(y, F(x)) \mathrm{~~~whenever~~} (x,y) \in U \times V.
\end{equation}
The infimum of $\kappa $ over all such combinations of $\kappa $, $U$ and $V$ is called the \textit{regularity modulus} of $F$ at $\bar{x}$ for $\bar{y}$ and is denoted by reg $(F; \bar{x}|\bar{y})$.\\

\note 
The absence of metric regularity is signalled by reg $(F; \bar{x}|\bar{y} ) = \infty $. 

\begin{eg} \textbf{(Metric Regularity for a Bounded Linear Operator)} \citep{artacho2007regular}\\
The concept of metric regularity goes back to the classical Banach open mapping principle. It is a well known result by Lyusternik and Graves \citep{durea2014introduction} that a stronger property than the usual openness can be deduced from the open mapping principle (which we already introduced it as linear openness, see Definition \ref{linear openness-def}).\\
For any bounded linear mapping $A$ between two Banach spaces $X$ and $Y$, the Banach open mapping theorem states the equivalence between:
\begin{enumerate} [topsep=-1ex, itemsep=0ex, partopsep=0ex, parsep=0ex, leftmargin = 7ex]
\item[(a)] $A$ is surjective;
\item[(b)] $A$ is an open map (at every point);
\item[(c)] there exists $ \kappa > 0 $ such that for every $ y \in Y $ there exists $ x \in X $ with 
$$ y = Ax \mathrm{~~and~~} \norm{x} \, \leq \, \kappa \norm{y}. $$
\end{enumerate}
Consider a linear function $ \smap{f}{n}{m} $ which is surjective; then, according to (c), there exists $ \kappa > 0 $  such that 
$$ d \big( 0, f^{-1} (y) \big) \, \leq \, \kappa \norm{y} \mfa y \in \R^m = \mathrm{\, rge} \,f. $$
Because of linearity, for any $ x \in \R^n $ we have $ f^{-1} (y - f(x) ) = f^{-1}(y) - x $, and then 
$$ d(x, f^{-1}(y) ) = d \big(0, f^{-1}(y - f(x)) \big) \, \leq \, \kappa \norm{y - f(x)} = \kappa \, d \big( y, f(x) \big). $$
This is a global form of the metric regularity condition for the function $f$. Hence, we showed that, in terms of metric regularity, the equivalence between (a) and (c) can
be stated as follows: $f$ is surjective if and only if it is metrically regular at any point $ (x, y) \in \R^n \times \R^m $.
\end{eg}

\rem \label{numerical MR}
\citep[p. 177]{implicit} Metric regularity is a valuable concept in its own right, especially for numerical purposes (see also \citep{dontchev2003radius}). For a general set-valued mapping $ F $ and a vector $y$, it gives an estimate for how far a point $x$ is from being a solution to the generalized equation $F(x) \ni y$ in terms of the \lq\lq residual\rq\rq
\footnote{
\textbf{Historical note}. 
The term \emph{residual function} seems to have an origin in the theory of Error Bounds, as explicitly discussed in \citep[chapter 6]{facchinei2003finite}.\\
Let $X$ be a given subset of $ \R^n$. We are interested in obtaining inequalities that bound the distance $d(x,X)$ from points $x \in S$ (a given set) to $X$ in terms
of a computable nonnegative-valued function $r : S \cup X \rightarrow \R_{+} $ whose zeros coincide with the elements of $X$; that is,
\begin{equation*}
r(x) = 0 \, \Leftrightarrow \, x \in X.
\end{equation*}
Given a residual function $r$ of the set $X$ and a subset (called \emph{test set}) $S$ of $\R^n$, we wish to establish the existence of positive scalars $c_1, c_2, \gamma_1$, and $\gamma_2$ such that
\begin{equation*}
c_1 r(x)^{\gamma_1} \leq \, d(x,X) \, \leq c_2 r(x)^{\gamma_2},~~~~~~ \forall x \in S.
\end{equation*}
Such an inequality is called an \emph{error bound} for the set $X$ with residual $r$ and with respect to $S$. \\
Since the residual function is a computable quantity whereas the exact distance function is not (because we do not know the elements of the set $X$), the former function can therefore be used as a surrogate of the latter function for computational purposes; such as in the design of solution algorithms for computing an (approximate) element of $X$.
}
$ d(y, F(x) )$.\\
To be specific, let $\bar{x}$ be a solution of the inclusion $ \bar{y} \in F(x) $, let $F $ be metrically regular at $ \bar{x} $ for $ \bar{y} $, and let  $ x_a $ and $ y_a  $ be approximations to $ \bar{x} $ and $ \bar{y} $, respectively. Then from (\ref{metric regularity}), one can deduce that the distance from $x_a$ to the set of solutions of the inclusion $ y_a \in F(x) $ ( i.e. $F^{-1}(y_a) $ ) is bounded by the constant $ \kappa $ times the residual $ d(y_a,F(x_a)) $.\\
In applications, the residual is typically easy to compute or estimate, whereas finding a solution might be considerably more difficult. Metric regularity says that there exists a solution to the inclusion $ y_a \in F(x) $ at distance from $x_a$ proportional to the residual. In particular, if we know the rate of convergence of the residual to zero, then we will obtain the rate of convergence of approximate solutions to an exact one. 

Apart from the importance of metric regularity mentioned in the previous remark, the next two theorems reveal another aspect of importance of this concept, which would be its relation with the Aubin property and linear openness.

\begin{thm} [\textbf{Equivalence of Metric Regularity and the Aubin Property of the Inverse}] \citep[p. 177]{implicit} \label{MR-AP}\\
A set-valued mapping $ \mmap{F}{n}{m} $ is metrically regular at $\bar{x}$ for $\bar{y}$ with a constant $ \kappa >0 $ if and only if its inverse $ \mmap{F^{-1}}{m}{n} $ has the Aubin property at $\bar{y}$ for $\bar{x}$ with constant  $ \kappa >0 $, i.e. there exist neighborhoods $ U $ of  $\bar{x}$ and  $V$ of  $\bar{y}$ such that
\begin{equation}\label{thmformula}
 e \, (F^{-1} (y') \cap U, F^{-1}(y) ) \leq \kappa ~ \| y' - y \| \mathrm{~~~for all~~} y',y \in V. 
\end{equation}
Thus, 
\begin{equation} \label{thmformula 5}
\mathrm{lip } (F^{-1}; \bar{y} | \bar{x} ) = \mathrm{reg } (F ; \bar{x} | \bar{y} ).
\end{equation}
\end{thm}
\begin{proof}
Clearly, the local closedness of the graph of $ F $ at $ \rfp{x}{y} $ is equivalent to the local closedness of the graph of $ F^{-1} $ at $ \rfp{y}{x} $. \\
Let $ \kappa > \mathrm{reg } (F; \bar{x}| \bar{y}) $; then there are positive constants $ a $ and $b$ such that (\ref{metric regularity}) holds with $ U = \B_a (\bar{x}), ~~ V = \B_b (\bar{y}) $ and with this $ \kappa $. Without loss of generality, assume $ b <  \frac{a}{\kappa} $ (to be more precise, one can take $b_0 < \min \, \{ ~b, ~ a/ \kappa ~ \}$ ). Choose $ y,y' \in \B_b (\bar{y}) $. If $ F^{-1} (y) \cap \B_a (\bar{x}) = \emptyset $, then $ d( \bar{x},F^{-1}(y)) > a $. 
But then the inequality (\ref{metric regularity}) with $ x = \bar{x} $ yields 
\begin{equation*}
a ~ \leq ~ d( \bar{x},F^{-1}(y)) ~ \leq ~ \kappa ~d(y, F(\bar{x})) ~ \leq ~  \kappa ~\| y - \bar{y}\|  ~ \leq ~\kappa b ~ < ~ a,
\end{equation*}
a contradiction. Hence there exists $ x \in F^{-1}(y) \cap \B_a (\bar{x}) $, and for any such $ x $ we have from (\ref{metric regularity}) that
\begin{equation}
d(x, F^{-1}(y')) ~ \leq~ \kappa ~ d(y',F(x)) ~ \leq ~ \kappa ~\| y - y' \|.
\end{equation} 
Taking the supremum with respect to $ x \in F^{-1}(y) \cap \B_a (\bar{x}) $ we obtain (\ref{thmformula}) with $ U = \B_a (\bar{x}) $ and  $ V = \B_b (\bar{y}) $, and therefore (considering the infimum in the definition of Lipschitz modulus)
\begin{equation} \label{thmformula3}
\mathrm{reg } (F; \bar{x} | \bar{y}) ~ \geq ~ \mathrm{lip }(F^{-1}; \bar{y}| \bar{x}).
\end{equation}
Conversely, suppose there are neighborhoods $ U $ of $ \bar{x}$ and $ V $ of $ \bar{y} $ along with a constant $ \kappa $ such that (\ref{thmformula}) holds. Take $U$ and $V$ smaller if necessary so that, according to Proposition (\ref{alternative Aubin}), we have
\begin{equation}\label{thmfomula2}
e(F^{-1}(y') \cap U, F^{-1}(y)) ~ \leq ~ \kappa ~ \| y' - y \| \mathrm{~~~~for all~} y' \in \R^m \mathrm{~and~} y \in V.
\end{equation}
Let $ x \in U $ and $ y \in V$. If $ F(x) \not = \emptyset $, then for any $ y' \in  F(x) $ we have $ x \in F^{-1}(y') \cap U $.
From (\ref{thmfomula2}), we obtain
\begin{equation}
d( x, F^{-1}(y)) ~ \leq ~ e(F^{-1}(y') \cap U, F^{-1}(y)) ~ \leq ~ \kappa ~ \| y - y' \|.
\end{equation}
This holds for any $ y' \in F(x) $, hence, by taking the infimum with respect to $ y' \in F(x) $ in the last expression we get 
\begin{equation}
d (x, F^{-1}(y)) ~ \leq ~ \kappa ~ d(y,F(x)).
\end{equation}
(If $ F(x) = \emptyset $, then because of the convention $ d(y, \emptyset ) = \infty $, this inequality holds automatically.) Hence, $ F $ is metrically regular at $\bar{x}$ for $\bar{y}$ with a constant $ \kappa $. Then we have $ \kappa ~ \geq ~ \mathrm{reg } (F; \bar{x} | \bar{y}) $ and hence $ \mathrm{reg } (F; \bar{x} | \bar{y}) ~ \leq ~  \mathrm{lip } (F^{-1}; \bar{y}| \bar{x}) $. This inequality together with (\ref{thmformula3}) results in (\ref{thmformula 5}).
\end{proof}

\note \citep[p. 180]{implicit}
From the equivalence of metric regularity of $F $ at $\bar{x}$ for $\bar{y}$ and the Aubin property of $F^{-1} $ at $\bar{y}$ for $\bar{x}$, and Proposition \ref{AP properties1}, we obtain that if a mapping $ F $ is metrically regular at $\bar{x}$ for $\bar{y}$, then $F$ is open at $\bar{x}$ for $\bar{y}$. \\
Metric regularity is actually equivalent to the stronger version of the openness property as mentioned in the following theorem: \\

\begin{thm} [\textbf{Equivalence of Linear Openness and Metric Regularity}]  \label{LO-MR relation}\hfill \\ \citep[p. 180]{implicit}
A set-valued mapping $ \mmap{F}{n}{m} $ is linearly open at $\bar{x}$ for $\bar{y}$ if and only if $ F $ is metrically regular at $\bar{x}$ for $\bar{y}$. In this case the infimum of $ \kappa $ for which (\ref{linear openness}) holds is equal to reg $(F; \bar{x} | \bar{y})$.
\end{thm}
\begin{proof}
Both properties require local closedness of the graph of $F$ at the reference point.\\
Let (\ref{linear openness}) hold. Choose $ y \in V $ and $ x' \in U $. Let $ y' \in F(x') $ (if there is no such $y'$ there is nothing to prove). Since $ y = y'+ \| y - y' \| \, w $ for some $ w \in \B$, denoting $r = \| y- y' \|$, for every $ \epsilon > 0$ we have $y \in (F(x')+r(1+\epsilon) \mathrm{~int} \, \B) \cap V $. From (\ref{linear openness}), there exists $ x \in F^{-1}(y) $ with $ \| x - x' \| \, \leq \, \kappa (1+ \epsilon )r = \, \kappa (1+ \epsilon ) \| y' - y \| $. 
Then 
\begin{equation*}
d(x',F^{-1} (y)) \, \leq  \, \kappa (1+ \epsilon ) \| y' - y \|. 
\end{equation*}
Taking infimum with respect to $ y' \in  F(x') $ on the right and passing to zero with $ \epsilon $ (since the left side does not depend on $ \epsilon $), we obtain that $ F $ is metrically regular at $\bar{x}$ for $\bar{y}$ with constant $ \kappa $.\\
For the converse, let $ x\in U $, $r >0$, and let $y' \in ( F(x) + r \mathrm{~int} \, \B) \cap V $. Then there exists $ y \in F(x)$ such that $ \| y - y' \| \, < \, r $. \\
If $ y=y' $ then $y' \in F(x) \subset F( x+ \kappa r \mathrm{~int} \, \B) $,which yields (\ref{linear openness}) with constant $ \kappa $.\\
Let $ y \not = y' $. From the equivalence of metric regularity of $ F $ and Aubin property of  $ F^{-1} $ (proved in Theorem \ref{MR-AP}), and by using the characterization of the Aubin property given in Proposition (\ref{alternative Aubin}) we obtain 
\begin{equation*}
d( x, F^{-1}(y')) \, \leq \, e \left ( F^{-1} (y) \cap U, F^{-1} (y') \right ) \, \leq \, \kappa \, \| y -y' \| \, < \, \kappa r .
\end{equation*}
Then there exists $ x' \in  F^{-1}(y') $ such that $ \| x -x' \| < \kappa r $. But then
\begin{equation*}
y' \in F(x')  \subset \, F ( x+ \kappa r \mathrm{~int} \,\B),
\end{equation*}  
which again yields (\ref{linear openness}) with constant $ \kappa $.
\end{proof}

\begin{eg}
Consider the function $\smap{f}{2}{}$, defined as:
\begin{figure}[ht]
  \begin{minipage}{.52\textwidth}
    \begin{eqnarray*}
	 	f(x_1, x_2) = \left \{ \begin{matrix}
		x_2 + x_1 ^2 & & ~ x_2 \geq 0, \\
		x_2 - x_1 ^2 & & ~ \mathrm{otherwise}.
		\end{matrix} \right.	
	\end{eqnarray*} 
  \end{minipage}
  \begin{minipage}{.45\textwidth}
    \centering
		\includegraphics[width=5.8cm]{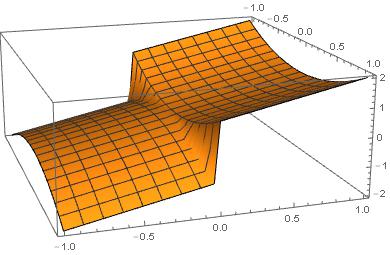}
  \end{minipage}
\caption{A map which is open but not linearly open at at $\bar{x} = (0, 0) $ for $\bar{y} = 0 $}  
\label{fig: Linear openness - eg}
\end{figure}

This function does not have the linear openness property at $\bar{x} = (0, 0) $ for $\bar{y} = 0 $, but it is open at that point. 
Indeed, one can consider the neighborhood $U$ of $\bar{x}$ as $U = \B_{r_1} (0) \times \B_{r_2} (0) $ for $r_1, r_2 > 0$, and observe that 
$$ \norm{f(x) - 0 } \leq r_2 + r_1 ^2 \mathrm{~~~~~~~for~any~~~}  x = (x_1, x_2) \in \B_{r_1} (0) \times \B_{r_2} (0). $$
Thus, $f(U)$ is a neighborhood of $0$, and $f$ is open at $\bar{x} = (0, 0) $ for $\bar{y} = 0 $. \\
To see that it is not linearly open at the reference point, in view of Theorem \ref{LO-MR relation}, we would assume by contradiction that it is metrically regular \at{x}{y}. 
Thus, there should exist a constant $ \kappa > 0 $, together with \NUV such that
$$ d \big(x, f^{-1}(y) \big) \, \leq \, \kappa \, d \big( y, f(x) \big) \mathrm{~~for~any~} x \in U, \, y \in V. $$
Assume that $ y = \bar{y} = 0 $, $ U = \B_r \big( (0, 0) \big) $, with a fixed $ r > 0 $, and the sequence of points $ x_n = ( \frac{1}{n}, 0) $, with $n \in \N$. Then, for $n$ large enough, $x_n \in U $, and we have $ d \big( y, f(x_n) \big) = \frac{1}{n^2} $. On the other hand, $ d \big(x_n, f^{-1}(y) \big) = d (\frac{1}{n}, 0) = \frac{1}{n} $. 
This implies that $ \kappa = n \longrightarrow \infty $, which is a contradiction.
\end{eg}

\defn \textbf{(Strong Metric Regularity)} \label{def. smr} \index{metric regularity ! strong metric regularity} \hfill \\ \citep[p. 194]{implicit}
A mapping $ \mmap{F}{n}{m} $ with $ \rfp{x}{y} \in \gph{F} $ whose inverse $ F^{-1} $ has a Lipschitz continuous single-valued localization around $\bar{y}$ for $\bar{x}$ will be called \emph{strongly metrically regular} at $\bar{x}$ for $\bar{y}$.

\note
The terminology of strong metric regularity offers a way of gaining new perspectives on earlier results by translating them into the language of metric regularity. Indeed, strong metric regularity is just metric regularity plus the existence of a single-valued localization of the inverse.

The following theorem, provides a tool to check strong metric regularity via metric regularity, in examples. Indeed, it suggests that having metric regularity in hand, it only suffices to check the localization of the inverse map for not being multi-valued. In view of Theorem \ref{MR-AP}, it also enlighten the relation between strong metric regularity and the Aubin property.

\begin{prop} [\textbf{Single-Valued Localizations and Metric Regularity}] \label{AP and SMR}\hfill \\ \citep[p. 192]{implicit}
For a mapping $ \mmap{F}{n}{m} $ and a pair $\rfp{x}{y} \in \gph{F}$, the following properties are equivalent:
\begin{enumerate}[topsep=-1ex,itemsep=0ex,partopsep=1ex,parsep=1ex, leftmargin = 7ex]
\item[(a)] $F^{-1} $ has a Lipschitz continuous single-valued localization $s$ around $\bar{y} $ for $\bar{x}$;
\item[(b)] $F$ is metrically regular at $\bar{x}$ for $\bar{y}$ and $F^{-1}$ has a localization \at{y}{x} that is nowhere multivalued.
\end{enumerate}
Indeed, in the circumstances of (b) the localization $s$ in (a) has $ \mathrm{lip~}(s; \bar{y}) = \mathrm{reg~} (F; \bar{x}| \bar{y}) $.
\end{prop}

There are certain situations in which, one can not distinguish between metric regularity and strong metric regularity; one of these situations is described in the following proposition.

\begin{prop} [\textbf{Strong Metric Regularity of Locally Monotone Mappings}] \citep[Theorem 3G.5, p. 195]{implicit} \\
If a mapping $ \mmap{F}{n}{m} $ that is locally monotone \at{x}{y} is metrically regular \at{x}{y}, then it must be strongly metrically regular \at{x}{y}.
\end{prop}

The \lq\lq one-variable \rq\rq \,version of metric regularity is defined as metric sub-regularity. Considering Theorem \ref{MR-AP}, and definition of calmness and isolated calmness, one may think of a relation between metric sub-regularity and the calmness of the inverse map. Such a relation actually exists and will be expressed afterwards.\\

\defn \textbf{(Metric Sub-Regularity)} \label{def. msr} \index{metric regularity ! metric sub-regularity} \hfill \\ \citep[p. 198]{implicit}
A mapping $ \mmap{F}{n}{m} $ is called \emph{metrically sub-regular} \at{x}{y} if $ \rfp{x}{y} \in \gph{F} $ and there exists $ \kappa \geq 0 $ along with \NUV such that
\begin{equation} \label{metric sub-regularity}
d ( x, F^{-1} (\bar{y})) ~ \leq ~ \kappa \, d (\bar{y}, F(x) \cap V )  \mfa x \in U .
\end{equation}
The infimum of all $\kappa$ for which (\ref{metric sub-regularity}) holds is the \emph{modulus} of metric sub-regularity, denoted by $ \mathrm{subreg\,} ( F ; \bar{x} | \bar{y}) $.

The absence of metric sub-regularity is signaled by $ \mathrm{subreg~} ( F ; \bar{x} | \bar{y}) = \infty $.\\
The main difference between metric sub-regularity and metric regularity is that the data input $\bar{y}$ is now fixed and not perturbed to a nearby $y$.\\
It is worth mentioning that, since $ d( \bar{y},F(x)) \leq d( \bar{y},F(x) \cap V) $, it is clear that sub-regularity is a weaker condition than metric regularity, and
$$ \mathrm{subreg~} ( F ; \bar{x} | \bar{y}) \leq \mathrm{reg~} ( F ; \bar{x} | \bar{y}). $$

\begin{thm}  [\textbf{Characterization by Calmness of the Inverse}] \citep[p. 198]{implicit} \label{Characterization by Calmness of the Inverse}\\
For a mapping $ \mmap{F}{n}{m} $, let $ \bar{y} \in F(\bar{x}) $. Then $F$ is metrically sub-regular \at{x}{y} if and only if its
inverse $ \mmap{F^{-1}}{m}{n} $ is calm \at{y}{x}, in which case
$$ \mathrm{clm~} (F^{-1}; \bar{y}|\bar{x}) = \mathrm{subreg~} (F; \bar{x}| \bar{y}). $$
\end{thm}
\begin{proof}
First assume that $F^{-1} $ is calm, i.e., there exist a constant $ \kappa > 0 $ and \NUV such that
\begin{equation} \label{calmness criteria}
F^{-1}(y) \cap U \subset F^{-1}(\bar{y}) + \kappa  \, \| y - \bar{y} \| \B \mfa y \in V.
\end{equation}
Let $ x \in U $. If $ F(x) \cap V = \emptyset $, then the right side of (\ref{metric sub-regularity}) is $\infty$ and we are done. If not, having $ x \in U $ and $ y \in F(x) \cap V $ is the same as having $ x \in F^{-1}(y) \cap U $ and $ y \in V $. For such $x$ and $y$, the inclusion in (\ref{calmness criteria}) requires the ball $ x+ \kappa \, \| y - \bar{y} \| \B $ to have non-empty intersection with $ F^{-1}( \bar{y}) $. Then $ d(x,F^{-1}( \bar{y})) \leq \kappa \|y - \bar{y} \| $. Thus, for any $ x \in U $, we must have
 $ d(x,F^{-1}( \bar{y})) \leq \inf_{y} \, \{ \, \kappa \, \| y - \bar{y} \| \, : \, y \in F(x) \cap V  \} $ which is (\ref{metric sub-regularity}). This shows that (\ref{calmness criteria}) implies (\ref{metric sub-regularity}) and that
 \begin{equation*}
 \inf \{ \kappa ~|~ U, \, V, \, \kappa \mathrm{~satisfying~} \eqref{calmness criteria} \} \geq \inf \{ \kappa ~|~ U, \, V, \, \kappa \mathrm{~satisfying~} \eqref{metric sub-regularity} \} 
 \end{equation*}
the latter being by definition $ \mathrm{subreg~} (F; \bar{x}|\bar{y}) $. \\
For the opposite direction, we have to demonstrate that if $ \mathrm{subreg~} (F; \bar{x}|\bar{y}) < \kappa < \infty $, then (\ref{calmness criteria}) holds for some choice of neighborhoods $U$ and $V$. Consider any $\kappa'$  with $ \mathrm{subreg~} (F; \bar{x}|\bar{y}) < \kappa' < \kappa $. 
For this $ \kappa' $ , there exist $U$ and $V$ such that 
$$ d( x, F^{-1}(\bar{y})) \leq \kappa'  d( \bar{y}, F(x) \cap V) \mfa  x \in U. $$
 Then we have $ d (x, F^{-1}( \bar{y} )) \leq \kappa  \| y - \bar{y} \| $ when $ x \in U $ and $ y \in F(x) \cap V $, or equivalently
  $ y \in V $ and $ x \in F^{-1}(y) \cap U $ .\\
Fix $ y \in V $. If $ y = \bar{y} $ there is nothing to prove; let $ y  \not = \bar{y} $. If $ x \in F^{-1} (y) \cap U $, then
 $ d(x,F^{-1}( \bar{y})) \leq \kappa'  \| y - \bar{y} \| < \kappa \| y - \bar{y} \| $. Then there must be a point $ x'  \in F^{-1}( \bar{y} ) $ having
  $ \| x' - x \| < \kappa \| y - \bar{y} \|$. Hence we have (\ref{calmness criteria}), as required, and the proof is complete.
\end{proof}

Note that, undesirably, the property of metric sub-regularity (and hence calmness) is not stable under smooth perturbation, even for convex multifunctions, as demonstrated by the following example.

\begin{eg} \textbf{(Instability of Metric Sub-Regularity under Perturbation)} \label{MSR-perturbing problem} \\ \citep[Example 2.1, p. 1441]{Gfrerer2011}
Consider the convex multifunction $ \mmap{G}{}{} $, defined as $ G(x) := \R_{+} $ for every $ x \in \R $. At the reference point $ \rfp{x}{y} = (0, 0) $, we would have 
$ G^{-1} (\bar{y}) = \{ x~|~ 0 \in G(x) \} = \R $. Then, $ d \big( x, G^{-1}(\bar{y}) \big) = d \big( \bar{y}, G(x) \big) = 0 $ holds for every $ x \in \R $, showing $G$ is metrically sub-regular at $ (0, 0) $. \\
On the other hand, let $ \smap{\varphi}{}{}_{+} $ denote any convex function differentiable at $0$ fulfilling $ \varphi (0) = \varphi' (0) = 0 $, and $ \varphi(t) > 0 $ for every $ t \neq 0 $. 
For arbitrarily fixed $ \epsilon > 0 $, let $ \widetilde{G} (x) := G(x) + \epsilon \, \varphi (x) $. Then $\widetilde{G}$ is still a convex multifunction, but it is no longer metrically sub-regular at $ (0, 0) $. In fact, 
$ \widetilde{G}^{-1}(\bar{y}) = \{ 0 \} $, and $ d \big( x, \widetilde{G}^{-1}(\bar{y}) \big) = d(x, 0) = |x| $. Also we can compute 
$  d \big( \bar{y}, \widetilde{G} (x) \big) = \epsilon \, \varphi(x) $. Since 
\begin{equation*}
\lim_{ x \rightarrow 0} \dfrac{ |x| }{ \varphi(x) } = \infty,
\end{equation*}
the metric sub-regularity condition \eqref{metric sub-regularity} does not hold for any finite $ \kappa$. \\
Of course we can simply take $ \varphi (x) = x^2 $, but as a more subtle choice for the function $\varphi$ we can take the function $ \varphi  \in C^{\infty} (\R) $ defined by $ \varphi(0) = \varphi' (0) = \varphi'' (0) = 0 $, $ \varphi'' (x) = e^{- x^{-2}} $ for every $ x \neq 0 $, which has the property that all derivatives vanish at $ 0 $, that is, 
$ \varphi^{(i)} (0) = 0 ~ \forall i \in \N $. \\
Then, any \lq\lq derivative-like\rq\rq \, tool offered for checking metric sub-regularity would be desired to have the property that the derivatives of the multifunctions $ G $ and $\widetilde{G}$ are different (so the metric sub-regularity of $G$ and the lack of this property for $\widetilde{G}$ could be detectable). But $ G $ and $ \widetilde{G} $ differ only by a $ C^{\infty} $-function, where all derivatives vanish at $\bar{x}$ and so the usual calculus rules (that is, the derivative of the sum being equal to the sum of derivatives) cannot be valid.
\end{eg}

\defn \textbf{(Strong Metric Sub-Regularity)} \label{def. smsr}\index{metric regularity ! strong metric sub-regularity} \hfill \\ \citep[p. 202]{implicit}
A mapping $\mmap{F}{n}{m}$ is said to be \emph{strongly metrically sub-regular} \at{x}{y} if $\rfp{x}{y} \in \gph{ F} $ and there is a constant $ \kappa \geq 0 $ along
with \NUV such that
\begin{equation} \label{smsr}
\| x - \bar{x} \| \leq \, \kappa \, d(\, \bar{y}, \,F(x) \cap V) ~ \mfa x \in U\footnote{
The strong metric sub-regularity property was first considered without giving it a name, by Rockafellar \citep{rockafellar1989proto}. The name \lq\lq strong metric sub-regularity\rq\rq \,was first used in \citep{dontchev2004regularity} where its equivalence with the isolated calmness was proved.
}.
\end{equation}

Clearly, the infimum of $\kappa$ for which (\ref{smsr}) holds is equal to $ \mathrm{subreg~} (F ; \bar{x} | \bar{y} ) $.

\begin{thm} [\textbf{Characterization by Isolated Calmness of the Inverse Map}] \label{Characterization by Isolated Calmness of the Inverse Map} \hfill 
\citep[p. 202]{implicit} A mapping $ \mmap{F}{n}{m} $ is strongly metrically sub-regular \at{x}{y} if and only if its inverse
$ F^{-1} $ has the isolated calmness property \at{y}{x}.\\
Specifically, for any $ \kappa > \mathrm{subreg~} (F; \bar{x}| \bar{y}) $ there exist \NUV such that 
\begin{equation} \label{relation 3}
 F^{-1}(y) \cap U \subset \, \bar{x} + \kappa \, \|y - \bar{y} \| \, \B \mathrm{~~~when~~} y \in V.
\end{equation}
 
Moreover, the infimum of all $\kappa$ such that the inclusion (\ref{relation 3}) holds for some neighborhoods $U$
and $V$ actually equals $\mathrm{subreg~} (F; \bar{x}| \bar{y})$.
\end{thm}
\begin{proof}
Assume first that $F$ is strongly metrically sub-regular \at{x}{y}. Let $ \kappa > \mathrm{subreg~} (F; \bar{x}| \bar{y}) $. Then there are neighborhoods $U$ for $\bar{x}$ and $V$ for $\bar{y}$ such that (\ref{smsr}) holds with the indicated $\kappa$. 
Consider any $ y \in V$. If $ F^{-1}(y) \cap U = \emptyset $, then (\ref{relation 3}) holds trivially. If not, let $ x \in F^{-1} (y) \cap U$. This entails $ y \in F(x) \cap V $, hence $ d( \, \bar{y},F(x) \cap V) \leq \|y - \bar{y}\| $ and consequently
 $ \|x - \bar{x}\| \leq \kappa \, \|y - \bar{y} \| $ by (\ref{smsr}). Thus, $ x \in \bar{x} + \kappa \| y - \bar{y} \| \,\B $, and we conclude that (\ref{relation 3}) holds. Also, we see that $ \mathrm{subreg~} (F; \bar{x}| \bar{y}) $ is not less than the infimum of all $ \kappa $ such that (\ref{relation 3}) holds for some choice of $U$ and $V$.\\
For the converse, suppose (\ref{relation 3}) holds for some $\kappa$ and neighborhoods $U$ and $V$.
Consider any $ x \in U$. If $ F(x) \cap V = \emptyset $ the right side of (\ref{smsr}) is $\infty$ and there is nothing more to prove. If not, for an arbitrary $ y \in  F(x) \cap V $ we have $x \in F^{-1}(y) \cap U $, and therefore $ x \in \bar{x} + \kappa \, \|y - \bar{y} \| \, \B $ by (\ref{relation 3}), which means
$ \|x - \bar{x}\| \leq \, \kappa \|y - \bar{y}\| $. This being true for all $ y \in F(x) \cap V$, we must have $ \|x - \bar{x}\| \leq \kappa d(\, \bar{y}, F(x) \cap V)$. Thus, (\ref{smsr}) holds, and
in particular we have $ \kappa \geq \mathrm{subreg~} (F; \bar{x}| \bar{y})$. Therefore, the infimum of $\kappa $ in (\ref{relation 3})
 equals $ \mathrm{subreg~} (F; \bar{x}| \bar{y})$.
\end{proof}

The relation of local stability properties in Section \ref{Local Stability Properties} and the metric regularities\footnote{
Referring to all four definitions of metric regularity (shown with MR in the figure, cf. Definition \ref{def. mr}), strong metric regularity (abbreviated as SMR, cf. Definition \ref{def. smr}), metric sub-regularity  (abbreviated as MSR, cf. Definition \ref{def. msr}), and strong metric sub-regularity (abbreviated as SMSR, cf. Definition \ref{def. smsr}). 
} of this section, is summarized in the following chart.
\begin{equation*}
\boxed{
\begin{matrix}
\mathrm{for ~~~~} F: & SMR & \subset & MR & \subset & MSR & \supset & SMSR \\
& & & \cong & & \cong  & & \cong \\
\mathrm{ ~~~for ~} F^{-1}: ~~&  & & Aubin \, Property & \subset & Calmness & \supset & Isolated \, Calmness 
\end{matrix}}
\end{equation*}
\\

\section{A Review on Perturbation Results} \label{A Review on Perturbation Results}
This section will contain results about perturbing a set-valued map that has one of the stability-like properties mentioned in previous sections, with a single-valued function.
The general idea is to provide conditions under which, the sum of two maps still has the stability-like property. In Chapter \ref{Chapter2}, we will give a physical interpretation to this problem.\\
In order to study the perturbation problem, we need to somehow approximate the nonsmooth function and replace it with another smooth function which is close enough to it. This would be the essence of the following three definitions.

\defn \textbf{(First-Order Approximations of Functions)}\index{function ! first-order approximation} \hfill \\ \citep[p. 39]{implicit}
Consider a function $ \smap{f}{n}{m} $ and a point $ \bar{x} \in \mathrm{int~dom} \, f $. A function $ \smap{h}{n}{m} $ with $ \bar{x} \in \mathrm{int~dom} \, h $ is a 
\emph{first-order approximation} to $f$ at $\bar{x}$ if  $ h(\bar{x}) = f (\bar{x}) $ and 
$$ \mathrm{clm\,}(e; \bar{x}) = 0 \mathrm{~~~~~~for~~} e(x) = f (x) - h(x),  $$
which can also be written as $ f(x) = h(x) + \, o (\|x - \bar{x}\| )$. \\
It is a \emph{strict first-order approximation} if the stronger condition holds that
$$ \mathrm{lip\,}(e; \bar{x}) = 0  \mathrm{~~~~~~for~~} e(x) = f (x) - h(x).  $$

\eg
When the function $ \smap{f}{n}{m} $ is differentiable around a point $\bar{x} \in \mathrm{int~dom} \, f $, one can define a first-order approximation as 
$ h(x) = f(\bar{x}) + \nabla f(\bar{x}) (x - \bar{x}) $.\\
First, observe that $h(\bar{x}) = f(\bar{x})$. Then, for every $ \epsilon > 0 $ one can find $ \delta > 0 $ such that
$$  \| f (x) - h(x) \|  = \| f(x) - f(\bar{x}) - \nabla f(\bar{x}) (x - \bar{x}) \| \leq \, \epsilon \, \|x - \bar{x}\| \mathrm{~~~for~every~} x \in \B_{\delta}(\bar{x}). $$

\defn \textbf{(Estimators)} \label{Estimators} \index{function ! estimator} \hfill \\ \citep[p. 41]{implicit}
Consider a function $ \smap{f}{n}{m} $ and a point $ \bar{x} \in \mathrm{int~dom} \, f $. A function $ \smap{h}{n}{m} $ with $ \bar{x} \in \mathrm{int~dom} \, h $ is an
\emph{estimator }of $f$ at $\bar{x}$ with constant $\mu$ if $ h(\bar{x}) = f (\bar{x}) $ and 
$$ \mathrm{clm}~(e; \bar{x}) \leq \mu < \infty \mathrm{~~~~~~for~~} e(x) = f (x) - h(x),  $$
which can also be written as $ \| f(x) - h(x) \| \leq \, \mu \, \|x - \bar{x} \| + \, o (\|x - \bar{x}\| )$. \\
It is a \emph{strict estimator} if the stronger condition holds that
$$ \mathrm{lip}~(e; \bar{x}) \leq \mu < \infty \mathrm{~~~~~~for~~} e(x) = f (x) - h(x).  $$

\note \label{use of linearization as an estimator}
\textbf{(a)} In this terminology, a first-order approximation is simply an estimator with constant $ \mu = 0 $. Through that, any result involving estimators can immediately be specialized to a result about first-order approximations.

\textbf{(b)} Although one can consider the idea of using estimators as a substitute for differentiability in some sense, estimators can be of interest even when differentiability is present. For instance, in the case of a function $ f $ that is strictly differentiable at $\bar{x}$ a strict estimator of $f$ at $\bar{x}$ with constant $\mu$ is furnished by
 $ h(x) = f (\bar{x}) + A (x - \bar{x}) $ for any matrix $A$ with $ \| \nabla f(\bar{x}) - A \| \leq \mu $ \footnote{
Throughout this thesis, we consider matrices as linear operators and thus, $ \| A \| $ is the operator norm of $A$ induced by the Euclidean norm,
\begin{equation*}
\norm{A} = \max_{\norm{x} \, \leq \, 1} \norm{Ax}.
\end{equation*}
 }. We will use this idea in Subsection \ref{Strong Metric (Sub-) Regularity with a Nonsmooth Single-valued Part}.

When dealing with multivariable functions, it is possible to consider one variable as a parameter (which is actually the case in Chapter \ref{Chapter4}), and thus, we need to find the approximation with respect to one variable, while the others are fixed. To do so, we introduce \emph{partial} estimators and approximations in the following way.

\defn \textbf{(Partial First-Order Estimators and Approximations)} \index{function ! partial estimator} \hfill \\ \citep[p. 49]{implicit}
For $ f : \R^d \times \R^n \longrightarrow \R^m $ and a point $ \rfp{p}{x} \in \mathrm{int~dom}\, f$, a function $ \smap{h}{n}{m} $ is said to be an \emph{estimator} of $f$
with respect to $x$ uniformly in $p$ at $ \rfp{p}{x}$ with constant $\mu$ if $  h(\bar{x}) = f (\bar{p},\bar{x}) $ and
\begin{equation*}
\widehat{\mathrm{clm}}_x \, (e; (\bar{p}, \bar{x})) \leq \mu < \infty \mathrm{~~~~~~for~~} e(p,x) = f (p,x) - h(x).
\end{equation*}
It is a \emph{strict estimator} in this sense if the stronger condition holds that
\begin{equation*}
 \widehat{\mathrm{lip}}_x \, (e; (\bar{p}, \bar{x})) \leq \mu < \infty \mathrm{~~~~~~for~~} e(p,x) = f (p,x) - h(x).
\end{equation*}
In which, $ \widehat{\mathrm{clm}}_x  $ is the uniform partial calmness modulus
\footnote{
\textbf{(Notation Review: Partial Calmness)} \index{calmness ! partial calmness}
A function $  f : \R^d \times \R^n \longrightarrow \R^m $ is said to be calm with respect
to $x$ at $\rfp{p}{x} \in \mathrm{dom} \, f $ when the function $\Phi$ with values $ \Phi(x) = f (\bar{p},x) $ is calm at $\bar{x}$.
Such calmness is said to be uniform in $p$ at $\rfp{p}{x}$ when there exists a constant $ \kappa \geq 0 $
and neighborhoods $Q$ of $\bar{p}$ and $U$ of $\bar{x}$ such that actually
$$ \| f (p,x) - f (p, \bar{x}) \| \leq \kappa \|x - \bar{x}\| \mfa (p,x) \in (Q \times U) \cap \mathrm{dom} \, f. $$
Correspondingly, the partial calmness modulus of $f$ with respect to $x$ at $\rfp{p}{x}$ is
denoted as $ \mathrm{clm}_x ( f ; \rfp{p}{x} ) $, while the uniform partial calmness modulus is
\begin{equation*}
\widehat{\mathrm{clm}}_x \, (f; (\bar{p}, \bar{x}))  := \limsup_{\mathclap{\substack{ x \rightarrow \bar{x}, p \rightarrow \bar{p}\\
																							 (p,x)\, \in \, \mathrm{dom }\, f, x \not = \bar{x} }}}
																							   \dfrac{\| f (p,x) - f (p, \bar{x}) \|}{\|x - \bar{x}\|}
\end{equation*}
} and similarly, $ \widehat{\mathrm{lip}}_x $ is the uniform partial Lipschitz modulus.\\
In the case of $ \mu = 0 $, such an estimator is called a \emph{partial first-order approximation}.

The first theorem of this section is about perturbing a strongly metrically regular set-valued map with a Lipschitz function. We would see under suitable conditions for regularity modulus and Lipschitz constant, the new map would remain strongly metrically regular. In order to state and prove this theorem, we use the idea presented in Proposition \ref{AP and SMR}, and go in two steps.\\
First, we would prove the single-valuedness (non-multivaluedness, to be more precise) of the new map in the following proposition. Next, we would prove the metric regularity of the new map in Theorem \ref{mr stability}. Then, we would sum up everything in Theorem \ref{Inverse Function Theorem with Strong Metric Regularity}.

\begin{prop} [\textbf{Stability of Single-Valuedness Under Perturbation}] \label{AP perturbed} \hfill \\ \citep[p. 193]{implicit}
Let $ \nu $ and $ \lambda $ be positive constants such that $ \nu \lambda < 1 $. Consider a mapping $ \mmap{F}{n}{m} $ and a pair $ \rfp{x}{y} \in \gph{F}$, 
such that $ F^{-1} $ has a Lipschitz continuous single-valued localization $s$ around $\bar{y}$ for $\bar{x}$ with $ \mathrm{lip\,} (s; \bar{y}) < \lambda $. 
Consider also a function $ \smap{g}{n}{m} $ with $ \bar{x} \in \mathrm{~int~ dom~} g $ and such that $ \mathrm{lip\,} (g; \bar{x}) < \nu $. \\
Then the mapping $ (g+F)^{-1} $ has a localization around $ g( \bar{x})+ \bar{y} $ for $\bar{x}$ which is nowhere multivalued.
\end{prop}

\begin{thm} [\textbf{Inverse Mapping Theorem with Metric Regularity}] \label{mr stability}\hfill \\ \citep[p. 184]{implicit}
Consider a mapping $ \mmap{F}{n}{m} $, a point $ \rfp{x}{y} \in \gph{F}$, and a function $ \smap{g}{n}{m} $ with $ \bar{x} \in \mathrm{int\,dom~} g $. 
Let $\kappa$ and $\mu $ be nonnegative constants such that
$$ \kappa \mu  < 1,~ \mathrm{reg\,} (F; \bar{x}| \bar{y}) \leq \kappa \mathrm{~~ and~~} \mathrm{lip\,} (g; \bar{x}) \leq \mu. $$
Then
$$ \mathrm{reg\,} (g + F; \bar{x}| g(\bar{x}) + \bar{y}) \leq \dfrac{\kappa}{1 - \kappa \mu } . $$
\end{thm}

\note
Although formally there is no inversion of a mapping in Theorem \ref{mr stability}, if this result is stated equivalently in terms of the Aubin property of the inverse mapping
$F^{-1}$, it fits then into the pattern of the inverse function theorem paradigm \footnote{
We can actually deduce the classical inverse function Theorem (\ref{classical inverse function theorem}) from Theorem \ref{mr stability}. Indeed, let $ \smap{f}{n}{n} $ be a smooth function around $\bar{x}$ and let $ \nabla f ( \bar{x}) $ be non-singular. \\
Then $ F(x) := Df( \bar{x}) (x - \bar{x}) $ is metrically regular everywhere and specially at $\bar{x}$ for $0$.
 Moreover, the function $ g(x) := f (x) - Df( \bar{x})(x - \bar{x}) $ is Lipschitz continuous at $\bar{x}$ with modulus $ \mathrm{lip~}(g; \bar{x}) = 0 $. \\
Thus, from Theorem \ref{mr stability}, we obtain that $ g+F = f$ is metrically regular at $ \bar{x}$ for $f ( \bar{x})$ with modulus
$  \mathrm{reg~} (f; \bar{x}) = \mathrm{reg~} (F; \bar{x}| 0) $. But then $ f $ must be open (cf. Theorem \ref{LO-MR relation}; in fact, $f$ is open at linear rate). So, the inverse map has a localization around $f(\bar{x})$ which is single-valued and Lipschitz continuous. }.\\
The theorem can also be viewed as a result concerning stability of metric regularity under perturbations by functions with small Lipschitz constants.
 
\begin{thm} [\textbf{Inverse Function Theorem with Strong Metric Regularity}] \label{Inverse Function Theorem with Strong Metric Regularity} \hfill \\ \citep[p. 194]{implicit} Let $ \kappa $ and $\mu $ be nonnegative constants such that $ \kappa \mu < 1 $. Consider a mapping $ \mmap{F}{n}{m} $ and any $ \rfp{x}{y} \in \gph{F}$ such that $F$ is strongly metrically regular \at{x}{y} with $ \mathrm{reg~} (F; \bar{x}| \bar{y}) \leq \kappa $ and a function $ \smap{g}{n}{m} $ with $ \bar{x} \in \mathrm{~int~dom~} g $ and $ \mathrm{lip~} (g; \bar{x}) \leq \mu $.\\
Then the mapping $ g+F $ is strongly metrically regular at $\bar{x}$ for $g(\bar{x}) + \bar{y}$. Moreover,
$$ \mathrm{reg~} (g + F; \bar{x}| g(\bar{x}) + \bar{y}) \leq \dfrac{\kappa}{1 - \kappa \mu } . $$
\end{thm}
\begin{proof} 
We try to prove the equivalent statement (that is (b) in Proposition \ref{AP and SMR}), so we have to show
\begin{enumerate}[topsep=-1pt,itemsep=0ex,partopsep=0ex,parsep=0ex]
\item[(1)] $ g + F $ is strongly metrically regular at $\bar{x}$ for $\bar{y} + g(\bar{x})$; and
\item[(2)] $ (g + F) ^{-1} $ has a localization at $\bar{y} + g(\bar{x})$ for $\bar{x}$ that is nowhere multivalued.
\end{enumerate}
 Our hypothesis that $F$ is strongly metrically regular \at{x}{y} implies that a
graphical localization of $F^{-1}$ around $ \rfp{y}{x} $ is single-valued near $\bar{y}$. Further, by fixing
$ \lambda > \kappa $ such that $ \lambda \mu < 1 $ and using Proposition \ref{AP and SMR}, we can get \NUV such that for every $ y\in V $ the set
$  F^{-1}(y) \cap U $ consists of exactly one
point, which we may denote by $s(y)$ and know that the function $ s : y  \mapsto F^{-1}(y) \cap U $ is
Lipschitz continuous on $V$ with Lipschitz constant $\lambda$ . Let $ \mu < \nu < \lambda ^{-1} $ and choose
a neighborhood $ U' \subset U $ of $\bar{x}$ on which $ g $ is Lipschitz continuous with constant $\nu$.
Applying Proposition \ref{AP perturbed}, we obtain that the mapping $ (g+F)^{-1} $ has a localization
around $ g( \bar{x})+ \bar{y} $ for $\bar{x}$ which is nowhere multivalued.\\
On the other hand, we know from Theorem \ref{mr stability} that for such $g $ the mapping $g+F$ is metrically regular at $ g( \bar{x})+ \bar{y} $ for $\bar{x}$. 
Applying Proposition \ref{AP and SMR} once more, we get that $ (g + F ) ^{-1} $ has a Lipschitz single-valued localization $s$ around $ g( \bar{x})+ \bar{y} $ for $\bar{x}$ with modulus
$$ \mathrm{lip~} (s; \bar{y} + g(\bar{x}) =  \mathrm{reg~} (g + F; \bar{x}| g(\bar{x}) + \bar{y}) < \dfrac{\lambda}{1- \lambda \nu} $$
As $\lambda $ and $\nu$ could be arbitrarily close to $\kappa$ and $\mu$, respectively, the modulus criteria is satisfied.
\end{proof}

Now, we will consider the situation that working with the function $f$ in the sum $ f + F$ is not easy for some reasons. So, we would like to use an approximation function $h$, instead. It would be important to know whether the regularity properties obtained for $ h + F $, contain any information about the regularity of $ f + F $ or not.\\
The following proposition is needed to prove the theorem afterwards, which will be one of the important tools in order to give a meaning to the usage and conditions of the auxiliary map $G_t$ in Chapter \ref{Chapter4}. First, we will prove the parametric version of this result, and then express the \lq\lq inverse function theorem\rq\rq \,version in Theorem \ref{Inverse Function Theorem for Set-Valued Mappings}.

\begin{prop} [\textbf{Contraction Mapping Principle for Composition}] \label{contraction mapping principle} \hfill \\ \citep[p. 86]{implicit} 
Consider a function $ \varphi : \R^d \times \R^n \longrightarrow \R^m $ and a point $ \rfp{p}{x} \in \mathrm{int~dom~} \varphi $ and let the scalars
$ \nu \geq 0, b \geq 0, a > 0 $, and the set $ Q \subset \R^d $ be such that $ \bar{p} \in Q $ and
\begin{equation}\label{thm condition3}
\left\{\begin{matrix}
\| \varphi (p,x') - \varphi (p,x) \| \leq \nu \| x - x' \| & \mfa x',x \in \B_a (\bar{x}) \mathrm{~and~} p \in Q,\\ 
\| \varphi (p, \bar{x}) - \varphi \rfp{p}{x} \| \leq b ~~~~~~~~~~~ & \mfa p \in Q ~~~~~~~~~~~~~~~~~~~~~~~~
\end{matrix}\right.
\end{equation}
Consider also a set-valued mapping $  \mmap{M}{m}{n} $ with $ \rfp{y}{x} \in \gph{M} $ in which $ \bar{y} :=\varphi \rfp{p}{x} $, such that for each
$ y \in \B_{\nu a+b} (\bar{y}) $ the set $ M(y) \cap \B_a(\bar{x}) $ consists of exactly one point, denoted by $r(y)$, and suppose that the function
\begin{equation}
r : y \mapsto M(y) \cap \B_a (\bar{x}) \mathrm{~~for~} y \in \B_{ \nu a+b }( \bar{y})
\end{equation}
is Lipschitz continuous on $\B_{ \nu a+b }( \bar{y})$ with a Lipschitz constant $\lambda$ . In addition, suppose that the following relations hold
\begin{enumerate}[topsep=-1pt,itemsep=-1ex,partopsep=0ex,parsep=0ex]
\item[(a)] $\lambda \nu < 1$;
\item[(b)] $ \lambda \nu a + \lambda b \leq a $. 
\end{enumerate}
Then for each $ p \in Q $ the set $ \{ x \in \B_a(\bar{x}) \, | \, x \in M(\varphi (p,x)) \} $ consists of exactly one point, and the associated function
\begin{equation}
s : p \mapsto \{ x \, | \, x = M(\varphi (p,x)) \cap  \B_a(\bar{x}) \} \mathrm{~~for~} p \in Q
\end{equation}
satisfies
\begin{equation}\label{estimate}
\| \, s(p') - s(p) \, \| \leq \dfrac{\lambda}{1 - \lambda \nu} ~ \|\, \varphi (p', s(p)) - \varphi (p, s(p)) \, \| \mfa p', p \in Q.
\end{equation}
\end{prop}

\begin{thm} [\textbf{Robinson Theorem Extended Beyond Differentiability}]\label{robinson meaning} \hfill \\ \citep[p. 86]{implicit}
For $ f : \R^d \times \R^n \longrightarrow \R^m $ and $ \mmap{F}{n}{m} $, consider a generalized equation $ f(p,x) + F(x) \ni 0 $ with the solution mapping defined as 
$ S: t \mapsto \{ x \, | \, f(p,x) + F(x) \ni 0   \}$. Let $\bar{p}$ and $\bar{x}$ be such that $ \bar{x} \in S(\bar{p}) $. Assume that:\\
(a) $ f (\cdot, \bar{x}) $ is continuous at $\bar{p}$, and $h$ is a strict estimator of $f$ with respect to $x$ uniformly in $p$ at $\rfp{p}{x}$ with constant $\mu$;\\
(b) the inverse $G^{-1} $ of the mapping $ G = h + F $ , for which $ G(\bar{x}) \ni 0 $, has a Lipschitz continuous single-valued localization 
$\sigma$ around $0$ for $\bar{x}$ with $ \mathrm{lip~} (\sigma ; 0) \leq \kappa $ for a constant $\kappa$ such that $ \kappa \mu < 1 $. \\
Then $S$ has a single-valued localization $s$ around $\bar{p}$ for $\bar{x}$ which is continuous at $\bar{p}$, and moreover for every $ \epsilon > 0 $ there is a neighborhood $Q$ of $\bar{p}$ such that
\begin{equation}
\| \, s(p') - s(p) \, \| \leq \dfrac{\kappa + \epsilon}{1 - \kappa \mu} ~ \| \, f (p', s(p)) - f (p, s(p)) \, \| \mfa p', p \in Q.
\end{equation}
\end{thm}

\begin{proof}
For an arbitrary $ \epsilon > 0 $, choose any $ \lambda > \mathrm{~lip} (\sigma ;0) $ and $ \nu >\mu $ such that $ \lambda \nu < 1 $ and
\begin{equation*}
 \dfrac{\lambda}{1 - \lambda \nu} \leq \dfrac{\kappa + \epsilon}{1 - \kappa \mu} ,
\end{equation*}
as is possible under the assumption that $ \kappa \mu < 1$ . Let $a, b$ and $c$ be positive numbers such that
\begin{equation} \label{thm property1}
\begin{matrix}
 ~~ \| \sigma (y) - \sigma (y') \| \leq \lambda \, \| y - y' \| ~~~~~ ~ &  \mathrm{~for~} y,y' \in \B_{\nu a+b } (0),  \\ 
 ~~ \| e(p,x') - e(p,x) \| \leq \nu \, \| x - x' \| ~ & \mathrm{~for~} x,x' \in \B_a( \bar{x}), \, p \in  \B_c ( \bar{p}),  
\end{matrix}
\end{equation}
where $ e(p,x) = f (p,x) - h(x) $, the first inequality is guaranteed by Lipschitz continuity of  $\sigma$, and the second inequality comes from the definition of estimator.\\
Continuity of $ f(\cdot,\bar{x}) $ at $ \bar{p} $ implies that
\begin{equation} \label{thm property2}
\begin{matrix}
 ~~ \| f (p,\bar{x}) - f (\bar{p}, \bar{x}) \| \leq b  ~~~~~~~~~~~~~~ &~~~ \mathrm{~for~} p \in  \B_c ( \bar{p}).
\end{matrix}
\end{equation}
Take $b$ smaller if necessary so that $ b \lambda < a (1 - \lambda \nu ) $, and accordingly adjust $c$ to ensure having the last inequality. \\
Now apply Proposition \ref{contraction mapping principle} with $ \varphi = - e, Q = \B_c (\bar{p}), M = G^{-1} = (h + F) ^{-1}, \bar{y} = 0$, and $ r = \sigma $. 
The first condition in (\ref{thm condition3}) would be the second statement in (\ref{thm property1}). The second condition in (\ref{thm condition3}) will be
$$ \| e (p, \bar{x}) - e \rfp{p}{x} \| = \| f(p, \bar{x}) - h(\bar{x}) - f \rfp{p}{x} + h(\bar{x})  \|  \leq b $$
which holds because of (\ref{thm property2}). \\
Extra conditions of (a) and (b) of the proposition hold true trivially. Thus, one obtains that for any $ p \in \B_c (\bar{p}) $, the map
 $$ p \mapsto s^* (p) := \{ x \, | \, x = [ (h + F) ^{-1} ( - e(p,x) ) ] \cap \B_a (\bar{x}) \} $$ 
 is single-valued with the estimate mentioned in (\ref{estimate}). The only point to check out is that whether $s^*$ is related to the solution mapping $S$ or not, that is, if we have $ s^*(p) \in S(p) $ for any $ p \in \B_c (\bar{p}) $.\\
Let $ x = s^*(p) $. Then $ x =   (h + F) ^{-1} ( - e(p,x) )  \cap \B_a (\bar{x}) $.
\begin{equation*}
\begin{split}
~~ x \in (h + F) ^{-1} ( - e(p,x) ) ~ & \Leftrightarrow ~  - e(p,x) \in (h + F) (x) ~   \Leftrightarrow ~ h(x) - f(p,x) \in h(x) + F(x) \\
&   \Leftrightarrow ~ 0 \in f(p,x) + F(x) ~ \Leftrightarrow ~ x \in S(p)
\end{split}
\end{equation*}
So, $s^*$ is a single-valued localization for $S$ around $\bar{p}$ for $\bar{x}$. The estimate in (\ref{estimate}) also implies the continuity of $s^*$ at $\bar{p}$, since 
$ f(\cdot,\bar{x}) $ is assumed to be continuous at $\bar{p}$ by assumption (a). 
\end{proof}

The inverse function version of Theorem \ref{robinson meaning} has the following simpler form. Let us note that in view of Definition \ref{def. smr}, this theorem is nothing more than a translation of Theorem \ref{Inverse Function Theorem with Strong Metric Regularity} from the strong metric regularity term to the Lipschitz continuity (for the localized inverse map).

\begin{thm} [\textbf{Inverse Function Theorem for Set-Valued Mappings}] \label{Inverse Function Theorem for Set-Valued Mappings} \hfill \\ \citep[p. 89]{implicit} 
Consider a mapping $ \mmap{G}{n}{n} $ with $ G( \bar{x} ) \ni \bar{y} $ and suppose that $G^{-1}$ has a Lipschitz continuous single-valued localization $\sigma$ around $\bar{y}$ for $\bar{x}$ with lip $(\sigma ; \bar{y}) \leq \kappa $ for a constant $\kappa$. 
Let $ \smap{g}{n}{n} $ be Lipschitz continuous around $\bar{x}$ with Lipschitz
constant $\mu$ such that $ \kappa \mu < 1 $. \\
Then the mapping $ (g + G)^{-1} $ has a Lipschitz
continuous single-valued localization around $ \bar{y} + g( \bar{x}) $ for $\bar{x}$ with Lipschitz constant $ \kappa / ( 1 - \kappa \mu ) $.
\end{thm}

The last Theorem of this section would discuss the perturbation effect on metric regularity and strong metric regularity in terms of the graph of the new map. It would be an important tool in our study and we may refer to it several times in Chapter \ref{Chapter4}.\\
We need the following statement to prove the theorem afterwards. \\

\begin{thm} [\textbf{Contraction Mapping Principle for Set-Valued Mappings}] \label{Contraction Mapping Principle for Set-Valued Mappings}\hfill \\ \citep[p. 313]{implicit}
Let $X$ be a complete metric space with metric $\rho$ , and consider a set-valued mapping $ \Phi : X \rightrightarrows X $ and a point $ \bar{x} \in X$. Suppose that there exist scalars $ a > 0 $ and $ \lambda \in (0,1) $ such
that the set $ \gph{\Phi} \cap \big( \B_a (\bar{x}) \times \B_a(\bar{x}) \big) $ is closed and
\begin{align*}
& \mathrm{(a)~} d (\bar{x}, \Phi(\bar{x})) <  a ( 1 - \lambda ); \\
& \mathrm{(b)~} e \big( \Phi(u) \cap \B_a (\bar{x}), \Phi(v) \big) \, \leq \, \lambda \, \rho(u,v) \mfa u,v \in \B_a (\bar{x}).~~~~~~~~~~~~~~~~~~~~~~
\end{align*}
Then $\Phi$ has a fixed point in $\B_a (\bar{x})$ ; that is, there exists $ x \in \B_a (\bar{x}) $ such that $ x \in \Phi(x) $.\\
\end{thm}

\begin{thm} [\textbf{Perturbed [Strong] Metric Regularity}] \citep[p. 325]{implicit} \label{theorem 5G.3}\\
Let $ X,Y $ be Banach spaces. Consider a mapping $ F : X  \rightrightarrows Y $ and a point $ \rfp{x}{y} \in \gph{F} $ at which $F$
is metrically regular, \big[that is, there exist positive constants $ a, b $, and a nonnegative $ \kappa $ such that
the set $ \gph{F} \cap (\B_a (\bar{x}) \times \B_b (\bar{y})) $ is closed and 
$$ d(x,F^{-1}(y)) \, \leq \, \kappa d(y,F(x)) ~~~~~\mfa (x,y) \in  \B_a (\bar{x}) \times \B_b (\bar{y}). \big] $$
Let $ \mu > 0 $ be such that $ \kappa \mu < 1 $ and let $ \kappa' > \kappa / (1 - \kappa \mu) $. \\
Then for every positive $\alpha$ and $\beta$ such that
\begin{equation} \label{5G.3-1}
\alpha \leq a / 2, ~~ 2 \mu \alpha + 2 \beta \leq b, ~~ \mathrm{~and~~~} 2 \kappa' \beta \leq \alpha
\end{equation}
and for every function $ g : X  \longrightarrow Y $ satisfying
\begin{equation}\label{thm condition1}
\norm{g(\bar{x})} \leq \beta
\end{equation}
and
\begin{equation}\label{thm condition2}
\norm{g(x) - g(x') } \leq \mu \norm{x - x'} \mathrm{~~~for~every~~} x, x' \in \B_{2 \alpha} (\bar{x}),
\end{equation}
the mapping $g+F$ has the following property: \\
for every $y, y' \in \B_{\beta} (\bar{y}) $ and every $ x \in (g+F)^{-1}(y) \cap \B_{\alpha} (\bar{x}) $ there exists $x' \in (g+F)^{-1}(y') $ such that
\begin{equation} \label{5G.3-4}
 \norm{x - x'} \leq \kappa' \norm{y - y'}. 
\end{equation}
In addition, if the mapping $F$ is strongly metrically regular \at{x}{y}; 
\big[that is, the
mapping $ y \mapsto F^{-1}(y) \cap \B_a (\bar{x}) $ is single-valued and Lipschitz continuous on $ \B_b (\bar{y}) $ with a Lipschitz constant $\kappa$\big],
then for $ \mu, \kappa', \alpha $, and $\beta$ as above and any function
$g$ satisfying (\ref{thm condition1}) and (\ref{thm condition2}), the mapping $ y \mapsto (g+F)^{-1}(y) \cap \B_{\alpha} (\bar{x}) $ is a Lipschitz
continuous function on $\B_{\beta} (\bar{y})$ with a Lipschitz constant $\kappa'$. 
\end{thm}
\begin{proof} We consider two different cases.\\
\textbf{Case 1. }\emph{with metric regularity};\\
Choose $\mu$ and $\kappa'$ as required and then $\alpha$ and $\beta$ to satisfy (\ref{5G.3-1}). For any $x \in \B_{2 \alpha} (\bar{x}) $ and $ y \in  B_{\beta} (\bar{y}) $, using (\ref{thm condition1}), (\ref{thm condition2}) and the triangle inequality, we obtain
\begin{equation} \label{5G.3-2}
\begin{split}
\norm{ - g(x) + y - \bar{y} } & \, \leq \, \norm{g(\bar{x}) } + \norm{ g(\bar{x}) - g(x) } + \norm{y - \bar{y} } \\
& \, \leq \, \beta + \mu \norm{x - \bar{x} } + \beta \, \leq \, 2 \beta + 2 \mu \alpha  \, \leq \, b,
\end{split}
\end{equation} 
where the last inequality follows from the second inequality in (\ref{5G.3-1}). Fix $ y' \in \B_{\beta} (\bar{y}) $ and consider the mapping
$$  \B_{\alpha} (\bar{x}) \ni x \longmapsto \Phi_{y'} (x) := F^{-1} \big( - g(x) + y' \big). $$
Let $ y \in \B_{\beta} (\bar{y}), y \neq y' $ and let $ x \in (g+F)^{-1} (y) \cap \B_{\alpha} (\bar{x}) $. We will apply Theorem (\ref{Contraction Mapping Principle for Set-Valued Mappings}) with the complete metric space $ X $ identified with the closed ball $ \B_{\alpha} (\bar{x}) $ to show that there is a fixed point $ x' \in \Phi_{y'} (x') $ in the closed ball centered at $x$ with radius
\begin{equation} \label{5G.3-3}
r : = \kappa' \norm{ y - y' }
\end{equation}
From the third inequality in (\ref{5G.3-1}), we obtain $ r \, \leq \, \kappa' (2 \beta) \, \leq \, \alpha$. Hence, from the first inequality in (\ref{5G.3-1}) we get $ \B_r (x) \subset \B_a(\bar{x}) $. \\
Let $ (x_n, z_n ) \in \gph{\Phi}_y' \cap \big( \B_r (x) \times \B_r (x) \big) $ and $ (x_n, z_n) \to (\tilde{x}, \tilde{z}) $. \\
From (\ref{5G.3-2}), 
$ \norm{ - g(x_n) + y' - \bar{y} } \leq b $;
also note that $ \norm{z_n - \bar{x} } \leq a $. 
Using closedness of $ \gph{F} \cap (\B_a (\bar{x}) \times \B_b (\bar{y})) $ (by metric regularity assumption) and passing to the limit we obtain that 
$ (\tilde{x}, \tilde{z}) \in \gph{\Phi}_y' \cap \big( \B_r (x) \times \B_r (x) \big) $, hence this set is closed.\\
Since $ y \in g(x)+F(x) $ and $ (x,y) $ satisfies (\ref{5G.3-2}), from the assumed metric regularity of $F$ we have
\begin{align*}
d \big(x, \Phi_{y'} (x) \big) & = d \big( x, F^{-1} ( -g(x)+y') \big) \, \leq \, \kappa \, d \big( -g(x)+y', F(x)) \\
 & = \kappa \, d \big( y', g(x)+F(x) \big) \, \leq \, \kappa \, \norm{ y - y' } \\
 & < \kappa' \norm{y - y' } (1 - \kappa \mu ) = r (1 - \kappa \mu ).
\end{align*}
For any $ u, v \in \B_r (x) $, using (\ref{thm condition2}), we have
\begin{equation*}
\begin{split}
e \Big( \Phi_{y'} (u) \cap \B_r(x), \, \Phi_{y'} (v) \Big) & \, \leq \, e \Big( F^{-1} (- g(u) +y') \cap \B_a (\bar{x}), \, F^{-1} (-g(v)+y') \Big) \\
& \, \leq \, \kappa \, \norm{ g(u) - g(v) } , \leq \, \kappa \mu \norm{ u- v}.
\end{split}
\end{equation*}
Applying Theorem \ref{Contraction Mapping Principle for Set-Valued Mappings} to the mapping $ \Phi_{y'} $, with $\bar{x}$ identified with $x$ and constants $ a = r $ and 
$ \lambda  = \kappa \mu $, we obtain the existence of a fixed point $ x' \in \Phi_{y'} (x') $, which is equivalent to $ x' \in (g+F)^{-1} (y') $, within distance $r$ given by (\ref{5G.3-3}) from $x$. This proves (\ref{5G.3-4}). 

\textbf{Case 2. }\emph{with strong metric regularity};\\
For the second part of the theorem, suppose that $ y \mapsto s(y) := F^{-1} (y) \cap \B_a( \bar{x}) $ is a Lipschitz continuous function on $ \B_b (\bar{y}) $ with a Lipschitz constant $\kappa$. Choose $ \mu, \kappa', \alpha $ and $\beta$ as in the statement and let $g$ satisfy (\ref{thm condition1}) and (\ref{thm condition2}). \\
For any $ y \in \B_{\beta} (\bar{y}) $, since $ \bar{x} \in (g+F)^{-1}( \bar{y} + g(\bar{x})) \cap \B_{\alpha} (\bar{x}) $, from (\ref{5G.3-4}) we obtain that there exists
$ x \in (g+F)^{-1}(y) $ such that
$$ \norm{ x - \bar{x} } \, \leq \, \kappa' \, \norm{ y - \bar{y} - g(\bar{x}) }. $$
Since $ \norm{ y - \bar{y} - g(\bar{x}) } \leq 2 \beta $, by (\ref{5G.3-1}) we get $ \norm{ x - \bar{x} } \leq \alpha $, i.e., 
$ (g+F)^{-1}(y) \cap \B_{\alpha} (\bar{x}) \neq \emptyset $. Hence the domain of the mapping $ (g+F)^{-1}(y) \cap \B_{\alpha} (\bar{x}) $ contains 
$ \B_{\beta} (\bar{y}) $. \\
If $ x \in (g+F)^{-1} (y) \cap \B_{\alpha} (\bar{x}) $, then 
$$ x \in F^{-1} ( y - g(x) ) \cap \B_{\alpha} (\bar{x}) \subset F^{-1} (y - g(x)) \cap \B_{a} (\bar{x}) = s(y - g(x)) $$
since $ y - g(x) \in \B_b(\bar{y}) $ according to (\ref{5G.3-2}). Hence,
\begin{equation} \label{5G.3-5}
F^{-1} ( y - g(x) ) \cap \B_{\alpha} (\bar{x}) = s(y - g(x)) = x.
\end{equation}
Let $ y, y' \in \B_{\beta} (\bar{y}) $, and define a map $ \sigma $ on $\B_{\beta} (\bar{y}) $ with 
$$ y \mapsto \sigma (y) := (g+F)^{-1}(y) \cap \B_{\alpha} (\bar{x}). $$
Utilizing the equality $\sigma (y) = s( -g( \sigma(y)) + y ) $ which comes from (\ref{5G.3-5}), we have
\begin{equation} \label{5G.3-6}
\begin{split}
\norm{ \sigma(y) - \sigma (y') } & = \norm{ s( -g( \sigma(y)) + y ) - s( -g( \sigma(y')) + y' ) }\\
\, & \leq \, \kappa \, \norm{ g( \sigma(y)) - g( \sigma(y')) } + \kappa \norm{ y - y' } \\
\, & \leq \, \kappa \mu \norm{ \sigma (y) - \sigma (y') } +  \kappa \norm{ y - y' }.
\end{split}
\end{equation}
If $ y =y' $, taking into account that $ \kappa \mu <1  $ we obtain that $ \sigma(y)$ must be equal to $ \sigma(y')$. 
Hence, the mapping $y \mapsto \sigma (y) := (g+F)^{-1}(y) \cap \B_{\alpha} (\bar{x})$ is single-valued. 
From (\ref{5G.3-6}) and (\ref{5G.3-4}) this function satisfies
$$ \norm{ \sigma (y) - \sigma (y') } \, \leq \, \kappa' \norm{ y - y'}. $$
The proof is complete.
\end{proof}

\note
If $g(\bar{x}) \not =0 $ then $\rfp{x}{y} $ may be not in the graph of $g+F$ and we can not claim that $g+F$ is (strongly) metrically regular \at{x}{y}. This
could be handled easily by choosing a new function $\widetilde{g}$ with $ \widetilde{g}(x) = g(x) - g(\bar{x})$.

%% file: Chapters/Chapter02.tex
\chapter{An Introduction to the Electrical Problem} 

\label{Chapter2} 
\lhead{Chapter 2. \emph{An Introduction to the Electrical Problem}} 

\epigraph{\emph{Everything should be made as simple as possible, but not simpler.}}{\textit{Albert Einstein, a.e. \\ (almost every-word!)}}

In this chapter, we will provide a short review on the theory of electrical circuits, and a brief introduction to some electronic components (in Section \ref{Introduction to the Theory of Electrical Circuits}). Then, in 
Section \ref{Formulating the problem} a general form of generalized equations that would be considered in the rest of the thesis as an outcome of modelling process would be obtained. \\
Also, it provides physical explanation for the importance of studying the small perturbations of source signals. Another important role of this section is to provide a meaning for the stability-like properties introduced in Chapter \ref{Chapter1}: in the static case, the question about sensitivity of the circuit to small perturbations of the power source, is translated to the question whether the solution mapping of the obtained generalized equation has some stability-like properties or not.\\
This chapter ends with a review on the examples that we will discuss in depth in Chapter \ref{Chapter3}. Here, we just provide the mathematical model and general discription of solutions. 

\section{Introduction to the Theory of Electrical Circuits} \label{Introduction to the Theory of Electrical Circuits}

In this section, we will introduce the most basic concepts in the circuit theory in order to provide the appropriate language for formulation of the physical problem, which would be presented in the next section. \\
Mostly, the information about the electrical components and circuit theory are drawn from \citep{desoer1969basic}, while for the electronic details we refer to \citep{Sedra}, and \citep{kuphaldt2009}. Although the way we present these physical concepts is a bit different from their approach and purpose, it would be a relaxing information to know where to refer in case of any ambiguity in the physical meaning, or for the sake of curiosity.  \\
We hope that this short note legitimizes why we do not give the detail referencing in this chapter, unlike the rest of this thesis. 

\begin{figure}[ht]
  \centering
\includegraphics[width=8.8cm]{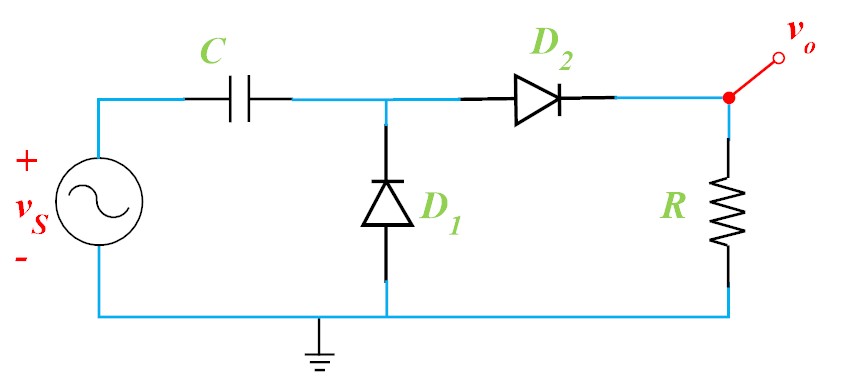}
	\caption{A simple electrical circuit}
	\label{fig: circuit}
\end{figure}

\paragraph*{Circuits}
Without going into technical details, one can think of an electrical circuit (like the one shown in Figure \ref{fig: circuit}) as a bunch of \lq\lq components\rq\rq \,connected together with wires. 
The special forms of connecting these components and the possible choice of different components in each combination, is due to the specific duty that each circuit is going to serve\footnote{
In other words, one has a goal and then tries to design a circuit, that is a combination of components and how to connect them, that provides the desired aim. It is referred to as \emph{circuit design}, in the literature. 
The other approach would be that given a particular circuit, one tries to use electrical rules and some tools to understand the goal of of that circuit. It is referred to as \emph{circuit analysis}, or \emph{synthesis}, in the literature.\\
}.

The word \emph{component}, in this context, refers to different materials which exhibit a particular electrical behaviour under an electromagnetic force. 
The language of circuit theory has two variables: 
\begin{itemize}[topsep=-1ex,itemsep=-1ex,partopsep=1ex,parsep=1ex]
\item current, shown with $I$ or $i$, and measured in amperes;
\item voltage, indicated with $V$ or $v$, and measured in volts\footnote{
The difference between small and capital letter symbols is that whether they indicate the quantities changing with time (known as \emph{alternative current} in the literature, with the abbreviation AC) or constant quantities with respect to time (known as \emph{direct current} in the literature, with the abbreviation DC), respectively.
}. 
\end{itemize}

Thus, for describing a component, we would look at the current passing through it, and the electric potential difference between its terminals, that is voltage over it. The behaviour of the component under different voltages dropped over its terminal, or equivalently, under various currents passing through it is usually described with a graph in the $ i - v $ plane, and referred to as the $ i - v $ \emph{characteristic} of the component.\\ 
In order to avoid confusion, a subscript indicates which component in the circuit we are talking about (for example, $v_s$ in the above circuit indicates the voltage of Source).\\
In circuit theory, we do not know nor are interested in the physical phenomenon causing such a behaviour, and consider a component as a black box that could be described with a map, that is, a relation between current and voltage. 
\begin{figure}[ht]
 \centering
\includegraphics[width=4.8cm]{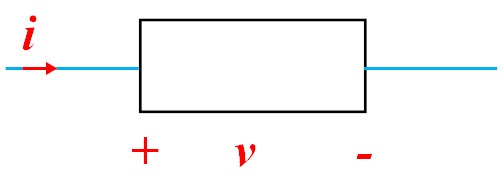}
	\caption{General black box viewpoint of components}
	\label{}
\end{figure}

In the sequel, we will use \lq\lq electrical\rq\rq \,and \lq\lq electronic\rq\rq \,quite very often in referring to circuits, components, or different phrases. It is better to say few words about it now. An \emph{electronic component} is usually made of semiconductors, which are materials that exhibit electrical behaviour somewhere between that of insulators and that of conductors\footnote{
Conductors present very low resistance to the flow of current, whereas insulators conduct very little current even when a large potential difference is applied. A semiconductor exhibits intermediate conductivity because it has more available charge carriers than an insulator but fewer than a conductor. \\
}. 
Examples of semiconductors are silicon and germanium, and Diodes and BJT Transistors are examples of electronic components.\\
An \emph{electronic circuit}, is simply a circuit with at least one electronic component. So, if one wants to be with mathematical precision in these engineering concepts, he/she can say the set of electronic circuits is a subset of electrical circuits.

\paragraph*{Components} We will now review those fundamental components of circuits that would be used in different examples of this thesis.
\begin{enumerate}[topsep=-1ex,itemsep=0ex,partopsep=1ex,parsep=1ex, leftmargin= 3ex]
\item[] \textbf{Resistor} \index{circuit components ! resistor}\\
In the following figure you can see the schematic,  and $i - v$ characteristic of the most common two-terminal electrical component, which has a linear map (i. e. $v_R = R \, i_R$). 

\begin{figure}[ht]
  \centering
    \includegraphics[width=9.8cm]{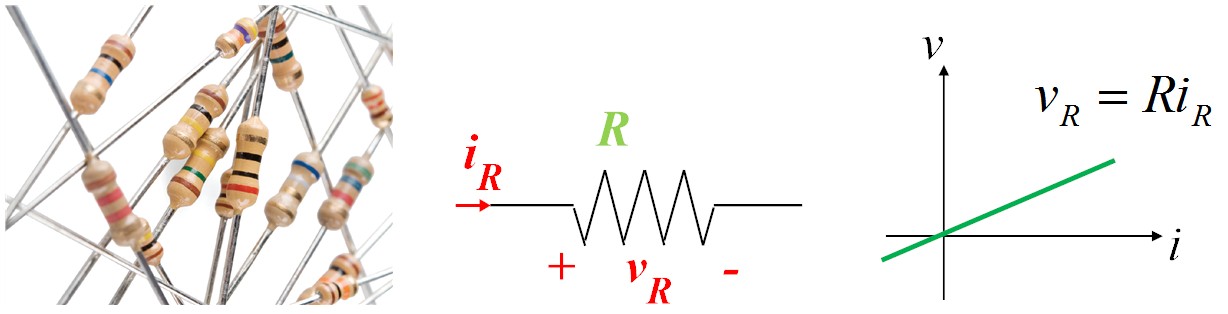}
    \caption{Circuit components: Resistor}
    \label{fig: resistor}
\end{figure}

In electronic circuits, resistors are used to reduce current flow, adjust signal levels, to divide voltages, to bias active elements\footnote{
In electronics, \emph{biasing} means establishing predetermined voltages or currents at various points of an electronic circuit in order to provide proper operating conditions in electronic components. Many electronic devices like transistors whose function is processing time-varying (AC) signals also require a steady (DC) current or voltage to operate correctly, a \emph{bias}.
}, and to model a typical consumer, among other uses.

\item[] \textbf{Capacitor} \index{circuit components ! capacitor}\\
A capacitor is a passive two-terminal electrical component that stores electrical energy in an electric field. \emph{Capacitance} is defined as the ratio of the electric charge $Q$ on each conductor to the potential difference $V$ between them. Unlike a resistor, an ideal capacitor does not dissipate energy.\\

\begin{figure}[ht]
  \centering
\includegraphics[width=11.8cm]{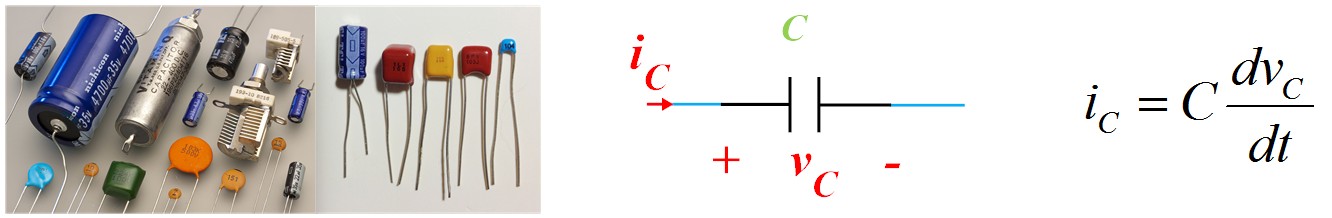}
	\caption{Circuit components: Capacitor}
	\label{fig: capacitor}
\end{figure}

Capacitors are widely used in electronic circuits for blocking direct current while allowing alternating current to pass. In electric power transmission systems, they stabilize voltage and power flow. The property of energy storage in capacitors was exploited as dynamic memory in early digital computers.

\item[] \textbf{Inductor} \index{circuit components ! inductor}\\
Inductors are components that store electrical energy in a magnetic field when electric current is flowing through it. An inductor typically consists of an electric conductor, such as a wire, that is wound into a coil, and it is characterized by its inductance, which is the ratio of the voltage to the rate of change of current. An \lq\lq ideal inductor\rq\rq \,has inductance, but no resistance or capacitance, and does not dissipate or radiate energy.

\begin{figure}[ht]
  \centering
\includegraphics[width=11.8cm]{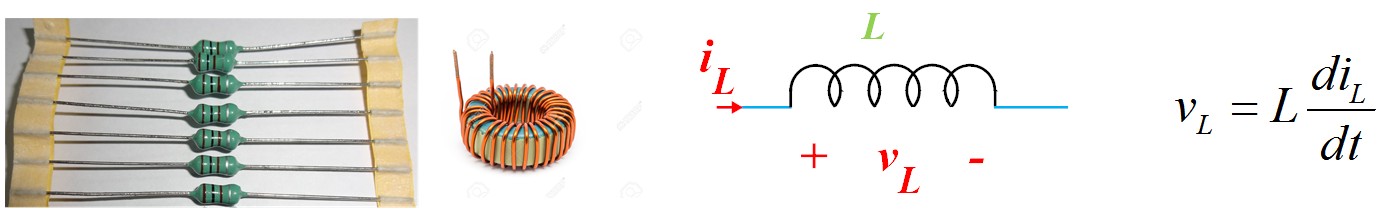}
	\caption{Circuit components: Inductor}
	\label{fig: inductor}
\end{figure}

Inductors are widely used in alternating current (AC) electronic equipment, particularly in radio equipment. They are used to block AC while allowing DC to pass; inductors designed for this purpose are called \emph{chokes}. They are also used in electronic filters to separate signals of different frequencies, and in combination with capacitors to make tuned circuits, used to tune radio and TV receivers.

\item[] \textbf{Voltage Source} \index{circuit components ! voltage source}\\
A voltage source is a two terminal device which can maintain a fixed voltage drop across its terminals. An ideal voltage source can maintain the fixed voltage independent of the load resistance or the output current. However, a real-world voltage source cannot supply unlimited current. 

\begin{figure}[ht]
  \centering
\includegraphics[width=12.8cm]{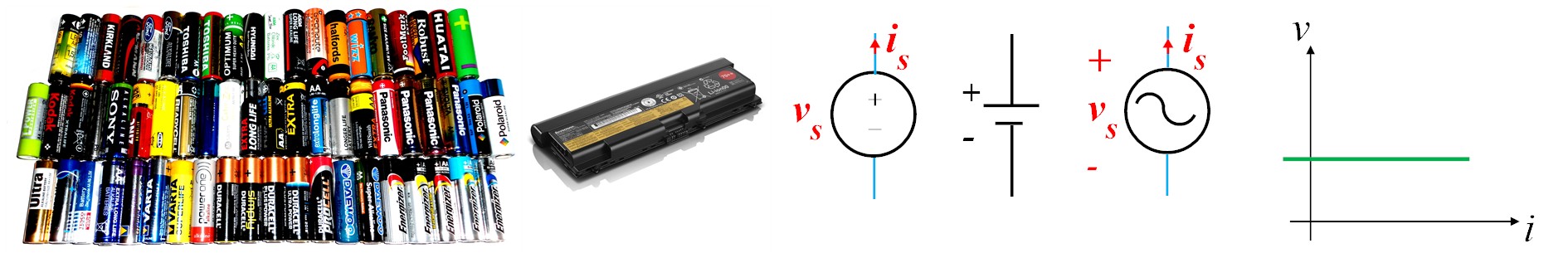}
	\caption{Circuit components: Voltage Source}
	\label{fig: voltage-source}
\end{figure}

Figure \ref{fig: voltage-source} shows the different schematics for DC and AC voltage sources. 
Real-world sources of electrical energy, such as batteries, generators, and power systems, can be modelled for analysis purposes as a combination of an ideal voltage source and additional combinations of impedance elements.\\
Most sources of electrical energy (mains electricity\footnote{
\textbf{\emph{Mains electricity}} is the general-purpose alternating-current (AC) electric power supply. Mains electricity is the form of electrical power that is delivered to homes and businesses, and it is the form of electrical power that consumers use when they plug kitchen appliances, televisions and electric lamps into wall sockets.
}, a battery, etc.) are best modelled as voltage sources, though in theory we can introduce current sources as components, too.

\item[] \textbf{Current Source} \index{circuit components ! current source}\\
A current source is an electronic circuit that delivers or absorbs an electric current which is independent of the voltage across it. A current source is the dual of a voltage source. The term, \emph{constant-current sink}, is sometimes used for sources fed from a negative voltage supply. 

\begin{figure}[ht]
  \centering
\includegraphics[width=12.8cm]{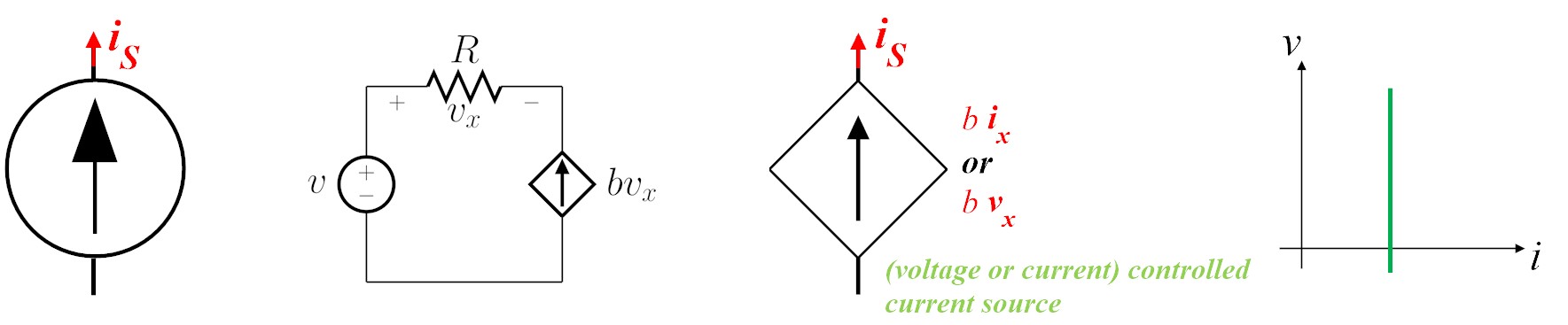}
	\caption{Circuit components: Current Source}
	\label{fig: current-source}
\end{figure}

If the current through an ideal current source can be specified independently of any other variable in a circuit, it is called an \emph{independent} current source. Conversely, if the current through an ideal current source is determined by some other voltage or current in a circuit, it is called a \emph{dependent} or \emph{controlled current source}.\\ 
The current value could be controlled by the voltage over another component, or the current of another branch in the circuit. Soon, we would see a usage of current-controlled current sources in modelling of transistors.

\item[] \textbf{Diode} \label{Diode} \index{circuit components ! diode}\\
A diode is an electrical device allowing current to move through it in one direction with far greater ease than in the other. Diode behaviour is analogous to the behaviour of a hydraulic device called \emph{check valve} (a check valve allows fluid flow through it in only one direction, also known as one-way valve). The most common kind of diode in modern circuit design is the semiconductor diode, although other diode technologies exist.

\begin{figure}[ht]
  \centering
\includegraphics[width=11.8cm]{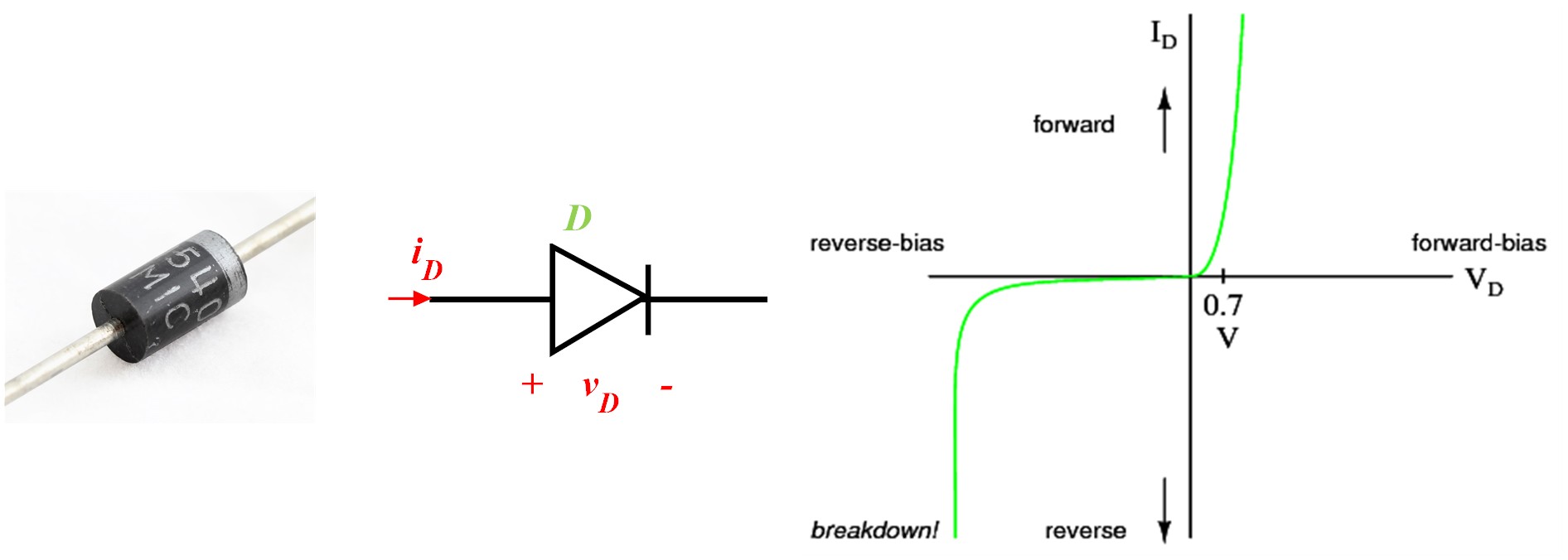}
	\caption{Circuit components: Diode}
	\label{fig: diode}
\end{figure}

When placed in a simple battery-lamp circuit, the diode will either allow or prevent current through the lamp, depending on the polarity of the applied voltage (Figure \ref{fig: diode-working}). The essential difference between forward-bias and reverse-bias is the polarity of the voltage dropped across the diode. 

\begin{figure}[ht]
  \centering
\includegraphics[width=8.8cm]{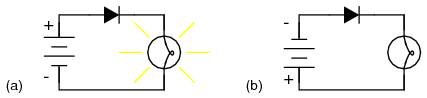}
	\caption{\footnotesize Diode operation: (a) Current flow is permitted; the diode is forward biased. (b) Current flow is prohibited; the diode is reversed biased.}
	\label{fig: diode-working}
\end{figure}

Thus, if the \emph{\lq\lq ideal diode\rq\rq} \,is reverse biased, the current flowing through it is zero. This ideal diode starts conducting at $ 0 ~volts $ and for any positive voltage an infinite current flows and the diode acts like a short circuit. The following figure (Figure \ref{Ideal Diode}), shows a simple circuit using this component with the 
$i - v $ characteristic of this diode and the two working models obtained from that behaviour. 

\begin{figure}[ht]
	\centering
		\includegraphics[width=13.1cm]{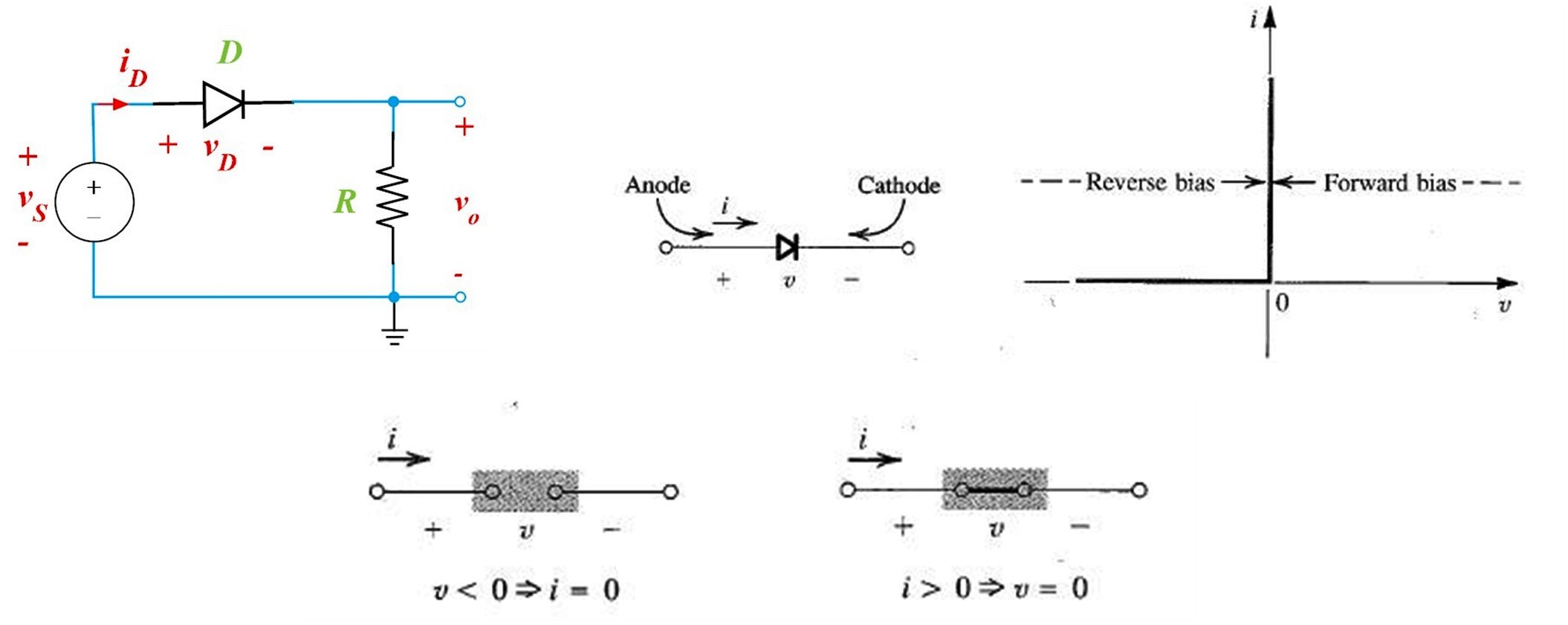}
		\caption{\footnotesize A circuit with Ideal Diode, its $i-v$ characteristic, and equivalent circuits in the reverse and forward directions.}
		\label{Ideal Diode}
\end{figure}

The first step toward practical model for the diode could be obtained by adding a voltage source to an ideal diode in series, compensating the small voltage drop over the component. Next, one can consider a small resistance, again added in series to recover the small slope in the $ i - v $ characteristic of Figure \ref{fig: diode}. The final step would be to consider the breakdown voltage, a fact that plays an important role in Zener diodes.

\paragraph*{Zener Diode} \index{circuit components ! Zener diode} \hfill \\
Unfortunately, when normal rectifying diodes \lq\lq breakdown\rq\rq, they usually do so, destructively. However, it is possible to build a special type of diode that can handle breakdown without failing completely. This type of diode is called a \emph{Zener diode}.
\begin{figure}[ht]
  \centering
\includegraphics[width=10.8cm]{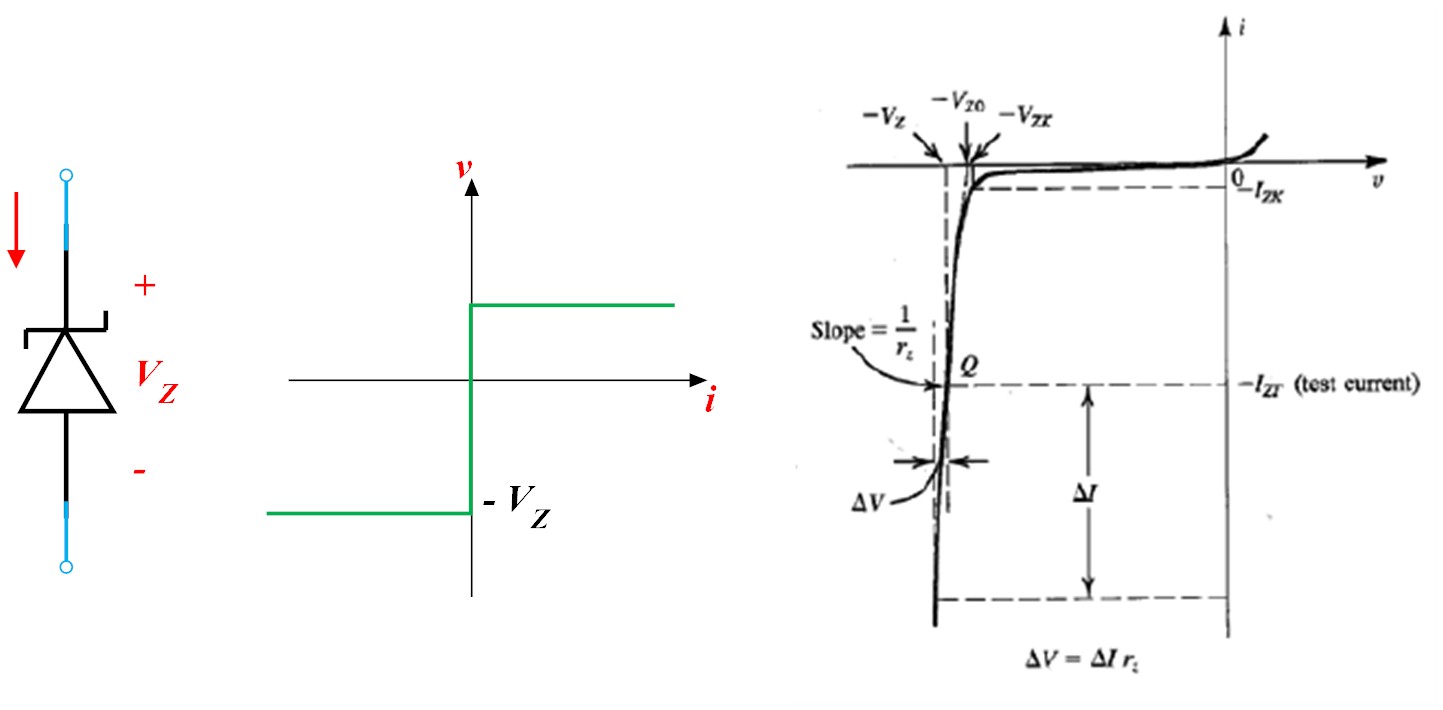}
	\caption{Circuit components: Zener Diode}
	\label{fig: Zener Diode}
\end{figure}

When forward-biased, Zener diodes behave much the same as standard rectifying diodes: they have a forward voltage drop which follows the \lq\lq diode equation\rq\rq \,and is about $ 0.7~volts $. In reverse-bias mode, they do not conduct until the applied voltage reaches or exceeds the so-called \emph{Zener voltage} ($V_Z$ in the figure), at which point the diode is able to conduct substantial current, and in doing so will try to limit the voltage dropped across it to that Zener voltage point. So as long as the power dissipated by this reverse current does not exceed the thermal limits of the diode, it will not be harmed.\\
Zener diodes are manufactured with Zener voltages ranging anywhere from a few volts to hundreds of volts. This Zener voltage changes slightly with temperature, and like common carbon-composition resistor values, may be anywhere from $5 \% $ to $10 \% $ in error from the manufacturer's specifications.

\paragraph*{LED} \index{circuit components ! LED} \hfill \\
Some semiconductor junctions, composed of special chemical combinations, emit radiant energy within the spectrum of visible light as the electrons change energy levels. Simply put, these junctions glow when forward biased. A diode intentionally designed to glow like a lamp is called a \emph{light-emitting diode}, or \emph{LED}.

\begin{figure}[ht]
  \centering
\includegraphics[width=7.8cm]{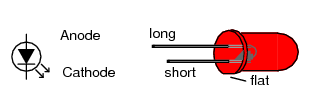}
	\caption{Circuit components: LED}
	\label{fig: LED-schematic}
\end{figure}

Despite the different schematic symbol of LEDs shown in Figure \ref{fig: LED-schematic}, there is no other difference (from the mathematical point of view) in analysing a circuit with standard rectifying diodes or with an LED. \\
It is worth mentioning that LEDs are more sensitive to voltage change, and it is reflected in the light they emit. Thus, for the design problems one should try to stabilize the voltage drop over them (cf. Example \ref{Driving an LED with AC, on the importance of understanding the behaviour of the circuit}). 

\paragraph*{Shockley Diode} \index{circuit components ! Shockley diode} \hfill \\
Shockley diodes are four-layer \emph{pnpn} diodes, which were one of the first semiconductor devices invented. The mathematical (and also physical) interesting point about these family is that their $ i - v $ characteristic exhibits \emph{hysteresis}, the property whereby a system fails to return to its original state after some cause of state change has been removed. 

\begin{figure}[ht]
  \centering
\includegraphics[width=10.8cm]{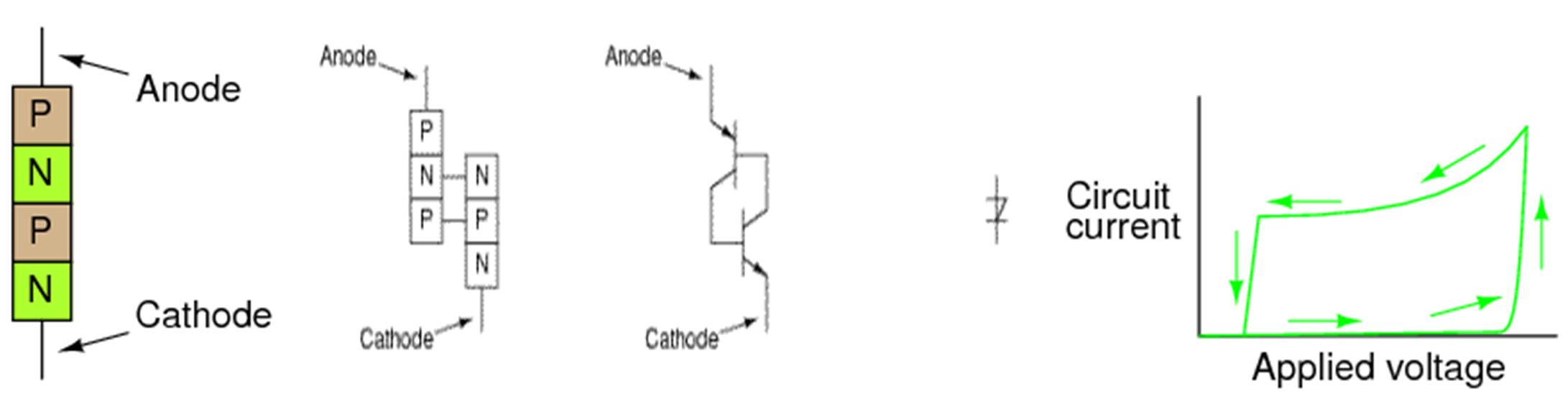}
	\caption{Circuit components: {\footnotesize Shockley diode. from left to right: schematic, physical diagram, equivalent schematic, circuit symbol, and $ i - v $ characteristic } }
	\label{fig: Shockley diode}
\end{figure}

Figure \ref{fig: Shockley diode} also gives the idea that how does a Shockley diodes behaves as a pair of interconnected \emph{pnp} and \emph{npn} transistors. We only refer to two components with this physical structure:

\subparagraph*{DIAC} \index{circuit components ! DIAC} \hfill \\
Like all diodes, Shockley diodes are unidirectional devices; that is, these only conduct current in one direction. If bidirectional (AC) operation is desired, two Shockley diodes may be joined in parallel facing different directions to form a DIAC. The term DIAC is an acronym of \lq\lq diode for alternating current\rq\rq. \\
When breakdown occurs, the diode enters a region of negative dynamic resistance, leading to a decrease in the voltage drop across the diode and, usually, a sharp increase in current through the diode. The diode remains in conduction until the current through it drops below a value characteristic for the device, called the holding current, IH. Below this value, the diode switches back to its high-resistance, non-conducting state. This behaviour is bidirectional, meaning typically the same for both directions of current.

\begin{figure}[ht]
  \centering
\includegraphics[width=10.8cm]{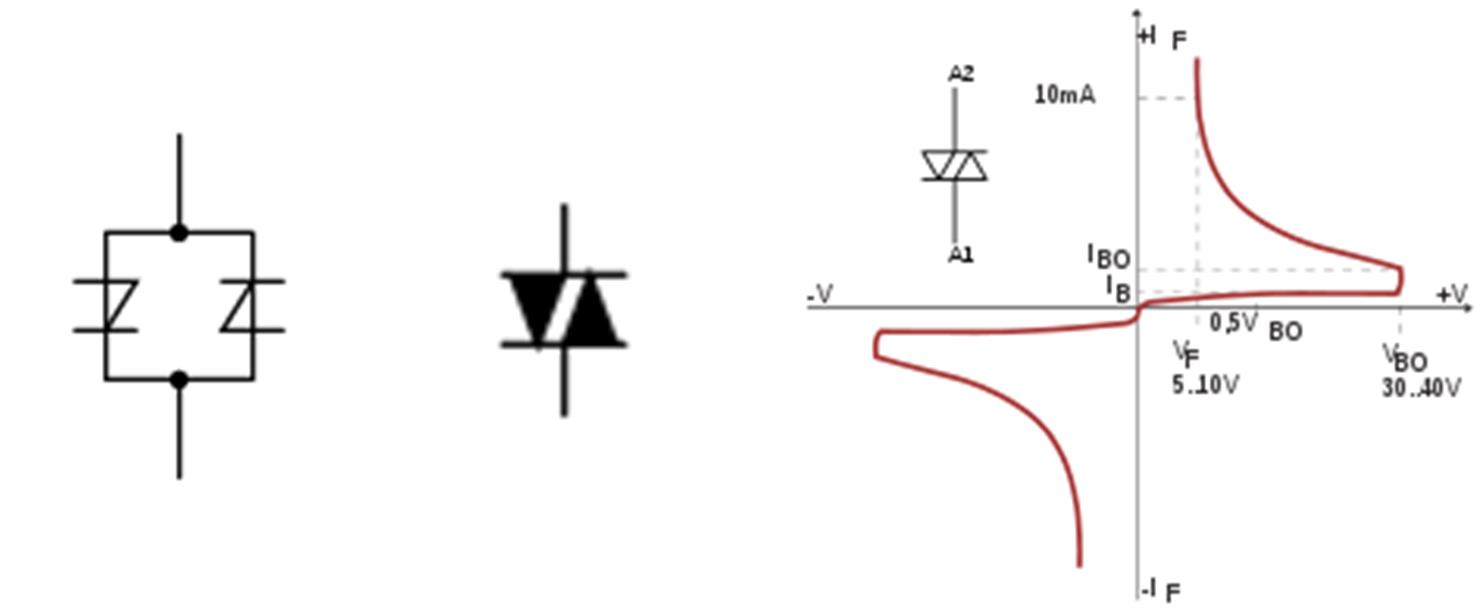}
	\caption{Circuit components: {\footnotesize DIAC. from left to right: equivalent circuit, schematic symbol, and $ i - v $ characteristic} }
	\label{fig: DIAC - symbol}
\end{figure}

DIACs are also called \lq\lq symmetrical trigger diodes\rq\rq \,due to the symmetry of their characteristic curve. Since DIACs are bidirectional devices, their terminals are not labelled as anode and cathode but as A1 and A2 or main teminal MT1 and MT2. DIACs are widely used in light dimmers, starter circuits for fluorescent lamps, and  in conjunction with TRIACs to equalise their switching characteristics. 

\subparagraph*{Silicon-Controlled Rectifier (SCR)} \index{circuit components ! SCR} \hfill \\
In order to expand the usefulness of Shockley diodes, one can equip them with another means of latching\footnote{
This term is used to describe the Shockley diode \lq\lq on\rq\rq \,state. 
To get a Shockley diode to latch, the applied voltage must be increased until break-over is attained.
}. 
In doing so, each of the \emph{npn} and \emph{pnp} junctions/transistors becomes true amplifying device, and we refer to these components as silicon-controlled rectifiers, or SCRs. The progression from Shockley diode to SCR is achieved with one small addition, actually nothing more than a third wire connection to the existing \emph{pnpn} structure. This extra terminal is called the gate, and it is used to trigger the device into conduction (latch it) by the application of a small voltage.

\begin{figure}[ht]
  \centering
\includegraphics[width=9.8cm]{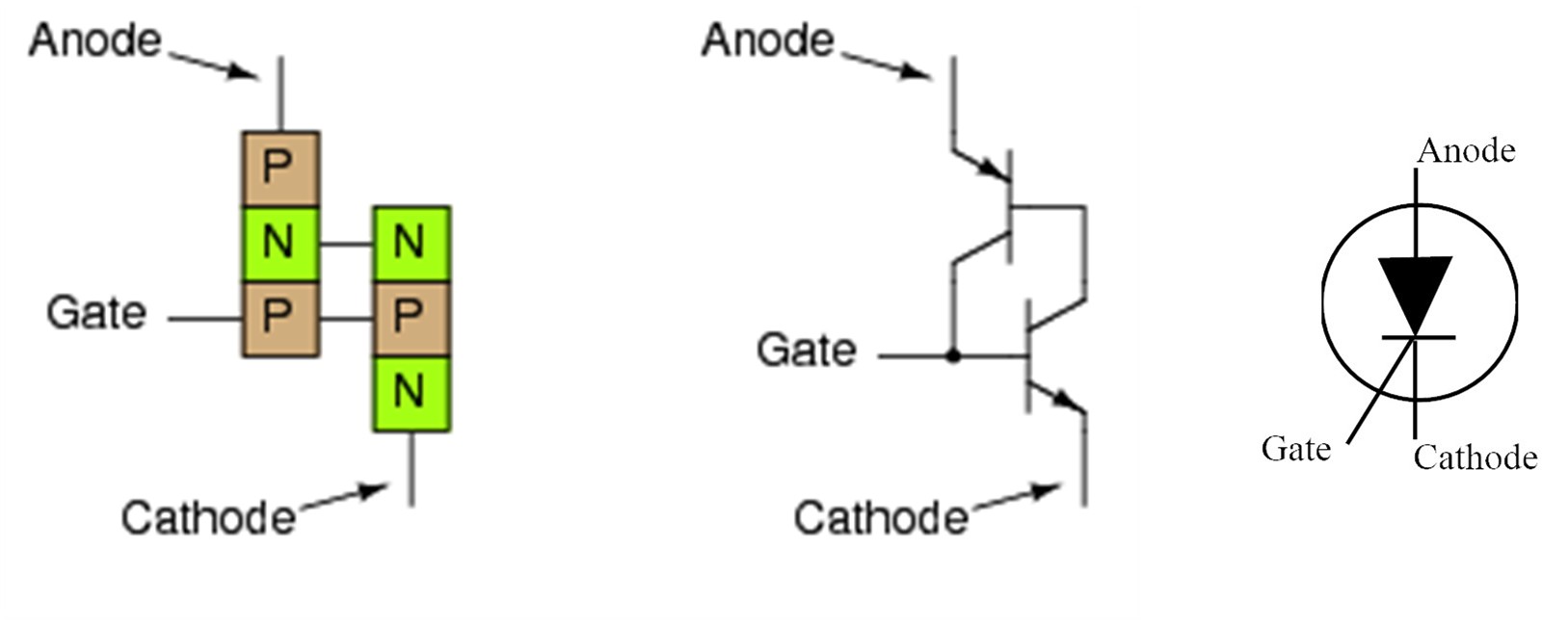}
	\caption{Circuit components: {\footnotesize SCR. from left to right: physical diagram, equivalent schematic, and circuit symbol} }
	\label{fig: SCR}
\end{figure}

SCRs are mainly used in devices where the control of high power, possibly coupled with high voltage, is demanded. Their operation makes them suitable for use in medium- to high-voltage AC power control applications, such as lamp dimming, regulators and motor control.

SCRs and similar devices are used for rectification of high-power AC in high-voltage direct-current power transmission

\item[] \textbf{Transistor} \index{circuit components ! transistor} \\
A transistor is a device composed of semiconductor material usually with at least three terminals for connection to an external circuit. The transistor is the fundamental building block of modern electronic devices, and is ubiquitous in modern electronic systems; ranging from signal amplification to the design of digital logic and memory circuits.\\
The basic principle involved is the use of the voltage between two terminals to control the current flowing in the third terminal, however, the material and technology used to produce transistors are very different, and thus circuit symbols, abbreviations and details of current-voltage relations vary from type to type. 
In this thesis, we only introduce a Bipolar\footnote{
Bipolar transistors are called bipolar because the main flow of electrons through them takes place in two types of semiconductor material: P and N, as the main current goes from emitter to collector (or vice versa). In other words, two types of charge carriers - electrons and holes - comprise this main current through the transistor.
} 
Junction Transistor (BJT).

\begin{figure}[ht]
  \centering
\includegraphics[width=9.8cm]{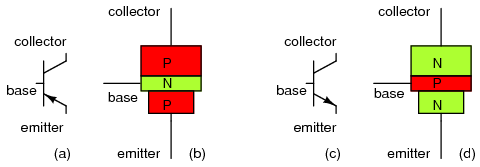}
	\caption{Circuit components: BJT transistor: (a) PNP schematic symbol, (b) physical layout (c) NPN symbol, (d) layout.}
	\label{fig: BJT-schematic}
\end{figure}

The BJT consists of two \emph{pn} junctions, the \emph{emitter-base} junction (\textbf{EBJ}) and the \emph{collector-base} junction (\textbf{CBJ}). Depending on the bias condition (forward or reverse) of each of these junctions, four different modes of operation of the BJT are obtained.\\
In order to understand and analysis the behaviour of a BJT in a circuit, one should study the voltage dropped over each junction and the current flow obtained. Figure \ref{fig: npn-BJT-active-mode-bias} describes the current flow in an \emph{npn} transistor biased to operate in the \emph{active mode} (that is, when EBJ is in forward bias, and CBJ is in reverse bias).

\begin{figure}[ht]
  \centering
\includegraphics[width=9.8cm]{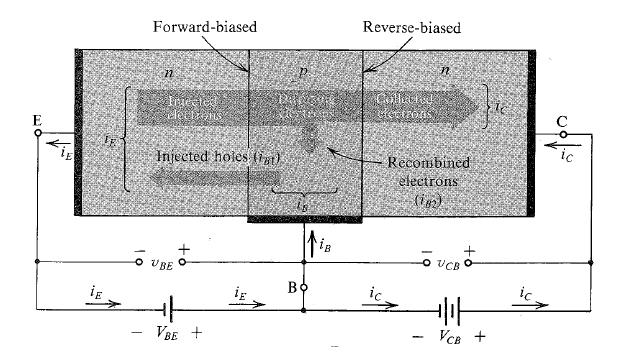}
	\caption{An \emph{npn} BJT biased in active mode}
	\label{fig: npn-BJT-active-mode-bias}
\end{figure}

Among different models that has been suggested for analysing BJT transistor circuits, Ebers and Moll, two early workers in the area, have shown that the following composite model (cf. Figure \ref{fig: EB-model-npn}) can be used to predict the operation of the BJT in \emph{all} of its possible modes.

\begin{figure}[ht]
  \centering
\includegraphics[width=9.8cm]{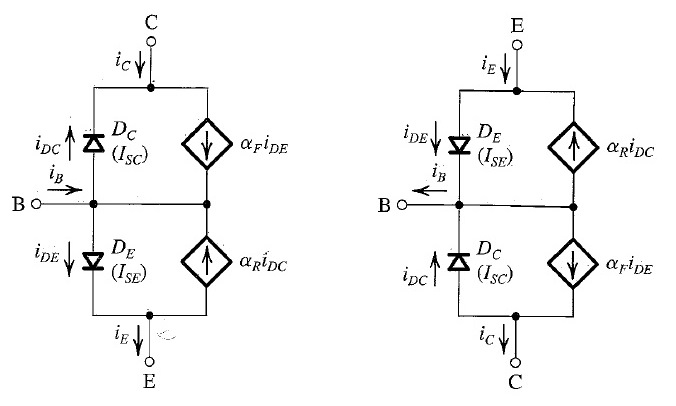}
	\caption{The Ebers-Moll (EM) model of the \emph{npn} transistor (left), and the \emph{pnp} transistor (right)}
	\label{fig: EB-model-npn}
\end{figure}

Thus, theoretically, one can replace a transistor with two diodes and two current-controlled current sources. The diodes in Ebers-Moll model are assumed to be ideal diodes. 
\end{enumerate}

\rem
Before proceeding to the next part of the circuit theory, we want to provide an interpretation of the set-valued map representing the diodes $ i -v $ characteristic. 
If we compare the $ i -v $ characteristic of a resistor with the characteristic of a practical diode (Figures \ref{fig: resistor} and \ref{fig: Zener Diode}), we can obtain two important facts:
\begin{enumerate} [topsep=-1ex, itemsep=1ex, partopsep=0ex, parsep=0ex, leftmargin = 5ex]
\item[\textbf{(a)}] the characteristic of a diode is not enough to identify the exact value of the voltage over it or the current passing thorough it by knowing the other one; while for the resistor at any point of the $ i -v $ characteristic, if you have in hand the value of $v_R$, you can obtain the value of $i_R$, and vice versa. However, the characteristic of the diode gives the idea of what could be the possible values of $v_D$ or $i_D$, and what range of values are not allowed.\\
In other words, the set-valued characteristic talks about the possibility of the values a diode can take as $v_D$ for a certain value $i_D$, and vice versa.

\item[\textbf{(b)}] in a specific circuit, if you want to change the current passing thorough a resistor $R_1$, while the voltage over it should be fixed at a certain value $v_1$, you have no option but to change the resistor; while for the diodes the possibility of accepting different voltages lets us use the same diode. In other words, the voltage and current of a resistor could be determined by its own characteristic, and for the diodes it depends also on the other components of the circuit. Thus, for a range of different components and different combinations, one can still use the same diode, that is, the diode could be matched with a \lq\lq set\rq\rq \,of other components. 
\end{enumerate}

Electrical components could be divided into two groups: active, and passive.\\
The \emph{active} components produce the electrical energy (again described in terms of voltage and current), while the \emph{passive} components are users of this energy. 
Resistors, and capacitors are examples of passive components, while voltage sources are considered as active elements. Transistors could be biased to work in an active mode, too, so that they would amplify the AC signal. \\
We need to consider this division while writing the Kirchhoff's circuit laws. 

\paragraph*{Kirchhoff's Laws} To analyse a circuit, we need to know the $i - v$ characteristics of all the components in the circuit, but that is not enough. In addition, we need some rules to describe how the current passes through different components and branches, and how to measure the drop of potential differences over each component in a specific circuit.  These rules (which are a reduced version of the Maxwell's electromagnetic equations), are known as \emph{Kirchhoff's laws}:
\begin{itemize}
\item[$\circ$] \emph{Kirchhoff's current law (\textbf{KCL})}: \index{KCL, Kirchhoff's current law} \\
at any node (junction) in an electrical circuit, the sum of currents flowing into that node is equal to the sum of currents flowing out of that node;
equivalently,\\
the algebraic sum of currents in a network of conductors meeting at a point is zero.
\item[$\circ$] \emph{Kirchhoff's voltage law (\textbf{KVL})}: \index{KVL, Kirchhoff's voltage law} \\
the directed sum of the electrical potential differences (voltage) around any closed network is zero;
equivalently,\\
the sum of electrical potential drops in a closed loop is equal to the total electromotive force (emf) available in that loop.
\end{itemize}

The word \lq\lq directed\rq\rq \,in the statement of KVL, means that the active and passive components have different voltage signs in the sum. It is mostly optional to choose, but both can not be positive (or negative) at the same time. Another way to consider this fact is that the product of current and voltage of passive component should be positive, while the $i \, v$ of an active component is negative, or vice versa. \\
For example, note that the conventional direction for voltages and currents in Figure \ref{fig: voltage-source} imply $ i_s \times v_s < 0 $, meaning that voltage sources are active components, while for a resistor as shown in Figure \ref{fig: resistor}, we have 
$ i_R \times v_R > 0  $. \\
Solving or analysing a circuit means to use the the equations obtained from KVL, KCL, and $i - v$ characteristics of components, to find the current passing through and the voltage over each component of the circuit. This usually ends up to a system of $n$ equations, and $n$ variables.

\begin{eg}
Figure \ref{divisor-circuit1} shows a simple circuit with a voltage source as the model for electrical power producer, a resistor $ R_B = 1.4 \, k \Omega $ as the bias resistor, and a consumer modelled as a load resistor $R_L = 5.6 \, k \Omega $. Our goal is to obtain the percentage of the produced voltage that is delivered to the consumer.

\begin{figure}[ht]
  \centering
\includegraphics[width=5.8cm]{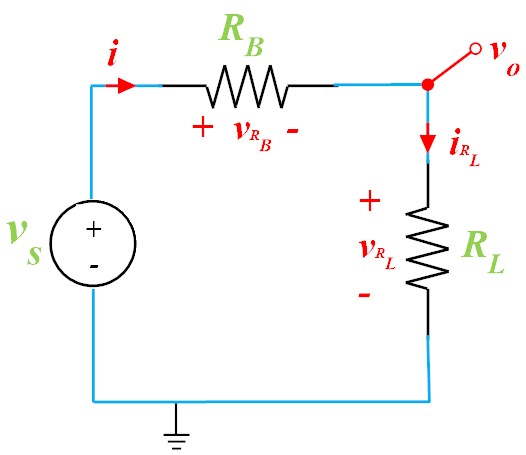}
	\caption{\footnotesize A simple circuit that divides voltage between a bias resistor and a load}
	\label{divisor-circuit1}
\end{figure}
\end{eg}

Since there is only one loop in this circuit, KCL implies $i_{R_B} = i_{R_L} = i $. For using KVL, consider a clockwise closed path starting from the voltage source. 
$$ - v_s + v_{R_B} + v_{R_L} = 0. $$
Now by using the component relations $  v_{R_B} =  R_B \,  i_{R_B} $, and $  v_{R_L} =  R_L \,  i_{R_L} $, one can calculate the current value in the circuit as  
$ i = \frac{v_s}{ R_B + R_L }$. Using the $i-v$ relation of the resistor once more, yields
$$ v_o = v_{R_L} = \dfrac{R_L}{R_B + R_L } \, v_s. $$
Thus, with the given values for $R_B$, and $R_L$, we would have a $20 \% $ loss in delivering the voltage produced by $v_s$ to the load. 

\section{Formulating the Problem} \label{Formulating the problem}
\begin{figure}[ht]
  \centering
\includegraphics[width=7.1cm]{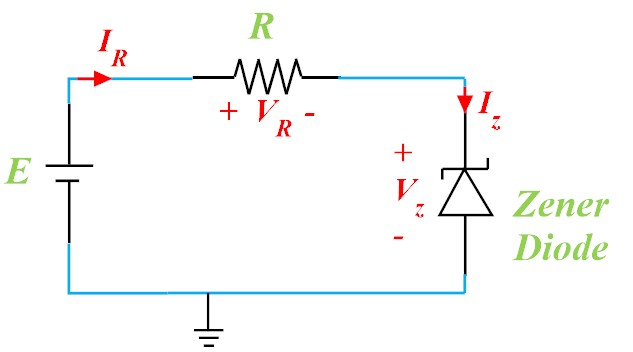}
$ ~~~~~ $
\includegraphics[width=3.8cm]{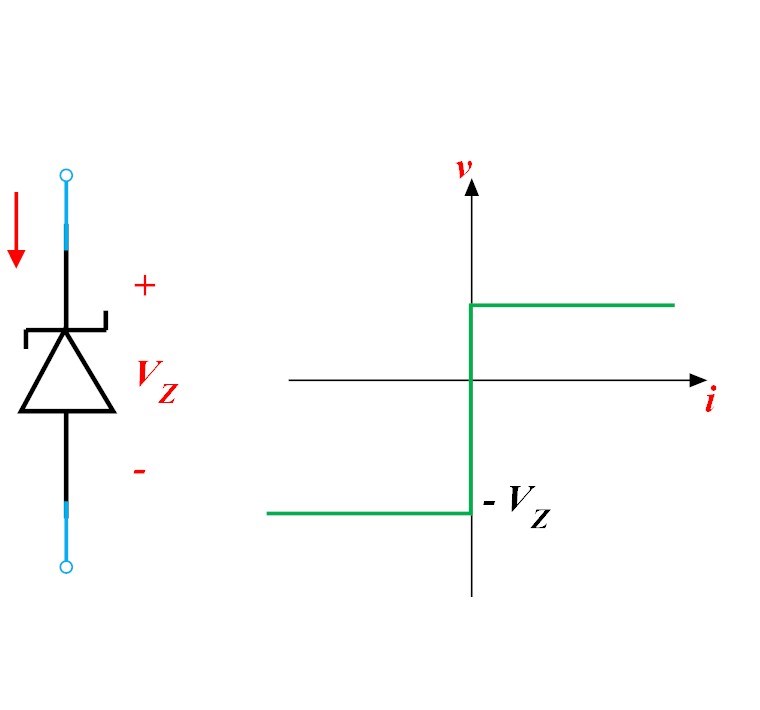}
	\caption{\footnotesize A simple electrical circuit, the schematic and   $i - v$ characteristic of Zener Diode}
	\label{Zener Diode simple circuit}
\end{figure}

Based on the components in the circuit, and the exactness of the solution required, there is a huge theory and lots of work done around it till now, in electrical engineering literature. We are not exactly interested in this topic, but we will use the setting of circuit theory and the rules mentioned in previous section to formulate our problem.\\
We start with a simple circuit as shown below (Figure \ref{Zener Diode simple circuit}) which involves an electronic component named \emph{Zener Diode}. Our aim is to find a relation between the given input voltage ($ E $, here) and the variable $ I $, current, in the circuit.
\begin{equation}\label{formulation1}
\left.\begin{matrix}
\mathrm{KVL: }  ~ -E + V_R + V_z = 0~\\ 
\mathrm{KCL: } ~~~~ I_R = I_z = I ~~~~~~~~\\ 
~~\, V_R = R \,   I_R \\ 
~~\, V_z \in F(I_z)
\end{matrix}\right\}~~\Rightarrow
~~~E \in R \, I ~ + F(I),
\end{equation}
in which $F(I) $ describes the relation between $I_z$ and $V_z$ as a set-valued map. \\
To match our mathematical setting, we change $E$ to $p$, in order to indicate that it is a \emph{parameter}; and $I$ to $z$, to show that $z$ is the variable, thus we get
\begin{equation}
p \in f(z) + F (z),
\end{equation}
in which $f(z) = R \, z$, in this particular example. For a given $p$, we are interested in the \emph{solution mapping}\index{solution mapping} defined as
\begin{equation}\label{solution mapping2}
S(p) = \left \{  z \in \R ~|~ p \in f(z) + F(z) \right \}.
\end{equation}
It is also useful to introduce another notation
\begin{equation}
\boxed{ \Phi (z) := f(z) + F(z) , \mathrm{~~~~~so~~~~} S(p) = \left \{  z \in \R ~|~ p \in \Phi (z)  \right \} = \Phi^{-1} (p).} 
\end{equation}

\begin{figure}[ht]
  \centering
\includegraphics[width=4.8cm]{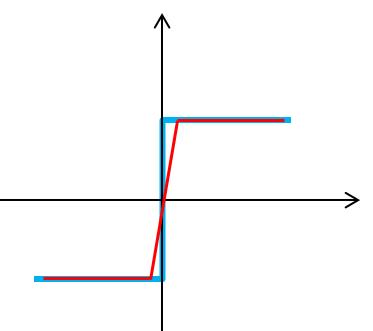}
	\caption{\footnotesize Approximation of the Diode's $i - v$ characteristic}
	\label{approximation of a set-valued map}
\end{figure}

Before proceeding our study with this model, we need to justify its appropriateness.\\
One may wonder why this form of generalized equation does not appear in electronic books. This depends on the way we try to model the behaviour of the diode. 
When one is using the ideal model for diodes, usually the two different working modes of the diode (that is, when the diode is on and acting like a wire, and when it is off and acting like a gap or open connection) are studied separately. \\
When one is only interested in the forward bias, the $ i - v $ relation is usually approximated by an exponential term, or some approximations of this exponential map. \\
Since the $i - v$ characteristic of a Zener Diode (and diodes in general) fails from a function-form description only in one point of its domain, those books approximate the given set-valued map (and similar other maps for other type of diodes) with a single-valued map (look at Figure \ref{approximation of a set-valued map}). So the later inclusion in the right hand side of (\ref{formulation1}) becomes equality, and they end up in an equation. 

Moreover, we are not going to solve the circuit equations, which is mostly the aim in electronics. We want to study the behaviour of the solution mapping in a neighborhood of a certain point (that is a local study) and it needs precise models concerning each point $(p, z)$ in the input-output relation.\\
Regarding different electronic components in the circuit and various circuits, the obtained solution mapping (\ref{solution mapping2}) might change a bit, but we will come back to this difficulty in few pages, after stating the general purpose of introducing these type of problems.\\
We are going to answer this question that what will happen if the voltage of source changes from $\bar{p}$ to $p'$. 
There are three specific reasons to validate the importance of this question:
\begin{enumerate}[topsep=0ex,itemsep=0ex,partopsep=1ex,parsep=1ex] 
\item[$\bullet$] \textbf{Failure in Precise Measurements}\\
When dealing with real world applications, one should always be aware of inexact measurements. Any producer provides an error percentage or tolerance range for his products. For example, when you have a resistor with the following colored bands : Brown, Green, Red, and Gold, then its resistance is $15 \times 10^2 \, \pm 5 \% $. \\
Thus, either the signal source is a single component with a tolerance range, or it is a circuit itself, made of several components, you would have a deviation from the exact value $\bar{p}$. \\
Also when you are designing a circuit, you find the optimal value of your DC voltage source $10.53~volts$, for example. But in the market, you cannot have such a battery or combination of batteries that provide this precise value. Thus, you are forced to use a physical component close to your calculations, but not exactly the same. 
\item[$\bullet$] \textbf{Process of Ageing}\\
Consider a simple chemical pile as the voltage source. As soon as you start to use it in a circuit, the chemical ionization process that provides the electrical energy runs inside the pile. After a certain time, you would notice that the battery is dead and you need to change it or charge it.\\
During this period, although the battery level was at the appropriate level, but it was decreasing slowly to get to the unacceptable level. In other words, there were small changes in its precise value.
\item[$\bullet$] \textbf{Thermal Effect}\\
Especially when the signal source is an electronic circuit itself, one should consider that semiconductors are very sensitive to change of temperature. \\
For example, when an output of a  circuit with diodes or transistors is used as the signal source for the next circuit (like the cascade structure in amplifiers), the behaviour of these components change with temperature and thus, the outcome would have a little tolerance around the precise expected value.
\end{enumerate}

We would now try to provide an interpretation of the local stability properties of the solution mapping in terms of the circuit parameters. Let us note that based on the problem one may face during the design process, one of these properties would fit better to his/her demands.

\begin{figure}[ht]
  \centering
    \includegraphics[width=0.68\textwidth]{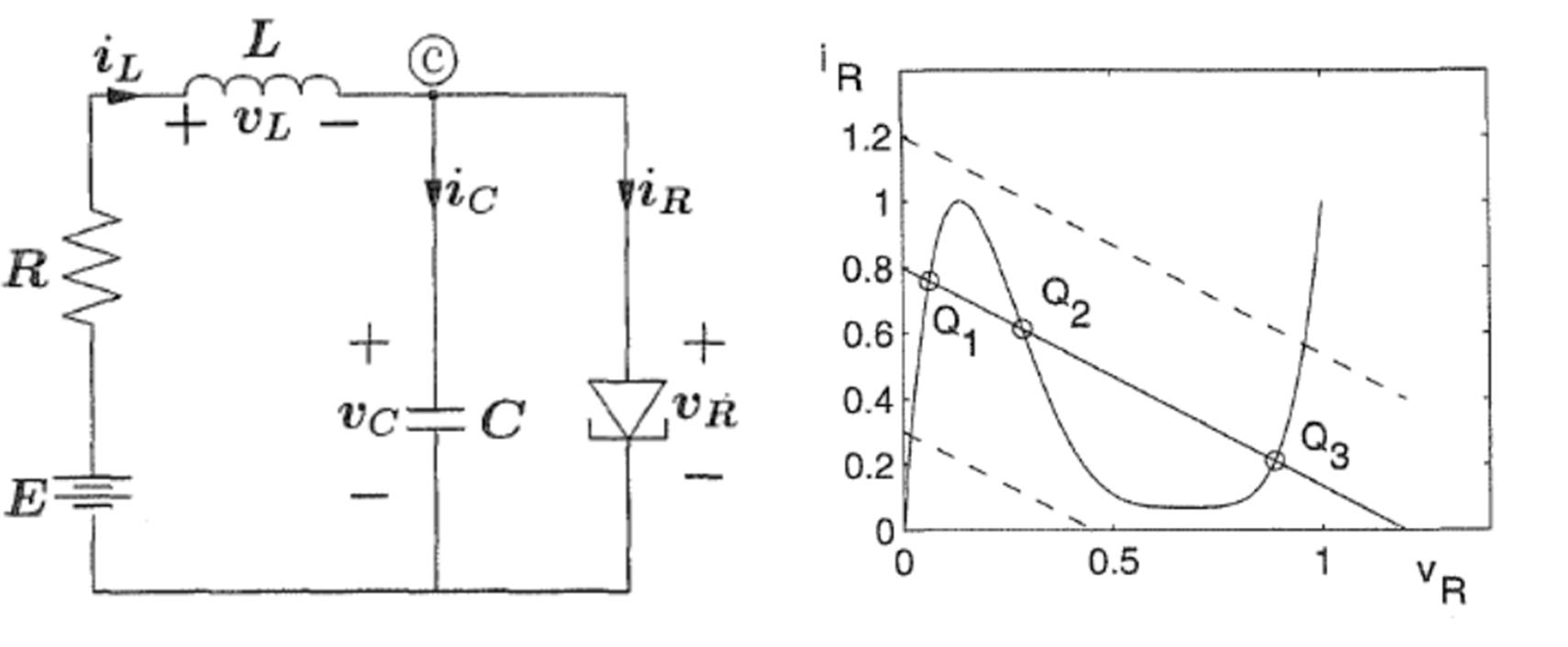}
    \caption{\footnotesize A circuit with Tunnel diode (left) and the graphical analysis of it (right)}
    \label{operatingpoint2}
\end{figure}

\paragraph*{Stability Formulation} Suppose for a given $\bar{p}$, we know the previous current of \emph{operating point}\footnote{ 
In the graphical analysis of the circuit, we plot two maps on the same $i - v$ plain: the $i-v $ characteristic of the diode, and the ordered (with respect to $i$) equation of the circuit gained from KVL. The  solution can then be obtained as the equilibrium point, that is the coordinates of the intersection point of the two graphs. This point is called ``operating point'' (look at the figure below). It is easy to see that it depends both on the inner structure of the diode and the rest of the circuit.\\
\begin{center}
  \includegraphics[width=0.48\textwidth]{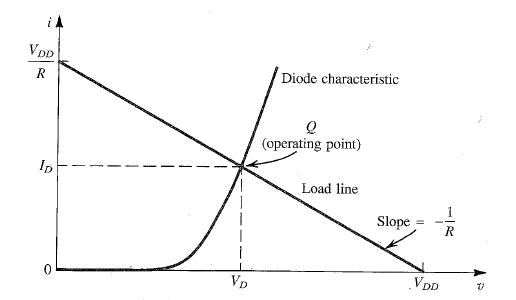}
\end{center}
}, 
say $\bar{z}$. We also assume that the input change is small, that is, in mathematical terms, $p ' \in \B_r (\bar{p}) =: V$ for some small $r >0$.\\
In fact, we are interested in those circuits that keep the small input-change, small. More precisely, the distance between a $z \in S(p')$ and $\bar{z}$, is controlled by the distance $ \|~p' - \bar{p}~\|$. Thus we wish (and search for) $S$ having the following property
\begin{equation} \label{IC-formulation}
 \|~z - \bar{z}~\| \, \leq \, \kappa ~ \|~p' - \bar{p}~\|~~\mathrm{~~when~~~}~z \in S(p') \cap U,~~p' \in V
\end{equation}
where $U$ is a neighborhood of $\bar{z}$, $V$ is a neighborhood of $\bar{p}$, and $\kappa \geq 0 $ is a constant.\\
We defined this property as \emph{\textbf{isolated calmness}} in Chapter 1 (cf. Definition \ref{IC definition}). \\
The reason we considered the intersection $ S(p') \cap U $ in the above formulation, is that while it is possible to have different values in $S(p')$, we are just interested in quarantining the existence of a $z$ near enough to $\bar{z}$. $S$ may be not single-valued at $p'$, but we do not care (in other words, the control is not about that); what matters is that it must have a value in $U$.\\
The very next question about equation \eqref{IC-formulation} would be whether $S$ is single-valued at $\bar{p}$ or not. Actually, single-valuedness is implicitly assumed in that formulation, but one could be in a situation described in Figure \ref{operatingpoint2} involving a Tunnel diode.

The answer is that it is not a necessary assumption in general, though in some particular cases it might be a request. Thus, we get a slight modification of the previous formulation and ask for the existence of $\kappa \geq 0$, neighborhoods $U$ of $\bar{z}$, and $V$ of $\bar{p}$ such that
\begin{equation}
e(S(p') \cap U, S(\bar{p})) \leq \kappa ~ \| p'- \bar{p} \| \mathrm{~~~~for~ all~} p' \in V
\end{equation}
where $\bar{z} \in S(\bar{p}) $, and $ e(A, B)  $ is the \emph{excess} of $ A $ beyond $ B $. 
We defined this property as \emph{\textbf{calmness}} in Chapter 1  (cf. Definition \ref{C definition}). \\
If we are investigating a general local property of the solution mapping and the point $\rfp{p}{z} \in \gph{S} $ does not play a crucial role in our study, we would be interested in the \lq\lq \emph{two-variable}\rq\rq \,version of the previous condition, that is
\begin{equation}
e(\,S(p') \cap U, S(p) \,) \leq \kappa ~ \| p' - p \| \mathrm{~~~for~all~~~~} p', p \in V := \B_r (\bar{p}),
\end{equation}
We defined this as \emph{\textbf{Aubin property}} in Chapter 1 (cf. Definition \ref{AP}). 

\paragraph*{Regularity Formulation} Up to now, our construction were built under the assumption that the explicit form of the solution mapping $S(\cdot)$ is in hand and so we can easily calculate values like $S(p')$. \\
It is not true in general. All we are sure we can get from the circuit is $ \Phi = S^{-1} $. It is not always simple (or even possible) to derive the formula of $S(\cdot)$, so we need to provide the proper formulation of the desired properties with $\Phi$.\\
Two slightly different interpretations of the problem will result in the two different approaches described below.

\textbf{Approach 1.} One can take almost the same procedure. Starting from $\bar{z}$, and $\bar{p} \in \Phi(\bar{z}) $, then consider a point $ z \in U =: \B_r (\bar{z}) $, calculate $ \Phi(z)$, and ask for a control $\kappa \geq 0$ satisfying the following condition:
\begin{equation}
d ( z, \bar{z}) ~ \leq ~ \kappa \, d (\bar{p}, \Phi(z) \, \cap V)  \mfa z \in U .
\end{equation}
This property was introduced as \emph{\textbf{strong metric sub-regularity}} in Chapter 1 (cf. Definition \ref{def. smsr}). \\
If one releases the single-valuedness condition $ \Phi(\bar{z}) = \left \{ \bar{p} \right \} $, gets
\begin{equation} 
d ( z, \Phi^{-1} (\bar{p})) ~ \leq ~ \kappa \, d (\bar{p}, \Phi(z) \cap V )  \mfa z \in U .
\end{equation}
which was defined as \emph{\textbf{metric sub-regularity}} in Chapter 1 (cf. Definition \ref{def. msr}). \\
Then, we can do a general local study around an arbitrary point $\rfp{z}{p} \in \gph{\Phi}$ by using a \emph{two variable} condition 
\begin{equation}
d \left (  z,\Phi^{-1} (p) \right ) \leq \kappa~ d \left ( p, \Phi(z) \right ) \mathrm{~~~whenever~~} (z,p) \in U \times V.
\end{equation}
which was defined as \emph{\textbf{metric regularity }} in Chapter 1 (cf. Definition \ref{def. mr}). 

\textbf{Approach 2.} One can consider a different approach toward the problem, by starting from a point $\rfp{z}{p} \in \gph{\Phi}$, and taking an arbitrary $z' \in U =: \B_{r_1} (\bar{z}) $, calculating $\Phi(z')$, taking a $p \in V =: \B_{r_2} (\bar{p})$, and then asking for the chance of having the output distance $ d \left (  z',\Phi^{-1} (p) \right ) $ being controlled by the input distance $ d \left ( p, \Phi(z') \right ) $, that is,
\begin{equation}
d \left (  z',\Phi^{-1} (p) \right ) \leq \kappa~ d \left ( p, \Phi(z') \right ) \mathrm{~~~whenever~~} (z',p) \in U \times V,
\end{equation}
which is again \textbf{metric regularity}. Then make slight modifications by considering the additional condition of single-valuedness at $\bar{z}$ (that is, 
$ \Phi(\bar{z}) = \left \{ \bar{p} \right \} $, called \emph{\textbf{strong metric regularity}}, cf. Definition \ref{def. smr}). Or considering the one-variable version of the above condition (introduced as \textbf{metric sub-regularity }).\\
 
\note \textbf{(a)} The numerical advantage of regularity formulation was already discussed in Remark \ref{numerical MR}.\\
\textbf{(b)} One should pay attention to the fact that $z'$, and $p$, are chosen independently. Indeed, $ p \in \Phi( z')$ ends up at an obvious situation $ 0 \leq 0$. \\
\textbf{(c)} We have already seen that there is some relations between the two mentioned formulations (cf. Section \ref{section1.3}):
\begin{equation*}
\begin{matrix}
\mathrm{ for~ } \Phi: & SMR & \subset & MR & \subset & MSR & \supset & SMSR \\
& & & \cong & & \cong  & & \cong  \\
\mathrm{ for~ } S: &  & & Aubin \, Property & \subset & Calmness & \supset & Isolated \, Calmness 
\end{matrix}
\end{equation*}

\paragraph*{Toward Complexity} We want to introduce two techniques that play an important role in our future calculations. They would also annihilate the idea that our modelling is so much simplified and will not be able to cover the real cases in practical complicated circuits, though we may not refer to such circuits here.
\begin{itemize}[topsep=0ex,itemsep=0ex,partopsep=1ex,parsep=1ex, leftmargin=4ex]
\item[$\bullet$] \textbf{Cascading} \\
In some circuits, one can ``slice'' the circuit to smaller ones, study each separately and sum them up at the end. This is a frequently used technique in designing Amplifiers, and other circuits with very special limitations. 
Here we used \lq\lq cascading\rq\rq \,as a general term for objects connected serially that start a chain reaction\footnote{
In fact, the title \lq\lq cascading\rq\rq \,was drawn from \emph{cascade amplifier circuits}, which are circuits made of different stages, each serves a special improvement of the current or voltage gain in a serial connection. \\
}.
For a better understanding of the situation look at the working Diagram \ref{cascade system}, in which the network is divided into three subnets:

\begin{figure}[ht]
	\centering
		\includegraphics[width=5.1cm]{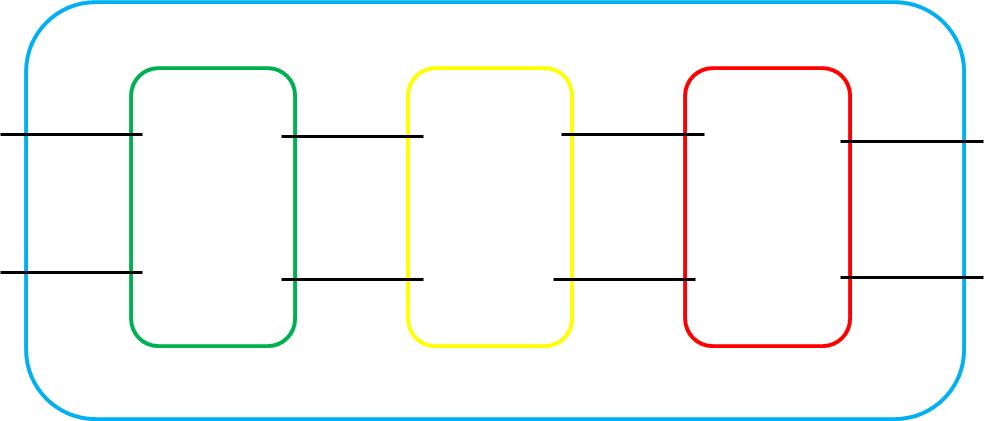}
		\caption{\footnotesize Cascading}
		\label{cascade system}
\end{figure}

The \lq\lq slicing\rq\rq \,process should be done by considering the working diagram and how components are connected to each other for their functionality. Thus,
one should be careful to avoid cutting the feedback loops\footnote{
Feedback loop is a design method that returns part of the output of the system into input, in order to enable the system to adjust its performance to meet a desired output response. The operational diagram of a typical feedback loop could be seen in the following figure

\begin{center}
		\includegraphics[width=4.1cm]{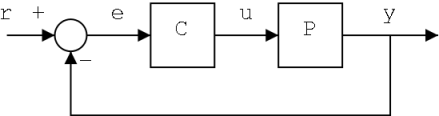}
\end{center}
}, 
or separating a current controlled source from the control branch, ...  and such mistakes.\\
The calculus rues provided in Subsection \ref{Calculus Rules}, especially the composition Propositions \ref{special chain rules for coderivatives} and \ref{special chain rules for coderivatives2} will be the mathematical tools which give us the permission to use this technique while studying the stability properties.

\begin{eg} 
In the following circuit (Figure \ref{cascading}) we can study the stability of the two simpler circuits on the right.

\begin{figure}[ht]
	\centering
	$~$
		\includegraphics[width=6.1cm]{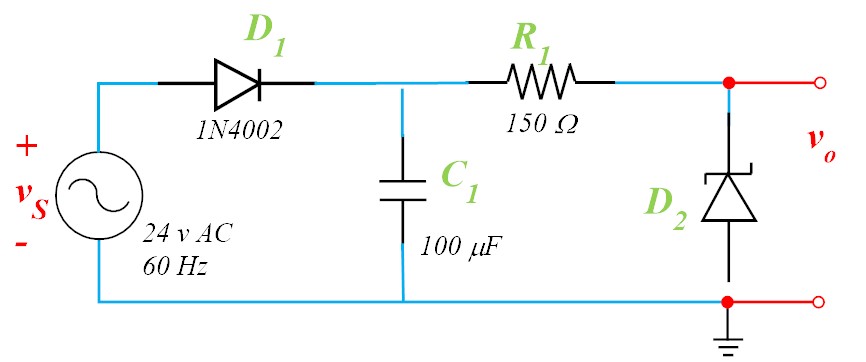}
		$ ~~~~~$
		\includegraphics[width=6.8cm]{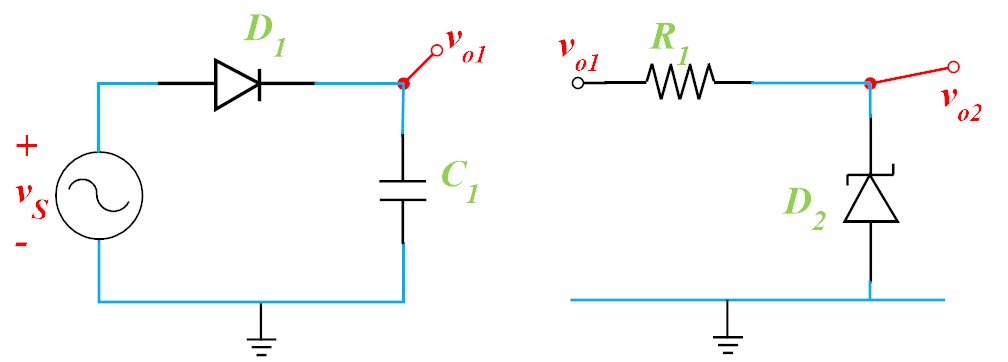}
		\caption{\footnotesize Simplifying the study of the circuit by splitting the layers}
		\label{cascading}
\end{figure}
\end{eg}

Indeed, in the right hand-side figure, we can study the stability of the \emph{transfer functions} $ \dfrac{v_{o1}}{v_s} $ and $\dfrac{v_{o2}}{v_{o1}}$ separately, considering the output of the first sub-circuit as an independent source $v_{o1}$ of the second one. Then use a composition theorem to obtain the stability of the transfer function $\dfrac{v_{o}}{v_s}$.

\item[$\bullet$] \textbf{Multivalued Simplification}\\
As we usually end up to a term like $f(z) \, + \, F(z)$, we try to make the set-valued map as simple as possible. Such a simplification may not have a real physical meaning, but it would be a correct mathematical operation.
To be more clear, consider the following figure that shows the use of this technique for a DIAC. \\
Although the pointwise sum of two graphs on the right would be equal to the graph on the left, the first graph (blue one, representing a single-valued relation) is not an $i-v$ of a component. One can think of the possibility of making a complicated circuit that has such an $i-v$ characteristic. Though it is possible in theory, and thus our technique is valid, in practice, it might be contrary to the idea of \lq\lq simplification\rq\rq \,we were trying to follow.

\begin{figure}[ht]
	\centering
		$~$
		\includegraphics[width=4.1cm]{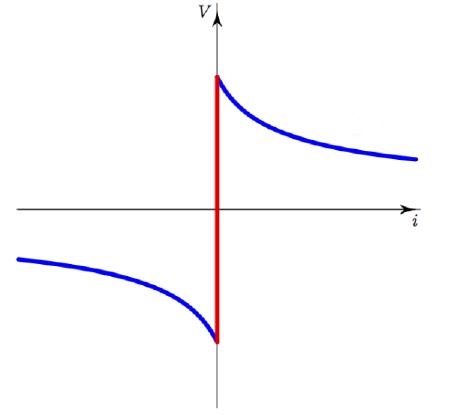}
		$ ~~~~~~~~$
		\includegraphics[width=7.8cm]{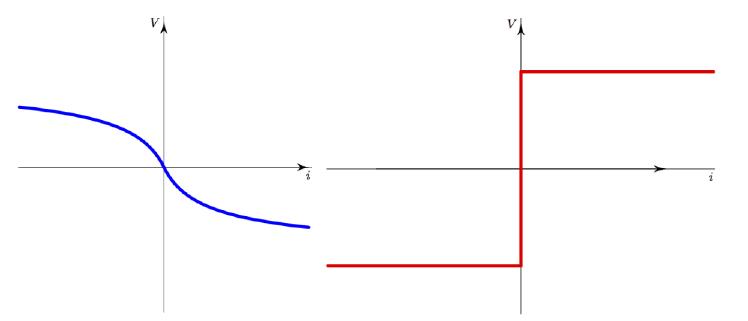}
		\caption{\footnotesize A simplified $i-v$ characteristic of a DIAC and the equivalent maps}
		\label{Multivalued Simplification example}
\end{figure}

This technique will be very useful, and essential somehow when calculating the tangent and normal cones, the tools we will provide for our study in next chapter.
\end{itemize}

Now we are ready to go back to the question we left before: the effect of changing components or their order on the final form of generalized equation that would be obtained. 
We try to provide the answer of this question with some examples.

\begin{eg}[Increasing Components and Loops] \label{Increasing Components and Loops}
In Figure \ref{loops}, we can see a circuit with two loops, each one containing a Diode. This will be a good example about how the setting will change in $\R^n$.
\begin{figure}[ht]
	\centering
		\includegraphics[width=5.1cm]{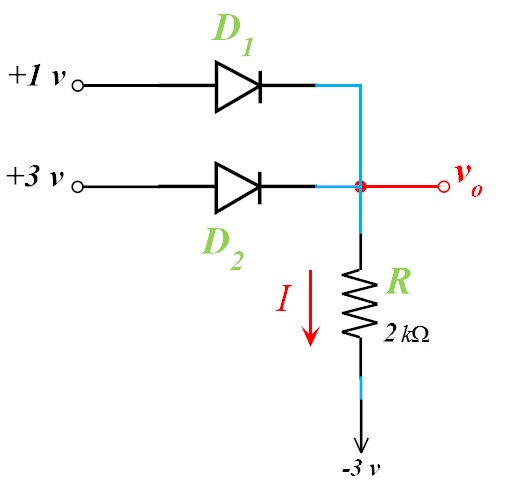}
		$ ~~~~~~~$
		\includegraphics[width=7.8cm]{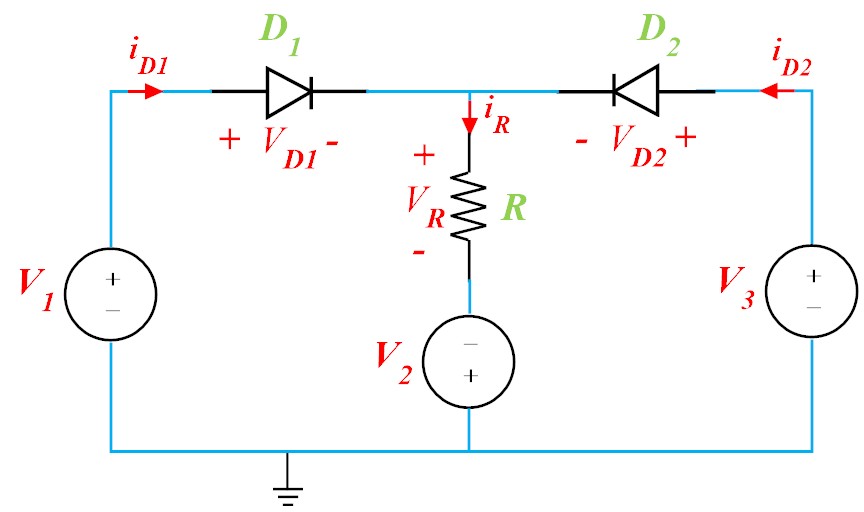}
		\caption{A circuit with two loops}
		\label{loops}
\end{figure}

Suppose the diodes do not have the same characteristic. One can write the following:
\begin{equation}
\left.\begin{matrix}
\mathrm{KVL1: }  -V_1  + V_{D_1}+ V_R - V_2 = 0\\ 
\mathrm{KVL2: } -V_3  + V_{D_2}+ V_R - V_2 = 0\\ 
\mathrm{KCL: }~~~~~~ i_{D_1} + i _{D_2} = i_R~~~~~~~~~\\
~~\, V_R = R \,   i_R ~~\\ 
~~~~\, V_{D_1} \in F_1( i_{D_1})\\
~~~~\, V_{D_2} \in F_2( i_{D_2})\\
\end{matrix}\right\}~~\Rightarrow
\begin{matrix}
  (V_1  + V_2) \in  R \,   ( i_{D_1} + i _{D_2})~ + F_1(i_{D_1})\\ 
   (V_3  + V_2) \in  R \,  ( i_{D_1} + i _{D_2}) ~ + F_2(i_{D_2})
\end{matrix}
\end{equation}
So we get the following unified form
\begin{equation}
\begin{pmatrix}
p_1 \\ p_2
\end{pmatrix}
\in 
\begin{pmatrix}
R & R \\ R & R
\end{pmatrix}
\begin{pmatrix}
i_{D_1} \\ i_{D_2}
\end{pmatrix}
+
\begin{pmatrix}
1 & 0 \\ 0 & 1
\end{pmatrix}
F \left(
\begin{pmatrix}
1 & 0 \\ 0 & 1
\end{pmatrix}
 \begin{pmatrix}
i_{D_1} \\ i_{D_2}
\end{pmatrix}
\right)
\end{equation}

in which $ p_1 =  V_1  + V_2 $, $p_2 =  V_3  + V_2$, and $\mmap{F}{2}{2}$ is defined by 
\begin{equation*}
F
\begin{pmatrix}
z_1 \\ z_2
\end{pmatrix}
:= 
\begin{Bmatrix}
(x, y) ~|~ x \in F_1 (z_1), y \in F_2 (z_2)
\end{Bmatrix}
= \prod_{i=1}^{2} F_i (z_i)
\end{equation*}
by letting $p := (p_1, \, p_2) ^{T}$, and $z := (z_1, \, z_2) ^T$, we get to a general form
\begin{equation}
p \in f(z) \, + B F(Cz)
\end{equation} 
In this example, $
f(z) = A z = 
\begin{pmatrix}
R & R \\ R & R
\end{pmatrix} z $,
$ B = C = 
\begin{pmatrix}
1 & 0 \\ 0 & 1
\end{pmatrix}$.\\
Most of our results till now is about this form of generalized equation. We correct the introduced $\Phi$ and the solution mapping in the following way
\begin{equation}
\boxed{
\begin{matrix}
\Phi (z) : =f(z) \, + B F(Cz) \\
S(p) := \left \{ z ~|~ p \in \Phi (z) \right \}~
\end{matrix}}
\end{equation}
\end{eg}

Thus, one can see increasing the components with the $i-v$ characteristic described as a set-valued map, in different loops will reformulate the problem in $\R^n$ with some correction coefficients to indicate the presence or absence of those set-valued components in different branches. We will study this type of problem in Chapter \ref{Chapter3}.\\
Now, we consider the change of power sources from constant (that is, DC sources) to time-varying ones (that is, AC sources).

\begin{eg}[Effect of AC Sources] \label{Effect of AC Sources}
In the following figure, we have an alternating voltage source, i.e. a voltage source depending on time.
\begin{figure}[ht]
	\centering
		\includegraphics[width=5.1cm]{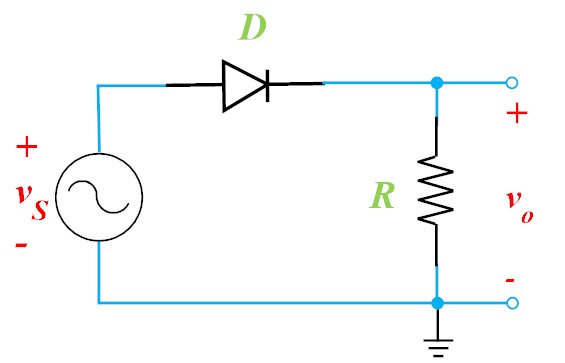}
		$ ~~~~~~~~~~$
		\caption{A circuit with AC voltage source}
		\label{ac}
\end{figure}

Mostly, such voltage sources would be of the form $v_s(t) = A \,\sin(t - t_0)$. To indicate the time-depending nature of this circuit, not only we use $(t)$ in front of the functions, but also use small letters to refer to a voltage or current. \\
We can easily obtain the generalized equation in this case as $ v_s(t) \in R z(t) \, + F( z(t)) $.\\
In complicated circuits we may have 
\begin{equation}\label{time-varying-case-first-formulation}
\begin{matrix}
\Phi (z(t)) : =f(z(t)) \, + B F(Cz(t)) \\
S(p) := \left \{ z ~|~ p \in \Phi (z)  \mfa t \in [a, \, b] ~\right \}.~
\end{matrix}
\end{equation}
\end{eg}
So this case may reduce to solving the previous case (in which $p$ was a fixed vector in $\R^n$) for any $t $ in an interval $ [a, b]$. We will study this type of problem in Chapter \ref{Chapter4}.

\begin{eg}[AC Source with a Capacitor; A Rectifier Circuit]
Using AC voltage sources will allow us to use passive components apart from resistors. In Figure \ref{DI}, we can see another simple circuit with a diode and a capacitor, and the voltage source is a sinusoid function. 
\begin{figure}[ht]
	\centering
		\includegraphics[width=5.1cm]{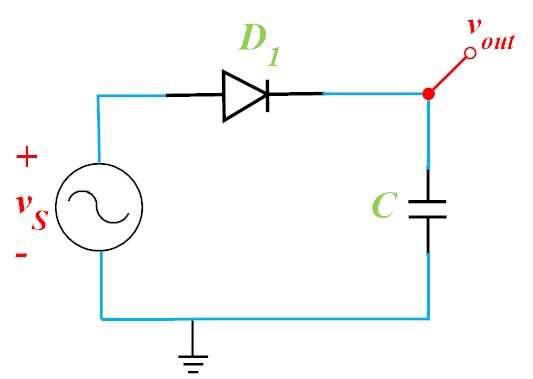}
		\caption{\footnotesize A simple AC to DC Rectifier circuit}
		\label{DI}
\end{figure}

Although we tried so much to avoid encountering capacitors and inductors till now, in almost any practical circuit one will find these components (mostly capacitors). So it is better to show in a simple circuit why we have ignored them up to now. \\
The $i - v$ relation of the capacitor is $ i_c = C \, \dfrac{dv_c}{dt} $, in which $C > 0 $ is a constant called capacitance, and the $i - v$ relation of the inductor is $ v_L = L \, \dfrac{di_L}{dt} $, in which $L> 0 $ is a constant called inductance.
\begin{equation}\label{dip}
\left.\begin{matrix}
\mathrm{KVL: }  -v_s + v_{D_1}+ v_c = 0\\ 
\mathrm{KCL: }~~~ i_{D_1} = i _c =: i~~~~~~~\\
~~i_c = C \, \dfrac{dv_c}{dt}\\ 
~~~~~~\, V_{D_1} \in F( i_{D_1})\\
\end{matrix}\right\}~~\Rightarrow
v_s(t) \in \frac{1}{C} \int_0^t i \, d\tau + F(i)
\end{equation}
\end{eg}

\subsection{Alternative Formulation} \label{Alternative Formulation, Pros and Cons}
In this subsection, we will shortly discuss another suggested formulation for components like diode, and transistor in the literature (see for example \citep{acary, cibulka-note, facchinei2003finite, facchinei2003finite2}). We would explain why the setting of generalized equation (GE) is more appropriate for our study and how it covers more circuits.\\
Although one can follow the process used in the previous section, the special $i-v$ characteristic of the ideal diode give us a chance to gain another formulation. Let us consider the circuit in Figure \ref{Zener Diode simple circuit} with ideal diode instead of a Zener diode and observe that at any point $ (i_D, v_D ) \in \gph{ F}$ (where $F$ is the set-valued map relating the current passing trough ideal diode to the voltage over it) one can write $v_D \, i_D = 0$. \\
So in this example one can write 
\begin{equation}
\langle {\, p - R z, z \,} \rangle = 0 ,
\end{equation}
where $p = V_I$, $z = i_D$, and $ \inp{\cdot}{ \cdot} $ indicates the scalar product in $\R^n$.
Following this method we will end up at what is called \emph{\textbf{Variational Inequality}}. \\
Not only for Ideal Diodes, but also when dealing with transistors as switches, we can formulate our problem as a VI. 

\section{Examples}
This section will contain some real world examples. Each example will display part of the mathematical formulation difficulties, and will shed light on different aspects of the problem. We will provide more details of these examples in the following chapters, after providing the necessary mathematical tools for this study.

\begin{figure}[ht]
  \centering
\includegraphics[width=10.8cm]{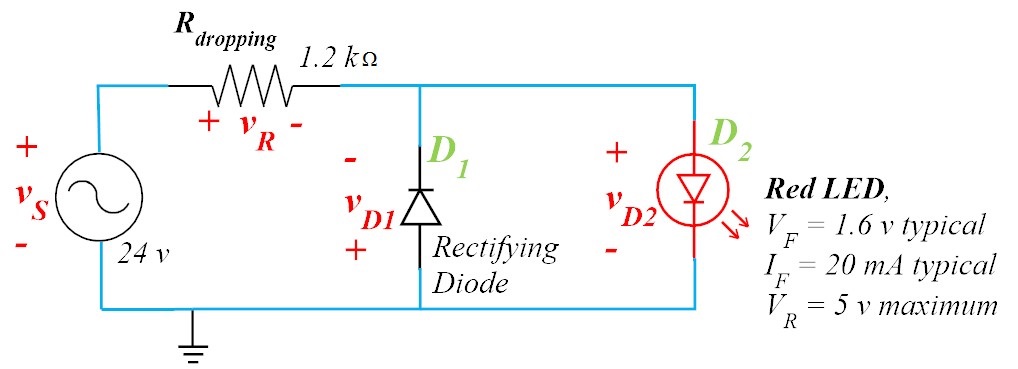}
	\caption{Driving an LED with AC}
	\label{fig: LED}
\end{figure}

\eg \textbf{(Driving an LED with AC, on the Importance of Understanding the Behaviour of the Circuit)} \citep{LED} \label{Driving an LED with AC, on the importance of understanding the behaviour of the circuit}\\
Figure \ref{fig: LED} depicts an LED paralled backward with a rectifying diode. 
Since LEDs are made of different chemical substances than silicon diodes, their forward voltage drops will be different. Typically, LEDs have much larger forward voltage drops than rectifying diodes, anywhere from about $1.6 \,volts$ to over $3 \,volts$, depending on the color. Typical operating current for a standard-sized LED is around $20 \,mA$. When operating an LED from a DC voltage source greater than the LED’s forward voltage, a series-connected ``dropping'' resistor must be included to prevent full source voltage from damaging the LED. \\
Also because of their unique chemical makeup, LEDs have much, much lower peak-inverse voltage (PIV) ratings than ordinary rectifying diodes. A typical LED might only be rated at 5 volts in reverse-bias mode. Therefore, when using alternating current to power an LED, connect a protective rectifying diode anti-parallel with the LED to prevent reverse breakdown every other half-cycle as in Figure \ref{fig: LED}.\\
Now, we try to give the appropriate mathematical model for this circuit; which one can guess such a model could not have ideal diode model, or same set-valued map for both diodes. Let $D_1$ indicates the protective rectifying diode, and refer to LED with $D_2$. KVL, KCL, and Diode Characteristics are as follows: 
\begin{eqnarray*}
\left \{ \begin{matrix}
- v_s + v_R - v_{D_1} = 0, \\
i_R = - i_{D_1} + i_{D_2}
\end{matrix} \right.
, \hspace*{2.5cm}
\left \{ \begin{matrix}
v_{D_1} = - v_{D_2},~~~~~~~~~~~~~~~\\
v_{D_j} \in F (i_{D_j} ), ~~~ j = 1,2.
\end{matrix} \right.
\end{eqnarray*}
Let 
$ x = \begin{pmatrix}
v_{D_1} \\ 
v_{D_1}
\end{pmatrix} $,
and 
$ u = \begin{pmatrix}
i_{D_1}\\ 
i_{D_2}
\end{pmatrix} $. 
Then from the above relations , we get
\begin{eqnarray*}
\begin{pmatrix}
x_1 \\ 
x_2
\end{pmatrix} 
= 
\begin{pmatrix}
-1 \\ 
+1
\end{pmatrix} v_s +
\begin{pmatrix}
-R & R \\ 
R & -R
\end{pmatrix} 
\begin{pmatrix}
u_1 \\ 
u_2
\end{pmatrix} \\ \\
\mathrm{with~} x_j \in F_j (u_j). ~~~~~~~~~~~~~~~~~~
\end{eqnarray*}
So, one can form the generalized equation as follows:
\begin{equation}
B + A U \in F (U),
\end{equation}
where $ B = 
\begin{pmatrix}
- v_s \\
+ v_s
\end{pmatrix}
$, $ A= 
\begin{pmatrix}
-R & ~R \\ 
~R & -R
\end{pmatrix}
$, and $ \mmap{F}{2}{2} $ is defined as 
$$ F \begin{pmatrix}
u_1 \\ 
u_2
\end{pmatrix} : = \prod_{i=1}^{2} F_i (u_i)
 $$
In Chapter \ref{Chapter3}, we would see why this formulation cannot be written as a variational inequality (see Example \ref{Driving an LED with AC-model}).

\begin{eg} \textbf{(A Simple Circuit with DIAC)} \label{DIAC-Chap2} \\
Let us consider the simple circuit in Figure \ref{fig: DIAC-circuit-chap02} with a DIAC, whose $ i - v $ characteristic for given $ V > 0 $, and $ a > 0 $ is as follows:
\begin{figure}[ht]
  \begin{minipage}{.42\textwidth}
    \begin{eqnarray*}
		F_D (i_D) :=
		\left \{ \begin{matrix}
		\dfrac{-V}{\sqrt{1 - \frac{2ai_D}{V}}} & & i_D < 0,\\[2em]
		[- V, V] & & i_D = 0,\\[0.5em]
		\dfrac{V}{\sqrt{1 + \frac{2a i_D}{V}}} & & i_D > 0.
		\end{matrix} \right.
	\end{eqnarray*} 
  \end{minipage}
  \begin{minipage}{.60\textwidth}
    \centering
		\includegraphics[width=8.8cm]{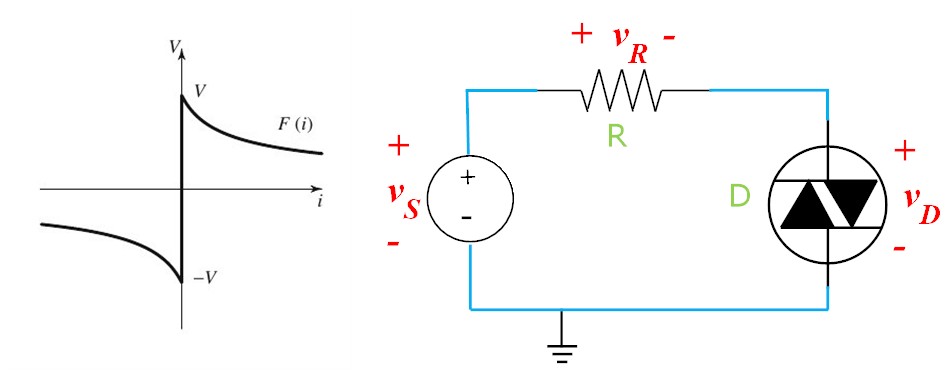}
  \end{minipage}
\caption{A simple circuit with DIAC}  
\label{fig: DIAC-circuit-chap02}
\end{figure}

Using KVL, KCL, and the characteristic relations we obtain that
\begin{eqnarray*}
\left \{ \begin{array}{l}
- v_s + v_R + v_{D} = 0, \\
i_R = i_{D} = i.
\end{array} \right.
\hspace*{2.5cm}
\left \{ \begin{array}{l}
v_{R} = R \, i_R\\
v_{D} \in F_D (i_{D} ).
\end{array} \right.
\end{eqnarray*}
By replacing $ p = v_s $, and $ z= i $, we get the generalized equation
\begin{equation} \label{DIAC-eg-chap02}
p \in f(z) + F (z),
\end{equation}
where $ f(z) = R z $, and $F(z) = F_D (z) $ for every $ z \in \R $. For each $ p \in \R $, one can look at the solution mapping associated to \eqref{DIAC-eg-chap02} as the intersection of the line 
$ \varphi (z) = p $ and the set-valued mapping $ \Phi = f + F $. This interpretation would provide a general idea about the solution sets $S(p)$. Figure \ref{fig: DIAC-circuit02-chap02}, shows the graph of $ \Phi$ for different values of parameters $R$ and $a$.

\begin{figure}[ht]
  \centering
\includegraphics[width=7.1cm]{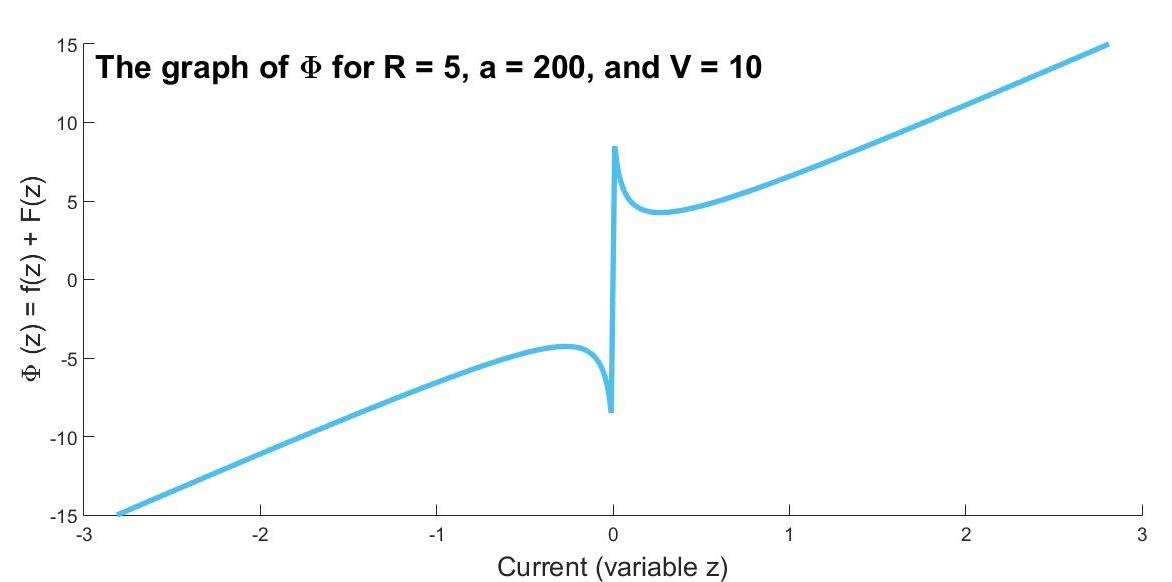}
$~~$
\includegraphics[width=7.1cm]{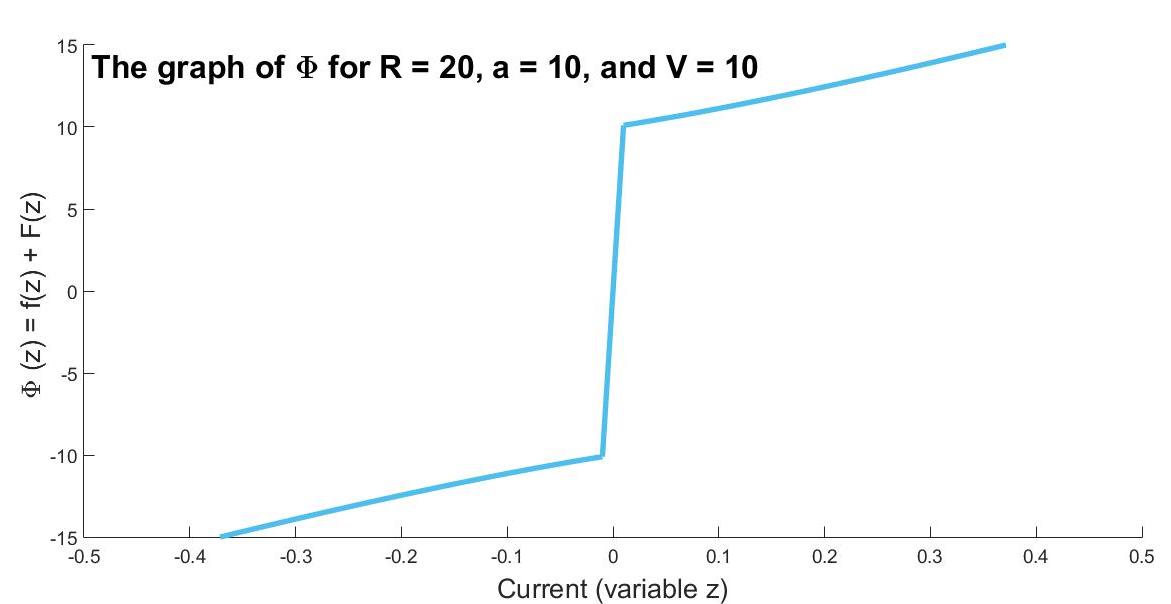}
	\caption{The graph of $\Phi$ for different values of parameters $R$, $a$ and $V$\\
	(left: $a = 200, R = 5, V = 10 $;~ right: $a = 10, R = 20, V = 10 $) }
	\label{fig: DIAC-circuit02-chap02}
\end{figure}

Note that for each $ p \in \R $ there is a solution $ z \in \R $ to the equation \eqref{DIAC-eg-chap02}. Indeed, as $ \Phi (z) = Rz + F(z)$ for $ z \in \R $, we have 
$ \mathrm{rge}\, \Phi = \R $. Moreover, this solution is unique if $ a \leq R $, since for a non-zero $ z \in \R $ one infers that 
\begin{equation*}
 \Phi'(z) = R - \dfrac{a}{\paren{ \dfrac{2 a |z|}{V} + 1 } ^{\frac{3}{2}}} = \dfrac{\, R \paren{ \dfrac{2 a |z|}{V} + 1 } ^{\frac{3}{2}} - a }{\paren{ \dfrac{2 a |z|}{V} + 1 } ^{\frac{3}{2}}} > \dfrac{R - a }{\paren{ \dfrac{2 a |z|}{V} + 1 } ^{\frac{3}{2}}} \, \geq \, 0.
\end{equation*}
In this case, the solution equals to zero if $ p \in [-V, V] $; it is positive if $ p > V $, and negative when $ p < - V $.\\
Clearly, we can split $\Phi$ into different functions $f$ and $F$ which produce the same inclusion \eqref{DIAC-eg-chap02}. Namely, from now on, we assume that
\begin{eqnarray*}
f(z) :=
\left \{ \begin{matrix}
R z - \dfrac{V}{\sqrt{1 - \frac{2az}{V}}} + V & & z < 0,\\
R z + \dfrac{V}{\sqrt{1 + \frac{2az}{V}}} - V & & z \geq 0,
\end{matrix} \right.
\mathrm{~~~~and~~~~}
F(z) :=
\left \{ \begin{matrix}
-V ~~~ & &  z < 0,\\
[- V, V] & & z = 0,\\
V ~~~ & & z > 0.
\end{matrix} \right.
\end{eqnarray*} 
This is the \emph{simplification idea} we discussed before in this chapter. Obviously, this new function $f$ does not represent an electronic component but the sum $ f + F $ and thus, the solution mapping remain unchanged. Hence, from the analytic point of view, we are allowed to do this change of functions in the sum.  
$$  f'(z) = R - \dfrac{a}{\paren{ \dfrac{2 a |z|}{V} + 1 } ^{3/2}} \mathrm{~~whenever~~} z \in \R. $$
Indeed, $ f'(z) = \Phi'(z) $ if $ z \neq 0 $, and $ f'(0) = R - a $ could be obtained from the following simple calculations:
\begin{equation*}
\lim_{h \longrightarrow 0^{+}} \dfrac{R h + \frac{V}{\sqrt{1 + \frac{2ah}{V}}} - V}{h} = \lim_{h \longrightarrow 0^{+}} \paren{R - \dfrac{a}{\paren{\frac{2 a h}{V} + 1 } ^{3/2}}} = R - a.
\end{equation*}
The case $ a > R $ could be discussed similarly, but the solution is not necessarily unique any more (see Figure \ref{fig: DIAC-circuit02-chap02}). In Chapter \ref{Chapter3}, we will examine the Aubin property and the isolated calmness of the solution mapping at various reference points (refer to Examples \ref{DIAC-example} and \ref{DIAC-example-IC}). \\
Apart from the electrical importance, this examples describes a simple situation in which the solution set $S(p)$ is not a singleton for some $ p \in \R $, even if we approximate the set-valued part of the sum $ f + F $ with a line from $ - \varepsilon $ to $ + \varepsilon $ (for a small $\varepsilon > 0$, and with a steep slope $\frac{V}{\varepsilon}$).\\
We can also observe the importance and ease of dealing with the metric regularity of $ \Phi $ instead of working with the Aubin property of $S$. Although they are describing the same thing, which is the behaviour of the system, it is much easier to obtain the graph of $\Phi$.
\end{eg}

\begin{eg} \textbf{(A Circuit with SCR and Zener Diode)} \label{SCR and Zener Diode-Chap2} \\
Consider the circuit in Figure \ref{fig: SCR and the Zener Diode-chap2} involving a load resistance $ R > 0 $, two bias sources $ E_2 > E_1 > 0 $, an input-signal source $ u $ with corresponding instantaneous current $i$, and two non-smooth elements: the SCR and the Zener Diode.

\begin{figure}[ht]
	\centering
		\includegraphics[width=14.8cm]{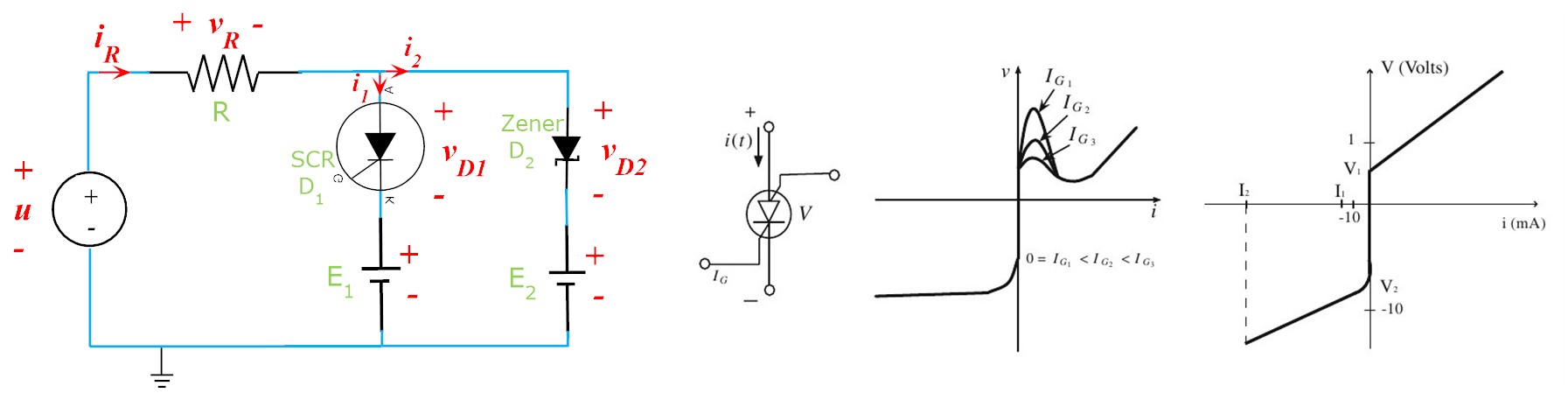}
	\caption{A circuit with SCR and Zener Diode}
	\label{fig: SCR and the Zener Diode-chap2}
\end{figure}

Suppose that $ V > 0, V_1 < 0, \alpha > 0, a > 0 $ and $ b > 0 $ are given numbers, and that $ \varphi : \R \longrightarrow (0, \infty) $ is a continuously differentiable function with $ \varphi(\alpha) < \varphi(0), \varphi'(0) > 0 $ and $ \varphi'(\alpha) > 0 $. Suppose that the $ i - v $ characteristics $ F_1 $ of SCR and $ F_2 $ of the Zener diode are defined by
\begin{eqnarray*}
F_1 (z) :=
\left \{ \begin{array}{lcl}
a z + V_1  & & z < 0, \\
$[$ V_{1}, \, \varphi(0) $]$  & &  z = 0, \\
\varphi(z) & & z \in [0, \alpha], \\
a (z - \alpha) + \varphi (\alpha) & & z > \alpha.
\end{array} \right.
~~~~
F_2(z) :=
\left \{ \begin{matrix}
b z - V & & z < 0,\\
[- V, V] & & z = 0,\\
b z + V & & z > 0.
\end{matrix} \right.
\end{eqnarray*}
Using KVL, KCL, and the characteristic relations we obtain that
\begin{eqnarray*}
\left \{ \begin{array}{l}
- u + v_R + v_{D_1} + E_1  = 0, \\
- u + v_R + v_{D_2} + E_2  = 0, \\
i_R = i_{D_1} + i_{D_2} .
\end{array} \right.
\hspace*{2.5cm}
\left \{ \begin{array}{l}
v_{R} = R \, i_R\\
v_{D_1} \in F_1 (i_{D_1} ), \\
v_{D_2} \in F_2 (i_{D_2} ).
\end{array} \right.
\end{eqnarray*}
By replacing $ p_1 = u - E_1 $, $ p_2 = u - E_2$, $ z_1 = i_{D_1} $, and $ z_2 = i_{D_2} $ we would have
\begin{equation} \label{SCR and Zener Diode-Chap2-eq01}
\left \{ \begin{array}{l}
p_1 \in R ( z_1 + z_2 ) + F_1 (z_1), \\
p_2 \in R ( z_1 + z_2 ) + F_2 (z_2),
\end{array} \right.
\end{equation}
which is in the form of the generalized equation
\begin{equation} \label{SCR-Zener-eg-chap02}
p \in f(z) + B F ( Cz ),
\end{equation}
with $ m = n = 2 $, $ B = C = I_2 $, and $ f(z) = A z $ for $ z \in \R^2 $, where 
$ A = \begin{pmatrix}
R & R \\
R & R
\end{pmatrix}$, and $ \mmap{F}{2}{2} $ is defined as 
\begin{equation*}
F \begin{pmatrix}
z_1 \\ 
z_2
\end{pmatrix} : = 
\begin{pmatrix}
F_1 (z_1) \\ 
F_2 (z_2)
\end{pmatrix}.
\end{equation*}
This example represents a more complicated situation with respect to the previous example. Not only there are two different set-valued maps in the model, but also any simplification for one of them, would destroy the symmetry of $f$. Also, despite the parallel connection of diodes, the role of resistor showed clearly in formulation \eqref{SCR and Zener Diode-Chap2-eq01} warns us not to study each diode (and so, each generalized equation) separately. \\
Later, in Example \ref{SCR with Zener-example} we will discuss the Aubin property of the solution mapping and we will see the dependence of this property to the various parameters of the circuit.
\end{eg}

\begin{eg} \textbf{(Sampling Gate)} \label{Sampling Gate-chap02}\\
Consider a particular circuit of the \emph{sample and hold circuits}\footnote{
In electronics, a sample and hold (S$\&$H) circuit is a device that captures (samples) the voltage of a continuously varying analog signal and freezes (holds) its value at a constant level for a specified minimum period of time. Sample and hold circuits and related peak detectors are the elementary analog memory devices. If the output is available during the sampling period, we would have a track and hold (T$\&$H) circuit.
} 
family in Figure \ref{fig: sampling-gate-chap02}, composed of four diodes $ D_1, D_2, D_3, D_4 $ which are controlled symmetrically by gate voltages $ +V_c $ and $ -V_c $, and the control resistors $ R_c > 0 $. The input signal is given by $ v_i $ and the output signal is defined by the voltage $ v_{\mathrm{out}} $ over the load resistor $ R_L > 0 $.
\begin{figure}[ht]
  \begin{minipage}{.40\textwidth}
  		\begin{eqnarray*}
		~ F_D (z) :=
		\left \{ \begin{matrix}
		V_{D1} ~~~~~~ & & z < 0,\\
		[V_{D1} , V_{D2}]  & &  z = 0,\\
		V_{D2} ~~~~~~ & & z > 0.
		\end{matrix} \right.
		\end{eqnarray*} 
		\vspace*{0.2cm}
  \end{minipage}
  \begin{minipage}{.60\textwidth}
    	\centering
		\includegraphics[width=9cm]{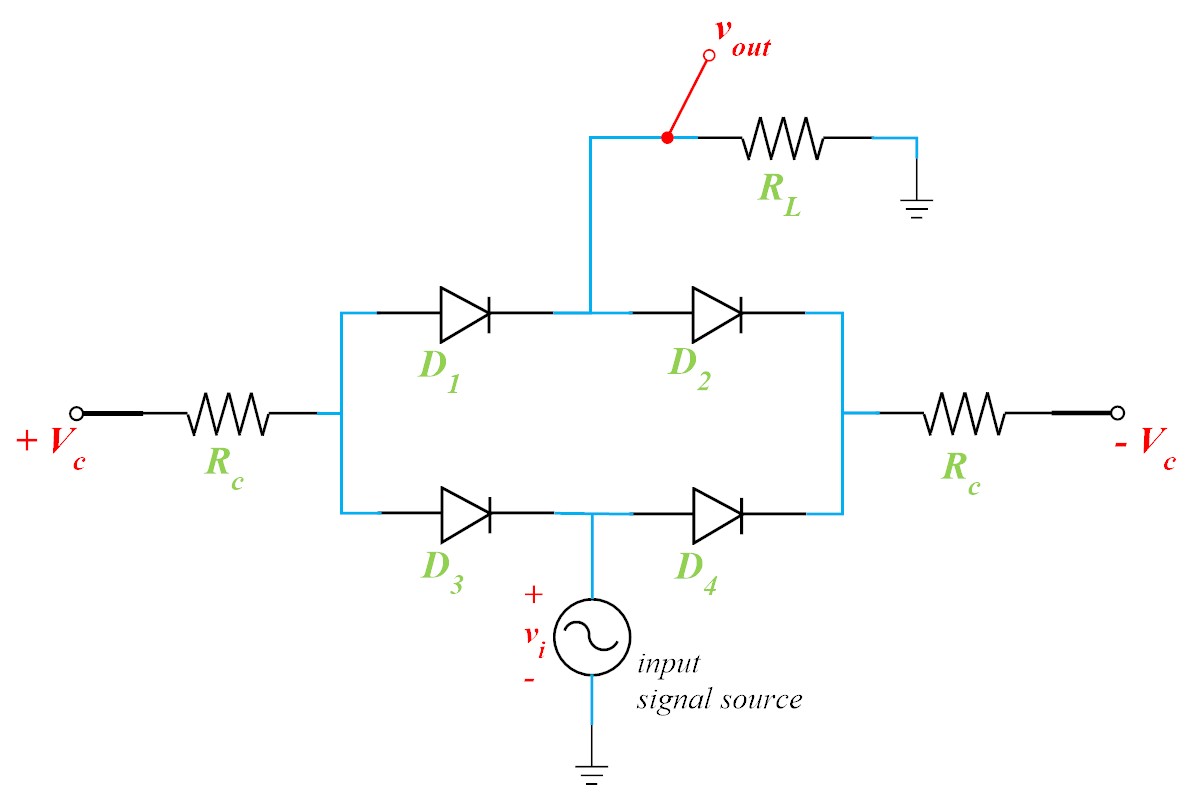}
  \end{minipage}
	\caption{Loops and variables definition in sampling gate circuit}
	\label{fig: sampling-gate-chap02}
\end{figure}

During \emph{sampling}, equal currents passing thorough resistors $ R_c $, flow into the diode bridge, and $ v_{in} $ is copied over $R_L$. During \emph{hold}, no current flows through the bridge. The reason we are interested in this circuit is that we have three independent loops (that is, $ n = 3$), and four diodes (that is, $ m = 4 $), while in all the previous examples we had $ m \, \leq \, n $. Assume that all the diodes have the same characteristics given for $ V_{D1} < 0 < V_{D2} $ by the mentioned set-valued map. \\
One can write the Kirchhoff's laws in the following form (Look at Figure \ref{fig: sampling-gate-02-chap02})
\begin{eqnarray*}
\left \{ \begin{array}{ll}
\mathrm{\scriptstyle{KVL1:}} & - v_i + v_4 - v_2 + R_L \, i_7 = 0, \\
\mathrm{\scriptstyle{KVL2:}} &  V_c - R_c \, i_6 - v_4 - v_3 - R_5 \, i_5 + V_c = 0, \\
\mathrm{\scriptstyle{KVL3:}} & v_1 + v_2 - v_4 - v_3 = 0,
\end{array} \right.
\hspace*{0.3cm}
\left \{ \begin{array}{ll}
\mathrm{\scriptstyle{KCL1:}} & i_5 = i_1 + i_3 = i_2 + i_4 = i_6, \\
\mathrm{\scriptstyle{KCL2:}} & i_1 - i_2 = i_4 - i_3 = i_7,
\end{array} \right.
\end{eqnarray*}
\begin{figure}[ht]
	\centering
		\includegraphics[width=0.84\textwidth]{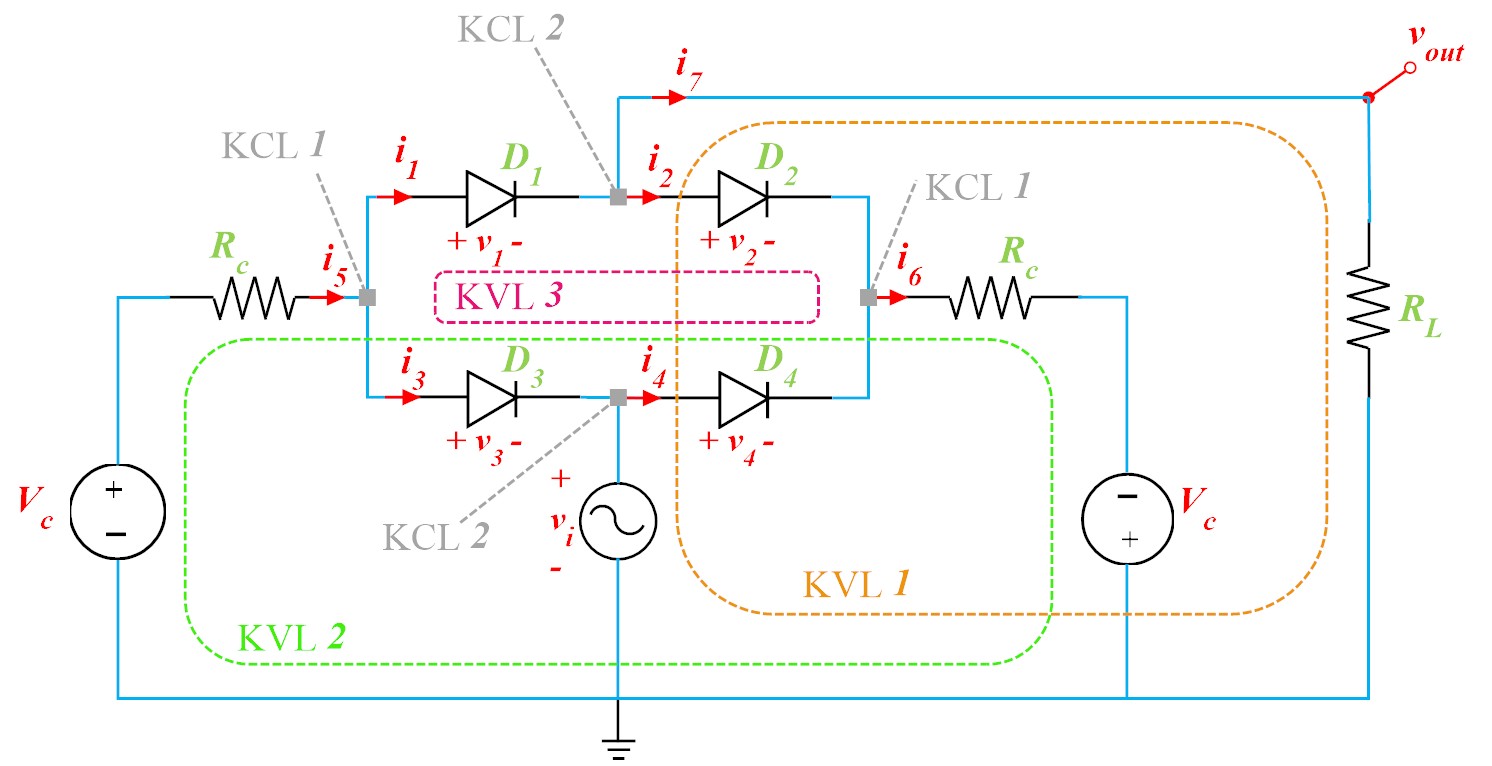}
	\caption{Sampling gate circuit}
	\label{fig: sampling-gate-02-chap02}
\end{figure} 

where $v_j $ is the voltage over the diode $D_j$ for $ j \in \{ 1, 2, 3, 4 \} $. Therefore, by letting $ z = (i_7, i_6, i_1 )^T$, $v = (v_1, v_2, v_3, v_4)^T$, and 
$ u = (v_i, \, 2 V_c)^T$, one can obtain the  linear equation 
$$ Du = A z + B v, $$
where
\begin{equation*}
D = \begin{pmatrix}  
1 & 0 \\
0 & 1 \\
0 & 0
\end{pmatrix}, ~
A = \begin{pmatrix}
R_L & 0 & 0 \\
0 & 2 \, R_c & 0 \\
0 & 0 & 0
\end{pmatrix},~ \mathrm{and~}
B = \begin{pmatrix}
0 & -1 & 0 & 1 \\
0 & 0 & 1 & 1 \\
1 & 1 & -1 & -1
\end{pmatrix}
\end{equation*}
Now, from the diode characteristics we have that $v_1 \in F_D (i_1)$, $v_2 \in F_D (i_2) = F_D (i_1 -  i_7)$, $v_3 \in F_D (i_3) = F_D (i_6 -  i_1)$, and $ v_4 \in F_D (i_4) = F_D (i_7+ i_6 - i_1) $. Hence, we obtain $ v \in F (Cz) $ where
\begin{equation*}
C = \begin{pmatrix}
0 & 0 & 1 \\
-1 & 0 & 1 \\
0 & 1 & -1 \\
1 & 1 & -1
\end{pmatrix} = B^T,~
\mathrm{and~}
F(y) = \begin{pmatrix}
F_D(y_1) \\
F_D(y_2) \\
F_D(y_3) \\
F_D(y_4)
\end{pmatrix} 
\mathrm{~for~} y \in \R^4.
\end{equation*}
So, we arrived at the generalized equation $ p \in f(z) + B F(Cz) $ with $n = 3$, $m = 4$, $ p := D u$, and $ f(z) := A z $. In Chapter \ref{Chapter3} we will investigate the isolated calmness property of the solution mapping at some reference points (ref. Example \ref{IC without A1}).

\end{eg}

%% file: Chapters/Chapter03-a.tex
\chapter{Tools of Variational Analysis for Studying the Local Stability of the Solution Mapping in the Static Case} 
\label{Chapter3} 

\lhead{Chapter 3. \emph{Variational Analysis and Stability of the Solution Mapping}} 
\epigraph{\emph{In mathematics the art of proposing a question must be held of higher value than solving it.}}{\textit{Georg Cantor}}

We start this chapter with introducing some tools from variational analysis. To be more precise, the concept of graphical derivative for set-valued maps is introduced in Section \ref{Variational Geometry}, and some calculus rules for these derivatives is provided in subsection \ref{Calculus Rules}. 
Then we provide two norm-like maps in Subsection \ref{Inner and Outer Norms}, which allow us to characterize the stability-like
properties of a set-valued map in terms of the \emph{inner} and \emph{outer norm} of their graphical derivatives.
Subsection \ref{Subdifferentials} will introduce \emph{subdifferentials} and \emph{generalized Jacobians} for functions that may not be differentiable. These concepts would be used later in this chapter to show a possible way to continue the stability study when dealing with a generalized equation with a non-smooth single-valued part. \\
The remaining sections of this chapter have a common structure. Each section will discuss one of the stability-like properties for the specific generalized equation $ f(\cdot) + B F(C \cdot) \ni p $, with the following assumptions:
\begin{enumerate}[topsep=-1ex, itemsep=-1ex, partopsep=1ex, parsep=1ex, leftmargin = 9ex]
\item[(A1)] $B$ is injective;
\item[(A2)] $\smap{f}{n}{n}$ is a continuously differentiable function; 
\item[(A3)] $\mmap{F}{m}{m}$ is a set-valued map with closed graph;
\item[(A4)] $C$ is surjective; and 
\item[(A5)] there are $ \mmap{F_i}{}{},~ i \in \{ 1,... \, , m \} $ such that $ F(x) = \prod\limits_{i=1}^{m} F_i(x_i) $ whenever $ x = (x_1, ..., x_m)^T \in \R^m $.
\end{enumerate}
We try to use the norm characterization obtained for an arbitrary set-valued mapping in Subsection \ref{Inner and Outer Norms} for each property, and calculus rules of Subsection \ref{Calculus Rules} to go step by step toward a pointwise easy-to-check criteria for the local stability of the solution mapping.\\
The main theorem of each section uses only assumptions (A1) - (A3). \\
Among the first three assumptions, (A3) is not such a strong requirement and holds for the $i - v$ characteristics of semiconductors like diodes in our study. In Subsection \ref{Isolated Calmness Without Injectivity Assumption} we try to consider the case where $B$ is not injective, and use the following condition instead of (A1) to provide some similar statements for isolated calmness in this case.
\begin{itemize}[nolistsep]
\item[]$(\widetilde{\mathrm{A1}})$ Suppose that there is $ \bar{v} \in F(C\bar{z}) $ such that
\begin{equation*}
\bar{p} = f(\bar{z}) + B \, \bar{v} ~\mathrm{~~and~~}~ \overline{ \bigcup_{t > 0} \dfrac{ \mathrm{rge~} F_C - \bar{v} } {t} } \,  \bigcap \mathrm{~ker~} B = \{ 0_{\R^m} \}.
\end{equation*}
\end{itemize}
The process is almost the same as before. \\
In Subsection \ref{Strong Metric (Sub-) Regularity with a Nonsmooth Single-valued Part} we focus on functions $f$ which are not smooth enough to satisfy (A2). A possible approach would be considering the generalized Jacobians and replace (A2) with
\begin{itemize}[nolistsep]
\item[]$(\widetilde{\mathrm{A2}})$ $f$ is locally Lipschitz continuous on $\R^n$.
\end{itemize}
In Section \ref{Results in terms of Metric Regularity}, using the relations between the regularity terms and local stability properties of the inverse map expressed in Section \ref{section1.3}, we restate the results of previous sections in terms of metric regularities in Theorems \ref{SMR-Theorem 2}, and \ref{Summary of MR Criteria for GE in Case of a Nonsmooth Single-valued Part}. 

\section{Variational Geometry} \label{Variational Geometry}

In this section we will first introduce some concepts that let us define derivatives of a set-valued map in a proper way. From elementary calculus, we know that the derivative of a function $f$ at a certain point $\bar{x}$, could be interpreted as the slope of the tangent line to the graph of $f$ at the point $ ( \bar{x}, f(\bar{x}) ) $. This tangent line could be then used to approximate the function in a vicinity of the reference point. \\
The idea is quite the same for set-valued maps. We will define approximating cones to the graph of a set-valued map at a certain point, and use this \emph{graphical approximation} to introduce the derivatives. \\
Then, in Subsection \ref{Calculus Rules} some simple calculus rules for the newly defined derivatives will be given. Subsection \ref{Inner and Outer Norms}, will bridge the local stability concepts of Chapter \ref{Chapter1} and derivative concepts of this chapter. \\
Finally, we will look at the non differentiable (in the standard common sense) functions and introduce another graphical-based concept that could be considered as a replacement of the non-existing derivative, in Subsection \ref{Subdifferentials}.

\defn \textbf{(Tangent Cones)} \hfill \\ \citep[p. 132]{Schirotzek} \citep[p. 69]{Mordukhovich} \citep[p. 162]{Aubin} Let $\Omega$ be an arbitrary non-empty subset of $\R^d$ containing a point $\bar{x}$. 
\begin{enumerate} [topsep=-1ex, itemsep=-1ex, partopsep=0ex, parsep=0ex]
\item[\textbf{(a)}] The \emph{Bouligand-Severi tangent/contingent cone}\index{cone ! contingent cone} $ T(\bar{x};\Omega)$ to $\Omega$ at $\bar{x}$ contains those $u \in \R^d$ for which there are sequences $(t_k)_{k \in \N} $ in $(0, \infty)$ converging to $0$, and $(u_k)_{k \in \N} $ in $\R^d$ converging to $u$, such that $ \bar{x} +t_k u_k \in \Omega $ whenever $k \in \N$;
\begin{equation}
 T(\bar{x};\Omega) := \left \{ u \in \R^d ~|~ \exists \, t_k \downarrow 0, \, u_k \rightarrow u \mathrm{~~~with~~} \bar{x} + t_k u_k \in \Omega \, \right \}.
\end{equation}

\item[\textbf{(b)}] the \emph{Bouligand paratingent cone}\index{cone ! paratingent cone} $\widetilde{T}(\bar{x};\Omega)$ to $\Omega$ at $\bar{x}$ contains those $u \in \R^d$ for which there are sequences 
$ (t_k)_{k \in \N}$ in $(0,\infty)$ converging to $0$, $(u_k)_{k \in \N} $ in $\R^d$ converging to $u$, and $(x_k)_{k \in \N}$ in $\Omega$ converging to $\bar{x}$, such that $x_k + t_k u_k \in \Omega$ whenever $k \in \N$;
 \footnote{One should be aware that this definition is different from \emph{Clarke tangent cone}, \index{Clarke ! Clarke tangent cone}\index{cone ! Clarke tangent cone}shown by $ T_C (\bar{x}; \Omega) $ and defined as 
\begin{equation*}
 T_C (\bar{x};\Omega) := \left \{ u \in \R^d ~|~  \boldsymbol {\forall} \, x_k \overset{\Omega}{\longrightarrow} \bar{x}, \,  \boldsymbol {\forall} \,  t_k \downarrow 0, \, \exists \, u_k \rightarrow u, \mathrm{~~~with~~} x_k + t_k u_k \in \Omega ~ \forall k \, \right \}.
\end{equation*}
}
\begin{equation}
 \widetilde{T}(\bar{x};\Omega) := \left \{ u \in \R^d ~|~ \exists \, t_k \downarrow 0, \, u_k \rightarrow u, \, x_k \overset{\Omega}{\longrightarrow} \bar{x} \mathrm{~~~with~~} x_k + t_k u_k \in \Omega \, \right \}. \\
\end{equation}
\end{enumerate}

After reviewing these definitions in the following graphical example, we will provide some general properties of tangent cones.

\begin{figure}[ht]
	\centering
		\includegraphics[width=0.30\textwidth]{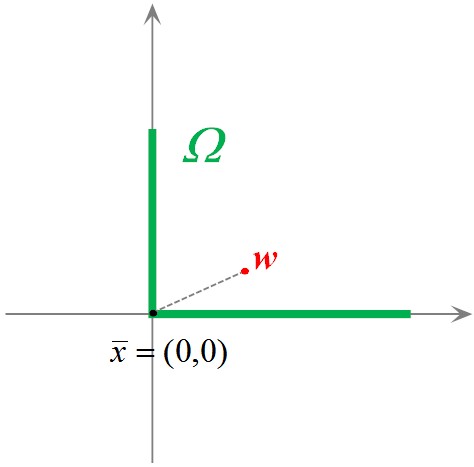}
	\caption{Calculating contingent and paratingent cones for a set}
	\label{fig: contingent and paratingent cones}
\end{figure}

\eg \label{tangent cones examples}
We try to calculate the contingent and paratingent cones for two different sets in $\R^2$ at the reference point $ (0,0) $. \\ 
\textbf{(a)} Consider the set $\Omega = \left \{ (x,y) \in \R^2 _{+} ~|~ xy = 0  \right \} $, as shown in Figure \ref{fig: contingent and paratingent cones}. 

To calculate the contingent cone, we must find points $ u \in \R^2 $, such that there exist sequences $ t_k \downarrow 0 $, and $ \, u_k \rightarrow u $ with 
$  t_k u_k \in \Omega $. 
Thus, any point $ u \in \Omega $ would belong to the contingent cone (for example let $ u_k = u $, and $ t_k = \frac{1}{k} $). Any other point like $w$ shown in the figure does not belong to 
$T \big( (0, 0); \Omega \big) $, since for any sequence $(w_k) \rightarrow w$, no matter how small the $t_k > 0 $ would be, $ t_k w_k$ would be out of $\Omega $. Hence,
\begin{equation*}
\begin{small}
T \Big( (0, 0); \Omega \Big) = \R_{+} \begin{pmatrix}
0 \\ 1 
\end{pmatrix} \, \bigcup \, \R_{+} \begin{pmatrix}
1 \\ 0
\end{pmatrix}.
\end{small}
\end{equation*}

In order to compute $ \widetilde{T} \big( (0,0); \Omega \big) $, first observe that by letting $ x_k = x $, we get that
$ T \big( (0, 0); \Omega \big)  \subset \widetilde{T} \big( (0,0); \Omega \big) $. We argue the set 
$ {\tiny \R_{+} \begin{pmatrix}
-1 \\ 0
\end{pmatrix}} $ also belongs to the paratingent cone at $ (0, 0)$. \\
Let $ u = ( \lambda, 0 )^T $ with $ \lambda < 0 $ be an arbitrary point. Then, by assuming $ u_k = u $, $ t_k = \frac{1}{ - \lambda k } $, and $ x_k = (\frac{2}{k}, 0 )^T $ for every $ k \in \N $, we get $ x_k + t_k \, u_k = (\frac{1}{k}, 0)^T \in \Omega $. \\
A similar argument would prove that $ {\tiny \R_{+} \begin{pmatrix}
0 \\ -1
\end{pmatrix} \subset  \widetilde{T} \big( (0,0); \Omega \big) }$. Hence,
\begin{equation*}
\begin{small}
\widetilde{T} \Big( (0, 0); \Omega \Big) = \R \begin{pmatrix}
0 \\ 1 
\end{pmatrix} \, \bigcup \, \R \begin{pmatrix}
1 \\ 0
\end{pmatrix}.
\end{small}
\end{equation*}
One may consider if $\Omega$ contains some points between the two lines, the situation would be different; in fact, it is true. 

\begin{figure}[ht]
	\centering
		$~$
		\includegraphics[width=0.45\textwidth]{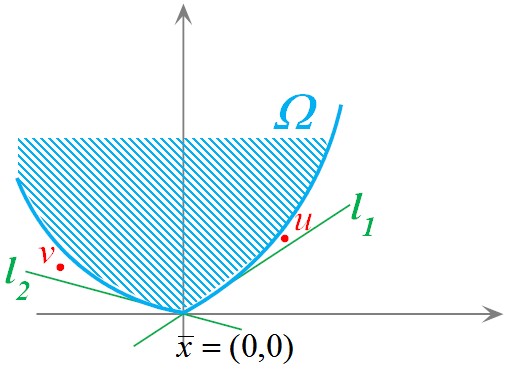}
		$~~~~~~~~$
		\includegraphics[width=0.45\textwidth]{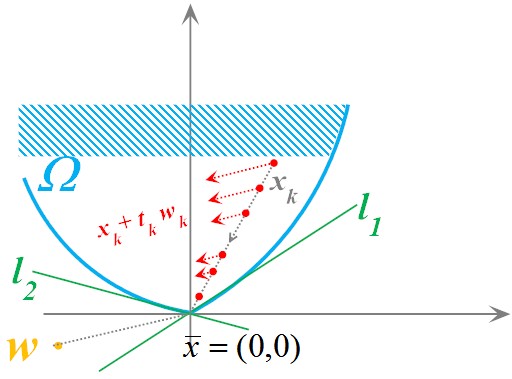}
	\caption{Calculating contingent and paratingent cones for a convex set}
	\label{fig: contingent and paratingent cones2}
\end{figure}

\textbf{(b)} This time we consider the convex set $\Omega$ in Figure \ref{fig: contingent and paratingent cones2} (left). 

Consider $ l_1 $, and $l_2$ to be the \lq\lq tangent lines\rq\rq \,to $\Omega$ at $\bar{x} =  (0, 0)$. Not only all points in $\Omega$ are inside the contingent cone, but also points like $u$, and $v$ (shown in the figure) belong to $T(\bar{x}; \Omega)$. The convexity of $ \Omega $, enables us to find a sequence $u_k \rightarrow u$, and a sequence of numbers $t_k $ such that $t_k u_k $ is small enough to fit inside $\Omega$, near $\bar{x}$. \\
Thus, the contingent cone to $\Omega$ at $(0,0)$ would be the whole area between line segments $l_1$, and $l_2$. Assuming the slope of these lines to be $m_1$, and $m_2$, respectively, we can write
\begin{equation*}
\begin{footnotesize}
T \Big( (0, 0); \Omega \Big) = \R_{+} \begin{pmatrix}
1 \\ m_1 
\end{pmatrix} \, \bigcup \, \R_{+} \begin{pmatrix}
1 \\ m_2
\end{pmatrix} \, \bigcup \, \mathrm{~cone~} \begin{Bmatrix}
\begin{pmatrix}
1 \\ m_1 
\end{pmatrix}, \,
\begin{pmatrix}
1 \\ m_2
\end{pmatrix} 
\end{Bmatrix}
\end{footnotesize}.
\end{equation*}
It is expected to assume that the paratingent cone will include the negative part of this cone, too. What is more, the paratingent cone is even bigger. \\
For example, consider the point $ w $ (shown in the figure). Then, one can take a sequence $ x_k \in \mathrm{\,int\,} \Omega $ converging to $ (0,0) $ slower than the convergent of $ w_k t_k $ to $ (0,0) $ (a possible choice is shown in Figure \ref{fig: contingent and paratingent cones2} (right)). Hence, $ x_k + t_k w_k $ remains inside the $\Omega$, and this reasoning reveals that 
$ \widetilde{T} \big( (0, 0); \Omega \big) = \R^2 $.

\begin{prop} [\textbf{Some Properties of the Contingent Cone}] \citep[p. 232]{Schirotzek} \\
Let $\Omega$ be a non-empty subset of $\R^d$, and consider a point $\bar{x} \in \Omega$. Then the following properties for the contingent cone hold:
\begin{itemize}[topsep=-1ex, itemsep=-1ex, partopsep=1ex, parsep=1ex, leftmargin=7ex]
\item[\textbf{(a)}]
$T(\bar{x}; \Omega)$ is a cone, containing the zero element
\item[\textbf{(b)}]
$T(\bar{x}; \Omega)$ is closed;
\item[\textbf{(c)}]
if $U$ is a neighborhood of $\bar{x}$, then $T(\bar{x}; \Omega) = T(\bar{x}; \, \Omega \cap U )$;
\item[\textbf{(d)}]
if $\Omega$ is convex, then $ T(\bar{x}; \Omega) = \mathrm{cl} \, (\R_+ (\Omega - \bar{x}) )  $.
\end{itemize}
\end{prop}

Statement \textbf{(c)} of the above proposition means that the contingent cone depends on the local properties of $\Omega$ near $\bar{x}$ only.\\
It is time to define the \lq\lq dual\rq\rq \,concept for tangent cones. 

\defn \textbf{(Normal Cones)} \label{normal cones} \hfill \\ \citep[p. 4]{Mordukhovich}\citep[p. 229]{implicit}
Let $\Omega $ be a non-empty subset of $ \R^n $ and consider a point $ \bar{x} \in \Omega $, at which $\Omega$ is locally closed.
\begin{enumerate}[topsep=0pt,itemsep=-1ex,partopsep=1ex,parsep=1ex]
\item[\textbf{(a) }] The \emph{Fr\'echet/regular normal cone}\index{cone ! regular normal cone} to $\Omega$ at $\bar{x}$ is the set  
\begin{equation}
\widehat{N}(\bar{x}; \Omega) := \left \{ v \in \R^n ~|~ \limsup_{x \overset{\Omega}{\longrightarrow}\bar{x}} \dfrac{\inp{v}{x - \bar{x}}}{\left \| x - \bar{x} \right \|} \leq 0 \right \} 
\end{equation}
in which $x \overset{\Omega}{\longrightarrow}\bar{x}$ means $x \longrightarrow \bar{x}$ with $x \in \Omega$.\\

\item[\textbf{(b) }] The \emph{Mordukhovich/general/limiting normal cone}\index{cone ! limiting normal cone} $ N(\bar{x};\Omega)$ to $\Omega$ at $\bar{x}$ contains those $ v \in \R^n$ for which there are sequences $ (x_k)_{k \in \N}$ in $\Omega$ converging to $\bar{x}$, and $(v_k)_{k \in \N}$ in $\R^n$ converging to $v$ such that $v_k \in \widehat{N}(x_k;\Omega)$ for each $k \in \N$;
\begin{equation}
N(\bar{x};\Omega) := \left \{ v \in \R^n ~|~ \exists \, v_k \rightarrow v, ~ x_k \overset{\Omega}{\longrightarrow}\bar{x} \mathrm{~~~with~~}  v_k \in \widehat{N}(x_k;\Omega) \right \}.
\end{equation}
\end{enumerate}

\rem \label{Clarke regularity} 
\textbf{(a)} We define $ \widehat{N}(x;\Omega) = \emptyset $ for any $ x \not \in \Omega$. \\
\textbf{(b)} Very often, the limit process in the definition of the general normal cone $N(\bar{x};\Omega)$ is superfluous; no additional vectors $v$ are produced in that manner, and one merely has $N(\bar{x};\Omega) = \widehat{N}(\bar{x};\Omega) $. This circumstance is termed as the \emph{Clarke regularity} 
\footnote{
\textbf{Definition (Clarke regularity of sets).} \index{Clarke ! Clarke regularity of sets}A set $C \subset \R^n $ is \emph{regular} at one of its points $\bar{x}$ in the sense of \emph{Clarke} if it is locally closed at $\bar{x}$ and every normal
vector to $ C $ at $\bar{x}$ is a regular normal vector, i.e., $N(\bar{x}; C) = \widehat{N}(\bar{x}; C)$
\citep[p. 199]{Rockafellar1998variational}.} of $\Omega$ at $\bar{x}$. \\
Anyway, $N(\bar{x};\Omega)$  is always a closed cone.

\eg \label{normals example}
Let us consider the set $\Omega = \left \{ (x,y) \in \R^2 _{+} ~|~ xy = 0  \right \} $. In Example \ref{tangent cones examples} we calculated the contingent and paratingent cones at $\bar{x} =  (0,0)$. Now we take a look at the normal cones to $\Omega$ at $\bar{x}$.
\begin{figure}[ht]
  \centering
$~$
\includegraphics[width=5.7cm]{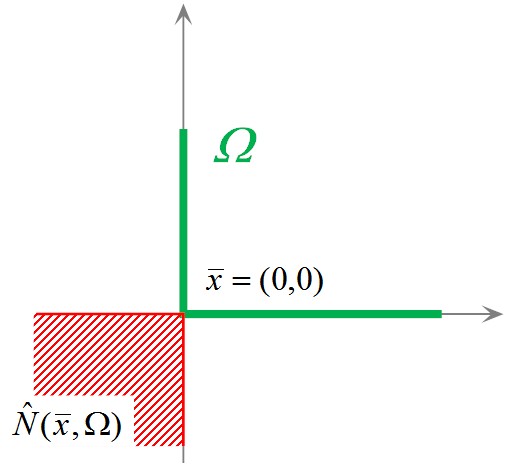}
$ ~~~~~ $
\includegraphics[width=6.8cm]{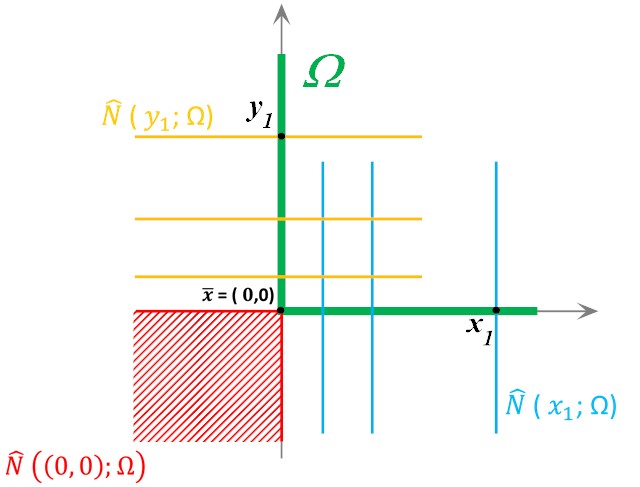}
	\caption{\footnotesize Regular and limiting normal cones to $\Omega$ at $\bar{x} = (0,0)$}
	\label{fig: normal cones}
\end{figure}
\begin{enumerate}[topsep=0ex, itemsep=-1ex, partopsep=1ex, parsep=1ex]
\item[\textbf{(a)}] In the left image of Figure \ref{fig: normal cones}, all vectors $v \in \R^2$ for which $\inp{v}{x} \leq 0$ for those $x \in \Omega$ that are near to $(0,0)$, are shown with red. Thus,
\begin{equation*}
\begin{footnotesize}
\widehat{N} \Big( (0, 0); \Omega \Big) = \R_{+} \begin{pmatrix}
0 \\ -1 
\end{pmatrix} \, \bigcup \, \R_{+} \begin{pmatrix}
-1 \\ 0
\end{pmatrix} \, \bigcup \, \mathrm{~cone~} \begin{Bmatrix}
\begin{pmatrix}
0 \\ -1
\end{pmatrix}, \,
\begin{pmatrix}
-1 \\ 0
\end{pmatrix} 
\end{Bmatrix}
\end{footnotesize}.
\end{equation*}
\item[\textbf{(b)}] In the right hand side of Figure \ref{fig: normal cones}, one can choose $x_k = \bar{x} ~~ (\forall k)$, to get $\widehat{N}(x_k;\Omega) = \widehat{N}((0,0);\Omega) $.
Another option is moving along the horizontal line in $\Omega$, that will produce blue lines which converge to the vertical line passing through origin, as $ k \rightarrow \infty $.\\
Moving along the vertical line in $\Omega$, produces normals shown as orange lines which converge to the horizontal line passing through origin. Any other sequence $x_k \in \Omega$, that converges to $\bar{x}$, will produce a normal cone subset of the above ones. Hence,
\begin{equation*}
\begin{footnotesize}
N \Big( (0, 0); \Omega \Big) = \R \begin{pmatrix}
\, 0 \\  1 
\end{pmatrix} \, \bigcup \, \R \begin{pmatrix}
\, 1 \\  0
\end{pmatrix} \, \bigcup \, \mathrm{~cone~} \begin{Bmatrix}
\begin{pmatrix}
0 \\ -1
\end{pmatrix}, \,
\begin{pmatrix}
-1 \\ 0
\end{pmatrix} 
\end{Bmatrix}
\end{footnotesize}.
\end{equation*}
\end{enumerate}

\eg
Given a differentiable convex function $\smap{h}{m}{}$, consider the set $\Omega := \{ u \in \R^m ~|~ h(u) \leq 0 \}$. Then
\begin{equation*}
N(u; \Omega) = 
\left\{\begin{matrix}
~\{ 0 \} ~~~~~~~~~~~~~~~~~~~~~ & \mathrm{~~~if~} h(u) < 0; \\ 
\{ \lambda \nabla h(u) ~:~ \lambda \geq 0 \}  &  \mathrm{~~~if~} h(u) = 0; \\ 
~ \emptyset ~~~~~~~~~~~~~~~~~~~~~ &  \mathrm{~~~if~} h(u) > 0.
\end{matrix}\right.
\end{equation*} 
Indeed, $\Omega$ is a convex set and one can easily deduce the direction of the normal vector at a point in the boundary of $\Omega$ from $\nabla h(\cdot) $ (since $h$ is differentiable). The other two relations are obvious. 

\begin{prop} [\textbf{Some Properties of the Fr\'echet Normal Cone}] \citep[p. 17]{Mordukhovich} \citep[p. 35]{Schirotzek}
Let $\Omega \subset \R^n $ be a non-empty set, and $\bar{x} \in \Omega$. Then, 
\begin{enumerate} [topsep=0.2pt,itemsep=-1ex,partopsep=1ex,parsep=1ex]
\item[\textbf{(a)}] the Fr\'echet normal cone $\widehat{N}(\bar{x};\Omega)$ is a non-empty closed convex cone; 
\item[\textbf{(b)}] the Fr\'echet normal cone to $\Omega$ at $\bar{x}$ is dual to the contingent cone to $\Omega$ at this point, that is,
\begin{equation}
\widehat{N}(\bar{x};\Omega) = T(\bar{x}; \Omega)^{\circ}
:= \left \{  x^* \in \R^n ~|~ \inp{x^*}{v} \leq 0 ~~ \mathrm{~~whenever~} v \in T(\bar{x}; \Omega) \, \right \}
\footnote{ 
For a set $A $, subset of a normed vector space $X$, the (negative) \emph{polar cone} of $A$ is defined as
\begin{equation*}
A^{\circ} := \left \{  x^* \in X^* ~|~ \inp{x^*}{x} \leq 0 ~~ \forall x \in A \right \},
\end{equation*}
while the \emph{dual cone} to $K$, a non-empty closed convex cone in $\R^m$, is defined as
\begin{equation*}
K^* := \{ p \in \R^m ~|~ \inp{p}{v} \geq 0 \mathrm{~~~for~all~} v \in K \}.
\end{equation*}
For more information, look at \citep[p. 35]{Schirotzek}.
}
\end{equation}
\item[\textbf{(c)}] the Mordukhovich normal cone is not a convex set in general (like the case in Example \ref{normals example}), and it is not dual to any tangent cone.
\end{enumerate}
\end{prop}

\begin{eg} \label{Driving an LED with AC-model}
Consider the circuit in Figure \ref{fig: LED} with an LED. In Example \ref{Driving an LED with AC, on the importance of understanding the behaviour of the circuit}, we observed that the generalized equation has the following form:
\begin{equation}
B + A U \in F (U),
\end{equation}
where $ B = 
\begin{pmatrix}
- v_s \\
+ v_s
\end{pmatrix}
$, $ A= 
\begin{pmatrix}
-R & ~R \\ 
~R & -R
\end{pmatrix}
$, and $ \mmap{F}{2}{2} $ is defined as 
$$ F \begin{pmatrix}
u_1 \\ 
u_2
\end{pmatrix} : = 
\begin{pmatrix}
F_1 (u_1) \\ 
F_2 (u_2)
\end{pmatrix} = 
\begin{pmatrix}
F_1 (u_1) \\ 
N( u_2; \R^{+})
\end{pmatrix} $$
Since the matrix $A$ is not invertible and $F_1$ could not be expressed as a normal cone to a set, One can see it is not possible to write this generalized equation as a VI (Variational Inequality). \\
In other words, although one can change the role of voltage and current in the ideal model of diodes using the following fact 
\begin{equation}
\boxed{ v_{\mathrm{\genfrac{}{}{0pt}{}{ideal}{diode}}} \in N \Big(i_{\mathrm{\genfrac{}{}{0pt}{}{ideal}{diode}}}; \,  \R^{+} \Big) ~ \Longleftrightarrow ~ - i_{\mathrm{\genfrac{}{}{0pt}{}{ideal}{diode}}} \in N \Big( - v_{\mathrm{\genfrac{}{}{0pt}{}{ideal}{diode}}}; \,  \R^{+} \Big) },
\end{equation}
in the practical model, or Zener diode case, one can only have the normal cone inclusion for the voltage, that is
\begin{equation}
\boxed{ i_{\mathrm{\genfrac{}{}{0pt}{}{non~ideal}{diode}}}\in N \Big(v_{\mathrm{\genfrac{}{}{0pt}{}{non~ideal}{diode}}}; \, [- \alpha, \beta] \Big) }.
\end{equation}
Thus, in order to use the VI setting, one needs to choose the variables in a way that everything reaches at $ i_D = f(v_s, v_D, v_R, \cdots) $ and $ i_D \in F(v_D) $.
\end{eg}

Before defining the \lq\lq derivatives\rq\rq, let us review few statements concerning product and sum operations with normal and tangent cones. Later, in Subsection \ref{Calculus Rules}, we would use these propositions to provide calculus rules.
 
\begin{prop} [\textbf{Tangents and Normals to Product Sets}] \citep[p. 227]{Rockafellar1998variational} \label{product rule}\\
With $ \R^n $ expressed as $\R^{n_1} \times \dots \times \R^{n_m} $, write $ x \in \R^n$ as $(x_1, ..., x_m)$ with components $x_i \in \R^{n_i} $. If $ C = C_1 \times \dots \times C_m$ for closed sets $C_i \in \R^{n_i} $, then at any $ \bar{x} = ( \bar{x}_1, \dots , \bar{x}_m) $ with $ \bar{x}_i \in C_i $ one has
\begin{eqnarray}
N(\bar{x}; C) = N(\bar{x}_1; C_1) \times \cdots \times N(\bar{x}_m; C_m), \\
\widehat{N}(\bar{x}; C) = \widehat{N}(\bar{x}_1; C_1) \times \cdots \times \widehat{N}(\bar{x}_m; C_m), \\
T(\bar{x}; C)  \,  \subset  \, T(\bar{x}_1; C_1) \times \cdots \times T(\bar{x}_m; C_m), \, \\
\widetilde{T} (\bar{x}; C)  \,  =  \, \widetilde{T}(\bar{x}_1; C_1) \times \cdots \times \widetilde{T}(\bar{x}_m; C_m). \,
\end{eqnarray}
Furthermore, if $C$ is a regular set, then the inclusion for $T(\bar{x}; C) $ becomes an equality like the others.
\end{prop}

\begin{eg} \citep[p. 203]{Rockafellar1998variational} \textbf{(Tangents and Normals to Convex Sets)}\\
For a  convex set $ C \subset \R^n $ and any point $ \bar{x} \in C$, one can observe that the generalized normal cone $N(\bar{x};\Omega)$ agrees with the normal cone in the sense of convex analysis. That is,
\begin{eqnarray*}
\begin{split}
& N(\bar{x}; C) =  \widehat{N}(\bar{x}; C) = \Big\{ v ~|~ \inp{v}{x - \bar{x}} \leq 0 \mfa x \in C \Big\}, \mathrm{and~~} \\
& T(\bar{x}; C) =  \mathrm{cl} \, \Big\{ w ~|~ \exists \lambda >0 \mathrm{~with~} \bar{x} + \lambda w \in C \Big\}.~~~~~~~~~~~~~~~~~~~
\end{split}
\end{eqnarray*}
Furthermore, $C$ is regular at $\bar{x}$ as long as $C$ is locally closed at $\bar{x}$.\\
\end{eg}

\eg \citep[p. 204]{Rockafellar1998variational}\textbf{ (Tangents and Normals to Boxes)}\\
Suppose $ C = C_1 \times \dots \times C_n$, where each $C_j$ is a closed interval in $ \R $ (not necessarily bounded, perhaps just consisting of a single number). 
Then $C$ is regular at every one of its points $ \bar{x} = (\bar{x}_1, . . . , \bar{x}_n) $. Its tangent cone has the form
\begin{eqnarray*}
T(\bar{x}; C)  \,  =  \, T(\bar{x}_1; C_1) \times ... \times T(\bar{x}_m; C_m), \mathrm{~where~~~~~~~~~~~~~~~~~~~~~~~~~~}\\
T(\bar{x}_j; C_j) = 
\left\{ \begin{matrix}
(- \infty, 0] ~ & \mathrm{~~if~} \bar{x}_j \mathrm{~is~the~right~endpoint~of~} C_j, \\
[0, \infty)~~~ & \mathrm{~~if~} \bar{x}_j \mathrm{~is~the~left~endpoint~of~} C_j,~ \\
( - \infty, \infty) & \mathrm{~~if~} \bar{x}_j \mathrm{~is~an~interior~point~of~} C_j,~ \\
[0, 0] ~~~~& \mathrm{~~if~} C_j \mathrm{~is~a~singleton~}~~~~~~~~~~~~~~~~\\
\end{matrix} \right.
\end{eqnarray*}
while its normal cones have the form
\begin{eqnarray*}
N(\bar{x}; C)  \,  =  \, N(\bar{x}_1; C_1) \times ... \times N(\bar{x}_m; C_m), \mathrm{~where~~~~~~~~~~~~~~~~~~~~~~~~~~}\\
N(\bar{x}_j; C_j) = 
\left\{ \begin{matrix}
[0, \infty)~~~ & \mathrm{~~if~} \bar{x}_j \mathrm{~is~the~right~endpoint~of~} C_j, \\
(- \infty, 0] ~ & \mathrm{~~if~} \bar{x}_j \mathrm{~is~the~left~endpoint~of~} C_j,~ \\
[0, 0] ~~~~ & \mathrm{~~if~} \bar{x}_j \mathrm{~is~an~interior~point~of~} C_j,~ \\
( - \infty, \infty) & \mathrm{~~~if~} C_j \mathrm{~is~a~singleton~}~~~~~~~~~~~~~~~~~\\
\end{matrix} \right.
\end{eqnarray*}
In particular, $C$ is a closed convex set. The formulas in previous example relative to a tangent vector $ w = (w_1, . . . , w_m) $ or a normal vector 
$ v = (v_1, . . . , v_n) $ translate directly into the indicated requirements on the signs of the components $ w_j $ and $v_j$. 

\begin{prop}  [\textbf{Tangents and Normals under Set Addition}] \citep[p. 230]{Rockafellar1998variational} \\
Let $ C = C_1+ \dots +C_m $ for closed sets $ C_i \subset \R^n $. Then, at any point $ \bar{x} \in C $, one has
\begin{eqnarray*}
T(\bar{x}; C)  ~~ \supset \bigcup_{\bar{x}_1 + ... + \bar{x}_m \, = \, \bar{x} \atop \bar{x}_i \in C_i } \Big[ \, T(\bar{x}_1; C_1) + \, \cdots \, + T(\bar{x}_m; C_m) \, \Big], \\
\widehat{N}(\bar{x}; C) \subset \bigcap_{\bar{x}_1 + ... + \bar{x}_m \, = \, \bar{x} \atop \bar{x}_i \in C_i } \Big[ \, \widehat{N}(\bar{x}_1; C_1) \cap \, \cdots \, \cap \widehat{N}(\bar{x}_m; C_m) \, \Big].
\end{eqnarray*}
\end{prop}

\rem
The principal difference between tangential and normal approximations is that the former constructions (tangents) provide local approximations of sets in primal spaces, while the latter ones (normals) are defined in dual spaces carrying \lq\lq dual\rq\rq \,information for the study of the local behaviour. \\
Being applied to epigraphs of extended-real-valued functions and graphs of set-valued mappings, tangential approximations will generate corresponding \emph{directional derivatives/subderivatives} of functions and \emph{graphical derivatives} of mappings, while normal approximations will relate to \emph{subdifferentials} and \emph{coderivatives}, respectively, as we would introduce them soon.

\defn \textbf{(Graphical Derivative)} \index{graphical derivative} \\ \citep[p. 215]{implicit} 
For a mapping $ \mmap{F}{n}{m} $ and a pair $ (x,y) $ with $y \in F(x)$, the \emph{graphical derivative} of $F$ at $x$ for $y$ is the mapping $ \mmap{DF(x|y)}{n}{m}$ whose graph is the tangent cone $T (\, (x,y); \gph{F} \, )$ to $\gph {F}$ at $(x,y)$:
\begin{equation}
v \in DF(x|y)(u) \Leftrightarrow  (u,v) \in T (\, (x,y); \gph{F} \, ).
\end{equation}
Thus, $v \in DF(x|y)(u)$ if and only if there exist sequences $u_k \rightarrow u,~ v_k \rightarrow v$ and $ \tau_k \downarrow 0$ such that $ y+ \tau_k v_k \in F(x + \tau_k u_k)$ for all $k$.

\defn \textbf{(Coderivative)}\label{coderivative} \index{coderivative} \\ \citep[p. 229]{implicit} 
For a mapping $ \mmap{F}{n}{m} $ and a pair $(x,y) \in \gph{F}$ at which $\gph{F}$ is locally closed, the \emph{coderivative }of $F$ at $x$ for $y$ is the mapping $ \mmap{D^* F(x|y)}{m}{n} $ defined by
\begin{equation}
 w \in D^{*} F(x|y)(z) \Leftrightarrow (w,-z) \in N (\, (x,y); \gph {F} \, ).
\end{equation}

\eg  \textbf{(Coderivatives of Differentiable Mappings)} \label{coderivative of strict differentiable maps}\hfill \\  \citep[p. 45]{Mordukhovich}
Let $ \smap{f}{n}{m} $ be strictly differentiable at $\bar{x}$. Then 
$$ D^{*} f(\bar{x})(u) = \{  \nabla f(\bar{x}) ^T u \} \mfa u \in \R^m. $$ 
It is sufficient to show that for an arbitrary $u \in \R^m$ and any $ w \in D^{*} f(\bar{x})(u) $, we have $w = \nabla f(\bar{x}) ^{T} u$. By Definitions \ref{coderivative} and \ref{normal cones}, we have sequences 
$ \epsilon_k \downarrow 0, x_k \rightarrow \bar{x} $, and $ (w_k, u_k) \rightarrow (w,u) $ such that
$$ \inp{w_k}{x - x_k} - \inp{u_k}{f(x) - f(x_k)} \, \leq \, \epsilon_k \big( \norm{x - x_k} + \norm{f(x) - f(x_k)} \big) $$
for all $x$ close enough to $x_k$ and all $k \in \N$. It follows from Definition \ref{strict differentiability} of strict differentiability that for any sequence $ \tau_j \downarrow 0 $ as $ j \rightarrow \infty $ there is a sequence of neighborhoods $ U_j $ of $\bar{x}$ with
$$ \norm{ f (x') - f (x) - \nabla f (\bar{x})(x' - x) } \leq \tau_j \norm{x' - x } \mfa x, x' \in U_j, ~ j \in \N. $$
This allows us to select a subsequence $ \{ k_ j \} $ of $ \N $ such that for all $x, x_{k_j} \in U_{k_j}, ~ j \in \N$, 
\begin{eqnarray*}
\begin{split}
\inp{u_{k_j}}{f(x) - f(x_{k_j})} & = & \inp{u_{k_j}}{\nabla f(\bar{x})( x - x_{k_j}) + o \big( x - \bar{x}, x_{k_j} - \bar{x} \big)} ~~~~~ \\
 & = & \inp{\nabla f(\bar{x})^{T} u_{k_j}}{x - x_{k_j}} + o \big( x - \bar{x}, x_{k_j} - \bar{x} \big) ~ u_{k_j}
\end{split}
\end{eqnarray*}
where $o(\cdot , \cdot) $ is defined in a natural way, that is 
$ \displaystyle \lim_{(x, \, x_{k_j}) \rightarrow \rfp{x}{x}} \frac{o \big( x - \bar{x}, x_{k_j} - \bar{x} \big)}{\norm{ \big( x - \bar{x}, x_{k_j} - \bar{x} \big)}} = 0 $.
\begin{eqnarray*}
\begin{split}
\inp{w_{k_j} - \nabla f(\bar{x}) ^{T} u_{k_j}}{x - x_{k_j}} & \leq \, \epsilon_{k_j} \big( \norm{x - x_{k_j}} + \norm{f(x) - f(x_{k_j})} \big) + o \big( x - \bar{x}, x_{k_j} - \bar{x} \big) u_{k_j} \\
& \leq \, \epsilon_{k_j} \big( \norm{x - x_{k_j}} + l  \norm{x - x_{k_j}} ) + \tau_{j} \, \norm{u_{k_j}} \, \norm{x - x_{k_j}} \\
&\leq \, \widetilde{\epsilon}_j \norm{ x - x_{k_j} } 
\end{split}
\end{eqnarray*}
for all $ x, x_{k_j} \in U_{k_j}, ~ j \in \N$, where $ \widetilde{\epsilon}_j := ( 1 + l ) ( \epsilon_{k_j} + \tau_j \norm{u_{k_j} } ) $ and $l$ denotes a Lipschitz constant of $f$ around $\bar{x}$ (existence of which is guaranteed by the strict differentiability of $f$ around $\bar{x}$).
The latter implies that
$$ \norm{w_{k_j} - \nabla f(\bar{x}) ^{T} u_{k_j} } \leq  \widetilde{\epsilon}_j \mathrm{~~~~for~large~} j \in \N, $$
which gives $ w = \nabla f (\bar{x})^{T} u $. \\
The inverse inclusion, which would be equivalent to show 
$$( \nabla f(\bar{x}) ^T u, - u ) \in N \big( (\bar{x}, f(\bar{x})); \gph{f} \big)  \mathrm{~~~for~any~~}  u \in \R^m, $$ 
is a direct consequence of the strict differentiability of $f$ around $\bar{x}$ and the definition of limiting normal cone.

\begin{rem} \label{derivative meaning for functions}
One may wonder the competence of \textit{``derivative\rq\rq} in the last two definitions. Let us see what happens when $F$ is single-valued, thus reducing to a function $ \smap{f}{n}{m}$. Suppose that $f$ is strictly differentiable at $x$; then, for $y = f(x)$, the graphical derivative $Df(x|y)$ is of course the linear mapping $Df(x)$ from $\R^n$ to $\R^m$ with matrix $ \nabla f(x)$. In contrast, the coderivative $D^* f(x|y)$ comes out as the adjoint linear mapping $Df(x) ^ *$ from $\R^m$ to $\R^n$ with matrix $ \nabla f(x) ^T$.
\end{rem}

\defn \textbf{(Strict Graphical Derivative)} \index{strict graphical derivative} \\ \citep[p. 238]{implicit} 
For a mapping $ \mmap{F}{n}{m} $ the \emph{strict graphical derivative} mapping $ \widetilde{D} F \rfpp{x}{y} $ \at{x}{y}, where $ \bar{y} \in F(\bar{x}) $, is defined as a mapping whose graph is the collection of vectors $(u,v)$ for which there exist sequences $ (x_k, y_k) \in \gph{ F}, (x_k, y_k) \rightarrow \rfp{x}{y} $, as well as $ \tau_k \downarrow 0$ and $(u_k, v_k) \rightarrow (u,v)$ such that $(x_k + \tau_k u_k, y_k + \tau_k v_k) \in \gph{ F} $. Thus,
\begin{equation}
v \in \widetilde{D}F(x|y)(u) \Leftrightarrow  (u,v) \in \widetilde{T} \big( (x,y); \gph{F} \big).
\end{equation}

\note
There are many different ways to define the tangential and normal approximation cones. We only defined those cones that tend to be used in derivatives definitions. One can refer to \cite[Chapter 11]{Schirotzek} to see other possibilities. \\
In a general Banach space setting, one has the possibility to define different versions of coderivatives (cf. \cite[Cahpter 1]{Mordukhovich}); but since they coincide in the case of finite dimensional spaces, we avoid introducing them. The one we defined here is referred to as \emph{Normal Coderivative} in \cite{Mordukhovich}.

\rem \citep[p. 216]{implicit}
Note that the graphical differentiation comes from an operation on graphs, and the graph of a mapping $F$ can be converted to the graph of its inverse $F^{-1}$ just by interchanging variables, thus, we immediately have the rule that
$$  D(F^{-1})(y|x) = DF(x|y)^{-1}. $$

\subsection{Calculus Rules} \label{Calculus Rules}
The very next step would be to provide some calculus rules for the derivation operators we introduced recently. 
Being interested in generalized equations (of the form $ f + F \ni 0 $), we will focus only on the rules dealing with the sum of a function and a set-valued map. For possible general calculus rules one can refer to \citep{Levy2004, uderzo2015strong, durea2012chain, durea2014chain, bianchi2015linear}. 

\begin{prop} [\textbf{Sum Rule for Graphical Derivatives}] \label{Sum Rule for Graphical Derivatives} \hfill \\ \citep[p. 216]{implicit}
For a function $\smap{f}{n}{m}$ which is differentiable at $x$, a set-valued mapping $\mmap{F}{n}{m}$ and any $ y \in F(x) $, one has
\begin{equation} \label{Sum Rule for Graphical Derivatives-01}
D (f +F) (x| f(x)+y) = Df(x)+DF(x|y). 
\end{equation}
\end{prop}

\begin{proof} 
If $v \in D(f +F)(x| f(x)+y)(u) $ there exist sequences $ \tau_k \downarrow 0$, $u_k \rightarrow u$ and $v_k \rightarrow v$ such that 
$( x, f(x) + y ) + \tau_k (u_k, v_k) \in \gph{(f + F)} $, that is,
$$ f(x)+y +\tau_k v_k - f(x+\tau_k u_k) \in F(x+\tau_k u_k) \mathrm{~for~every~} k. $$
By using the definition of the derivative for $f$, we get
$$ y + \tau_k (v_k - Df(x) u_k) + o(\tau_k) \in F(x + \tau_k u_k). $$
Hence, by the definition of the graphical derivative, $v - Df(x)u \in  DF(x|y)(u) $.\\
Conversely, if $ v - Df(x)u \in DF(x|y)(u) $ then there exist sequences $ \tau_k \downarrow 0$, $u_k \rightarrow u$ and $w_k \rightarrow v - Df(x)u$  such that 
$ y+\tau_k w_k \in  F(x +\tau_k u_k) $. By the differentiability of $f$, and letting $ v_k = w_k+Df(x) u_k $, we get
$$ y+ f(x) + \tau_k v_k + o (\tau_k) \in (f +F)(x + \tau_k u_k)  $$
which yields $ v \in D(f +F)(x| f(x)+y)(u) $.
\end{proof}

\eg \textbf{(Graphical Derivative for a Constraint System)}  \citep[p. 217]{implicit}
Consider a general constraint system of the form 
\begin{equation}\label{constraint system}
f (x) - D \ni y,
\end{equation}
for a function $ \smap{f}{n}{m} $, a set $ D \subset \R^m $ and a parameter vector $y$; and let $\bar{x}$ be a solution of (\ref{constraint system}) for $\bar{y}$ 
at which $f$ is differentiable. Then for the (set-valued) mapping
$$ G : x \mapsto f (x) - D, \mathrm{~~with~} \bar{y} \in G(\bar{x}), $$
one has
\begin{equation} \label{constraint system2}
DG( \bar{x}| \bar{y})(u) = Df( \bar{x}) u - T( f ( \bar{x}) - \bar{y}; D).
\end{equation}
First observe that for the constant mapping $ F : = - D $, the definition of the graphical derivative gives 
$$ v \in DF(x|z) (u) \Leftrightarrow  v \in T(z; -D) \Leftrightarrow  - v \in T(-z; D). $$
Now, using the sum rule \eqref{Sum Rule for Graphical Derivatives-01} at the reference point $ \rfp{x}{y}$ gives the relation (\ref{constraint system2}).\\
In the special case where $D = \R^s _{-} \times \{ 0 \}^{m-s} $ with $f =(f_1, \dots , f_m)$, the constraint system (\ref{constraint system}) with respect to $y =(y_1, \dots ,y_m)$ takes the form
\begin{equation*}
\left\{ \begin{matrix}
f_i(x)  \leq y_i & \mathrm{~~for~} i = 1, . . . , s, ~~~\\ 
f_i(x)  = y_i & ~~~\mathrm{~~for~} i = s+1, . . . ,m.
\end{matrix}\right.
\end{equation*}
The graphical derivative formula (\ref{constraint system2}) then says that a vector $v = (v_1, . . . ,v_m)$ is in $ DG(x|y)(u) $ if and only if
\begin{equation*}
\left\{ \begin{matrix}
\nabla f_i(x) u \leq v_i & \mathrm{~~for~} i = 1, . . . , s \mathrm{~~with~} f_i(x) = y_i,\\ 
\nabla f_i(x) u = v_i & ~~~\mathrm{~~for~} i = s+1, . . . ,m.~~~~~~~~~~~~~~~~~~
\end{matrix}\right.
\end{equation*}
For $ i = 1, \dots, s $, one can easily calculate that $ T( f_i(x) - y_i ; \, \R_{-}) = \R_{-} $, and then use the product law (\ref{product rule}) to obtain $ T( f(x) - y ; D ) $.

\begin{prop} [\textbf{Sum Rule for Strict Graphical Derivatives}] \label{Sum Rule for Strict Graphical Derivatives}
For a function $\smap{f}{n}{m}$ which is differentiable at $x$, a set-valued mapping $\mmap{F}{n}{m}$ and any $ y \in F(x) $, one has
$$ \widetilde{D} (f +F) (x| f(x)+y) = Df(x) + \widetilde{D}F(x|y). $$
\end{prop}

\begin{proof}
The proof is very similar to the proof of Proposition \ref{Sum Rule for Graphical Derivatives}. It is easy to check that for a differentiable function $f$ at a point $x \in \mathrm{dom} \, f$, the strict graphical derivative $ \widetilde{D}f \big( x|f(x) \big)$ is equal to its graphical derivative $Df \big( x|f(x) \big)$, both equal to the linear mapping $Df(x)$  with matrix $\nabla f(x)$.\\
We first prove that $ \widetilde{D} (f +F) (x| f(x)+y) \subset Df(x) + \widetilde{D}F(x|y) $.\\
Assume that $ v \in \widetilde{D} (f +F) (x| f(x)+y) (u) $. Then, $(u, v) \in \widetilde{T} \Big( (x, f(x)+y); \gph{(f+F)} \Big) $, which in turn, means that there are sequences 
$(t_k)_{k \in \N}$ in $(0, \infty)$, $(v_k)_{k \in \N}$ in $\R^d$, $(u_k)_{k \in \N}$ in $\R^d$, and $(x_k, w_k)_{k \in \N}$ in $\gph{(f+F)}$, converging to $0, v, u$,  and 
$ \big( x, f(x) + y \big) $ respectively, with $ d = n + m $ such that 
$$ \big( x_k + t_k u_k , w_k + t_k v_k \big) \in \gph{(f+F)}. $$
For each $k$, choose $y_k \in F(x_k) $ such that $w_k = f(x_k) + y_k $. Then, $ w_k + t_k v_k \in f(x_k + t_k u_k) + F(x_k + t_k u_k) $, yields
$ y_k + f(x_k) -  f(x_k + t_k u_k) + t_k v_k \in F(x_k + t_k u_k) $. The differentiability of $f$ implies
$ y_k + t_k \big( v_k - u_k D f(x_k) \big) + o (t_k) \in F(x_k + t_k u_k) $. By definition, this means $(u, v - Df(x) u) \in \widetilde{T} \Big( (x, y); \gph{F} \Big)$. 
So, $ v - Df(x) u \in \widetilde{D}F(x|y) (u) $.\\
On the other hand, assume $ v \in Df(x) u + \widetilde{D}F(x|y) (u) $. Then, $ (u, v - Df(x) u) \in \widetilde{T} \Big( (x, y); \gph{F} \Big) $. This implies the existence of sequences 
$ t_k \to 0 $, $ u_k \to u $, $ w_k \to v - Df(x) u $, and $(x_k, y_k) \to (x, y) $ such that $ y_k + t_k w_k \in F(x_k + t_k u_k) $. Let $ v_k = w_k + Df(x_k) u_k $. Using the differentiability of $f$ once more, we obtain
\begin{eqnarray*}
\begin{split}
 & y_k + t_k v_k - f(x_k + t_k u_k) + f(x_k) + o (t_k) \in F(x_k + t_k u_k),  \\
 & y_k + f(x_k) + t_k v_k + o (t_k) \in (f+F) (x_k + t_k u_k).
\end{split}
\end{eqnarray*}
By definition, this means $ (u, v) \in \widetilde{T} \Big( (x, y + f(x) ); \gph{(f+F)} \Big)  $. Thus, $ v $ belongs to $ \widetilde{D} (f +F) (x| y + f(x) ) (u) $, and the proof is complete. 
\end{proof}

\begin{prop} [\textbf{Sum Rule for Coderivatives}] \citep[p. 232]{implicit} \label{sum rule for coderivatve}\\
For a function $ \smap{f}{n}{m} $ which is strictly differentiable at $\bar{x}$ and a mapping $ \mmap{F}{n}{m} $ with $ \bar{y} \in F(\bar{x}) $, for all $u \in \R^m$ one obtains
\begin{equation}
D^{*} (f +F)(\bar{x}~|~ f(\bar{x})+ \bar{y})(u) = D f(\bar{x})^{*} u + D^{*} F(\bar{x}~|~\bar{y})(u).\footnote{
Once more, we want to draw the readers' attention to the distinction between Jacobian matrix and derivative linear mapping mentioned in \ref{notation}. This should explain the difference of notations regarding this theorem with some references. \\ }
\end{equation}
\end{prop}

\begin{proof}
Let $ w \in D^{*} (f +F)(\bar{x}~|~ f(\bar{x})+ \bar{y})(u) $. Then, $ (w, - u) \in N \big( \rfp{x}{z}; \gph{f+F} \big) $ where $ \bar{z} = f(\bar{x})+ \bar{y} $. By definition, there exist sequences $(w_k, u_k) \rightarrow (w, u) $, and $ (x_k, z_k) \xrightarrow{\tiny{\gph{f+F}}} (\bar{x},\bar{z}) $ with $z_k = f(x_k) + y_k$, such that for an arbitrary $ \tau > 0 $ one has
$$ \inp{ (w_k, - u_k) }{ (x - x_k, z - z_k) } \leq \tau \big( \norm{x - x_k} + \norm{z - z_k} \big) $$
for $x, z$ close enough to $x_k, z_k$, respectively. Thus, 
\begin{equation}\label{proof1}
\inp{ w_k }{ x - x_k } - \inp{ u_k}{z - z_k} \leq \tau \big( \norm{x - x_k} + \norm{z - z_k} \big).
\end{equation}
Now, we try to replace the terms involving $z, \bar{z}, z_k$ with known variables of the theorem in the above equation and use strict differentiability of $f$ (cf. Definition \ref{strict differentiability}).
\begin{eqnarray*}
\begin{split}
\inp{u_k}{z - z_k} & = & \inp{u_k}{f(x) + y - f(x_k) - y_k} = \inp{u_k}{y - y_k} + \inp{u_k}{f(x) - f(x_k)} \\
 & = &  \inp{u_k}{y - y_k} + \inp{u_k}{Df(\bar{x}) (x - x_k) + o \big( (x,x_k) - \rfp{x}{x} \big) } ~~~~~~~~~~~\\
 & = &  \inp{u_k}{y - y_k} +  \inp{u_k}{Df(\bar{x}) (x - x_k) } + o \big( (x - \bar{x}, x_k - \bar{x} ) \big) ~ u_k ~~~~~~~\\
 & = & \inp{u_k}{y - y_k} + \inp{ Df(\bar{x})^{*} u_k}{x - x_k} + o \big( (x - \bar{x}, x_k - \bar{x} ) \big) ~ u_k ~~~~~~~~
\end{split}
\end{eqnarray*}
For the term $\norm{z - z_k}$ on the left hand side of \eqref{proof1}, again by strict differentiability of $f$ (and maybe replacing $k \in \N$ with a subsequence $k_j$ of $\N$ such that for each $j \in \N$, the points $x, x_{k_j}$ belong to a neighborhood $U_j$ of $\bar{x}$ \footnote{
Details of this change is almost the same as Example \ref{coderivative of strict differentiable maps}, so we avoid repeating that.
}), one gets
\begin{eqnarray*}
\begin{split}
\norm{z - z_k} & \leq & \norm{y - y_k} + \norm{ f(x) - f( x_k)} ~~~~~~~~~~~~~~~~~~~~~\\
& \leq & \norm{y - y_k} + \epsilon \, \norm{x - x_k} + \norm{Df(\bar{x})} \norm{x - x_k}
\end{split}
\end{eqnarray*}
where $\epsilon > 0 $ could be chosen arbitrarily small. Putting all calculations in \eqref{proof1}, gives
\begin{equation*}
\inp{w_k - Df(\bar{x})^{*} u_k }{x - x_k} - \inp{u_k}{y - y_k} \leq \tau' \big( \norm{x - x_k} + \norm{y - y_k} \big),
\end{equation*}
which is by definition of coderivative, equal to $ w - Df(\bar{x})^{*} u \in D^{*} F(\bar{x}~|~\bar{y})(u) $ or $ w \in Df(\bar{x})^{*} u + D^{*} F(\bar{x}~|~\bar{y})(u)$.\\
Applying the currently proved inclusion to the sum $ ( f + F) + (- f ) $, we get 
\begin{eqnarray*}
\begin{split}
D^{*} (F) (\bar{x}|\bar{y}) (u) \, & = & D^{*}  \big(( f + F) + (- f ) \big)(\bar{x} | f(\bar{x}) + \bar{y} - f(\bar{x}) ) (u) ~~~~~~~~~~ \\
& \subset & D^{*} (f+F) (\bar{x}|f(\bar{x}) (u) + \bar{y})+ D^{*} (-f) (\bar{x}| - f(\bar{x})) (u) \\
& \subset & D^{*} (f+F) (\bar{x}|f(\bar{x}) (u) + \bar{y}) - D(f)^{*} (\bar{x}| f(\bar{x})) (u). ~~~~~
\end{split}
\end{eqnarray*}
which gives the opposite inclusion and thus establish the equality.
\end{proof}

We will end this subsection by two chain rules, both for coderivatives of the composition of a single-valued function with a set-valued map. We will need these special chain rules when dealing with calculation of $D^* (B F (Cz)) $, in which $B$ and $C$ are linear operators and $F$ is the set-valued map representing the diodes $i - v $ characteristic. \\
In this way, instead of referring to the definition for calculating the coderivative in each special example, where $B$, and $C$ could be different, we just calculate the coderivative of $F$ once, and then use the tools provided here and those in Sections \ref{Results about Aubin Property} and \ref{Results in terms of Metric Regularity} to obtain the coderivative of the composition. The first chain rule would be about the coderivative of $ f \circ G $.

\begin{prop} [\textbf{Special Chain Rules for Coderivatives}] \citep[p. 146]{Mordukhovich} \label{special chain rules for coderivatives} \\
For mappings $\mmap{G}{n}{d}$ and $ f : \R^n \times \R^d \rightarrow \R^m $, define
\begin{equation}
( f \circ G)(x) := f (x, G(x)) = \bigcup \big\{ f (x, y) ~|~ y \in G(x) \big\}.
\end{equation}
Given $ \bar{x} \in \mathrm{~dom~} G$, we assume that:
\begin{itemize}[topsep=-1ex, itemsep=0ex, partopsep=0ex, parsep=0ex]
\item[\textbf{(a)}] $ f (x, \cdot) \in L(\R^d, \R^m) $ around $\bar{x}$, i.e., it is a linear operator from $\R^d$ into $\R^m$. Moreover, $ f (\bar{x}, \cdot) $ is injective.
\item[\textbf{(b)}] The mapping $ x \to f (x, \cdot) $ from $\R^n$ into the operator space $ L(\R^d, \R^m) $ is strictly differentiable at $\bar{x}$.
\end{itemize}
Take any $\bar{y} \in G(\bar{x})$ and denote $ \bar{z} := f  \rfp{x}{y} $. Then one has
\begin{equation}\label{chain law}
D^{*} ( f \circ G)(\bar{x}|\bar{z})(z^*) = \nabla_x f \rfp{x}{y} ^{*} z^* + D^* G(\bar{x}|\bar{y}) \big( f (\bar{x}, \cdot )^* z^* \big) \mfa z^* \in \R^m.
\end{equation}
\end{prop}

As a counterpart of the previous proposition, the next chain rule would be about the coderivative of $ G \circ f $. Thanks to the special structure of $\R^m$, both propositions find a simplified presentation here, comparing to the a general Banach space setting in \citep{Mordukhovich, mordukhovich2004restrictive}, where they were first stated and proved. 

\begin{prop} [\textbf{Special Chain Rules for Coderivatives}] \label{special chain rules for coderivatives2} \hfill \\ \citep[Theorem 3.10, p. 2669]{mordukhovich2004restrictive}
Assume that $ \smap{f}{n}{m} $ is strictly differentiable at  $\bar{x} \in \mathrm{~dom} \, f$ and \emph{Restrictive Metrically Regular} (RMR)\footnote{
\textbf{Definition (RMR).} \index{metric regularity ! restrictive metrically regular (RMR)}A function $ \smap{f}{n}{m} $ is said to have the \emph{RMR property} around $\bar{x}$, or $f$ is RMR around this point, if the restrictive mapping $ f : \R^n \longrightarrow f(\R^n) $ between $\R^n$ and the metric space $f(\R^n) \subset \R^m $, whose metric is induced by
the norm on $\R^m$, is metrically regular around $\bar{x}$ in the sense of Definition \ref{def. mr}.\\
One can easily see, by the classical open mapping theorem, that for linear mappings $f$, the RMR property always holds since the subspace $ f(\R^n) $ is closed in $ \R^m$. However, the situation is much more complicated for nonlinear mappings when the RMR property may be violated even in the simplest cases as, for example, for 
$ f(x) = x^2 $ around $ \bar{x} = 0 \in \R $ (for more details, refer to \cite{mordukhovich2004restrictive}).
}
 around this point. Let $\mmap{G}{m}{d} $ be a set-valued mapping, such that $ \bar{z} \in G(f(\bar{x})) $. Then
$$ D^* (G \circ f) \rfpp{x}{z} = D f(\bar{x}) ^* \circ  D^* G|_f (f(\bar{x}) \, | \, \bar{z}), $$
where $G|_f $ is the restriction of $G$ on $f$ (cf. Notation \ref{restriction}).
\end{prop}

\subsection{Inner and Outer Norms} \label{Inner and Outer Norms}
This subsection starts with the introduction of positively homogeneous set-valued mappings; in particular, we will see that the derivatives defined at the beginning of this section, are positively homogeneous. Then, we provide two \emph{norm-like} operators which let us characterize, quantitatively, the behaviour of these maps (Proposition \ref{Norm Characterizations}). \\
The rest of this section would be devoted to provide general statements relating the local stability or (metric) regularity of a set valued map to the (inner or outer) norm of its derivative. 

\defn \textbf{(Positively Homogeneous Mappings)} \label{Positively Homogeneous Mappings} \index{positively homogeneous map} \hfill \\ \citep[p. 216]{implicit}
A mapping $\mmap{H}{n}{m} $ is called \emph{positively homogeneous} when $ \gph{H}$ is a cone, which is equivalent to $H$ satisfying
\begin{equation*}
0 \in H(0) \mathrm{~~~and~~~} H(\lambda x) = \lambda H(x) \mathrm{~for~} \lambda > 0.
\end{equation*}

\rem
The inverse of a positively homogeneous mapping is another positively homogeneous mapping. Linear mappings are positively homogeneous as a special case, their graphs being not just cones but linear subspaces.\\
Because graphical derivative mappings are positively homogeneous, the general properties of positively homogeneous mappings can be applied to them. Norm concepts are available in particular for capturing quantitative characteristics.\\

\defn \textbf{(Outer and Inner Norms)} \label{Outer and Inner Norms} \index{inner norm} \index{outer norm} \hfill \\ \citep[p. 218]{implicit}
For any positively homogeneous mapping $\mmap{H}{n}{m} $, the \emph{outer norm} and the \emph{inner norm} are defined, respectively, by
\begin{equation} \label{Outer and Inner Norms-relations}
|H|^{+} := \sup_{\norm{x} \leq 1} \, \sup_{y \in H(x)} \norm{y}, \mathrm{~~~and~~~}
|H|^{-} := \sup_{\norm{x} \leq 1} \, \inf_{y \in H(x)} \norm{y}
\end{equation}
with the convention $ \displaystyle \inf_{y \in \emptyset} \norm{y} = \infty $ and $ \displaystyle \sup_{ y \in \emptyset } \norm{y} = - \infty $.\\

\note \label{non-empty norm relations}
When $H$ is a linear mapping, both $|H|^{+}$ and $|H|^{-}$ reduce to the operator (matrix) norm $ \norm{H}$ associated with the Euclidean norm. However, it must be noted that neither $|H|^{+}$ nor $|H|^{-}$ satisfies the conditions in the definition of a true ``norm'', inasmuch as 
in general, the positively homogeneous mappings $ \mmap{H}{n}{m} $ do not even form a vector space under addition and scalar multiplication.\\
Consider the ``inverse elements of addition'' axiom. One may define $-H$ in a natural way, that is for any $ x \in \mathrm{~dom }\, H $, let $ (-H)(x) := \{ -y \in \R^m ~|~ y \in H(x) \} $. Then, of course $$ 0_{\R^m} \in \big( H + (-H) \big) (x) \mfa x \in \mathrm{~dom }\, (H \cap -H), $$
but $ H + (-H) \neq O $, the zero mapping, which is needed to be defined as $$  O(x) = \{ 0 \} \mfa x \in \R^n. $$
Nevertheless, some elementary rules are valid, like the following ones:
\begin{equation*}
\begin{matrix}
|\lambda H|^{+} = \norm{\lambda} |H|^{+}, & ~~~|H_1 + H_2|^{+} \leq |H_1|^{+} + |H_2|^{+}, &~~~ |H_2 \circ H_1 |^{+} \leq |H_2|^{+} \, |H_1|^{+}, \\
|\lambda H|^{-} = \norm{\lambda} |H|^{-}, & ~~~|H_1 + H_2|^{-} \leq |H_1|^{-} + |H_2|^{-}, & ~~~|H_2 \circ H_1 |^{-} \leq |H_2|^{-} \, |H_1|^{-}.
\end{matrix}
\end{equation*}
The following example shows that in general, some properties of the norm definition may not be satisfied as well.

\eg
Define the set-valued mapping $ \mmap{H}{}{}$ with $H(x) = \{ 0, \, x \} $. We show that $H$ is a positively homogeneous map with $ |H|^{-} = 0$, and so the inner norm is not a true norm, since $H$ is not the \emph{zero mapping}. \\
First, observe that $ H(0) = \{0 \} $. Letting $ \lambda > 0 $, we obtain
$$ H(\lambda x) = \{ 0, \lambda x \} = \lambda \, \{ 0, \, x \} = \lambda H(x). $$
Thus, $H$ is positively homogeneous; but
$$ |H|^{-} = \sup_{\norm{x} \leq 1} \, \inf_{y \in H(x)} \norm{y} = \sup_{\norm{x} \leq 1} \inf \{0, \norm{x} \} = 0 $$

\begin{note} 
For a positively homogeneous mapping $ \mmap{H}{n}{m} $, when $ \mathrm{dom} \, H = \R^n $, it would be a simple and immediate observation from the defining relation (\ref{Outer and Inner Norms-relations}) and its conventions concerning the empty set, that 
$$ |H|^{+} \, \geq \, |H|^{-}. $$
However, when $ \mathrm{dom} \, H \neq \R^n $, things are a bit tricky. For example, define $\mmap{H}{}{}$ as
\begin{equation*}
H(x) = \left \{
\begin{matrix}
0 & & \mathrm{for~} x \geq 0, \\
\emptyset & & \mathrm{otherwise}.
\end{matrix} \right.
\end{equation*}
Then, for $ x \in [-1, \,0) $, we get $ \displaystyle \sup_{y \, \in \, H(x)} \norm{y} = \sup_{ y \, \in \, \emptyset } \norm{y} = - \infty $. Thus, 
$$ |H|^{+} = \sup_{\norm{x} \leq 1} \, \sup_{y \, \in \, H(x)} \norm{y} = \sup \, \{ 0, - \infty \} = 0. $$
While $ \displaystyle  |H|^{-} = \sup_{\norm{x} \leq 1} \, \inf_{y \, \in \, H(x)} \norm{y} = \sup \, \{ 0, \, \infty \} = \infty $ (for more details refer to \citep{artacho2007}). \\
\end{note}

\eg 
It would be interesting to see the interpretations of the inner and outer norms, when $ H = A ^{-1} $ for a linear function $ \smap{A}{n}{m} $. Let the $ m \times n $ matrix for this linear function be denoted likewise by $A$, for simplicity. 

$\bullet$ If $ m < n$, we have $ A $ surjective (the associated matrix being of rank $m$) if and only if $ |A^{-1}|^{-}$ is finite, this expression being the norm of the right inverse of $ A$, that is, $ |A^{-1}|^{-} = \| A^T (A A^T) ^{-1} \| $. Then $ |A^{-1}|^{+} = \infty $.\\
This fact is a result of the convention and the definition of inner and outer norms. To be more clear, when $A$ is surjective, for any $y \in \R^m$, there exists $x \in \R^n $ such that $Ax = y$. That is, $A^{-1} $ as a set-valued mapping has domain equal to $ \R^m $ 
\begin{center}
$ \displaystyle  |H|^{-} = \sup_{\norm{y} \leq 1} \, \inf_{x \, \in \, A^{-1}(y)} < \infty $, since it is not the 
$ \displaystyle  \sup \inf_{x \, \in \, \emptyset} \norm{x} $.
\end{center}
And vice versa, when $  |H|^{-}  < \infty $, it means that $A^{-1} (y) \neq \emptyset$ for any $y \in \R^m$, that is, $A$ is surjective.\\
One can easily observe that $ H $ is not single-valued, by reminding that the system of linear equations $ A x = y $, with number of equations less than the variables (components of $x$). Then there are infinite solutions corresponding to free variables. Thus 
$$  |H|^{+} = \sup_{\norm{y} \leq 1} \, \sup_{x \, \in \, A^{-1}(y)} \norm{x} = \infty. $$
From linear algebra we know that since $ m < n $, and $A$ is surjective; it has a right inverse, that is, there exist a linear function 
$ A^{-1} _{ \mathrm{right} } : = A^T (A A^T) ^{-1} $ such that 
$ A A^{-1} _{ \mathrm{right} }  = I_m $. Hence, $ |A^{-1}|^{-} = \norm{ A^{-1} _{ \mathrm{right} } }$. 

$\bullet$ On the other hand, if $m > n$, we have $|A^{-1}|^{+} < \infty $ if and only if $A$ is injective (the associated matrix has rank $n$), and then $ |A^{-1}|^{+} = |(A^T A )^{-1} A^T | $ but $|A^{-1}|^{-} = \infty $. \\
Since $ \mathrm{dom }\, H \subsetneq \R^m $, then for some $y \in \R^m,~~ H(y) = \emptyset $, and so, by convention 
$$ |H|^{-} = \sup_{\norm{y} \leq 1} \, \inf_{x \, \in \, \emptyset} \norm{x} = \infty . $$
But these points are not important for the outer norm calculations (since 
$ \displaystyle \sup_{x \, \in \, \emptyset} \norm{x}= - \infty $). For the other points we have a single-valued map which is the left inverse of $A$, given by 
$ A^{-1} _{ \mathrm{left} } : = ( A^T A ) ^{-1} A^T $ and thus the assertion yields.

$\bullet$ For $m=n$, of course, both norms agree with the usual matrix norm $\|A^{-1}\|$, and the finiteness of this quantity is equivalent to non-singularity of $A$.\\

The inner and outer norms have interesting properties specially when dealing with sublinear mappings. We only mention the most needed ones here, however, the enthusiast reader can refer to \citep{artacho2007} for more information. \\

\begin{prop} [\textbf{Norm Characterizations}] \citep[p. 219]{implicit}\label{Norm Characterizations} \\
The inner norm of a positively homogeneous mapping $\mmap{H}{n}{m} $ satisfies 
\begin{equation} \label{Norm Characterizations-01}
| H |^{-} = \inf \big\{ \kappa > 0 ~|~ H(x) \cap \kappa \, \B \not = \emptyset \mfa x \in  \B \big\}.
\end{equation}
In parallel, the outer norm satisfies
\begin{equation} \label{Norm Characterizations-02}
|H|^{+} =  \inf \big\{ \kappa > 0 ~|~~ y \in H(x) \Rightarrow | y | \leq  \kappa |x| \big\} = \sup_{|y|=1 } \dfrac{1}{d(0, H^{-1}(y))} .
\end{equation}
If $H$ has closed graph, then furthermore
\begin{equation} \label{Norm Characterizations-03}
|H| ^{+} < \infty ~ \iff ~ H(0) = \{ 0 \} .
\end{equation}
If $H$ has closed and convex graph, then
\begin{equation} \label{Norm Characterizations-04}
|H|^{-} < \infty ~ \iff ~ \mathrm{dom} \, H = \R^n, 
\end{equation}
and in that case $ |H|^{-} < \infty $ if and only if $H$ is surjective.
\end{prop}

In order to prove the next theorem, we need a proposition which is a simple case of \citep[Theorem 4B.3]{implicit}.

\begin{prop}
For a set-valued mapping $\mmap{F}{n}{m}$ and a point $ \rfp{x}{y} \in \gph{F} $ such that $ \gph{F} $ is locally closed at it, and for any $ c \in (0, \infty)$ satisfying 
\begin{equation} \label{prop-01}
\limsup_{(x,y) \rightarrow \rfp{x}{y} \atop (x,y) \, \in \, \mathrm{\tiny {gph}} F } | DF(x|y) ^{-1}|^{-} < c,
\end{equation}
there are neighborhoods $V$ of $\bar{y}$ and $U$ of $\bar{x}$ such that
\begin{equation} \label{prop-02}
d(x, F^{-1}(y)) \, \leq \,  c \, d( y, F(x) ) \mathrm{~~for~} x \in U, \mathrm{~and~} y \in V.
\end{equation}
\end{prop}

\begin{proof}
As mentioned above, this is a special case of \citep[Theorem 4B.3, p. 222]{implicit} where $ G(p,x) = F(x) - p $, and $y$ takes the place of $p$; we have 
$ S(y) = F^{-1}(y) $ and $ d(0, G(y,x) ) = d(y, F(x)) $. The upper semicontinuity of $ y \longmapsto d(y, F(\bar{x})) $ is granted and the rest are simple modifications.
\end{proof}

We start derivative description statements with derivative criteria for metric regularity, in terms of graphical derivative, in the following theorem and coderivative, in the next theorem.

\begin{thm} [\textbf{Graphical Derivative Criterion for Metric Regularity}] \label{GDCMR} \hfill \\ \citep[p. 221]{implicit} 
For a mapping $\mmap{F}{n}{m} $ and a point $\rfp{x}{y} \in \gph{F}$ at which the $\gph{F} $ is locally closed, one has 
\begin{equation}\label{GDCMR-1}
\mathrm{reg} \, (F;\bar{x}|\bar{y}) = \limsup_{(x,y) \rightarrow \rfp{x}{y} \atop (x,y) \, \in \, \mathrm{\tiny {gph}} F } | DF(x|y) ^{-1}|^{-}
\end{equation}
Thus, $F$ is metrically regular \at{x}{y} if and only if the right side of (\ref{GDCMR-1}) is finite.
\end{thm}

\begin{proof}
For short, let $ d_{DF} $ denote the right side of (\ref{GDCMR-1}). \\
We will start by showing that $ \mathrm{reg} \, (F;\bar{x}|\bar{y}) \leq d_{DF} $. If $ d_{DF} = \infty $, there is nothing to prove. Let $ d_{DF} < c < \infty $. 
Applying the previous proposition, Condition (\ref{prop-02}), becomes the definition of metric regularity of $F$ \at{x}{y}, and therefore 
$ \mathrm{reg} \, (F;\bar{x}|\bar{y}) \leq  c $. 
Since $c$ can be taken arbitrarily close to $d_{DF} $ we conclude that $ \mathrm{reg} \, (F;\bar{x}|\bar{y}) \leq d_{DF} $.\\
Now, we turn  to demonstrating the opposite inequality,
\begin{equation} \label{GDCMR-01}
\mathrm{reg} \, (F;\bar{x}|\bar{y}) \geq d_{DF}.
\end{equation} 
If $ \mathrm{reg} \, (F;\bar{x}|\bar{y}) = \infty $ we are done. Suppose therefore that $F$ is metrically regular \at{x}{y} with respect to a constant $\kappa$ and neighborhoods $U$ for $\bar{x}$ and $V$ for $\bar{y}$. Then
\begin{equation} \label{GDCMR-02}
d(x', F^{-1}(y)) \, \leq \, \kappa  \, d(y, F(x')) \, \leq \, \kappa \, \norm{y - y' }, 
\end{equation}
whenever $ (x', y') \in \gph{F}, x' \in U, y \in V $.\\
We know from Proposition \ref{AP properties1} that $V$ can be chosen so small that $ F^{-1}(y) \cap U \neq \emptyset $ for every $ y \in V $. Pick any $ y' \in V $ and $ x' \in  F^{-1}(y) \cap U $, and let $ v \in \B $. Take a sequence $ \tau_k \downarrow 0 $ such that $ y_k := y' + \tau_k v \in V \mfa k $. By (\ref{GDCMR-02}) and the local closedness of $ \gph{F} $ at $ \rfp{x}{y} $ there exists $ x_k \in F^{-1} (y' + \tau_k v) $ such that
$$ \norm{x' - x_k}  = d(x',F^{-1}(y_k)) \, \leq \, \kappa \, \norm{y_k - y'} = \kappa \tau_k \norm{v}. $$
For $ u_k := (x_k - x') / \tau_k $ we obtain
\begin{equation} \label{GDCMR-03}
\norm{u_k} \, \leq \, \kappa \, \norm{v}.
\end{equation}
Thus, $u_k$ is bounded, so $u_{k_i} \to u $ for a subsequence $ k_i \to \infty $. Since $ (x_{k_i}, y' + \tau_{k_i} v ) \in \gph{F} $, we obtain 
$ (u,v) \in T \big( (x', y'); \gph{F} \big) $. Hence, by the definition of the graphical derivative, we have $ u \in DF^{-1} (y' | x') (v) = DF (x' | y')^{-1}(v) $. The bound (\ref{GDCMR-03}) guarantees that
$$ | DF(x|y) ^{-1}|^{-} \, \leq \, \kappa . $$
Since $ (x, y) \in \gph{F} $ is arbitrarily chosen near $ \rfp{x}{y}$, and $\kappa$ is independent of this choice, we conclude that (\ref{GDCMR-01}) holds and hence we have (\ref{GDCMR-1}).
\end{proof}

\begin{thm} [\textbf{Coderivative Criterion for Metric Regularity}] \citep[p. 232]{implicit} \label{CDCMR theorem} \\
For a mapping $\mmap{F}{n}{m} $ and a pair $\rfp{x}{y} \in \gph{F}$ at which $\gph{F} $ is locally closed, one has
\begin{equation} \label{CDCMR}
\mathrm{reg} \, (F; \bar{x}| \bar{y}) = | D^{*}F(\bar{x}| \bar{y})^{-1} |^{+} .
\end{equation}
Thus, $F$ is metrically regular \at{x}{y} if and only if the right side of (\ref{CDCMR}) is finite, which is equivalent to
$$ D^{*}F(\bar{x}| \bar{y})(u) \ni 0 ~ \Longrightarrow ~ u = 0. $$
\end{thm}

\begin{proof}
Indeed, the proof is a consequence of the following equality and the previous theorem (Theorem \ref{GDCMR}). \\
(\textbf{Basic Equality}) Let $ \mmap{F}{n}{m} $ be a set-valued map, let $ \bar{y} \in F(\bar{x}) $, and assume that $ \gph{F} $ is locally closed at $\rfp{x}{y}$. Then
\begin{equation}
\limsup_{(x,y) \rightarrow \rfp{x}{y} \atop (x,y) \, \in \, \mathrm{\tiny {gph}} F } | DF(x|y) ^{-1}|^{-} =  | D^{*}F(\bar{x}| \bar{y})^{-1} |^{+}.
\end{equation}
For the proof of this equality one can refer to \citep[Theorem 4C.3, p. 233]{implicit}.
\end{proof}

Considering Theorem \ref{MR-AP}, describing the equality of metric regularity of $F$ and Aubin property of $F^{-1}$ at the proper reference points, the following theorem would be another way of expressing the previous criteria.

\begin{thm} [\textbf{Graphical Derivative Criterion for Aubin Property}] \hfill \\ \citep[p. 222]{implicit} 
For a mapping $ \mmap{S}{m}{n} $ and a point $\rfp{y}{x} \in \gph{S}$ at which $ \gph{S} $ is locally closed, one has
\begin{equation}\label{GDCAP}
\mathrm{lip}(S;\bar{x}|\bar{y}) =  \limsup_{(y,x) \rightarrow \rfp{y}{x} \atop (y,x) \, \in \, \mathrm{\tiny {gph}} S } |DS (y |x )|^{-}.
\end{equation}
Thus, $S$ has the Aubin property \at{y}{x} if and only if the right side of (\ref{GDCAP}) is finite.
\end{thm}

Next, we provide a criterion for strong metric sub-regularity in terms of graphical derivative, and as a consequence of Theorem \ref{Characterization by Isolated Calmness of the Inverse Map}, a corollary about isolated calmness.

\begin{thm} [\textbf{Graphical Derivative Criterion for Strong Metric Sub-regularity}] \label{GDCSMSR} \citep[p. 246]{implicit} \\
A mapping $ \mmap{F}{n}{m} $ whose graph is locally closed at $ \rfp{x}{y} \in \gph{F} $ is strongly metrically sub-regular \at{x}{y} if and only if
\begin{equation} 
DF \rfpp{x}{y} ^{-1} (0) = \{ 0 \}, 
\end{equation}
this being equivalent to 
\begin{equation}
| D F \rfpp{x}{y} ^{-1} |^{+} < \infty, 
\end{equation}
and in that case 
\begin{equation}
\mathrm{subreg}\, (F; \bar{x}\, |\, \bar{y}) = |DF\rfpp{x}{y} ^{-1}|^{+}.
\end{equation}
\end{thm}

\begin{cor} [\textbf{Graphical Derivative Criterion for Isolated Calmness}] \label{GDCIC} \hfill \\ \citep[p. 246]{implicit}
For a mapping $ \mmap{S}{m}{n} $ and a point $ \rfp{y}{x} \in \gph{S} $ at which $ \gph{S} $ is locally closed, one has
\begin{equation}
\mathrm{clm}\, (S; \bar{y} \,| \, \bar{x} ) = | D S \rfpp{y}{x} |^{+}.
\end{equation}
\end{cor}

The last statement of this subsection would be a criterion for strong metric regularity in terms of strict graphical derivative. 

\begin{thm} [\textbf{Strict Graphical Derivative Criterion for Strong Metric Regularity}] \citep[p. 238]{implicit} \label{SGDCSMR} \\
Consider a set-valued mapping $ \mmap{F}{n}{m} $ and $ \rfp{x}{y} \in \gph{F} $. If $F$ is strongly metrically regular \at{x}{y}, then 
\begin{equation} \label{SGDCSMR-1}
| \, \widetilde{D} F \rfpp{x}{y} ^{-1} \, | ^{+} < \infty.
\end{equation}
On the other hand, if the graph of $F$ is locally closed at $\rfp{x}{y} $ and
\begin{equation} \label{SDGCSMR-2}
\bar{x} \, \in \, \liminf_{ y \to \bar{y}} F^{-1} (y),
\end{equation}
then condition \eqref{SGDCSMR-1} is also sufficient for strong metric regularity of $F$ \at{x}{y}. In this case the quantity on the left side of \eqref{SGDCSMR-1} equals to
$ \mathrm{reg} \, (F; \bar{x} \, | \, \bar{y} ) $.
\end{thm}

\begin{proof} 
Proposition \ref{AP and SMR} says that a mapping $F$ is strongly metrically regular \at{x}{y} if and only if it is metrically regular there and $ F^{-1} $ has a localization around $\bar{y}$ for $ \bar{x} $ which is nowhere multivalued. Furthermore, in this case for every $ c > \mathrm{reg} \, (F; \bar{x} \, | \, \bar{y}) $ there exists a neighborhood $V$ of $\bar{y} $ such that $ F^{-1} $ has a localization around $\bar{y}$ for $\bar{x}$ which is a Lipschitz continuous function on $V$ with constant $c$.\\
Let $F$ be strongly metrically regular \at{x}{y}, let $ c > \mathrm{reg} \, (F; \bar{x} \, | \, \bar{y}) $ and let $U$ and $V$ be open neighborhoods of $\bar{x}$ and $\bar{y}$, respectively, such that the localization $V \ni y \mapsto \varphi(y) := F^{-1} (y) \cap U $ is a Lipschitz continuous function on $V$ with a Lipschitz constant $c$. \\
We will show first that for every $ v \in \R^m $ the set $ \widetilde{D} F \rfpp{x}{y} ^{-1} (v) $ is non-empty. Let $ v \in \R^m $. Since 
$ \mathrm{dom} \, \varphi \supset V $, we can choose sequences $ \tau_k \downarrow 0 $ and $ u_k $ such that $ \bar{x} + \tau_k u_k = \varphi ( \bar{y} + \tau_k v) $ for large $k$. Then, from the Lipschitz continuity of $\varphi$ with Lipschitz constant $c$ we conclude that $ \norm{u_k} \leq c \norm{v} $, hence $ u_k $ has a cluster point $u$ which, by definition, is from $ \widetilde{D} F \rfpp{x}{y} ^{-1} (v) $. \\
Now, choose any $ v \in \R^m $ and $ u \in \widetilde{D} F \rfpp{x}{y} ^{-1} (v) $; then, there exist sequences $ (x_k, y_k) \in \gph{F} $, $ (x_k, y_k) \longrightarrow \rfp{x}{y} $, $ \tau_k \downarrow 0 $, $u_k \longrightarrow u $ and $ v_k \longrightarrow v $ such that $ y_k + \tau_k v_k \in V $, $ x_k = \varphi (y_k) $ and 
$ x_k + \tau_k u_k = \varphi (y_k + \tau_k v_k) $ for $k$ sufficiently large. But then, again from the Lipschitz continuity of $\varphi$ with Lipschitz constant $c$, we obtain that $ \norm{u_k} \leq c \norm{v_k} $. Passing to the limit we conclude that $ \norm{u_k} \leq c \norm{v} $, which implies that 
$$ |\widetilde{D} F \rfpp{x}{y} ^{-1} | ^{+} \, \leq \, c. $$
Hence (\ref{SGDCSMR-1}) is satisfied; moreover, the quantity on the left side of \eqref{SGDCSMR-1} is less than or equal to $ \mathrm{reg} \, (F; \bar{x} \, | \, \bar{y} )$.\\
To prove the second statement, we first show that $F^{-1}$ has a single-valued bounded localization, that is there exist a bounded neighborhood $U$ of $\bar{x}$ and a neighborhood $V$ of $\bar{y}$ such that $ V \ni y \mapsto F^{-1} (y) \cap U $ is single valued. \\
On the contrary, assume that for any bounded neighborhood $U$ of $\bar{x}$ and any neighborhood  $V$ of $\bar{y}$ the intersection $\gph{F^{-1}} \cap (V \times U) $ is the graph of a multivalued mapping. This means that there exist sequences $ \epsilon_k \downarrow 0, x_k \longrightarrow \bar{x}, x'_k \longrightarrow \bar{x}, x_k \neq x'_k \mfa k $ such that 
$$ F(x_k) \cap F(x'_k) \cap \B_{\epsilon_k} (\bar{y}) \neq \emptyset \mfa k. $$
Let $ t_k = \norm{x_k - x'_k} $ and let $ u_k = \frac{x_k - x'_k}{t_k} $. Then $ t_k \downarrow 0 $ and $ \norm{u_k} = 1$ for all $k$. Hence $ \{ u_k \} $ has a cluster point $ u \neq 0 $. Consider any $ y_k \in F(x_k) \cap F(x'_k) \cap \B_{\epsilon_k} (\bar{y}) $. Then, 
$ y_k + t_k 0 \in F ( x'_k + t_k u_k ) \mfa k $. By the definition of the strict graphical derivative, $ 0 \in \widetilde{D} F \rfpp{x}{y} (u) $. Hence 
$  |\widetilde{D} F \rfpp{x}{y} ^{-1} | ^{+} = \infty $, which contradicts \eqref{SGDCSMR-1}. Thus, there exist neighborhoods $U$ of $\bar{x}$ and $V$ of $\bar{y}$ such that $ \varphi (y) := F^{-1}(y) \cap U$ is at most single-valued on $V$, and $U$ is bounded. By assumption \eqref{SDGCSMR-2}, there exists a neighborhood $ V' \subset V $ of $\bar{y}$ such that 
$F^{-1}(y) \cap U \neq 0$ for any $ y \in V' $, hence $ V' \subset \mathrm{dom} \, \varphi $. Further, since $\gph{F}$ is locally closed at 
$\rfp{x}{y}$ and $\varphi$ is bounded, there exists an open neighborhood $V'' \subset V' $ of $\bar{y}$ such that $ \varphi $ is a continuous function on $V''$. \\
From the definition of the strict graphical derivative we obtain that the set-valued mapping $(x,y) \mapsto \widetilde{D} F (x \, | \, y) $ has closed graph. We claim that condition \eqref{SGDCSMR-1} implies that 
\begin{equation} \label{SGDCSMR-3}
\limsup_{(x,y) \rightarrow \rfp{x}{y} \atop (x,y) \, \in \, \mathrm{\tiny {gph}} F } | \widetilde{D} F (x \, | \, y) ^{-1} |^{+} < \infty.
\end{equation}
On the contrary, assume that there exist sequences $ (x_k, y_k) \in \gph{F} $ converging to $ \rfp{x}{y} $, $ v_k \in \B $ and 
$ u_k \in  \widetilde{D} F (x_k \, | \, y_k) ^{-1} (v_k) $ such that $ \norm{u_k} > k \norm{v_k} $.\\
\emph{Case 1.} There exists a subsequence $ v_{k_i} = 0 $ for all $k_i$. Since $ \gph{\widetilde{D} F (x_{k_i} \, | \, y_{k_i}) ^{-1}} $ is a cone, we may assume that $ \norm{u_{k_i} } = 1 $. Let $u$ be a cluster point of $u_{k_i}$. Then, passing to the limit we get 
$ 0 \neq u \in \widetilde{D} F \rfpp{x}{y} ^{-1} (0) $ which, combined with formula (\ref{Norm Characterizations-02}) in Proposition \ref{Norm Characterizations}, contradicts \eqref{SGDCSMR-1}.\\
\emph{Case 2.} For all large $k$, $ v_k \neq 0 $. Since $ \gph{\widetilde{D} F (x_k \, | \, y_k) ^{-1} } $ is a cone, we may assume that $ \norm{v_k} = 1 $. Then,
$ \displaystyle \lim_{k \to \infty } \norm{u_k}  = \infty $. Define
\begin{equation*}
w_k := \dfrac{u_k}{\norm{u_k}} \, \in \, \widetilde{D} F (x_k \, | \, y_k) ^{-1} \Bigg( \dfrac{v_k}{\norm{u_k}} \Bigg)
\end{equation*}
and let $w$ be a cluster point of $w_k$. Then, passing to the limit we obtain $0 \neq w \in \widetilde{D} F \rfpp{x}{y} ^{-1} (0) $ which, combined with formula (\ref{Norm Characterizations-02}), again contradicts \eqref{SGDCSMR-1}.\\
Hence \eqref{SGDCSMR-3} is satisfied. Therefore, there exists an open neighborhood $ \tilde{V} \subset V'' $ of $ \bar{y} $ such that 
$ | \widetilde{D} F (\varphi(y) \, | \, y) ^{-1} | ^{+} < \infty $ for all $ y \in \tilde{V} $. 
We will now prove that for every $ (x, y) \in \gph{F} $ near $\rfp{x}{y}$ and every $ v \in \R^m $ we have that $ D F (x \, | \, y) ^{-1} (v) \neq \emptyset $. 
Fix $ (x, y) \in \gph{F} \cap (U \times \tilde{V} ) $ and $ v \in \R^m $, and let $ t_k \downarrow 0 $; then there exist $ u_k \in \R^n $ such that 
$ x + t_k u_k = F^{-1} (y + t_k v) \cap U = \varphi (y + t_k v) $ for all large $k$ and we also have that $ t_k u_k \longrightarrow 0 $ by the continuity of $\varphi$. 
Assume that $ \norm{u_k} \longrightarrow \infty $ for some subsequence (which is denoted in the same way without loss of generality). Set $ \tau_k = t_k \norm{u_k} $ and  $ w_k = \frac{u_k}{\norm{u_k}} $. Then $ \tau_k \downarrow 0 $ and, for a further subsequence, $ w_k \longrightarrow w $ for some $w$ with $ \norm{w} = 1 $. 
Since $ (y + \frac{\tau_k 1}{\norm{u_k}} v, x + \tau_k w_k ) \in \gph{F^{-1}} $, we obtain that $ w \in  D F (x \, | \, y) ^{-1} (0) \subset \widetilde{D} F (x \, | \, y)  ^{-1} (0) $ for some $ w \neq 0 $.\\
Thus, $ | \widetilde{D} F (x \, | \, y) ^{-1} | ^{+} = \infty $ contradicting the choice of $\tilde{V}$. Hence the sequence $ \{ u_k \} $ cannot be unbounded and since 
$ y + t_k v \in F (x + t_k u_k) $ for all $k$, any cluster point $u$ of $ \{ u_k \} $ satisfies $ u \in D F (x \, | \, y) ^{-1} (v) $. Hence, $D F (x \, | \, y) ^{-1}$ is non-empty-valued. From this, Note \ref{non-empty norm relations}, and the inclusion $ D F (x \, | \, y) ^{-1} (v) \subset \widetilde{D} F (x \, | \, y) ^{-1} (v) $ we obtain
\begin{equation} \label{SGDCSMR-4}
| D F (x \, | \, y) ^{-1} |^{-} \, \leq \, | \widetilde{D} F (x \, | \, y) ^{-1} |^{+} .
\end{equation}
Putting together \eqref{SGDCSMR-3} and \eqref{SGDCSMR-4}, and utilizing the derivative criterion for metric regularity in Theorem \ref{GDCMR}, we obtain that $F$ is metrically regular \at{x}{y} with $ \mathrm{reg} \, (F; \bar{x} \, | \, \bar{y} ) $ bounded by the quantity on the left side of \eqref{SGDCSMR-3}. But since $ F^{-1} $ has a single-valued localization \at{y}{x}, we conclude that $F$ is strongly metrically regular \at{x}{y}. Moreover, 
$ \mathrm{reg} \, (F; \bar{x} \, | \, \bar{y} ) = | \widetilde{D} F (x \, | \, y) ^{-1} |^{+} $. The proof is complete.
\end{proof}

\subsection{Subdifferentials} \label{Subdifferentials}

In this subsection we will have a very short look into the theory of subdifferentials. We have two reasons for that; first of all, though it is undesired, it is still probable that we encounter problems in which the single-valued function $f$ is not smooth enough. Secondly, the subdifferential of a function is (in general) a set-valued map itself and it has the same idea as graphical derivatives we already introduced in this section.\\
To be more clear, as we already saw in Chapter \ref{Chapter2}, we need to deal with generalized equations of the form $ f + F \ni 0$, that is, sum of a function and a set-valued map. In the previous subsections of this chapter, we tried to provide a tool (based on derivatives of $f$, and $F$) for studying the local stability properties of the solution of the GE. \\
We would see in the next sections of this chapter that continuous differentiability or twice differentiability of $f$ is part of the assumptions of most of our theorems. Thus, one may get worried, or at least curious to know what will happen when $f$ does not satisfy such conditions. \\
Let us recall that the single valued function in our model not only represents the combination of all components of circuit apart from diodes and transistors, but also by using some simplification ideas we extract some parts of the graph of the set-valued map and add it to $f$ (refer to Figure \ref{Multivalued Simplification example} and explanations there). \\
The concept of subdifferential (for functionals, and subgradient for functions to $\R^n$) is an appropriate substitute to generalize the derivative concept to functions which are not differentiable. \\
According to \citep{Schirotzek}, the theory of the subdifferential and the conjugate of convex functionals as well as its various applications originated in the work of Moreau and Rockafellar in the early 1960s. The Rockafellar's \emph{Convex Analysis} \citep{rockafellar1970} is one of the classic texts on the subject in finite-dimensional spaces.

\defn \textbf{(Subdifferentials)} \label{subdifferentials} \index{function ! subdifferential} \\  \citep[p. 167, p. 288]{Schirotzek} 
Given a lower semicontinuous\footnote{
\textbf{Definition (lower semicontinuity).}\index{function ! lower semicontinuous} The functional $ f : \R^d \longrightarrow \overline{\R} $ is called \emph{lower semicontinuous} (\emph{l.s.c.}) at $ \bar{x} \in \mathrm{~dom~} f $ if either
$ f(\bar{x}) = - \infty $ or for every $ k < f(\bar{x}) $ there exists a neighborhood $ U $ of $\bar{x}$ such that
\begin{equation*}
k < f(x) \mathrm{~~~for~all~~} x \in \mathrm{~dom~} f \cap U.
\end{equation*}
\begin{center}
  \includegraphics[width=0.28\textwidth]{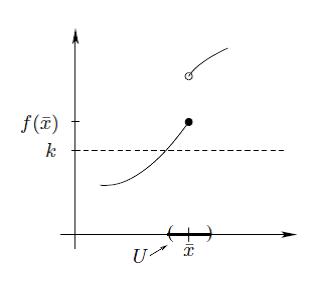}
\end{center}
The figure gives an idea of lower semicontinuity property \citep[p. 22]{Schirotzek}.\\
} 
function $ f : \R^d \longrightarrow \overline{\R} := \R \cup \{\pm \infty \} $ and a point $ \bar{x} \in \R^d $ such that $ f(\bar{x}) \in \R $,
\begin{enumerate}[topsep=0pt,itemsep=1ex,partopsep=1ex,parsep=1ex]
\item[\textbf{(a)}] the functional $f$ is said to be \emph{Fr\'{e}chet subdifferentiable}\index{Fr\'{e}chet subdifferentiable} (\emph{F-subdifferentiable}) at $\bar{x}$ if there exists 
$\xi \in \R^d $, the \emph{F-subderivative} of $f$ at $\bar{x}$, such that
\begin{equation}
\liminf_{h \to 0 } \dfrac{ f(\bar{x} + h ) - f(\bar{x}) - \inp{\xi}{h} }{\norm{h}} \, \geq \, 0.
\end{equation}
The set of all F-subderivatives of $f$ at $ \bar{x} $ is called \emph{Fr\'{e}chet subdifferential} (\emph{F-subdifferential}) of $f$ at $\bar{x}$, and is shown with $\partial_F f(\bar{x}) $.

\item[\textbf{(b)}] the \emph{limiting (Mordukhovich) subdifferential}\index{Mordukhovich subdifferential} of $f$ at $\bar{x}$ is the set $ \partial f(\bar{x}) $ containing all $ \xi \in \R^d $ such that there are sequences $ ( x_k )_{k \in \N} $ and $ ( \xi_k )_{ k \in \N} $ converging to $\bar{x}$ and $\xi$, respectively, with
\begin{equation*}
f(x_k) \to f(\bar{x}) \mathrm{~as~} k \to \infty , \mathrm{~~and~} \xi_k \in \partial_F f(x_k) \mathrm{~~for~each~~} k \in \N.
\end{equation*}

\item[\textbf{(c)}] the \emph{outer subdifferential}\index{outer subdifferential} of $f$ at $ \bar{x} $ is the set $ \partial_{>} f(\bar{x}) $, containing those $ \xi \in \R^d $ for which, there are sequences $ ( x_k )_{k \in \N} $ and $ ( \xi_k )_{ k \in \N} $ converging to $\bar{x}$ and $\xi$, respectively, with
\begin{equation*}
f(x_k) \downarrow f(\bar{x}) \mathrm{~as~} k \to \infty , \mathrm{~~and~} \xi_k \in \partial_F f(x_k) \mathrm{~~for~each~~} k \in \N.
\end{equation*}
\end{enumerate}

\note
If $ f(\bar{x}) $ is infinite, then all the above subdifferentials of $f$ at $\bar{x}$ are defined to be empty sets.\\
Note that nothing will change if for each $k \in \N$ one will take $ \xi_k $ from $ \partial f(x_k) $ instead of $\partial_F f(x_k)$ in the definition of the outer subdifferential.\\
Furthermore, note that the choice of any other (equivalent) norm on $ \R^d $ (instead of the usual one induced by the scalar product) does not affect the above subdifferential constructions.

When the function $f$ is in $\R^n$ we can follow the similar procedure and define generalized Jacobians.

\defn \textbf{(Generalized Jacobians)} \index{generalized Jacobian} \citep[p. 70]{clarke1990}
\begin{enumerate} [topsep=0pt,itemsep=1ex,partopsep=1ex,parsep=1ex]
\item[(a)] The \emph{Bouligand's limiting Jacobian} of a locally Lipschitz continuous function $ \smap{h}{l}{d} $ at a point $ \bar{u} \in \R^l $, is the (non-empty compact) set $ \partial_B h (\bar{u}) $, consisting of all matrices $ A \in \R^{d \times l} $ for which there is a sequence $ (u_n)_{ n \in \N } $ converging to $\bar{u}$ such that $h$ is differentiable at each $u_n$ and $ \nabla h (u_n) \longrightarrow A $ as $ n \to \infty$.
\item[(b)] The \emph{Clarke's generalized Jacobian}\index{Clarke ! Clarke's generalized Jacobian} of $h$ at $\bar{u}$, denoted by $ \partial h (\bar{u}) $, is the convex hull of $ \partial_B h (\bar{u}) $.
\end{enumerate}

\section{Results about Aubin Property} \label{Results about Aubin Property}
In this section and the rest of this chapter we would try to focus on the setting we find convenient for the study of electronic circuits in the static case in Chapter \ref{Chapter2}, that is 
\begin{equation} \label{inclusion}
\boxed{
\begin{matrix}
\Phi (z) : =f(z) \, + B F(Cz) \\
S(p) := \left \{ z ~|~ p \in \Phi (z) \right \}~
\end{matrix}}
\end{equation}
where $ p \in \R^n $  is a fixed vector, $ \smap{f}{n}{n} $ is a function, and $\mmap{F}{m}{m}$ is a set-valued map (with certain assumptions), and $B \in \R^{n \times m},~ C \in \R^{m \times n}$ are given matrices with $m \leq n$, unless otherwise is stated. \\
The problem is how we can provide some handy theorems in order to check the local stability properties of the of solution mapping in \eqref{inclusion}. We will use the derivative criteria of Subsection \ref{Inner and Outer Norms}, and the calculus rules of Subsection \ref{Calculus Rules}.\\
Let us provide the general description of the process for Aubin property here. The starting point would be a point-based description of the Aubin property. 

\begin{thm} [\textbf{Point-based Characterizations of Aubin Property}]\footnote{
This theorem has been proved by Mordukhovich (hence known as Mordukhovich criterion in the literature) in a more general setting, that is when $ F: X \rightrightarrows Y $ is a set-valued mapping between Asplund spaces. In that case, one needs to add another condition to $(b)$, and $(c)$ which is:\\
 $F$ is partially sequentially normally compact (PSNC) at $ \rfp{x}{y} $ \citep [Section 4.2]{Mordukhovich}. 
Since we work on $\R^m$, $F$ is automatically PSNC, and we avoid expressing such definitions and results. \\
Another important situation when the conditions of Theorem \ref{Mordukhovich criteria} can be essentially simplified and efficiently specified, concerns set-valued mappings with closed and convex graphs; see, for example, the Aubin property of convex-graph multifunctions \citep[p. 389]{Mordukhovich}.
} \label{Mordukhovich criteria}\\
Let $ \mmap{F}{n}{m} $ be a set-valued mapping with closed-graph around $ \rfp{x}{y} \in \gph{~F} $. Then the following properties are equivalent:
\begin{itemize}[topsep=-1ex,itemsep=-1ex,partopsep=1ex,parsep=1ex]
\item[\textbf{(a)}] $F$ has the Aubin property around $ \rfp{x}{y} $;
\item[\textbf{(b)}] $ | D^* F \rfp{x}{y} |^{+} < \infty $;
\item[\textbf{(c)}] $ D^* F \rfp{x}{y} (0) = \{ 0 \} $.
\end{itemize}
Moreover, in this case for the exact Lipschitzian bound of $ F $ around  $ \rfp{x}{y} $, one has
\begin{equation}
\mathrm{lip }F \rfp{x}{y} = | D^* F \rfp{x}{y} |^{+}.
\end{equation}
\end{thm}

\begin{proof}
One only needs to combine Theorem \ref{MR-AP} (concerning equivalence relation between metric regularity of $F^{-1}$ and the Aubin property of $F$) with the outer norm characterization in Proposition \ref{Norm Characterizations} and the coderivative criterion for metric regularity in Theorem \ref{CDCMR theorem}.
\end{proof}

\note
Considering the property $(c)$ in the above theorem, we try to go step by step from $ D^* S $ to $D^* F $ using some calculus rules and adapting them to our situation. Almost the same procedure would be done to obtain criteria for other local stability properties or for the metric regularity counterparts.

For the setting \eqref{inclusion}, let us define the mappings $ \mmap{Q}{n}{n} $, and $ \mmap{F_C}{m}{m}$, by 
\begin{eqnarray} \label{Q and F_C}
\begin{split}
& ~~~~~~~~~~~~ Q(z) := B F(Cz), \\
F_C(u) & :=
\left\{ \begin{matrix}
 F(u) &&  \mathrm{ ~if~} u = C z \mathrm{~~for~some~} z \in \R^n, \\
 \emptyset ~~ &&  \mathrm{otherwise. ~}~~~~
\end{matrix} \right.
\end{split}
\end{eqnarray}
We suppose that we have in hand a point $ \rfp{z}{p} \in \gph{~\Phi} $. Put $ \bar{v} := (B^T B) ^{-1} B^T ( \bar{p} - f(\bar{z}) ) $.
We may also need to refer to the following assumptions:

\fbox{\parbox{\textwidth -1ex}{
General Assumptions: \label{general assumptions}
\begin{enumerate}[topsep=-1ex, itemsep=-1ex, partopsep=1ex, parsep=1ex, leftmargin = 9ex]
\item[(A1)] $B$ is injective;
\item[(A2)] $f$ is continuously differentiable in $\R^n$; 
\item[(A3)] $F$ has closed graph;
\item[(A4)] $C$ is surjective; and 
\item[(A5)] there are $ \mmap{F_i}{}{},~ i \in \{ 1,... \, , m \} $ such that $ F(x) = \prod\limits_{i=1}^{m} F_i(x_i) $ whenever $ x = (x_1, ..., x_m)^T \in \R^m $.\\
\end{enumerate}}}

A short review of components descriptions (especially Diodes (jump to \ref{Diode})) of the last chapter, and the examples there (especially, Example \ref{Increasing Components and Loops}), reveals that conditions (A3), and (A5) hold true automatically. Assumption (A4) helps us to replace $F_C$ with $F$, without any further concern. A discussion about loosing the assumption (A1) will be given in Subsection \ref{Isolated Calmness Without Injectivity Assumption}.\\
The first proposition would provide a description for the coderivative of $\Phi$ in terms of $\nabla f$, and the coderivative of $F$.

\begin{prop} [\textbf{Coderivative of Sum for GE}] \citep[Proposition 3.1, p. 94]{adly2}\label{coderivative for GE}\\
Under the assumptions (A1) - (A3), for any $ \xi \in \R^n$ one has
\begin{equation*}
D^* \Phi \rfpp{z}{p} (\xi) = \nabla f(\bar{z})^T \xi + C^T D^*F_C (C\bar{z} \, | \, \bar{v}) (B^T \xi).
\end{equation*}
\end{prop}

\begin{proof}
Fix any $\xi \in \R^n$. As (A2) ensures the strict differentiability of $f$ at $\bar{z}$ (see  Remark \ref{strict differentiability properties}), Theorem \ref{sum rule for coderivatve} implies that 
$$ D^* \Phi \rfpp{z}{p} (\xi) = \nabla f(\bar{z})^T \xi + D^* Q \big( \bar{z} \, | \, \bar{p} - f(\bar{z}) \big) (\xi). $$
Define $ h : \R^n \times \R^m \to \R^n $ by $ h(z,v) := B v$, for each $ z \in \R^n,~ v \in \R^m $ and $ \mmap{G}{n}{m} $, by $ G(z) := F(Cz),~ z \in \R^n$.
Thus, defining $h \circ G$ as introduced in Proposition \ref{special chain rules for coderivatives}, gives
$$ h (z, G(z)) = \bigcup \big\{ Bv  : v \in G(z) \big\} = B G(z) = B F(Cz) = Q(z), \mathrm{~~~whenever~~} z \in \R^n. $$
As $\bar{p} - f(\bar{z}) \in B F(C \bar{z}) $, it is a simple observation that $\bar{v} = (B^T B) ^{-1} B^T ( \bar{p} - f(\bar{z}) ) \in G(\bar{z})$ and $h \rfp{z}{v} = B \bar{v} = \bar{p} - f(\bar{z}) $.\\
Obviously, $ z \to h(z, \cdot) $ is strictly differentiable at $\bar{z}$. 
Since $B$ is injective (A1), $\mathrm{rge} \, h(\bar{z}, \cdot ) = \R^n $; thus, using Proposition \ref{special chain rules for coderivatives} reveals that
\begin{eqnarray*}
\begin{split}
D^* Q \big( \bar{z} \, | \, \bar{p} - f(\bar{z}) \big) (\xi) & = D^* (h \circ G) \big( \bar{z} \, | \, \bar{p} - f(\bar{z}) \big) (\xi) \\
& = \nabla_z h \rfpp{z}{v} ^* (\xi) + D^* G \rfpp{z}{v} \big( h(\bar{z}, \cdot) ^* \xi \big) \\
& = D^* G \rfpp{z}{v} (B^T \xi).
\end{split}
\end{eqnarray*}
Finally, notice that $ F_C $ is actually the restriction of $F$ to the linear mapping $ C$ (that is, the function $\smap{g}{n}{m}$ defined by $g(z) := Cz $, cf. Notation \ref{restriction}), which easily satisfies the conditions of strict differentiability at $\bar{z}$ and RMR around this point. Thus, Proposition \ref{special chain rules for coderivatives2} yields that
\begin{eqnarray*}
\begin{split}
D^* G \rfpp{z}{v} (B^T \xi) & = D^* (F \circ g) \rfpp{z}{v} (B^T \xi) \\
& = \nabla g (\bar{z})^* \circ D^* (F|_g) ( g(\bar{z}) \, | \, \bar{v}) (B^T \xi) \\
& = C^T D^* F_C \big( C\bar{z} \, | \, \bar{v} \big) (B^T \xi).
\end{split}
\end{eqnarray*}
Combining the above expressions ends the proof.
\end{proof}

\begin{thm} [\textbf{Aubin Property Criterion for the Solution Mapping of GE}] \citep[Theorem 3.1, p. 95]{adly2} \label{theorem 3.1}
Under the assumptions (A1)–(A3), $ S $ has the Aubin property at $\rfp{p}{z}$ if and only if
\begin{equation}
0 \in \nabla f(\bar{z})^T \xi + C^T D^*F_C (C\bar{z} \, | \, \bar{v}) (B^T \xi) ~~~~\Rightarrow~~~ \xi = 0
\end{equation}
Moreover, its Lipschitz modulus is given by
\begin{equation}
\mathrm{lip } \big( S; \rfp{p}{z} \big) = \sup \left \{ \, \left \| \xi  \right \| \, : \left (   \nabla f(\bar{z})^T \xi + C^T D^* F_C (C\bar{z} \, | \, \bar{v}) (B^T \xi) \, \right ) \cap \B \neq \, \emptyset ~ \right \}.
\end{equation}
\end{thm}

\begin{proof}
First observe that since $ S = \Phi ^{-1}$, one has $ (\eta, - \xi) \in N \big( \rfp{z}{p}; \gph{\Phi} \big) $ if and only if $ (- \xi, \eta) \in N \big( \rfp{p}{z}; \gph{S} \big)   $ or, equivalently,
$$ \eta \in D^* \Phi \rfpp{z}{p} (\xi) \mathrm{~~if~and~only~if~~} - \xi \in D^* S \rfpp{p}{z} ( - \eta). $$
Having Proposition \ref{coderivative for GE} in hand, one applies the well-known Mordukhovich criterion in the finite-dimensional setting (Theorem \ref{Mordukhovich criteria})  to conclude the proof.
\end{proof}

Employing more assumptions, we may get the following corollaries.

\begin{cor} \citep[Corollary 3.1, p. 95]{adly2} \label{Corollary 3.1} 
Suppose that the assumptions (A1)–(A4) hold true. Then $S$ has the Aubin property at $\rfp{p}{z}$ if and only if
\begin{equation}\label{AP condition}
\left. \begin{matrix}
\big((CC^T)^{-1} C \nabla f(\bar{z})^T \xi , \, B^T \xi \big) \, \in - N \big( (C\bar{z},\bar{v}) ;\gph{F} \big) \\
\nabla f(\bar{z})^T \xi \in \mathrm{rge }C^T~~~~~~~~~~~~~~~~~~~~~~~~~~~~~~~~~~~~~~~~~
\end{matrix} \right\} ~~~~\Longrightarrow~~~ \xi = 0.
\end{equation}
\end{cor}

\begin{proof}
Since $C$ is surjective, then $F_C = F$ and $C C^ T\in \R^{m \times m} $ is non-singular. \\ 
First, let $ \xi \in \R^n $ be such that $ 0 \in \nabla f(\bar{z})^T \xi + C^T D^* F(C \bar{z} \, | \, \bar{v}) (B^T \xi)$. \\
Find $ w \in D^* F(C \bar{z} \, | \, \bar{v}) (B^T \xi) $ with $ \nabla f(\bar{z})^T \xi + C^T w = 0 $. Thus $ - (C C^T)^{-1} C \nabla f(\bar{z})^T \xi $ is in $ D^* F(C \bar{z} \, | \, \bar{v}) (B^T \xi) $ or equivalently, 
$ \big( - (C C^T)^{-1} C \nabla f(\bar{z})^T \xi, - B^T \xi \big) \in N \big( (C\bar{z},\bar{v}) ;\gph{F} \big) $. Clearly, we have $ \nabla f(\bar{z})^T \xi \in \mathrm{~rge~} C^T $ (for each $\xi$, one can define $ u \in \R^n $ as $ u = (C C^T)^{-1} C \nabla f(\bar{z})^T \xi $, thanks to surjectivity of $C$) and Theorem \ref{theorem 3.1} yields the rest.\\
On the other hand, pick any $ \xi \in \R^n $ with $\big((CC^T)^{-1} C \nabla f(\bar{z})^T \xi , \, B^T \xi \big) \, \in - N \big( (C\bar{z},\bar{v}) ;\gph{F} \big)$ and 
$ \nabla f(\bar{z})^T \xi \in \mathrm{~rge~} C^T $. The definition of the coderivative of $F$ implies that $ w :=  - (C C^T)^{-1} C \nabla f(\bar{z})^T \xi \in D^* F(C \bar{z} \, | \, \bar{v})(B^T \xi) $. Thus $ C C^T w = - C \nabla f(\bar{z})^T \xi $. This implies that $ C^T w + \nabla f(\bar{z})^T \xi \in \mathrm{~ker } C \, \cap \mathrm{~rge } C^T  = \{ 0 \}$. Therefore $ 0 \in \nabla f(\bar{z})^T \xi + C^T D^* F(C \bar{z} \, | \, \bar{v}) (B^T \xi)$ and Theorem \ref{theorem 3.1} ends the proof.  
\end{proof}

\rem \label{remark coordinate-wise}
If (A5) also holds, then $ \gph{~F} =\prod\limits_{j=1}^m \gph{\, F_j} $, and Proposition \ref{product rule} (about tangents and normals to product sets) implies that $$ N \big( (C\bar{z},\bar{v}) ;\gph{F} \big) =\prod\limits_{j=1}^m N \Big( \big( (C\bar{z})_j,\bar{v}_j \big) ;\gph{F_j} \Big). $$ 
Thus, the first condition in (\ref{AP condition}) can be checked coordinate-wise.\\

In order to express the next corollary, we need to remind the following definition from linear algebra.
\defn \textbf{ (P-matrix) } \citep[p. 145]{ortega1970} \label{P-matrix} \index{P-matrix}
A matrix $ A \in \R^{n \times n} $ is called a \emph{P-matrix} if all its $k$-by-$k$ principal minor determinants are positive whenever $ k \in \{1, \dots ,n \} $.\\
Using the well known Sylvester's criterion for positive definite matrices\footnote{
In linear algebra, a symmetric $ n \times n $ real matrix $M$ is said to be \emph{positive definite} if the scalar $ z^T M z $ is positive for every non-zero column vector $z \in \R^n $. \\
To see a proof of Sylvester's necessary and sufficient criterion, refer to \citep{gilbert1991}.
}, 
one can obtain an easy to verify rule for P-matrices:
\begin{center}
$ A $ is a P-matrix if and only if for any non-zero $ x \in \R^n $,\\ 
there is $ j \in \{1, \dots, n \} $ such that $ x_j (A x)_j > 0 $.
\end{center}

\begin{cor} \citep[Corollary 1, p. 336]{adly} \label{AP2}
In addition to (A1) - (A5), assume that $ n = m $, that $B = C= I_n $, that $ \nabla f(\bar{z})^T $ is a P-matrix, and that for each $j \in \{1,2,...,n \}$, we have
\begin{equation*}
N \left ( (\bar{z}_j, \bar{v}_j); \, \gph {F_j} \, \right ) \subset \left \{ \, (a,b)^T \in  \R^2 \, : \, a b \leq 0 \right \},
\end{equation*}
Then $S$ has the Aubin property at $ \rfp{p}{z} \in \gph {~S} $.
\end{cor}

\begin{proof}
In view of Remark \ref{remark coordinate-wise}, the condition (\ref{AP condition}) says that $S$ has the Aubin property at $\rfp{p}{z} $ provided that $ \xi = 0_{\R^n} $ is the only point which satisfies
\begin{equation*}
( \nabla f(\bar{z})^T \xi , \xi ) \, \in - \prod\limits_{j=1}^n N \Big( (\bar{z}_j,\bar{v}_j ) ;\gph{F_j} \Big) \subset \Big[ \left \{ \, (a,b)^T \in  \R^2 \, : \, a b \leq 0 \right \} \Big] ^n.
\end{equation*}
Suppose on the contrary that $\xi$ is non-zero. Then for each $j \in \{1,2,...,n \}$, we must have $ \xi_j \big( \nabla f(\bar{z})^T \xi \big)_j \leq 0 $; which is a contradiction, since $\nabla f(\bar{z})^T$ is a P-matrix.
\end{proof}

\rem
Note that the inclusion in Corollary \ref{AP2} is satisfied when $F_j$ is a maximal monotone operator (ref. Definition \ref{Monotone Mappings}).\\
For simplicity, let $\mmap{F}{}{}$ be a maximal monotone operator, and $ (a, b)^T \in N \Big( (\bar{z},\bar{v}) ;\gph{F} \Big) $. Then, by definition of normal cones there exist sequences $ (a_k, b_k) \to (a, b) $, and $ (z_k, v_k) \in \gph{F} $ with $ (z_k, v_k) \to \rfp{z}{v} $, such that
$$ \limsup_{(z, v) \xrightarrow{\tiny{\gph{F}}} (z_k, v_k) } \dfrac{\inp{(a_k, b_k)}{(z, v) - (z_k, v_k)}}{ \norm{(z, v) - (z_k, v_k)} } \, \leq \, 0. $$
Thus, we must have $ a_k ( z - z_k ) + b_k ( v - v_k ) \leq 0 $.\\
From the maximal monotonicity of $F$ we have $ \inp{v - v_k}{ z - z_k} \geq 0 $ (ref. Remark \ref{maximal monotone maps properties}), which means $z - z_k$ and $v - v_k$ must have the same sign. There are three possibilities:
\begin{enumerate} [topsep=-1ex,itemsep=0ex,partopsep=1ex,parsep=1ex]
\item[\textbf{I.~~}] $ z - z_k = 0 $, then $v - v_k$ could be an arbitrary real number. Thus $ b_k = 0 $ and for any $ a_k \in \R $, the desired relation $ a_k b_k \leq 0 $ holds.  
\item[\textbf{II.~}] $z - z_k  > 0 $, then $ v - v_k \geq 0 $. So one can obtain $ a_k + b_k \frac{v - v_k}{z - z_k} \leq 0 $. Thus, 
$$ a_k b_k \, \leq  \,  \frac{- b_k (v - v_k) }{z - z_k} b_k \, \leq \, \dfrac{- b_k ^2  (v - v_k) }{z - z_k} \, \leq \, 0 $$
\item[\textbf{III.}] $z - z_k < 0 $, then $ v - v_k \leq 0 $. So one can obtain $ a_k  \geq - b_k \frac{v - v_k}{z - z_k} $.\\
One can easily observe that if $b_k \leq 0 $, then we must have $ a_k \geq 0$. On the other hand, if $ a_k \leq 0 $, then $b_k \geq 0 $ (omitting the case when $ v - v_k = 0 $, and we must have $ a_k = 0 $).\\
Thus we have proved $ (a, b)^T \in N \Big( (\bar{z},\bar{v}) ;\gph{F} \Big) $ implies $ a b \leq 0 $.
\end{enumerate}

\begin{eg} \textbf{(A Simple Circuit with DIAC)} \citep[p. 339]{adly} \label{DIAC-example} \\
Consider the circuit in Figure \ref{fig: DIAC-circuit-chap03} with a DIAC. Suppose that $ V > 0 $ and $ a > 0 $ are given, and the $ i-v $ characteristic of the DIAC is given by:
\begin{figure}[ht]
  \begin{minipage}{.35\textwidth}
    \begin{eqnarray*}
		F(z) :=
		\left \{ \begin{matrix}
		\dfrac{-V}{\sqrt{1 - \frac{2az}{V}}} & & z < 0,\\[2em]
		[- V, V] & & z = 0,\\[0.5em]
		\dfrac{V}{\sqrt{1 + \frac{2az}{V}}} & & z > 0.
		\end{matrix} \right.
	\end{eqnarray*} 
  \end{minipage}
  \begin{minipage}{.65\textwidth}
    \centering
		\includegraphics[width=0.90\textwidth]{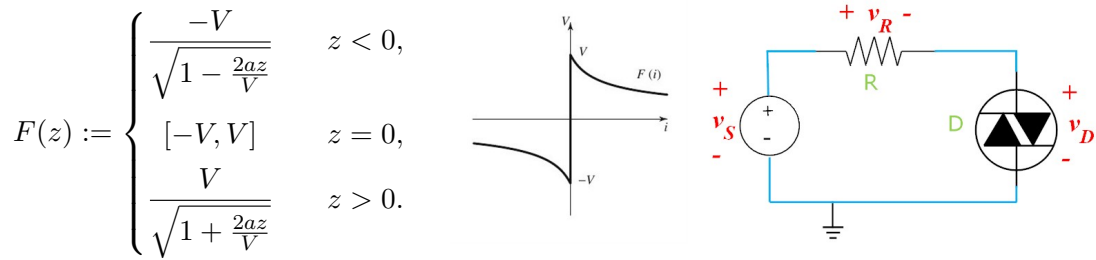}
  \end{minipage}
\caption{Calculating the Aubin property in a circuit with DIAC}  
\label{fig: DIAC-circuit-chap03}
\end{figure}

In Example \ref{DIAC-Chap2}, we discussed the effect of parameters $R$, and $a$ on the solution set at each point and changed the maps to 
\begin{eqnarray*}
f(z) :=
\left \{ \begin{matrix}
R z - \dfrac{V}{\sqrt{1 - \frac{2az}{V}}} + V & & z < 0,\\
R z + \dfrac{V}{\sqrt{1 + \frac{2az}{V}}} - V & & z \geq 0,
\end{matrix} \right.
\mathrm{~~~~and~~~~}
F(z) :=
\left \{ \begin{matrix}
-V ~~~ & &  z < 0,\\
[- V, V] & & z = 0,\\
V ~~~ & & z > 0.
\end{matrix} \right.
\end{eqnarray*} 
In this setting, $ B = C = I_1 $ and the assumptions (A1) - (A4) hold true. For a fixed $ \bar{p} \in \R $, let $\bar{z}$ be the corresponding solution to \eqref{inclusion}, and $\bar{v} = \bar{p} - f(\bar{z}) $. Now, condition \eqref{AP condition} in Corollary \ref{Corollary 3.1} reads as
$$  \big( f'(\bar{z}) \, \xi, \, \xi \big) \, \in - N \Big( (\bar{z},\bar{v}); \, \gph{F} \Big) \Longrightarrow~ \xi = 0, $$
where $ - N \big( (x, y); \, \gph{F} \big) $ is given by 
{\small
\begin{equation*}
~ \left \{ \begin{array}{lcl}
\R \begin{pmatrix}
0 \\ 1
\end{pmatrix}
& ~~ & x > 0,~ y = V, \mathrm{~~or~~} x < 0, ~y = - V, \\[2em]
\R \begin{pmatrix}
1 \\ 0
\end{pmatrix} 
& ~~ & x = 0,~ y \in (-V, \, V), \\[2em]
\R \begin{pmatrix}
0 \\ 1
\end{pmatrix} 
\, \bigcup \, \R \begin{pmatrix}
1 \\ 0
\end{pmatrix} 
\, \bigcup \, \mathrm{~cone} \, 
\begin{Bmatrix} 
	\begin{pmatrix}
	0 \\ -1
	\end{pmatrix}, 
	\begin{pmatrix}
	1 \\ 0
	\end{pmatrix}   
\end{Bmatrix}
& ~~ & x= 0, y = V, \\[2em]
\R \begin{pmatrix}
0 \\ 1
\end{pmatrix} 
\, \bigcup \, \R \begin{pmatrix}
1 \\ 0
\end{pmatrix} 
\, \bigcup \, \mathrm{~cone} \, 
\begin{Bmatrix} 
	\begin{pmatrix}
	0 \\ 1
	\end{pmatrix}, 
	\begin{pmatrix}
	-1 \\ 0
	\end{pmatrix}   
\end{Bmatrix}
& ~~ & x= 0, y = -V, \\[2em]
\emptyset & ~~~~~~~~ & \mathrm{otherwise.}
\end{array} \right.
\end{equation*}
}
Let us check the first rule of the normal cone, for example. The condition $\big( f'(\bar{z}) \xi, \, \xi \big) = (0, \lambda) $ for some $\lambda \in \R$, when $\bar{z} > 0$, and $\bar{p} - f(\bar{z}) = V $ or $\bar{z} < 0$, and $\bar{p} - f(\bar{z}) = -V $, should result in $ \xi = 0 $. 
The condition will fail if for $ \lambda \neq 0 $, we can obtain $ f'(\bar{z}) = 0$. \\
If $ a < R$, one knows that $ f'(z) > 0 ~ \forall  z \in \R $. If $a = R $, then $f'(z) = 0$ only happens at $z = 0$ which is out of the available range of $\bar{z}$. Thus, 
$ \lambda = 0 $ and the condition is satisfied. \\
In the third and forth rule of the normal cone,  $ \bar{z} = 0 $ is allowed and so, in case of $a = R$, for $\bar{p} \in \{V, -V \} $, the condition \eqref{AP condition} will be violated.\\
In the forth rule, let us check the criteria $ \big( f'(\bar{z}) \xi, \, \xi \big) \in \mathrm{~cone} \, 
\begin{Bmatrix} 
	\begin{pmatrix}
	0 \\ 1
	\end{pmatrix}, 
	\begin{pmatrix}
	-1 \\ 0
	\end{pmatrix}   
\end{Bmatrix} $ 
when $\bar{z} = 0$, and $ \bar{v} = \bar{p} - f(0) = - V $. For $ a < R $, since $ f'(0) = R - a > 0 $, the points $\big( f'(\bar{z}) \xi, \, \xi \big) $ define a line passing thorough origin in the first and third quadrant. Thus, the only possible intersection with the cone would be at $(0, 0)$ which imply $ \xi = 0 $. To sum up these observations, 
\begin{enumerate}[topsep=-1ex, itemsep=0ex, partopsep=1ex, parsep=1ex, leftmargin = 5ex]
\item[] if $ a < R $, then $S$ has the Aubin property at any point $ \rfp{p}{z} \in \gph{S}$;
\item[] if $ a = R $, then $S$ has the Aubin property at $ \rfp{p}{z} $ provided that $ \bar{p} \not \in \{-V, V \} $.
\end{enumerate}
The case $ a > R $ could be discussed similarly, but the solution is not unique any more. For example, if $ \bar{z} = 0 $, then $S$ has the Aubin property at $ \rfp{p}{z} $ when $ \bar{p} \in (-V, V) $, since the \lq\lq cone\rq\rq \,terms in the third and forth pieces of the normal cone calculation permit the $ \xi $ to be non-zero. 
\end{eg}

\begin{eg} \textbf{(A Circuit with SCR and Zener Diode)} \label{SCR with Zener-example} \citep[p. 344]{adly} \\
Consider the circuit in Figure \ref{fig: SCR and the Zener Diode} involving two non-smooth elements in parallel: the SCR and the Zener Diode. 

\begin{figure}[ht]
	\centering
		\includegraphics[width=0.95\textwidth]{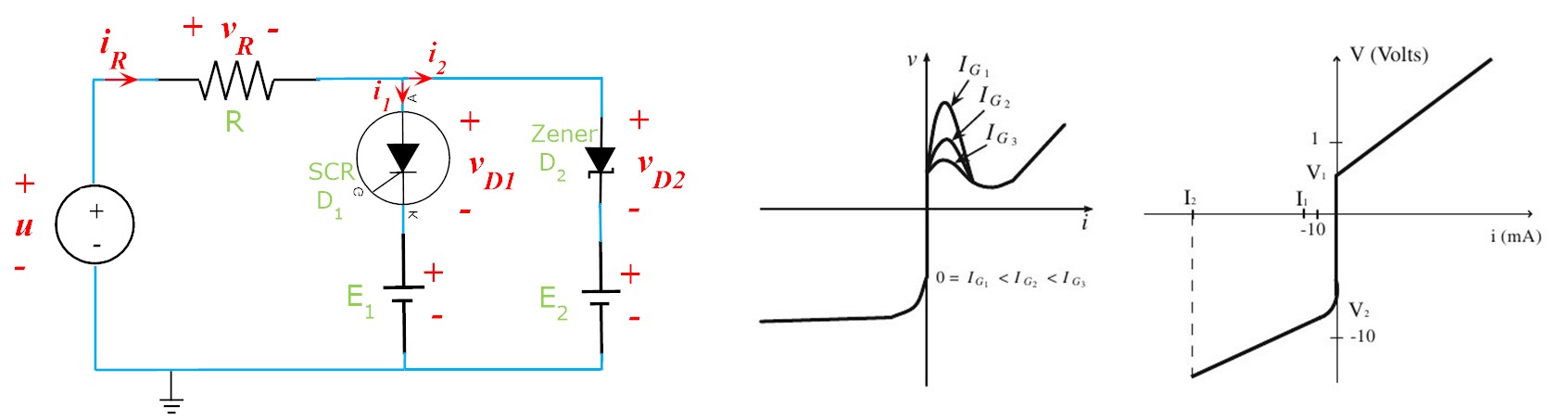}
	\caption{Calculating the Aubin property in a circuit with SCR and Zener Diode.}
	\label{fig: SCR and the Zener Diode}
\end{figure}

In Example \ref{SCR and Zener Diode-Chap2}, we obtained the mathematical model for this circuit in the form of \eqref{inclusion} with $ m = n = 2 $, $ B = C = I_2 $, 
$ f(z) = A z $ for $ z \in \R^2 $, where 
$ A = \begin{pmatrix}
R & R \\
R & R
\end{pmatrix}$, and $ \mmap{F}{2}{2} $ is defined as 
$ F \begin{pmatrix}
z_1 \\ 
z_2
\end{pmatrix} : = 
\begin{pmatrix}
F_1 (z_1) \\ 
F_2 (z_2)
\end{pmatrix}
$. \\
The $ i - v $ characteristics $ F_1 $ of SCR and $ F_2 $ of the Zener diode are defined as
\begin{eqnarray*}
F_1 (z) :=
\left \{ \begin{array}{lcl}
a z + V_1  & & z < 0, \\
$[$ V_{1}, \, \varphi(0) $]$  & &  z = 0, \\
\varphi(z) & & z \in [0, \alpha], \\
a (z - \alpha) + \varphi (\alpha) & & z > \alpha,
\end{array} \right.
~~~~
F_2(z) :=
\left \{ \begin{matrix}
b z - V & & z < 0,\\
[- V, V] & & z = 0,\\
b z + V & & z > 0,
\end{matrix} \right.
\end{eqnarray*} 
where $ V, \, -V_1, \, \alpha, \, a $, and $ b $ are positive constants, and $ \varphi : \R \longrightarrow (0, \infty) $ is a continuously differentiable function with 
$ \varphi(\alpha) < \varphi(0)$, $ \varphi'(0) > 0 $, and $ \varphi'(\alpha) > 0 $. \\
Given $ \bar{p} = (\bar{p}_1, \bar{p}_2) \in \R^2 $, denote by $ \bar{z} = (\bar{z}_1, \bar{z}_2)^T $ the corresponding solution to \eqref{inclusion}. 
Condition \eqref{AP condition} of Corollary \ref{Corollary 3.1} reads as 
$$ \big( I_2 \, A \, \xi , I_2 \, \xi \big) \in - N \Big( (I_2 \bar{z}, \bar{v}); \gph{F} \Big) \Longrightarrow \xi  = (0, 0), $$
where $ \bar{v} = \bar{p} - f(\bar{z}) $. In view of Remark \ref{remark coordinate-wise}, the necessary and sufficient condition for the Aubin property of $S$ at the reference point is
\begin{eqnarray}\label{SCR with Zener-example-eq01}
\left. \begin{array}{l}
\big( R (\xi_1 + \xi_2), \, \xi_1 \big) \in - N \Big( \big( \bar{z}_1, \, \bar{p}_1 - R (\bar{z}_1 + \bar{z}_2) \big) ; \gph{F_1} \Big)  \\[1em]
\big( R (\xi_1 + \xi_2), \, \xi_2 \big) \in - N \Big( \big( \bar{z}_2, \, \bar{p}_2 - R (\bar{z}_1 + \bar{z}_2) \big) ; \gph{F_2} \Big)  
\end{array} \right \} \Longrightarrow \xi_1 = \xi_2 = 0,
\end{eqnarray}
where the first normal cone, $ N \big( (x, y); \gph{F_1} \big) $, can be calculated as

{\footnotesize
\begin{eqnarray*}
~ \left \{ \begin{array}{lll}
\R \begin{pmatrix}
a \\ -1
\end{pmatrix}
&  & 
\begin{matrix}
x < 0,~ y = a x + V_1, \mathrm{~~or~~}~~~~ \\  x > \alpha, ~y = a(x - \alpha) + \varphi(\alpha)
\end{matrix}, \\[2em]
\R \begin{pmatrix}
1 \\ 0
\end{pmatrix} 
& ~~ & x = 0,~ y \in (V_1, \, \varphi(0)), \\[2em]
\R \begin{pmatrix}
\varphi'(x) \\ - 1
\end{pmatrix} 
& ~~ & x \in (0, \alpha),~ y = \varphi(x), \\[2em]
\R \begin{pmatrix}
a \\ -1
\end{pmatrix} 
\, \bigcup \, \R \begin{pmatrix}
1 \\ 0
\end{pmatrix} 
\, \bigcup \, \mathrm{~cone} \, 
\begin{Bmatrix} 
	\begin{pmatrix}
	a \\ -1
	\end{pmatrix}, 
	\begin{pmatrix}
	1 \\ 0
	\end{pmatrix}   
\end{Bmatrix}
& ~~ & x= 0, ~y = V_1, \\[2em]
\R \begin{pmatrix}
- \varphi'(0) \\  1
\end{pmatrix} 
\, \bigcup \, \R \begin{pmatrix}
-1 \\ 0
\end{pmatrix} 
\, \bigcup \, \mathrm{~cone} \, 
\begin{Bmatrix} 
	\begin{pmatrix}
	- \varphi'(0) \\  1
	\end{pmatrix}, 
	\begin{pmatrix}
	- 1 \\ 0
	\end{pmatrix}   
\end{Bmatrix}
& ~~ & x= 0,~ y = \varphi(0), \\[2em]
\R \begin{pmatrix}
a \\ -1
\end{pmatrix} 
\, \bigcup \, \R \begin{pmatrix}
\varphi'(\alpha) \\ -1
\end{pmatrix} 
\, \bigcup \, \mathrm{~cone} \, 
\begin{Bmatrix} 
	\begin{pmatrix}
	a \\ -1
	\end{pmatrix}, 
	\begin{pmatrix}
	\varphi'(\alpha) \\ -1
	\end{pmatrix}   
\end{Bmatrix}
& ~~~~~ & x= \alpha,~ y = \varphi(\alpha), \\[2em]
\emptyset & ~~ & \mathrm{otherwise.}
\end{array} \right.
\end{eqnarray*}
}

The second normal cone, $ N \big( (x, y); \gph{F_2} \big) $, could be obtained similarly  

{\footnotesize
\begin{eqnarray*}
~ \left \{ \begin{array}{lll}
\R \begin{pmatrix}
b \\ -1
\end{pmatrix}
&  & 
\begin{matrix}
x < 0,~ y = b x - V, \mathrm{~~or~~} \\  x > 0, ~y = b x + V~~~~~~~
\end{matrix}, \\[2em]
\R \begin{pmatrix}
1 \\ 0
\end{pmatrix} 
& ~~ & x = 0,~ y \in (-V, \, V), \\[2em]
\R \begin{pmatrix}
b \\ -1
\end{pmatrix} 
\, \bigcup \, \R \begin{pmatrix}
1 \\ 0
\end{pmatrix} 
\, \bigcup \, \mathrm{~cone} \, 
\begin{Bmatrix} 
	\begin{pmatrix}
	b \\ -1
	\end{pmatrix}, 
	\begin{pmatrix}
	1 \\ 0
	\end{pmatrix}   
\end{Bmatrix}
& ~~ & x= 0, ~ y = -V, \\[2em]
\R \begin{pmatrix}
-b \\ 1
\end{pmatrix} 
\, \bigcup \, \R \begin{pmatrix}
-1 \\ 0
\end{pmatrix} 
\, \bigcup \, \mathrm{~cone} \, 
\begin{Bmatrix} 
	\begin{pmatrix}
	-b \\ 1
	\end{pmatrix}, 
	\begin{pmatrix}
	-1 \\ 0
	\end{pmatrix}   
\end{Bmatrix}
& ~~ & x= 0,~ y = V, \\[2em]
\emptyset & ~\hspace*{1.5cm}~ & \mathrm{otherwise.}
\end{array} \right.
\end{eqnarray*} 
}

In order to simplify the process of checking condition \eqref{SCR with Zener-example-eq01}, we multiply the inclusions in (\ref{SCR with Zener-example-eq01}) by the matrices $ M_1 := \begin{pmatrix}
0 & 1 \\
\frac{1}{R} & -1 
\end{pmatrix} $, 
and 
$ M_2 := \begin{pmatrix}
\frac{1}{R} & -1 \\
0 & 1 
\end{pmatrix} $, 
respectively, to get the following equivalent condition
\begin{equation} \label{SCR with Zener-example-eq02}
M_1 \Big( N \big( (\bar{z}_1, \bar{v}_1); \gph{F_1} \big) \Big) \bigcap  M_2 \Big( N \big( (\bar{z}_2, \bar{v}_2); \gph{F_2} \big) \Big) = (0, 0) \footnote{
The reason is that 
\begin{equation*}
M_1 \begin{pmatrix}
R \, \xi_1 + R \, \xi_2 \\
R \, \xi_1
\end{pmatrix} = 
\begin{pmatrix}
 \xi_1 \\
 \xi_2
\end{pmatrix} =
M_2 \begin{pmatrix}
R \, \xi_1 + R \, \xi_2 \\
R \, \xi_2
\end{pmatrix}.
\end{equation*}
Therefore, the condition will hold true if and only if $ (0, 0) $ is the only point in the intersection of $ M_1 \Big( N \big( (\bar{z}_1, \bar{v}_1); \gph{F_1} \big) \Big) $ and 
$ M_2 \Big( N \big( (\bar{z}_2, \bar{v}_2); \gph{F_2} \big) \Big) $. 
}.
\end{equation}
Now, one infers that $M_1 \Big( N \big( (\bar{z}_1, \bar{v}_1); \gph{F_1} \big) \Big)$ is given by

{\footnotesize
\begin{eqnarray*}
~ \left \{ \begin{array}{lcl}
\R \begin{pmatrix}
-1 \\ \dfrac{a}{R} + 1 
\end{pmatrix}
& & 
\begin{matrix}
\bar{z}_1 < 0,~ \bar{v}_1 = a \bar{z}_1 + V_1,~\mathrm{or}~~~~~~~~ \\  \bar{z}_1 > \alpha, ~\bar{v}_1 = a(\bar{z}_1 - \alpha) + \varphi(\alpha),
\end{matrix} \\[2em]
\R \begin{pmatrix}
0 \\ \dfrac{1}{R}
\end{pmatrix} 
& & \bar{z}_1 = 0,~ \bar{v}_1 \in (V_1, \, \varphi(0)), \\[2em]
\R \begin{pmatrix}
-1 \\ \dfrac{\varphi'(x)}{R} + 1
\end{pmatrix} 
& & \bar{z}_1 \in (0, \alpha), ~ \bar{v}_1 = \varphi(\bar{z}_1), \\[2em]
\R \begin{pmatrix}
-1 \\ \dfrac{a}{R} + 1 
\end{pmatrix} 
\, \bigcup \, \R \begin{pmatrix}
0 \\ \dfrac{1}{R}
\end{pmatrix} 
\, \bigcup \, \mathrm{~cone} \, 
\begin{Bmatrix} 
	\begin{pmatrix}
	-1 \\ \dfrac{a}{R} + 1 
	\end{pmatrix}, 
	\begin{pmatrix}
	0 \\ \dfrac{1}{R}
	\end{pmatrix}   
\end{Bmatrix}
& & \bar{z}_1= 0, ~\bar{v}_1 = V_1, \\[2em]
\R \begin{pmatrix}
1 \\  \dfrac{- \varphi'(0)}{R} + 1
\end{pmatrix} 
\, \bigcup \, \R \begin{pmatrix}
0 \\ \dfrac{-1}{R}
\end{pmatrix} 
\, \bigcup \, \mathrm{~cone} \, 
\begin{Bmatrix} 
	\begin{pmatrix}
	1 \\  \dfrac{- \varphi'(0)}{R} + 1
	\end{pmatrix}, 
	\begin{pmatrix}
	0 \\ \dfrac{-1}{R}
	\end{pmatrix}   
\end{Bmatrix}
& & \bar{z}_1= 0, ~\bar{v}_1 = \varphi(0), \\[2em]
\R \begin{pmatrix}
-1 \\ \dfrac{a}{R} + 1
\end{pmatrix} 
\, \bigcup \, \R \begin{pmatrix}
-1 \\ \dfrac{\varphi'(\alpha)}{R} + 1
\end{pmatrix} 
\, \bigcup \, \mathrm{~cone} \, 
\begin{Bmatrix} 
	\begin{pmatrix}
	-1 \\ \dfrac{a}{R} + 1
	\end{pmatrix}, 
	\begin{pmatrix}
	-1 \\ \dfrac{\varphi'(\alpha)}{R} + 1
	\end{pmatrix}   
\end{Bmatrix}
& & \bar{z}_1= \alpha, ~ \bar{v}_1 = \varphi(\alpha), \\[2em]
\emptyset & &  \mathrm{otherwise.}
\end{array} \right.
\end{eqnarray*}
}

And $M_2 \Big( N \big( (\bar{z}_2, \bar{v}_2); \gph{F_2} \big) \Big)$ is given by\\

{\footnotesize
\begin{eqnarray*}
~ \left \{ \begin{array}{lll}
\R \begin{pmatrix}
\dfrac{b}{R} + 1  \\ -1
\end{pmatrix}
&  & 
\begin{matrix}
\bar{z}_2 < 0,~ \bar{v}_2 = b \bar{z}_2 - V, \mathrm{~~or~~} \\  \bar{z}_2 > 0, ~ \bar{v}_2 = b \bar{z}_2 + V~~~~~~~
\end{matrix}, \\[2em]
\R \begin{pmatrix}
\dfrac{1}{R}  \\ 0
\end{pmatrix} 
& ~~ & \bar{z}_2 = 0,~ \bar{v}_2 \in (-V, \, V), \\[2em]
\R \begin{pmatrix}
\dfrac{b}{R} + 1  \\ -1
\end{pmatrix} 
\, \bigcup \, \R \begin{pmatrix}
\dfrac{1}{R}  \\ 0
\end{pmatrix} 
\, \bigcup \, \mathrm{~cone} \, 
\begin{Bmatrix} 
	\begin{pmatrix}
	\dfrac{b}{R} + 1  \\ -1
	\end{pmatrix}, 
	\begin{pmatrix}
	\dfrac{1}{R}  \\ 0
	\end{pmatrix}   
\end{Bmatrix}
& ~~ & \bar{z}_2= 0,~ \bar{v}_2 = -V, \\[2em]
\R \begin{pmatrix}
\dfrac{-b}{R} - 1 \\ 1
\end{pmatrix} 
\, \bigcup \, \R \begin{pmatrix}
\dfrac{-1}{R} \\ 0
\end{pmatrix} 
\, \bigcup \, \mathrm{~cone} \, 
\begin{Bmatrix} 
	\begin{pmatrix}
	\dfrac{-b}{R} - 1 \\ 1
	\end{pmatrix}, 
	\begin{pmatrix}
	\dfrac{-1}{R} \\ 0
	\end{pmatrix}   
\end{Bmatrix}
& ~~ & \bar{z}_2= 0, ~\bar{v}_2 = V, \\[2em]
\emptyset & ~\hspace*{0.4cm}~ & \mathrm{otherwise.}
\end{array} \right.
\end{eqnarray*} 
}

Hence, one can easily observe that the Aubin property of the solution map depends explicitly on the exact values of the parameters of the circuit and the reference point. For example, if 
$ (\bar{z}_1, \bar{v}_1) $, and $ (\bar{z}_2, \bar{v}_2) $ reside in the first sub-domain of $M_1 \big( N \big( (\bar{z}_1, \bar{v}_1); \gph{F_1} \big) \big)$ and 
$M_2 \big( N \big( (\bar{z}_2, \bar{v}_2); \gph{F_2} \big) \big)$, respectively, then any point of the form 
$$ \lambda \begin{pmatrix}
-1 \\ \dfrac{a}{R} + 1 
\end{pmatrix} = \lambda'  \begin{pmatrix}
\dfrac{b}{R} + 1  \\ -1
\end{pmatrix}, $$
would belong to the intersection in \eqref{SCR with Zener-example-eq02} for $ \lambda, \lambda' \in \R $. The equality holds if and only if $ \lambda = \lambda' = 0 $. Thus, the solution mapping has the Aubin property. \\
When $ \bar{z} = (0, 0 ) $, condition \eqref{SCR with Zener-example-eq02} holds if $ \bar{p}_1 \in \big( V_1, \varphi (0) \big) $, and 
$ \bar{p}_2 \in [ -V, V ] $. Thus, the solution mapping has the Aubin property for these reference points. \\
For $ \bar{z}_1 \in (0, \alpha), \bar{p}_1 = \varphi(\bar{z}_1) $, and $ \bar{z}_2 = 0, \bar{p}_2  \in \{ - V, V \} $; one should check the intersection of the line 
$ \R \, \big(-1 , \, \frac{\varphi'(\bar{z}_1)}{R} + 1 \big)^T  $ with the third piece in $M_2 \Big( N \big( (\bar{z}_2, \bar{v}_2); \gph{F_2} \big) \Big)$, which is shown in Figure \ref{fig: SCR-with-Zener-chap03-AP}.
 
\begin{figure}[ht]
	\centering
		\includegraphics[width=0.65\textwidth]{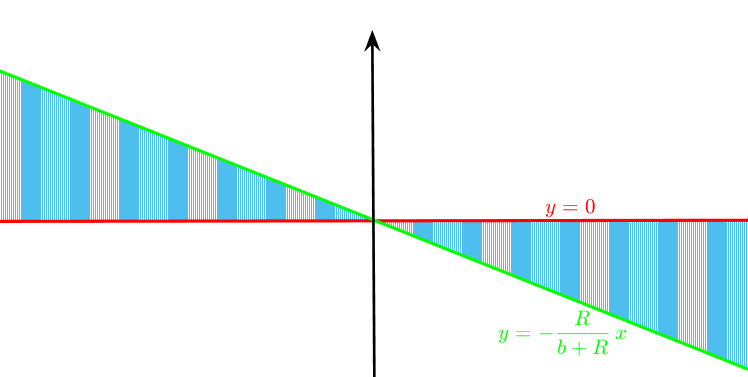}
	\caption{A part of the set $ M_2 \big( N \big( (\bar{z}_2, \bar{v}_2); \gph{F_2} \big) \big) $}
	\label{fig: SCR-with-Zener-chap03-AP}
\end{figure}

Clearly, the condition \eqref{SCR with Zener-example-eq02} is violated when $ 0 \, \leq \, \dfrac{\varphi'(\bar{z}_1)}{R} +1 \, \leq \, \dfrac{R}{b + R} $. Therefore, the solution mapping does not have the Aubin property if $ - R \, \leq \, \varphi'(\bar{z}_1) \, \leq \, \dfrac{- R \, b}{b + R} $.

\end{eg}

%% file: Chapters/Chapter03-b.tex
\section{Results about Isolated Calmness}

The starting point here would be a characterization of strong metric sub-regularity with graphical derivatives (see Theorem \ref{GDCSMSR} and Corollary \ref{GDCIC}). 
First we assume that the linear operator $B$ in the setting (\ref{inclusion}) is injective, and in the subsection afterwards, we change this assumption with a slightly different one. \\
Consider the mappings $Q$ and $F_C$ as defined in (\ref{Q and F_C}). In order to pass from the graphical derivative of $ Q $ to the graphical derivative of $ F_C $ (which in fact, needs a special chain rule for the graphical derivative of the composition mapping like the idea in Propositions \ref{special chain rules for coderivatives}, and \ref{special chain rules for coderivatives2} for coderivatives), we need the following lemma which provides a rule on how a linear operator alters a tangent cone relation.

\begin{lem} [\textbf{Matrix Operation over Tangent Cone Relations}] \citep[Lemma 4.1, p. 96]{adly2} \label{Lemma 4.1}
Let $E \in \R^{k\times d}$ be any matrix, let $G \in \R^{l \times d} $ be injective, and let $ \Gamma $ be a subset of rge $E$. Put $ \Xi := E^{-1}(\Gamma )$ and 
$ \Lambda := G(\Xi)$. For $\bar{x} \in \Lambda $, denote by $\bar{y}$ the (unique) point in $\Xi$ with $G \bar{y} = \bar{x}$. Then 
\begin{equation*}
T(\bar{x}; \Lambda ) = \big \{ u \in \R^l~:~ \exists \, w \in \R^d \mathrm{~~such~that~} u = Gw \mathrm{~and~} Ew \in T(E \bar{y}; \Gamma ) \, \big \}.
\end{equation*}
\end{lem}

\begin{proof}
We claim that 
\begin{equation} \label{lemma-01}
T(\bar{y}; \Xi) = \{ w \in \R^d: E w \in T(E \bar{y}; \Gamma ) \}.
\end{equation}
 First, take any $ w \in T(\bar{y}; \Xi) $. Find $ (t_n)_{n \in \N} $ in $ (0, \infty ) $ converging to $0$, and $ (w_n)_{ n \in \N } $ in $\R^d$ converging $w$, such that $ \bar{y} + t_n w_n \in \Xi $  whenever $n \in \N$. Then we have that 
$$ E \bar{y} + t_n E w_n = E(\bar{y} + t_n w_n) \in \Gamma \mathrm{~~for~each~~} n \in \N. $$
 Hence $ Ew \in T(E \bar{y}; \Gamma ) $. On the other hand, let $ w \in \R^d $ be such that $ Ew \in T(E \bar{y}; \Gamma) $. Pick $ (t_n)_{n \in \N} $ in $ (0, \infty ) $ and $ (v_n)_{ n \in \N } $ in $\R^k$ converging to $0$ and $ Ew $, respectively, such that $ E \bar{y} + t_n v_n \in \Gamma $ whenever $ n \in \N $. 
 As $ \Gamma \subset \mathrm{rge \,} E $, and $ \mathrm{rge \,} E$ is a closed subspace of $\R^k$ (which is a direct result of the continuity of $E$ as a bounded linear operator on $\R^k$), one infers that $ v_n \in \mathrm{~rge \, } E $ for each $ n \in \N $. \\
Therefore, by Banach open mapping theorem there is $ (w_n)_{n \in \N} $ in $\R^d$ converging to $w$ such that $ E w_n= v_n $ for each $n \in \N$\footnote{
\label{OMT}
Indeed, the restricted map $ \widetilde{E}: ( \R^d / \mathrm{ker}\, E ) \longrightarrow \mathrm{rge} \, E $, induced by $E$, is a bijective continuous linear operator, hence from the open mapping theorem, $ \widetilde{E}^{-1} $ is continuous, too. Then, $ E w_n= v_n \longrightarrow Ew $ implies that $ w_n \longrightarrow w $.\\
}.
Thus, for an arbitrary index $n$, we have $ E(\bar{y} + t_n w_n) \in \Gamma $, hence $ \bar{y} + t_n w_n \in  E^{-1}( \Gamma ) = \Xi $. So $ w \in T( \bar{y}; \Xi) $. The claim is proved. \\
Now we prove that 
\begin{equation} \label{lemma-02}
T( \bar{x}; \Lambda ) = \{ G w : w \in T(\bar{y}; \Xi) \}, 
\end{equation}
using exactly the same steps as in the proof of Lemma 4.1 in \citep{adly3}\footnote{
Actually the Lemma 4.1 of \citep{adly3} is a similar statement with two more conditions: $E$ is assumed to be surjective and $\Gamma$ is a closed subset of $\R^n$. \\
But these two conditions are only used in the first part of their proof (referring to \citep[Exercise 6.7]{Rockafellar1998variational}) which is equal to the claim above, proved without these assumptions. Thus, this lemma could be considered as a generalization of that result.
}. \\
Let us first prove the inclusion $ \{ G w : w \in T(\bar{y}; \Xi) \} \subset T( \bar{x}; \Lambda ) $.\\
Consider $ w \in T(\bar{y}; \Xi) $. By definition of the tangent cone, there are sequences $ t_k \downarrow 0 $, and $ w_k \to w $ such that $ \bar{y} + t_k w_k \in \Xi $. Clearly, we have 
\begin{equation*}
\bar{x} + t_k G w_k = G \bar{y} + t_k G w_k \in G( \Xi ) = \Lambda
\end{equation*}
Thus, $ G w \in T( \bar{x}; \Lambda ) $. \\
It remains to prove that $ T( \bar{x}; \Lambda ) \subset \{ G w : w \in T(\bar{y}; \Xi) \} $, which can be conducted by a similar reasoning.
Let $ x \in T( \bar{x}; \Lambda ) $. By definition, there are sequences $x_k \to x $, and $ \lambda_k \downarrow 0 $ such that $ \bar{x} + \lambda_k x_k \in \Lambda $. Hence, there is a sequence $w_k$ in $\Xi$ such that for each $k$, one has $ \bar{x} + \lambda_k x_k  = G w_k $.
Thus, $ x_k = G \dfrac{w_k - \bar{y}}{\lambda_k}$.\\
Consider the sequence $( h_k )$ defined as
$$ h_k = \dfrac{w_k - \bar{y}}{\lambda_k}. $$
For each $k$, let $\lambda'_k = \lambda_k $ and observe that $\lambda'_k h_k + \bar{y} = w_k \in \Xi $, and $  \lambda'_k \downarrow 0 $. In order to show that $h_k $ has a convergent subsequence, having Bolzano-Weierstrass Theorem in mind, it only suffices to prove the boundedness of $h_k$.\\
Assume by contradiction that $ h_k $ is unbounded. Hence there exists a subsequence (still denoted by $ h_k $) such that $ \norm{ h_k } \to + \infty $ as $ k \to + \infty $. By passing to a subsequence if necessary, we have
$$ \dfrac{h_k}{\norm{h_k}} \longrightarrow s \mathrm{~~with~} \norm{s} = 1. $$ 
It follows that $ \displaystyle Gs = \lim_{ k \to \infty } \dfrac{x_k}{\norm{h_k}}  = 0 $, which contradicts the injectivity of $G$. Thus, the sequence $h_k$ possesses a convergent subsequence with the limit $h_0$ in $ T(\bar{y}; \Xi) $. We can also conclude that $ x = G h_0 \in G T(\bar{y}; \Xi) $. Hence, (\ref{lemma-02}) is proved. \\
Combining (\ref{lemma-01}) and (\ref{lemma-02}) yield the assertion.
\end{proof}

The above Lemma could be written in the following symbolic representation now, which helps us remember what is happening there. 
\begin{eqnarray}
\begin{split}
Gw \in T \big( \bar{x}; G(E^{-1} (\Gamma)) \big) & \,  \Longleftrightarrow \, G^{-1} Gw \in T \big( G^{-1} \bar{x}; G^{-1} G (E^{-1} (\Gamma ) \big) \\
& \,  \Longleftrightarrow \, w \in T \big( \bar{y}; E^{-1} (\Gamma) \big) \\
& \,  \Longleftrightarrow \, Ew \in T \big( E\bar{y}; \Gamma \big).
\end{split}
\end{eqnarray}

\begin{prop} [\textbf{Graphical Derivative of Sum for GE}] \citep[Proposition 4.1, p. 96]{adly2} \label{Proposition 4.1}
Under the assumptions (A1) - (A2), for any $b \in \R^n$ one has
\begin{equation*}
D \Phi \rfpp{z}{p} (b) = \nabla f(\bar{z}) b + B \, D F_C (C\bar{z}|\, \bar{v}) (Cb).
\end{equation*}
\end{prop}

\begin{proof}
Fix any $ b \in \R^n $. By Proposition \ref{Sum Rule for Graphical Derivatives} and Remark \ref{derivative meaning for functions}, we have 
$$ D \Phi \rfpp{z}{p} (b) = \nabla f(\bar{z}) b + D Q \big( \bar{z} | \, \bar{p} - f (\bar{z})\big) (b). $$
Observe that by the definition of $Q$ as $Q(z) := B F(Cz) $, 
\begin{equation*}
\begin{split}
\gph{Q} =   
\Bigg \{ ~&
\begin{pmatrix}
u \\ 
v
\end{pmatrix}
\in \R^{2n} ~:~ \exists 
\begin{pmatrix}
b\\ 
c
\end{pmatrix}  
\in \R^n \times \R^m 
\mathrm{~such~that~} \\
&
\begin{pmatrix}
u\\ 
v
\end{pmatrix}
= G 
\begin{pmatrix}
b\\ 
c
\end{pmatrix}  
\mathrm{~and~}
E 
\begin{pmatrix}
b\\ 
c
\end{pmatrix}  
\in \gph{F_C} ~~~~ ~~~~ \Bigg \},
\end{split} 
\end{equation*}
with 
\begin{equation*}
G := 
\begin{pmatrix}
I_n & 0 \\
0 & B
\end{pmatrix},~~ 
E := 
\begin{pmatrix}
C & 0 \\
0 & I_m
\end{pmatrix}.
\end{equation*}
As $B$ is injective, so is $G$. Using Lemma \ref{Lemma 4.1} with $ k := 2m, l := 2n, d := n + m $, $ \Gamma := \gph{F_C}, \bar{x} := \big( \bar{z}, \bar{p} - f (\bar{z}) \big) ^T $, and $ \bar{y} := \rfp{z}{v} ^T $ reveals that
\begin{equation*}
T \Big( \big( \bar{z}, \bar{p} - f (\bar{z}) \big); \, \gph{Q} \Big) =
\begin{Bmatrix}
\begin{pmatrix}
b \\
Bc
\end{pmatrix}
:
\begin{pmatrix}
Cb \\
c
\end{pmatrix}
\in T \big(  (C \bar{z}, \bar{v} ) ; \, \gph{F_C} \big)
\end{Bmatrix}.
\end{equation*}
This means that $ D Q \big( \bar{z}| \, \bar{p} - f (\bar{z})\big) (b) = B D F_C ( C\bar{z}| \, \bar{v})(Cb) $. The assertion is proved.
\end{proof}

\begin{thm} [\textbf{Isolated Calmness Criterion for the Solution Mapping of GE}] \citep[Theorem 4.1, p. 97]{adly2} \label{Theorem 4.1}
Under the assumptions (A1) - (A3), $S$ has the isolated calmness property at $ \rfp{p}{z} $ if and only if
\begin{equation}
0 \in \, \nabla f(\bar{z}) b + B \, D F_C (C\bar{z}| \, \bar{v}) (Cb)~~ \Rightarrow ~~b = 0.
\end{equation}
Moreover, its calmness modulus is given by
\begin{equation} \label{calmness modulus}
\mathrm{clm } \big(S; \rfp{p}{z} \big) = \sup \big \{ \, \left \| b \right \| \, : \left (   \nabla f(\bar{z}) b + B \, D F_C (C\bar{z}| \, \bar{v}) (Cb) \, \right ) \cap \B \neq \, \emptyset ~ \big \}.
\end{equation}
\end{thm}

\begin{proof}
Note that $ x \in  D \Phi \rfpp{z}{p} (y) $ if and only if $ y \in D S \rfpp{p}{z} (x) $. \\
Indeed, 
$ (x, y) \in T \big( \rfp{p}{z}; \gph{S} \big) $ implies the existence of sequences $ t_k \downarrow 0, (u_k , v_k ) \longrightarrow (x, y) $ such that 
$ \bar{z} + t_k v_k \in S ( \bar{p} + t_k u_k ) $.
Then, $ \bar{p} + t_k u_k \in \Phi (\bar{z} + t_k v_k) $, by definition. Thus, $ (y, x) \in T \big( \rfp{z}{p}; \gph{\Phi} \big) $ and vice versa.\\
It would be enough to combine Corollary \ref{GDCIC}, Theorem \ref{GDCSMSR}, and Proposition \ref{Proposition 4.1}, to conclude the proof. We just want to explain a bit, the modulus formula here.
\begin{eqnarray*}
\begin{split}
\mathrm{clm\,} \big(S; \rfp{p}{z} \big)  & = | D S \rfpp{p}{z} |^{+} : = \sup_{\norm{x} \, \leq \, 1} ~ \sup_{y \, \in \, D S \rfpp{p}{z} (x)} \norm{y} \\
& = \sup_{\norm{x} \, \leq \, 1} ~ \sup_{x \, \in \, D \Phi \rfpp{z}{p} (y)} \norm{y} \\
& = \sup_{\norm{x} \, \leq \, 1} ~ \sup \, \big \{ \, \norm{y} ~:~ x \in \nabla f(\bar{z}) y + B \, D F_C (C\bar{z}| \, \bar{v}) (Cy) \big \} \\
& = \sup \, \big \{ \, \norm{y} ~:~ \norm{ \nabla f(\bar{z}) y + B \, D F_C (C\bar{z}| \, \bar{v}) (Cy) } \, \leq \, 1 \big \} \\
& = \sup \, \big \{ \, \norm{y} ~:~ \Big( \nabla f(\bar{z}) y + B \, D F_C (C\bar{z}| \, \bar{v}) (Cy) \Big) \cap \B \neq \emptyset \big \}.
\end{split}
\end{eqnarray*}
The equations show the step by step process of combining the mentioned results and definition of outer norm to obtain the claimed formula.
\end{proof}

Again, imposing the additional assumption that $C$ is surjective we can go one step further and get the following statement. \\

\begin{cor}  \citep[Corollary 4.1, p. 97]{adly2} \label{Corollary 4.1}
Suppose that the assumptions (A1) - (A4) (cf. Note \ref{general assumptions}) hold true. Then $S$ has the isolated calmness property at $ \rfp{p}{z} $ if and only if
\begin{equation}\label{IC condition}
\left. \begin{matrix}
(Cb,  -(B^T B)^{-1} B^T \nabla f(\bar{z}) \, b\, ) \, \in  \, T \big( (C\bar{z},\bar{v}) ;\, \gph{F} \big) \\
\nabla f(\bar{z}) b \in \mathrm{rge } \, B~~~~~~~~~~~~~~~~~~~~~~~~~~~~~~~~~~~~~~~~~~
\end{matrix} \right\} ~~~~\Longrightarrow~~~ b= 0.
\end{equation}
\end{cor}

\begin{proof}
Indeed, if $C$ is surjective, then $F_C = F $. Note that (A1) ensures that $ B^T B \in \R^{m \times m} $ is non-singular. \\
First, let $ b \in \R^n $ be such that $ 0 \in \nabla f(\bar{z}) b + B D F (C \bar{z}, \bar{v})(Cb) $. Then, find a point $ w \in D F (C \bar{z}, \bar{v})(Cb) $ with
 $ \nabla f(\bar{z}) b + Bw = 0 $. Thus, $ -(B^T B) ^{-1} B^T \nabla f(\bar{z}) b $ is in $ D F (C \bar{z}, \bar{v})(Cb) $. Clearly, we have $  \nabla f(\bar{z}) b \in \mathrm{~rge} \,B $ and the definition of the contingent derivative of $F$ yields the rest. \\
On the other hand, pick any $ b \in \R^n $ with $ (C b, - (B^T B) ^{-1} B^T \nabla f(\bar{z}) b ) \in T \big( (C \bar{z}, \bar{v}) ; \gph{F} \big) $ and 
$  \nabla f(\bar{z}) b \in \mathrm{rge\,} B $. 
The definition of the contingent derivative yields that 
$$ w := - (B^T B)^{-1} B^T \nabla f(\bar{z})b \in D F( C \bar{z}, \bar{v})(Cb). $$ 
Thus, we have 
$ B^T B w = - B^T \nabla f( \bar{z} ) b $. So $ B w + \nabla f( \bar{z} ) b \in \mathrm{ker }\, B^T \cap \mathrm{ rge }\, B = \{ 0 \} $. Therefore, 
$ 0 \in \nabla f( \bar{z} ) b + B D F(C \bar{z},  \bar{v})(Cb) $. Using Theorem \ref{Theorem 4.1} ends the proof.
\end{proof}

\rem \label{remark-IC-condition}
If (A5) is also satisfied and $\gph{F_j} $ is \emph{Clarke regular} at $ ((C\bar{z})_j, \bar{v}_j) $ for each $ j \in \{1,...,m \} $ 
 (cf. Remark \ref{Clarke regularity}), then Proposition \ref{product rule} implies that 
 $$ T ((C\bar{z},\bar{v}) ;\gph{F} ) = \prod_{j=1} ^m T (((C\bar{z})_j,\bar{v}_j ) ;\gph{F}_j ). $$ 
Hence the first condition in (\ref{IC condition}) can be checked coordinate-wise. \\

\begin{eg} \textbf{(A Simple Circuit with DIAC)} \citep[p. 339]{adly} \label{DIAC-example-IC}\\
Consider the circuit in Figure \ref{fig: DIAC-circuit-chap03} with a DIAC. We discussed the Aubin property of the solution mapping at different points of its graph in Example \ref{DIAC-example}. Here we want to investigate the isolated calmness property. \\
If $ a < R $, the uniqueness of the solution (ref. Example \ref{DIAC-Chap2}) implies that S has also the isolated calmness property at any $ \rfp{p}{z} \in \gph{S}$. 
In order to use Corollary \ref{Corollary 4.1}, we observe that the assumptions (A1) - (A4) hold true and the condition (\ref{IC condition}) has the form
\begin{equation} \label{DIAC-example-IC-eq01}
\big( b, - f'(\bar{z})b \big) \in T \Big( \big( \bar{z}, \bar{p} - f(\bar{z}) \big); \gph{F} \Big) \, \Longrightarrow \, b = 0,
\end{equation}
with $ T \big( (x,y); \gph{F} \big) $  given by 

{\small
\begin{equation*}
~ \left \{ \begin{array}{lcl}
\R \begin{pmatrix}
1 \\ 0
\end{pmatrix}
& ~~ & x > 0,~ y = V, \mathrm{~~or~~} x < 0, ~y = - V, \\[2em]
\R \begin{pmatrix}
0 \\ 1
\end{pmatrix} 
& ~~ & x = 0,~ y \in (-V, \, V), \\[2em]
\R_{+} \begin{pmatrix}
1 \\ 0
\end{pmatrix} 
\, \bigcup \, \R_{+} \begin{pmatrix}
0 \\ -1
\end{pmatrix} 
& ~~ & x= 0, y = V, \\[2em]
\R_{+} \begin{pmatrix}
-1 \\ 0
\end{pmatrix} 
\, \bigcup \, \R_{+} \begin{pmatrix}
0 \\ 1
\end{pmatrix} 
& ~\hspace*{2cm} ~ & x= 0, y = -V, \\[2em]
\emptyset & ~~~~~~~~ & \mathrm{otherwise.}
\end{array} \right.
\end{equation*}
}

From the third and forth pieces one can not guarantee that $ b= 0 $ when $ a = R $. Hence, if $ a = R $, $S$ has the isolated calmness property at $ \rfp{p}{z} $ provided that $ \bar{p} \not \in \{-V, V \} $. \\
The case $ a > R $ is more interesting, since the solution is no more unique.
Simple calculations show that it is possible to have $f'(\bar{z}) = 0 $ for some $ \bar{z} \neq 0 $. So, the solution mapping is not isolated calm when 
$ \bar{z} < 0 $, or $ \bar{z} > 0 $.  \\
Let $ \bar{z} = 0 $. Then, $ f'(0) = R - a \neq 0 $, and \eqref{DIAC-example-IC-eq01} reveals that $S$ has the isolated calmness property at $ \rfp{p}{z} $ when $ \bar{p} \in [-V, V] $ (while we have already observed that it does not have the Aubin property if $ \bar{p} \in \{-V, V \} $). Let us compute the calmness modulus provided by Theorem \ref{Theorem 4.1} in this case, with simplifying assumptions $ R = V= 1$. Formula \eqref{calmness modulus} reads as
$$
\mathrm{clm \,} \big( S; (-1, 0) \big) = \sup 
	\begin{Bmatrix}
	\norm{b} ~:~ \Big( (1 - a ) b + D\, F(0 \, | -1) (b) \Big) \cap \B \neq \emptyset 
	\end{Bmatrix}.
$$
Since 
{\footnotesize
$  D\, F(0 \, | -1) (b) = 
\begin{Bmatrix}
v \in \R ~|~ 
\begin{pmatrix}
b \\ v 
\end{pmatrix} \in 
\begin{pmatrix}
- \lambda \\ 0 
\end{pmatrix}
\bigcup
\begin{pmatrix}
0 \\ \lambda 
\end{pmatrix}
\mathrm{~~for~some~} \lambda \in \R_{+}
\end{Bmatrix} $}, 
one obtains that the supremum of $b$ should be taken over the constraint $ -1 \, \leq \, (1-a) b + v \, \leq \, 1 $, where the points $ (b, v)^T $ belong to the previously mentioned sets. Hence,
\begin{equation*}
\mathrm{clm \,} \big( S; (-1, 0) \big) = \left \{
\begin{array}{lcl}
\dfrac{1}{| 1 - a |} &  & a \neq 1, \\[1em]
\infty & & \mathrm{otherwise.}
\end{array} \right.
\end{equation*}
\end{eg}

\vspace*{0.2cm}
\begin{eg} \textbf{(A Circuit with SCR and Zener Diode)} \citep[p. 342]{adly}\\
Consider the circuit in Figure \ref{fig: SCR and the Zener Diode} with a SCR and Zener diode. We want to investigate the isolated calmness property of the solution mapping at different points of its graph (in Example \ref{SCR with Zener-example}, we discussed the Aubin property). 
The condition \eqref{IC condition} in view of Remark \ref{remark-IC-condition} reads as
\begin{eqnarray}\label{SCR with Zener-example-IC-eq01}
\left. \begin{array}{l}
\big( \beta_1 , \, - R (\beta_1 + \beta_2) \big) \in T \Big( \big( \bar{z}_1, \, \bar{p}_1 - R (\bar{z}_1 + \bar{z}_2) \big) ; \gph{F_1} \Big)  \\[1em]
\big(  \beta_2 , \, - R (\beta_1 + \beta_2) \big) \in T \Big( \big( \bar{z}_2, \, \bar{p}_2 - R (\bar{z}_1 + \bar{z}_2) \big) ; \gph{F_2} \Big)  
\end{array} \right \} \Longrightarrow \beta_1 = \beta_2 = 0,
\end{eqnarray}
where the first tangent cone, $T \big( (x, y); \gph{F_1} \big)$, can be calculated as

{\footnotesize
\begin{eqnarray*}
~ \left \{ \begin{array}{lll}
\R \begin{pmatrix}
1 \\ a
\end{pmatrix}
&  & 
\begin{matrix}
x < 0,~ y = a x + V_1, \mathrm{~~or~~}~~~~ \\  x > \alpha, ~y = a(x - \alpha) + \varphi(\alpha)
\end{matrix}, \\[2em]
\R \begin{pmatrix}
0 \\ 1
\end{pmatrix} 
& ~~ & x = 0,~ y \in (V_1, \, \varphi(0)), \\[2em]
\R \begin{pmatrix}
1 \\ \varphi'(x) 
\end{pmatrix} 
& ~~ & x \in (0, \alpha),~ y = \varphi(x), \\[2em]
\R_{+} \begin{pmatrix}
-1 \\ -a
\end{pmatrix} 
\, \bigcup \, \R_{+} \begin{pmatrix}
0 \\ 1
\end{pmatrix} 
& ~~ & x= 0, ~y = V_1, \\[2em]
\R_{+} \begin{pmatrix}
1 \\ \varphi'(0) 
\end{pmatrix} 
\, \bigcup \, \R_{+} \begin{pmatrix}
0 \\ -1
\end{pmatrix} 
& ~~ & x= 0,~ y = \varphi(0), \\[2em]
\R_{+} \begin{pmatrix}
1 \\ a
\end{pmatrix} 
\, \bigcup \, \R_{+} \begin{pmatrix}
-1 \\ - \varphi'(\alpha) 
\end{pmatrix} 
& ~\hspace*{4cm}~~~ & x= \alpha,~ y = \varphi(\alpha), \\[2em]
\emptyset & ~~ & \mathrm{otherwise.}
\end{array} \right.
\end{eqnarray*}
}

and the second tangent cone, $T \big( (x, y); \gph{F_2} \big)$, can be computed similarly to obtain

{\footnotesize
\begin{eqnarray*}
~ \left \{ \begin{array}{lll}
\R \begin{pmatrix}
1 \\ b
\end{pmatrix}
&  & 
\begin{matrix}
x < 0,~ y = b x - V, \mathrm{~~or~~} \\  x > 0, ~y = b x + V~~~~~~~
\end{matrix}, \\[2em]
\R \begin{pmatrix}
0 \\ 1
\end{pmatrix} 
& ~~ & x = 0,~ y \in (-V, \, V), \\[2em]
\R_{+} \begin{pmatrix}
-1 \\ -b
\end{pmatrix} 
\, \bigcup \, \R_{+} \begin{pmatrix}
0 \\ 1
\end{pmatrix} 
& ~~ & x= 0,~ y = -V, \\[2em]
\R_{+} \begin{pmatrix}
1 \\ b
\end{pmatrix} 
\, \bigcup \, \R_{+} \begin{pmatrix}
0 \\ -1
\end{pmatrix} 
& ~~ & x= 0,~ y = V, \\[2em]
\emptyset & ~\hspace*{4.5cm}~ & \mathrm{otherwise.}
\end{array} \right.
\end{eqnarray*} 
}

In order to facilitate the process of checking the inclusions in (\ref{SCR with Zener-example-IC-eq01}), we multiply them by the matrices 
$ M_3:= \begin{pmatrix}
1 & 0 \\ 
-1 & \frac{-1}{R}
\end{pmatrix} $ 
and 
$ M_4:= \begin{pmatrix}
-1 & \frac{-1}{R} \\
1 & 0
\end{pmatrix} $, 
respectively, to get the following equivalent condition
\begin{equation*}
M_3 \Big( T \big( ( \bar{z}_1, \, \bar{v}_1); \gph{F_1} \big) \Big) \bigcap M_4 \Big( T \big( ( \bar{z}_2, \, \bar{v}_2); \gph{F_2} \big) \Big) = 0.
\end{equation*}
Hence, $ M_3 \Big( T \big( ( \bar{z}_1, \, \bar{v}_1); \gph{F_1} \big) \Big) $ is given by

{\footnotesize
\begin{eqnarray*}
~ \left \{ \begin{array}{lcl}
\R \begin{pmatrix}
1 \\[0.5em] \dfrac{-a}{R} - 1 
\end{pmatrix}
& & 
\begin{matrix}
\bar{z}_1 < 0,~ \bar{v}_1 = a \bar{z}_1 + V_1,~\mathrm{or}~~~~~~~~ \\  \bar{z}_1 > \alpha, ~\bar{v}_1 = a(\bar{z}_1 - \alpha) + \varphi(\alpha),
\end{matrix} \\[2.5em]
\R \begin{pmatrix}
0 \\[0.5em] \dfrac{-1}{R}
\end{pmatrix} 
& & \bar{z}_1 = 0,~ \bar{v}_1 \in (V_1, \, \varphi(0)), \\[2.5em]
\R \begin{pmatrix}
1 \\[0.5em] \dfrac{-\varphi'(\bar{z}_1)}{R} - 1
\end{pmatrix} 
& & \bar{z}_1 \in (0, \alpha), ~ \bar{v}_1 = \varphi(\bar{z}_1), \\[2.5em]
\R_{+} \begin{pmatrix}
-1 \\[0.5em] \dfrac{a}{R} + 1 
\end{pmatrix} 
\, \bigcup \, \R_{+} \begin{pmatrix}
0 \\[0.5em] \dfrac{-1}{R}
\end{pmatrix} 
& & \bar{z}_1= 0, ~\bar{v}_1 = V_1, \\[2.5em]
\R_{+} \begin{pmatrix}
1 \\[0.5em]  \dfrac{- \varphi'(0)}{R} - 1
\end{pmatrix} 
\, \bigcup \, \R_{+} \begin{pmatrix}
0 \\[0.5em] \dfrac{1}{R}
\end{pmatrix} 
& & \bar{z}_1= 0, ~\bar{v}_1 = \varphi(0), \\[2.5em]
\R_{+} \begin{pmatrix}
1 \\[0.5em] \dfrac{-a}{R} - 1
\end{pmatrix} 
\, \bigcup \, \R_{+} \begin{pmatrix}
-1 \\[0.5em] \dfrac{\varphi'(\alpha)}{R} + 1
\end{pmatrix} 
& & \bar{z}_1= \alpha, ~ \bar{v}_1 = \varphi(\alpha), \\[2.5em]
\emptyset & \hspace*{2.5cm} &  \mathrm{otherwise.}
\end{array} \right.
\end{eqnarray*}
}
And $M_4 \Big( T \big( (\bar{z}_2, \bar{v}_2); \gph{F_2} \big) \Big)$ is given by
{\footnotesize
\begin{eqnarray*}
~ \left \{ \begin{array}{lcl}
\R \begin{pmatrix}
\dfrac{-b}{R} - 1  \\[0.7em] 1
\end{pmatrix}
&  & 
\begin{matrix}
\bar{z}_2 < 0,~ \bar{v}_2 = b \bar{z}_2 - V, \mathrm{~~or~~} \\  \bar{z}_2 > 0, ~ \bar{v}_2 = b \bar{z}_2 + V~~~~~~~
\end{matrix}, \\[2.5em]
\R \begin{pmatrix}
\dfrac{-1}{R}  \\[0.7em] 0
\end{pmatrix} 
&  & \bar{z}_2 = 0,~ \bar{v}_2 \in (-V, \, V), \\[2.5em]
\R_{+} \begin{pmatrix}
\dfrac{b}{R} + 1  \\[0.7em] -1
\end{pmatrix} 
\, \bigcup \, \R_{+} \begin{pmatrix}
\dfrac{-1}{R}  \\[0.7em] 0
\end{pmatrix} 
&  & \bar{z}_2= 0,~ \bar{v}_2 = -V, \\[2.5em]
\R_{+} \begin{pmatrix}
\dfrac{-b}{R} - 1 \\[0.7em] 1
\end{pmatrix} 
\, \bigcup \, \R_{+} \begin{pmatrix}
\dfrac{1}{R} \\[0.7em] 0
\end{pmatrix} 
&  & \bar{z}_2= 0, ~\bar{v}_2 = V, \\[2.5em]
\emptyset & ~\hspace*{3cm}~ & \mathrm{otherwise.}
\end{array} \right.
\end{eqnarray*} 
}

One can easily check each sub-domain of $ M_3 \Big( T \big( (\bar{z}_1, \bar{v}_1); \gph{F_1} \big) \Big) $ with sub-domains of 
$ M_4 \Big( T \big( (\bar{z}_2, \bar{v}_2); \gph{F_2} \big) \Big) $ to see if they have a non-zero intersection. The result is that the isolated calmness of the solution map depends on the parameters of a particular circuit, and the reference point. \\
For example, in the case of 
$ \bar{z}_1 < 0,~ \bar{v}_1 = a \bar{z}_1 + V_1 $ or $ \bar{z}_1 > \alpha, ~\bar{v}_1 = a(\bar{z}_1 - \alpha) + \varphi(\alpha) $, and for any 
$(\bar{z}_2, \bar{v}_2) \in \gph{F_2} $, the solution map has the isolated calmness property. While for $\bar{z}_1 \in (0, \alpha), ~ \bar{v}_1 = \varphi(\bar{z}_1)$, the condition \eqref{SCR with Zener-example-IC-eq01} will be violated for $ \bar{z}_2 \neq 0 $ if $ \varphi'(\bar{z}_1) = \dfrac{- \, R \, b}{ R + b }$; and for 
$ \bar{z}_2 = 0,~ \bar{v}_2 \in (-V, \, V) $ if $  \varphi'(\bar{z}_1) = - R $; and for $ \bar{z}_2 = 0,~ \bar{v}_2 \in \{-V, \, V\} $ if 
$  \varphi'(\bar{z}_1) \in \Big \{ - R, \, \dfrac{- \, R \, b}{ R + b } \Big \} $.
\end{eg}

\subsection{Isolated Calmness Without Injectivity Assumption} \label{Isolated Calmness Without Injectivity Assumption}
Up to now, we have assumed that $ m \leq n$. However, one can have $m > n$ in some applications. In such a case, one cannot expect the assumptions (A1) and (A4) (cf. Note \ref{general assumptions}) to hold true. So till the end of this section, $ m, n \in \N$ are not related to each other by an inequality. 
Also, in order to avoid confusion during the calculations of this subsection, we refer to the zero vector of different spaces with the space indicated as a subscript (like $ 0_{\R^d} $). We assume the following assumption (instead of (A1)):
\begin{itemize}[nolistsep]
\item[]$(\widetilde{\mathrm{A1}})$ Suppose that there is $ \bar{v} \in F(C\bar{z}) $ such that
\begin{equation*}
\bar{p} = f(\bar{z}) + B \, \bar{v} ~\mathrm{~~and~~}~ \overline{ \bigcup_{t > 0} \dfrac{ \mathrm{rge~} F_C - \bar{v} } {t} } \,  \bigcap \mathrm{~ker~} B = \{ 0_{\R^m} \}.
\end{equation*}
\end{itemize}
Of course, if $\bar{v}$ is an interior point of rge $F_C$ then $(\widetilde{\mathrm{A1}})$ reduces to (A1). However, when $ \bar{v} $ is in the boundary of $ \mathrm{rge~} F_C $, this is not true. For example, consider $ \mathrm{rge\,} F_C = (-3, -1] \subset \R $. Then, for $ \bar{v} = -2 $, we have 
$ \mathrm{rge\,} F_C  - \bar{v} = (-1, 1] $, and 
$$ \overline{ \bigcup_{t > 0} \dfrac{ \mathrm{rge~} F_C - \bar{v} } {t} } \,  \bigcap \mathrm{~ker~} B = \R ~ \cap \mathrm{~ker\,} B, $$
which implies $(\widetilde{\mathrm{A1}})$  being equivalent with (A1). While for $ \bar{v} = -1 $, we have $ \mathrm{rge\,} F_C  - \bar{v} = (-2, 0] $, and 
$$ \overline{ \bigcup_{t > 0} \dfrac{ \mathrm{rge~} F_C - \bar{v} } {t} } \,  \bigcap \mathrm{~ker~} B = \R_{-} ~ \cap \mathrm{~ker\,} B, $$
which is not necessarily equal to $ \mathrm{ker\,} B $.

\begin{lem}  \citep[Lemma 2, p. 353]{adly} \label{Lemma2}
Let $ G \in \R^{l \times d} $, let $ \Xi \subset \R^d $, and set $ \Lambda = G(\Xi )$. Suppose that $\bar{x} \in \Lambda $ and $\bar{y} \in \Xi $ are such that
\begin{equation} \label{lem-03} 
G(\bar{y}) = \bar{x}~~\mathrm{~and~}~~ \overline{ \bigcup_{t > 0} \dfrac{ \Xi - \bar{y} } {t} } \,  \bigcap \mathrm{~ker~} G = \{ 0_{\R^d} \}.
\end{equation}
Then
\begin{equation*}
T ( \bar{x}; \, \Lambda ) = G \left ( T ( \bar{y}; \, \Xi ) \right ).
\end{equation*}
\end{lem}

\begin{proof}
To prove that $G \big(T(\bar{y}; \, \Xi) \big) \subset T(\bar{x}; \, \Lambda ) $, pick any $ w \in G \big(T(\bar{y}; \, \Xi) \big)$. Find $ v \in T( \bar{y}; \, \Xi) $ with 
$G(v) = w $. Thus there is $(t_n)_{n \in \N} $ in $ (0, \infty) $ converging to $0$ and $(v_n)_{n \in \N} $ in $\R^d$ converging to $v$ such that $ \bar{y} + t_n v_ n \in \Xi $ whenever $ n \in \N $. For each $ n \in \N $, put $w_n := G(v_n) $. Clearly, $ (w_n)_{n \in \N } $ converges to $w$. Moreover,
$$ x + t_n w_n = G( \bar{y} + t_n v_n) \in G(\Xi) = \Lambda \mathrm{~~whenever~} n \in \N. $$
So $ w \in T(\bar{x}; \, \Lambda) $. \\
To see the opposite inclusion, pick any $ w \in T(\bar{x}; \, \Lambda) $. Find $ (t_n)_{ n \in \N} $ in $ (0, \infty) $ converging to $0$ and $(w_n)_{n \in \N}$ in $\R^l$ converging to $w$ such that
$$ \bar{x} + t_n w_n \in \Lambda \mathrm{~or~} G(\bar{y}) + t_n w_n \in G(\Xi) \mathrm{~~for~each~} n \in \N. $$
For each $n \in \N$, find $ v_n \in \Xi $ such that $ w_n = G \big( ( v_n - \bar{y}) / t_n \big) $, and put $ u_n = (v_n - \bar{y}) / t_n $. We claim that $ (u_n)_{n \in \N} $ is bounded. \\
Suppose on the contrary that this is not the case. Passing to a subsequence if necessary, we may assume that $ (u_n / \norm{u_n} )_{n \in \N} $ converges, to 
$ \bar{u} \in S_{\R^d} $, the unit sphere in $\R^d$. For each $n \in \N$, one has that
\begin{equation*}
\dfrac{u_n}{\norm{u_n}} = \dfrac{v_n - \bar{y}}{\norm{v_n - \bar{y}}} \mathrm{~~~and~~~} \dfrac{w_n}{\norm{u_n}} = G \bigg( \dfrac{u_n}{\norm{u_n}} \bigg).
\end{equation*}
Passing to the limit for $ n \rightarrow \infty $ one gets, that 
$$ \bar{u} \in \overline{ \bigcup_{t > 0} \dfrac{ \Xi - \bar{y} } {t} } \mathrm{~~~~and~~~~} G(\bar{u}) = \lim_{ n \to \infty } \dfrac{w_n}{\norm{u_n}} = 0_{\R^l},$$ 
which, by condition (\ref{lem-03}) implies $\bar{u} = 0_{\R^d} $, a contradiction since $ S_{\R^d} \ni \bar{u} \neq 0_{\R^d} $. \\
Having the claim in hand, one may assume without loss of generality that $(u_n)_{n \in \N} $ converges to some $ u \in  \R^d $. For each $n \in \N$, we have $ \bar{y} + t_n u_n = v_n \in \Xi $, therefore $ u \in T( \bar{y}; \, \Xi ) $. Moreover, $ w = G(u)$. The proof is finished.
\end{proof}

\begin{prop}  \citep[Proposition 1, p. 354]{adly} \label{Proposition 1}
Assume that ($\widetilde{\mathrm{A1}}$) and (A2) hold true. Then
\begin{equation*}
D \Phi \rfpp{z}{p} (b) = \nabla f(\bar{z}) b + B \, D F_C (C\bar{z}| \, \bar{v}) (Cb) ~~\mathrm{~~whenever~~}~ b \in \R^n.
\end{equation*}
\end{prop}

\begin{proof}
Fix any $ b \in \R^n $. By Proposition \ref{Sum Rule for Graphical Derivatives} and Remark \ref{derivative meaning for functions}, we have
$$ D \Phi \rfpp{z}{p} (b) = \nabla f(\bar{z}) b + D Q \big( \bar{z}| \, \bar{p} - f(\bar{z}) \big) (b). $$
Define $ \mmap{H}{n}{m} $ by $ H(z) := F (C z) $, for $ z \in \R^n $. Then, 
\begin{equation*}
\gph{Q} = G ( \gph{H}) \mathrm{~~~with~~~} G : = 
\begin{pmatrix}
I_n & 0 \\
0 & B
\end{pmatrix}
\end{equation*}
Then $ \mathrm{ker} \, G = \{ 0_{\R^n} \} \times \mathrm{ker} \, B $ and $ \gph{H} \subset \R^n \times \mathrm{rge} \, F_C $. Therefore
\begin{equation*}
\overline{ \bigcup_{t > 0} \dfrac{ \gph{H} - \rfp{z}{v} } {t} } \, \bigcap \mathrm{ker \,} G \, \subset \, \{ 0_{\R^n} \} \times \Bigg( 
\overline{ \bigcup_{t > 0} \dfrac{ \mathrm{rge \,} F_C - \bar{v} } {t} } \, \bigcap \mathrm{ker \,} B  \Bigg) = \{ 0_{\R^{n + m}} \}.
\end{equation*}
Applying Lemma \ref{Lemma2} with $ l := 2n$, $d := n + m $, $ \Xi := \gph{H}$, $\bar{x} = ( \bar{z}, \bar{p} - f(\bar{z})) $, and $ \bar{y} = \rfp{z}{v} $ reveals that
\begin{equation*}
T\Big( (\bar{z}, \bar{p} - f(\bar{z})); \gph{Q} \Big) = 
\begin{Bmatrix}
\begin{pmatrix}
b \\
B c
\end{pmatrix}
~ : ~
\begin{pmatrix}
b\\
c
\end{pmatrix}
\in T \big( \rfp{z}{v}; \gph{H} \big)
\end{Bmatrix}.
\end{equation*}
This means that $ D Q \big( \bar{z}, \bar{p} - f(\bar{z}) \big)(b) = B D H \rfp{z}{v} (b) $. Moreover,
\begin{equation*}
\gph{H} = 
\begin{Bmatrix}
\begin{pmatrix}
z \\
v
\end{pmatrix}
~ \in \R^n \times \R^m ~ : ~ E
\begin{pmatrix}
z\\
v
\end{pmatrix}
\in \gph{F_C}
\end{Bmatrix}
\mathrm{~~~with~~~} E := 
\begin{pmatrix}
C & 0 \\
0 & I_m
\end{pmatrix}.
\end{equation*}
Finally, using the claim in the first part of the proof of Lemma \ref{Lemma 4.1} with $ k := 2m$, $d := n + m$, $ \Gamma := \gph{F_C} $, and $ \bar{y} = \rfp{z}{v} $ reveals that
\begin{equation*}
T\Big( \rfp{z}{v} ; \gph{H} \Big) = 
\begin{Bmatrix}
\begin{pmatrix}
b \\
c
\end{pmatrix}
~ \in \R^n \times \R^m ~ : ~ E
\begin{pmatrix}
b\\
c
\end{pmatrix}
\in T \big( (C \bar{z}, \bar{v} ); \gph{F_C} \big)
\end{Bmatrix}.
\end{equation*}
Hence $ D H \rfpp{z}{v} (b) = D F_C (C \bar{z}| \, \bar{v})(Cb) $. Combining the above facts ends the proof.
\end{proof}

\begin{thm} \citep[Theorem 3, p. 355]{adly} \label{IC2}
Suppose that ($\widetilde{\mathrm{A1}}$), (A2), and (A3) hold true. Then $S$ has the isolated calmness property at $ \rfp{p}{z} $ if and only if
\begin{equation*}
0 \in \, \nabla f(\bar{z}) b + B \, D F_C (C\bar{z}| \, \bar{v}) (Cb)~~ \Rightarrow ~~b = 0_{\R^n}.
\end{equation*}
\end{thm}

\begin{proof}
The proof is similar to the proof of Theorem \ref{Theorem 4.1}, and could be concluded by combining Corollary \ref{GDCIC}, Theorem \ref{GDCSMSR}, and Proposition \ref{Proposition 1} above; noticing that $ x \in D \Phi \rfpp{z}{p}(y) $ if and only if $ y \in D S \rfpp{p}{z} (x)$.
\end{proof}

To use the above statement, one has to know the range of the matrix $C$. Sometimes the following sufficient condition may be useful.\\

\begin{cor} \citep[Corollary 2, p. 355]{adly} \label{IC without A1-cor}
Under the assumptions of Theorem \ref{IC2}, $S$ has the isolated calmness property at $\rfp{p}{z}$, if $ b = 0_{\R^n} $ is the only point that
\begin{equation} \label{IC without A1-cor-formula}
\nabla f(\bar{z}) b + B w = 0 ~~\mathrm{~and~}~~(Cb, w) \in T \Big( (C\bar{z}, \bar{v}); \, \gph{F} \Big) ~~~\mathrm{for~some~~}~~ w \in \R^m.
\end{equation}
\end{cor}

\begin{proof}
Clearly, $ \gph{F_C} \subset  \gph{F} $, so $ T \Big( (C \bar{z}, \bar{v}); \gph{F_C} \Big) \subset T \Big( (C \bar{z}, \bar{v}); \gph{F} \Big) $. Take any $ b \in \R^n $ such that $ 0 \in \nabla f(\bar{z}) b + B D F_C (C \bar{z}, \bar{v}) (Cb) $. Find $ w \in \R^m $ such that 
$ 0 = \nabla f(\bar{z}) b + B w $ and $ (Cb, w) \in T \Big( (C \bar{z}, \bar{v}); \gph{F_C} \Big) $. Thus $ b = 0_{\R^n} $.
\end{proof}

Let us finish this subsection with an example of a circuit that exhibits the failure of the assumption (A1), yet with the tools provided here under the assumption ($\widetilde{\mathrm{A1}}$), we can investigate the isolated calmness property.

\begin{eg} \textbf{(Sampling Gate)} \label{IC without A1} \citep[page 355]{adly}\\
Consider the circuit in Figure \ref{fig: sampling-gate}, composed of four diodes $ D_1, D_2, D_3, D_4 $ which are controlled symmetrically by gate voltages $ +V_c $ and $ -V_c $, and the control resistors $ R_c > 0 $. Assume that all the diodes have the same characteristics given for $ V_{D1} < 0 < V_{D2} $ by
\begin{figure}[ht]
  \begin{minipage}{.35\textwidth}
  		\begin{eqnarray*}
		~~ F_D (z) :=
		\left \{ \begin{matrix}
		V_{D1} ~~~~~~ & & z < 0,\\
		[V_{D1} , V_{D2}]  & &  z = 0,\\
		V_{D2} ~~~~~~ & & z > 0.
		\end{matrix} \right.
		\end{eqnarray*} 
		\vspace*{0.2cm}
  \end{minipage}
  \begin{minipage}{.65\textwidth}
    	\centering
    	$~~$
		\includegraphics[width=0.91\textwidth]{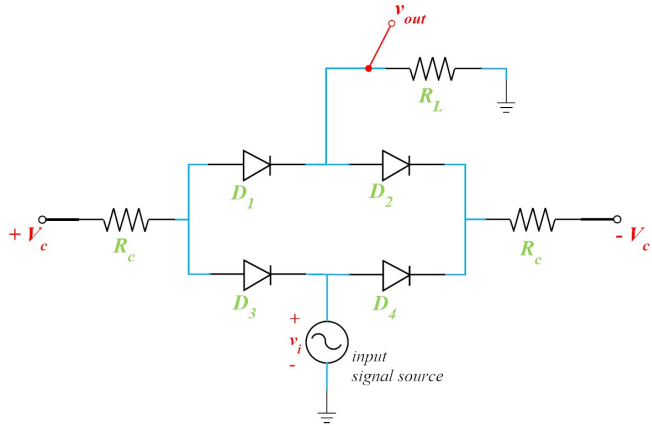}
	\end{minipage}
	\caption{Calculating isolated calmness in a circuit without the (A1) assumption}
	\label{fig: sampling-gate}
\end{figure}

In Example \ref{Sampling Gate-chap02}, we obtained the mathematical model for this circuit in the form of \eqref{inclusion} with $ m = 4 $, $ n = 3 $, 
$ p = (v_i, \, 2 V_c, \, 0 )^T $, $ C = B^T $, 
$ f(z) = A z $ for $ z \in \R^3 $, where 
\begin{equation*}
A = \begin{pmatrix}
R_L & 0 & 0 \\
0 & 2 \, R_c & 0 \\
0 & 0 & 0
\end{pmatrix},~ 
B = \begin{pmatrix}
0 & -1 & 0 & 1 \\
0 & 0 & 1 & 1 \\
1 & 1 & -1 & -1
\end{pmatrix},
\mathrm{~and~}
F(y) = 
{\small \begin{pmatrix}
F_D(y_1) \\
F_D(y_2) \\
F_D(y_3) \\
F_D(y_4)
\end{pmatrix} }
\mathrm{~for~} y \in \R^4.
\end{equation*}
Since $B$ is not invertible (in fact, $ \mathrm{ker\,} B = \big \{ \lambda \, (-1, \, 1, \, -1, \, 1)^T ~|~ \lambda \in \R \big \} $), we can not use Theorem \ref{Theorem 4.1} or Corollary \ref{Corollary 4.1}. However, for the points 
$ \rfp{z}{v} \in \gph{F_C} $  with $ \bar{v} $ in the boundary of $ \mathrm{rge\,} F_C $, we can check whether the assumptions of Theorem \ref{IC2} are satisfied and then, use Corrolary \ref{IC without A1-cor}. 
Let $ \bar{p} =(0, \, 2 V_{D1}, \, 0 )^T $ and $ \bar{z} = (0, 0, 0)^T $. Put $ \bar{v} = (V_{D1}, \, V_{D1}, \, V_{D1}, \, V_{D1})^T $. Then 
$$ \rfp{z}{p} \in \gph{\Phi}, ~~~~~~ \bar{p} = A \bar{z} + B \bar{v} \mathrm{~~~~ and~~~~} \bar{v} \in F(C \bar{z}). $$ 
Moreover, $ \mathrm{rge\,} F_C \, \subset \,  \mathrm{rge\,} F = [V_{D1}, \, V_{D2}]^4 $. Therefore
\begin{equation*} 
\overline{ \bigcup_{t > 0} \dfrac{ \mathrm{rge~} F_C - \bar{v} } {t} } \,  \bigcap \mathrm{~ker~} B ~ \subset ~ \R_{+}^4 ~ \bigcap ~ 
\R \, \Big \{  (-1, \, 1, \, -1, \, 1)^T \Big \} = \{ 0_{\R^4} \}. 
\end{equation*}
Hence $(\widetilde{\mathrm{A1}})$ holds. 
In order to check criteria \eqref{IC without A1-cor-formula}, let $ b \in \R^3 $ be such that
$$ A b + B w = 0_{\R^3} \mathrm{~~and~~} (Cb, w) \in T \big( (0, \bar{v}); \, \gph{F} \big) \mathrm{~~for~some~} w \in \R^4. $$
Clearly, $ T \big( (0, \bar{v}); \, \gph{F} \big) $ equals to
\begin{equation*}
\prod_{j = 1}^{4} T \big( (0, V_{D1}); \, \gph{F_D} \big) = 
\begin{bmatrix}
	\begin{Bmatrix}
		\begin{pmatrix}
		x \\ y
		\end{pmatrix} \in \R^2 ~|~ x \leq 0, \, y \geq 0, \, x y = 0  
	\end{Bmatrix}
\end{bmatrix}^4.
\end{equation*}
Therefore, one should have $ \inp{Cb}{w} = 0 $. But 
\begin{eqnarray*}
\begin{split}
 \inp{Cb}{w} & = \inp{B^T \, b}{w} = \inp{b}{Bw} = \inp{b}{-Ab} = 
 \left \langle 
 {\footnotesize \begin{pmatrix}
 b_1\\b_2\\b_3
 \end{pmatrix},~
\begin{pmatrix}
- R_L \, b_1 \\ - 2 \, R_c \, b_2 \\ 0
 \end{pmatrix}} \right \rangle  \\
& = - R_L \, b_1 ^2 - 2 \, R_c \, b_2 ^2,
\end{split}
\end{eqnarray*}
which implies that $ b_1 = b_2 = 0 $. Thus, $ C b = \big( b_3, \, b_3, \, - b_3, \, - b_3 \big)^T $. From the tangent cone relation, we know that all the coordinates of this vector should be non-positive, so $ b_3 = 0 $.  Applying Corollary \ref{IC without A1-cor} one concludes that $S$ has the isolated calmness property at the reference point.
\end{eg}

\section{Results about Calmness} \label{Results about Calmness}
Although we know that Aubin property or isolated calmness implies calmness, up to the results obtained so far, the criteria guaranteeing the calmness are more complicated than those for the isolated calmness or Aubin property. Hence, one should always consider the latter property first, and employ more sophisticated tools of this section only in case of their failure. \\ 
One reason for this difficulty is the lack of stability under perturbation for calmness. To be more clear, 
suppose that $S$ has the Aubin property at $\rfp{p}{z} $. Then, there exists $ \gamma > 0 $ such that for every function $ \smap{\tilde{f}}{n}{n} $, every point 
$ ( \tilde{p}, \tilde{z}) \in \B_{\gamma} ( \bar{p}) \times \B_{\gamma} ( \bar{z}) $ with 
\begin{eqnarray*}
\begin{split}
& \tilde{p} \in \tilde{f}(\tilde{z}) + B F ( C \tilde{z} ),~~ \norm{ \tilde{f}(\tilde{z}) - f( \tilde{z}) } \, \leq \, \gamma, \mathrm{~and~} \\
& \norm{ \Big( \tilde{f}(z') - f( z') \Big) - \Big( \tilde{f}(z'') - f( z'')  \Big)} \, \leq \, \gamma \norm{ z' - z'' } \mathrm{~~whenever~~} z', z'' \in  \B_{\gamma} ( \bar{z}), 
\end{split}
\end{eqnarray*}
the mapping $\widetilde{S} := ( \tilde{f} + BF ( C \cdot ) )^{-1} $ has the Aubin property at $( \tilde{p}, \tilde{z})$ (cf. \citep[Theorem 1]{dontchev2009}). A similar statement can be said about isolated calmness (cf. \citep{dontchev1995characterizations}), but not for calmness. \\
In fact, in Example \ref{MSR-perturbing problem} we have observed that 
we cannot expect to find an equivalent characterization of metric sub-regularity by means of a derivative-like object which is really computable in general, because any such derivative would have the property that the derivatives of the multifunctions $ G $ and $\widetilde{G}$ of that example are different. But it is possible that $ G $ and $ \widetilde{G} $ differ only by a $ C^{\infty} $-function, where all derivatives vanish at $\bar{x}$ and so the usual calculus rules cannot be valid. \\
We would provide a sufficient condition (Proposition \ref{Proposition 5.1}) for calmness and an estimate for calmness modulus in this section. In order to prove that proposition, we need some tools connecting the calmness of the solution mapping in the setting (\ref{inclusion}) to the outer subdifferential of an auxiliary function.

\begin{lem} \citep[Lemma 1, p. 441]{henrion2005} \label{Lemma 1} 
Let $ X, U, Y $ be normed spaces. Consider multifunctions 
$ M : X \rightrightarrows Y $ defined by 
\begin{equation*}
M(y) := \big \{ x \in X ~:~ g(x) + y \in \Lambda \big \},
\end{equation*}
where $ g : X \longrightarrow Y $ and $ \Lambda \subseteq Y $ is a closed subset; and $ M^{*} : U \rightrightarrows X $ defined on the basis of some locally Lipschitzian (with respect to the product topology) function $ h : U \times X \longrightarrow Y $ by means of
\begin{equation*}
M^{*} (u) :=  \{ x \in X ~|~ h(u, x) \in \Lambda \}.
\end{equation*}
Assume that $ h \rfp{u}{x} \in \Lambda $ for some $ \bar{u} \in U $ and $ \bar{x} \in X $. Then, $ M^{*} $ is calm at $ \rfp{u}{x} $ provided that $M$ is calm at
 $ (0, \bar{x}) $  with $ g(x) := h( \bar{u}, x)  \mfa x \in X $.
\end{lem}

\begin{proof} 
The local Lipschitz continuity of $h$ and the calmness of $M$ yield constants $ K, L, \epsilon > 0 $ such that
\begin{eqnarray*}
\begin{split}
\norm{ h(u', x) - h(u'', x) } & \, \leq \, K \norm{ u' - u'' } ~ \mfa u', u'' \in \B_{\epsilon} ( \bar{u} ), \mfa x \in \B_{\epsilon} ( \bar{x} ), \\
d \big( x, M(0) \big) & \, \leq \, L \norm{y} ~~~~~ \mfa y \in \B_{\epsilon} ( 0 ), \mfa x \in \B_{\epsilon} ( \bar{x} ) \cap M(y).
\end{split}
\end{eqnarray*}
Choose $ \epsilon ' $  such that $ 0 < \epsilon' \leq \epsilon $ and 
$$ \norm{ h(u, x) - h( \bar{u}, x) } \, \leq  \epsilon \mfa (u,x) \in \B_{\epsilon '} ( \bar{u} ) \times \B_{\epsilon '} ( \bar{x} ) . $$
Let $ x \in M^{*} (u) \cap \B_{\epsilon'} ( \bar{x} ) $, and $ u \in \B_{\epsilon ' } ( \bar{u} ) $, be arbitrary. 
Then, $ x \in M \big( h(u, x) - g(x) \big) \cap \B_{\epsilon '} ( \bar{x} ) $, by definition of $M$ and $ M^{*} $. It follows the calmness of $ M^{*} $ at $ \rfp{u}{x} $:
$$ d \big( x, M^{*} ( \bar{u}) \big) = d \big( x, M(0) \big) \, \leq \,  L \norm{h(u, x) - g(x)} \, \leq \, L K  \norm{ u - \bar{u} }. $$
\end{proof}

\begin{thm} \citep[Theorem 2.1, p. 203]{ioffe2008metric} \label{Theorem 2.1} \\
Let $ f : \R^n \longrightarrow \overline{\R} $ be lower semicontinuous in a neighborhood of $ \bar{x} \in \R^n $ and $ f (\bar{x}) = 0 $. Fix a $ \bar{\gamma} > 0 $ and consider the following properties:
\begin{enumerate}[topsep=-1ex,itemsep=-1ex,partopsep=1ex,parsep=1ex]
\item[(a)] For any $ \gamma < \bar{\gamma} $ there is a $ \delta > 0 $ such that $ d (x, [ f \leq 0]) \leq \gamma^{-1} f^{+} (x) $ if $ \norm{ x - \bar{x} } < \delta $;
\item[(b)] For any $ \gamma < \bar{\gamma} $ there is a $ \delta > 0 $ such that $ \liminf \norm{h_k} ^{-1} f (x_k + h_k) \geq \gamma $ whenever
$ f (x_k) \leq 0, \norm{x_k - \bar{x} }\leq \delta $ and $ h_k \to 0, h_k \in N_F ([  f \leq 0], x_k) \setminus \{ 0 \} $;
\item[(c)] For any $ \gamma < \bar{\gamma} $ there is a $ \delta > 0 $ such that $ \norm{ x^* } \geq \gamma $ if $ x^* \in \partial f (x) $ for some $x$ satisfying 
$ \norm{ x - \bar{x} } < 2 \delta $ and $ 0 < f (x) < \delta \gamma $;
\item[(d)] $ \norm{ x^* } \geq \bar{\gamma} $ if $ x^* \in \partial_{>} f(\bar{x}) $;
\item[(e)] For any $ \gamma < \bar{\gamma} $ there is a $ \delta > 0 $ such that $ | \nabla f | (x) \geq \gamma $ if $ \norm{ x - \bar{x} } < 2 \delta $ and $ 0 < f (x) < \delta \gamma $.
\end{enumerate}
Then, (e) $\Rightarrow$ (d) $\Leftrightarrow$ (c) $\Rightarrow$ (a) $\Leftarrow$ (b).\\
\end{thm}

\begin{prop} [\textbf{Calmness Sufficient Criterion for the Solution Mapping of GE}] \citep[Proposition 5.1, p. 98]{adly2} \label{Proposition 5.1}\\
Suppose that the assumptions (A1) - (A3) are satisfied. Put $ \Lambda = \gph {~F} \times \mathrm{rge~} B$ and define the functions\\
 $ ~~~ g : \R^n \rightarrow \R^m \times \R^m \times \R^n ~$ by $~ g(z) := \left ( Cz, \, (B^T B)^{-1} B^T ( \bar{p} - f(z)), \,  \bar{p} - f(z) \, \right ) $, and\\
 $ ~~~ \smap{h}{n}{+} ~$ by $~ h(z) := d \, (g(z), \, \Lambda ) ~~$ for each $~ z \in \R^n $.\\
 Then $S$ has the calmness property at $\rfp{p}{z}$, provided that 
\begin{equation} \label{calmness criterion} 
0 \, \not \in \, \partial_> h(\bar{z}).
\end{equation} 
\end{prop}

\begin{proof} 
The inclusion (\ref{inclusion}) says that, for each $ p \in \R^n $, one has
\begin{equation*}
S(p) = \Big \{ z \in \R^n ~:~ \Big(C z, (B^T B)^{-1} B^T \big( p - f(z) \big), p - f(z) \Big) \in \Lambda \Big \}.
\end{equation*}
Indeed, fix any $ p \in \R^n $. First, pick any $ z \in \R^n $ with $ (z, p) \in \gph{\Phi} $. Inclusion (\ref{inclusion}) reveals that $ p - f(z) \in \mathrm{\, rge \,} B $. We have already mentioned that $ B^T B \in \R^{m \times m} $ is non-singular. Hence, (\ref{inclusion}) implies that $ \Big(C z, (B^T B)^{-1} B^T ( p - f(z)) \Big) \in \gph{F} $. On the other hand, fix an arbitrary $ z \in \R^n $ with 
$$ \Big(C z, (B^T B)^{-1} B^T \big( p - f(z) \big), p - f(z) \Big) \in \Lambda. $$
Then $ B^T (p - f(z) ) = B^T B w $ for some $ w \in F ( C z ) $. As $ p - f(z) \in \mathrm{\, rge \,} B $, we have $ p - f(z) - B w \in \mathrm{\,ker \,}B^T \cap \mathrm{\, rge \, } B = \{ 0 \} $. Therefore $ (z, p) \in \gph{\Phi} $.\\
Define the mapping $ M : \R^m \times \R^m \times \R^n \rightrightarrows \R^n $ as follows:
\begin{equation*}
M(y) := \big \{ z \in \R^n ~:~ g(z) + y \in \Lambda \big \}, ~~~~ y \in \R^m \times \R^m \times \R^n. 
\end{equation*}
 As $\Lambda$ is closed and the mapping
 \begin{equation}
 \R^n \times \R^n \ni (p, z) \mapsto \Big(C z, (B^T B)^{-1} B^T \big( p - f(z) \big), p - f(z) \Big) ~ \in \R^m \times \R^m \times \R^n
 \end{equation}
is continuously differentiable (hence strictly differentiable and therefore locally Lipschitz continuous), Lemma \ref{Lemma 1} says that $S$ is calm at $\rfp{p}{z} $, provided that so is $M$ at $ (0, \bar{z}) $. As observed in \citep[p. 438]{henrion2005}, $M$ is calm at $ (0, \bar{z}) $ if and only if there are $ L > 0$ and $ \epsilon > 0 $ such that
\begin{equation*}
d \big(z, g^{-1}( \Lambda ) \big) \, \leq \,  L \, d \big( g(z), \Lambda \big) \mathrm{~~~whenever~~} z \in \B ( \bar{z}, \epsilon ).
\end{equation*}
Since $g$ is locally Lipschitz continuous, so is $h$. Moreover, $ [h \leq 0] := \{ u \in \R^n ~:~ h(u) \leq 0 \} = g^{-1} (\Lambda) $. Theorem \ref{Theorem 2.1} reveals that $M$ is calm at $ (0, \bar{z}) $, provided that (\ref{calmness criterion}) holds. The proposition is proved.
\end{proof}

\begin{rem} \textbf{(Estimate for Calmness Modulus)} \label{Estimate for Calmness Modulus} \\
Clearly, if the matrix $B$ is surjective, it suffices to consider $ \Lambda = \gph{~F} $, and $ g(z) = \left( Cz,B^{-1} (\bar{p} - f(z)) \right ),~ z \in \R^n$. Furthermore, 
Theorem \ref{Theorem 2.1} also gives an upper estimate of the corresponding calmness modulus.\\ 
Namely, if there is $ \gamma > 0 $ such that $ \left \| \xi \right \| \geq \gamma \mathrm{~for~each~} \xi \in \partial_> h(\bar{z}) $, then $S$ is calm at $\rfp{p}{z}$ with modulus not exceeding $ K/ \gamma $, where $K > 0$ is the Lipschitz constant of the map $\varphi$ at $\rfp{p}{z}$. 
$\varphi$ is defined in the proof of the Proposition \ref{Proposition 5.1} as 
\begin{equation*}
\R^n \times \R^n \ni ~ (p,z) \mapsto  \left ( Cz, \, (B^T B)^{-1} B^T (p - f(z)), \,  p - f(z) \, \right ) \in \R^m \times \R^m  \times \R^n
\end{equation*}
To be more precise, there is $ \epsilon > 0$ such that
\begin{equation*}
\left \| \varphi (p',z) - \varphi (p'',z) \right \| \, \leq \, K \, \left \| p' - p'' \right \| ~~\mathrm{~whenever~~} p', p'' \in \B_{\epsilon}(\bar{p}), ~ z \in \B_{\epsilon}(\bar{z}).
\end{equation*}
Without loss of generality, we may assume that the norm on $ \R^m \times \R^m  \times \R^n $ for each point $ (u, v, w) $ is given by
$ \norm{ (u, v, w) } = \norm{u}_m + \norm{v}_m + \norm{w}_n $, where $ \norm{\cdot}_k $ is the Euclidean norm on $ \R^k $.
\end{rem}

\begin{eg} \citep[page 99]{adly2} \\
Suppose that $ \smap{f}{}{} $ is given by $ f(z) := z ~$ for every $ z \in \R $, and $ \mmap{F}{}{} $ is defined as:
\begin{figure}[ht]
  \begin{minipage}{.41\textwidth}
  		\begin{eqnarray*}
			~~~ F (z) :=
			\left \{ \begin{matrix}
			-1 - z ~ & & z < 0,\\
			[-1, 1]  & &  z = 0,\\
			1 - z ~ & & z > 0.
			\end{matrix} \right. 
		\end{eqnarray*}
		\vspace*{0.2cm}
  \end{minipage}
  \begin{minipage}{.60\textwidth}
    	\centering
		\includegraphics[width=0.85\textwidth]{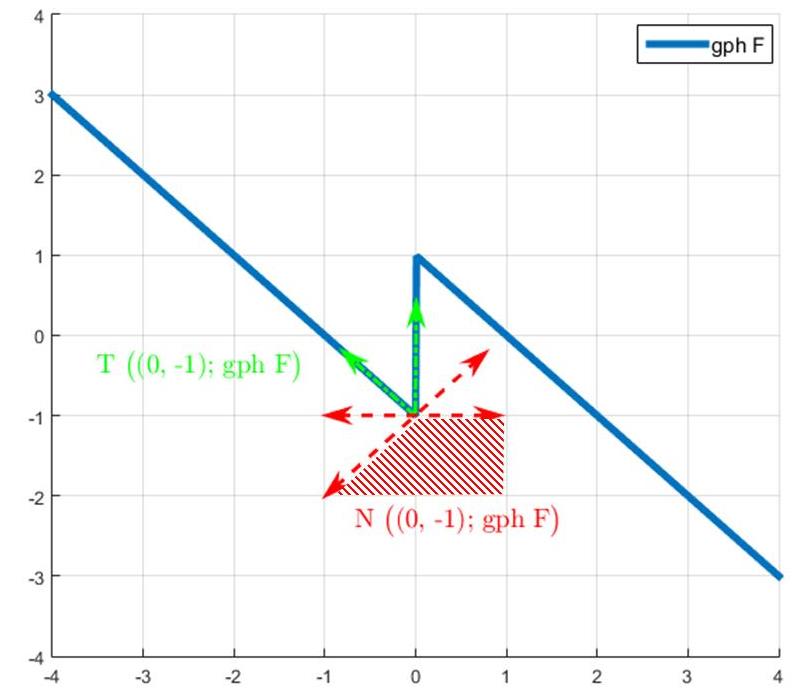}
	\end{minipage}
	\caption{Graph of $F$ with normal and tangent cones at $ (0, -1) $.}
	\label{fig: calmness-example-01}
\end{figure}
Let $ \bar{z} := 0 $ and $ \bar{p} := -1 $. It is easy to verify that both criteria (\ref{AP condition}) and (\ref{IC condition}) are violated. 
In fact, $\bar{v} = \bar{p} - f ( \bar{z} ) = -1 $, $ \nabla f(\bar{z}) = 1 $, and \eqref{AP condition} reads as 
$$ ( \xi, \, \xi) \in - N \Big( (0, -1); \, \gph{F} \Big) \Longrightarrow \, \xi = 0, $$
which is not true, since points other than $ (0, 0) $ from the identity line belong to the normal cone (see Figure \ref{fig: calmness-example-01}). 
Indeed, 
{\footnotesize 
$ \displaystyle N \Big( (0, -1); \, \gph{F} \Big) = {\scriptsize \R \begin{pmatrix}
1 \\ 0 
\end{pmatrix} \, \bigcup \,  \R 
\begin{pmatrix}
- 1 \\ - 1
\end{pmatrix} \, \bigcup \mathrm{~cone~} 
\begin{Bmatrix}
 \begin{pmatrix}
1 \\ 0 
\end{pmatrix}, \, \begin{pmatrix}
- 1 \\ - 1
\end{pmatrix} 
\end{Bmatrix}}
 $}. \\
The condition \eqref{IC condition} reads as
$$ (b, \, -b) \in T \Big( (0, -1); \, \gph{F} \Big) \Longrightarrow \, b = 0, $$
which is not true, since points other than $ (0, 0) $ from the bisector of second quadrant belong to the tangent cone (see Figure \ref{fig: calmness-example-01}). 
Indeed, 
{\footnotesize
$ \displaystyle T \Big( (0, -1); \, \gph{F} \Big) = {\scriptsize \R_{+} \begin{pmatrix}
0 \\ 1 
\end{pmatrix} \, \bigcup \, \R_{+} 
\begin{pmatrix}
- 1 \\  1
\end{pmatrix} } 
 $}. 
So, $S$ has neither the Aubin nor the isolated calmness property at the reference point. Now we want to check the sufficient criteria of Proposition \ref{Proposition 5.1} for the calmness property. \\
Note that, $ g(z) = \big( z, \, -1 -z, \, -1 -z \big) $ in this example and $ \Lambda = \gph{F} \times \R $. Thus,
$$ h(z) := d \big( g(z), \Lambda \big) = d \Big( (z, \, -1 -z, \, -1 -z), \, \gph{F} \times \R \Big) = d \big( (z, \, -1 -z), \, \gph{F} \big).  $$
Let us compute the outer subdifferential of $h$ at $ z = 0 $. Recall that, the set $\partial_{>} h(0) $ contains those points $ \xi \in \R $ such that there are sequences 
$ (z_k)_{k \in \N } $ and $ (\xi_k)_{ k \in \N } $ converging to $ 0 $ and $\xi$, respectively, with $ h(z_k) \downarrow h(0) = 0 $ as $ k \longrightarrow \infty $ and 
$ \xi_k \in \partial_F h(z_k) $ for each $ k \in \N $ (ref. Definition \ref{subdifferentials}). \\
If $ z < 0 $, then the point $ (z, -1 - z) \in \gph{F} $, thus $ h(z) = 0 $, hence we just need to consider $ (z_k)_{k \in \N} $ in $ (0, \infty ) $. 
Moreover, there is $ \delta > 0 $ such that, for each $ z \in (0, \delta) $, the unique nearest point in $\gph{F}$ to the point  $(z, \, -1-z)$  is $ (0, -1) $. Therefore,
$$ h(z) = d \big( (z, \, -1 -z), \, \gph{F} \big) = d \big( (z, \, -1 -z), \, (0, -1)  \big) = \sqrt{2} \, z \mathrm{~~~for~~}  z \in (0, \delta). $$
Hence, 
\begin{equation*}
\displaystyle \partial_F h(z) : = \Big \{ \xi \in \R ~|~ \liminf_{ 0 \, \neq \, r \, \rightarrow \, 0 } \dfrac{h(z+r) - h(z) - \inp{\xi}{r}}{\norm{r}} \, \geq \, 0 \Big \}= 
\big \{ \sqrt{2} \big \}, 
\end{equation*}
for any point $z \in (0, \delta) $. Let $ (z_k)_{k \in \N} $ be any sequence converging to $ 0 $ such that $ h(z_k) \downarrow h(0) = 0 $ as $ k \longrightarrow \infty $. We may assume without any loss of generality that $ z_k \in (0, \delta) $ for each $ k \in \N $. Thus $ \partial_F h(z_k) = \big \{ \sqrt{2} \big \} $, hence 
$ \partial_{>} h(0) = \big \{ \sqrt{2} \big \} $. Therefore, in view of Remark \ref{Estimate for Calmness Modulus}, $ S $ is calm at the reference point with modulus not exceeding $ \dfrac{2}{\sqrt{2}} $.
\end{eg}

\begin{eg} \textbf{(A Simple Circuit with DIAC)} \citep[page 100]{adly2} \\
Let us go back to Examples \ref{DIAC-example} and \ref{DIAC-example-IC}, and this time, check the calmness property in the circuit of 
Figure \ref{fig: DIAC-circuit-chap03} with a DIAC, at the reference point $ \rfp{p}{z} = (-1, 0)  $.\\
In order to avoid non-necessary details and to ease the computations, let us assume that $R = 1$, $V = 1$, and the $ i- v $ characteristic of the DIAC is given by a set-valued mapping 
$ \mmap{\Psi}{}{} $, where $ \Psi(0) := [ -1, 1] $, and $\Psi$ is single-valued and continuously differentiable on $ \R \setminus \{0\}$; its graph is symmetric with respect to the origin, and $ -1 < \Psi(z) < -1 - z $ whenever $ z \in ( - \infty, 0) $; and finally $ \Psi(0_{-}) = -1 $ and $\Psi'(0_{-}) = - a $ for some $a > 0$. \\
Putting $ F = \partial | \cdot | $, one infers that there is a continuously differentiable odd function $ \smap{\psi}{}{} $ such that $ \Psi = \psi + F $. Moreover, $ \psi(0) = 0, 
\psi' (0) = - a $, and $ 0 < \psi(z) < -z $ whenever $ z < 0 $ (see Figure \ref{Multivalued Simplification example}). 
Therefore, we arrive at \eqref{inclusion} with $ m = n = 1 $, $ B = C = I $, and $f$ defined for each $z \in \R $ by $ f(z) := z + \psi(z) $. \\
We have already seen that if $ a = 1 $, then $S$ has neither the Aubin nor the isolated calmness property at $(-1, 0)$. To see whether Proposition \ref{Proposition 5.1} is helpful or not, note that assumptions (A1) - (A3) are satisfied, $ \Lambda = \gph{F} \times \R $, $ g(z) = \big( z, -1 - f(z), -1 - f(z) \big) $, and 
$ h(z) = d \Big( \big( z, -1 -z - \psi(z) \big), \gph{F} \Big) $. \\
When $ z < 0 $, $ \psi (z) > 0 $ and thus, $ -1 -z - \psi(z) < -1 - z $ and by assumption $ \psi(z) < - z $. So, the nearest point in $ \gph{F} $ to $ -1 - z - \psi(z) $ for 
$ z < 0 $ is $ (z, -1) $. Therefore, 
$$ h(z) = \norm{ \big( z, -1 -z - \psi(z) \big), (z, -1) } = - z - \psi(z). $$
Hence, $ \partial_F h(z) = \{ - 1 - \psi' (z) \} $ whenever $ z < 0 $. For each $ k \in \N $, put $ z_k := \dfrac{-1}{k} $ and $ \xi_k := - 1 - \psi' (z_k) $. 
Then both $ (z_k)_{k \in \N} $ and $(\xi_k)_{k \in \N} $ converge to $0$ (in fact, $ \xi_k \longrightarrow -1 - \psi'(0) = - 1 + a = 0 $). 
For any $ k \in \N $, we have $ \xi_k \in \partial_F h(z_k) $, $ h(z_k) = - z_k - \psi(z_k) > 0 $, and $ h(z_k) \downarrow h(0) = 0 $. Thus $ 0 \in \partial_{>} h(0) $. 
Hence, the sufficient condition \eqref{calmness criterion} is not satisfied in this case, and Proposition \ref{Proposition 5.1} does not provide any information about the calmness property of $S$ at the reference point $ (-1, 0 ) $. \\
However, $S$ is not calm at $ (-1, 0)$. Indeed, suppose on the contrary that it is. As $ S(-1) = \{ 0 \} $, the point $0$ is an isolated point of $ S(-1) $, thus $S$ would have the isolated calmness property at $ (-1, 0)$, which is not the case as seen before.
\end{eg}

\section{Results in terms of Metric Regularity} \label{Results in terms of Metric Regularity}
Considering the relation between local stability properties and regularity definitions, one may be able to restate the last three sections' results, in terms of (strong) metric (sub-) regularity. All the new things here are about strong metric regularity, for which we need some calculus rules regarding strict graphical derivatives. \\
We will first provide such rules and then sum up all the results in this area (translated into regularity paradigm) in a single theorem to ease the access (cf. Theorem \ref{SMR-Theorem 2}). \\
The last subsection of this chapter, will discuss the case of a non-differentiable single-valued part (in the sum $f + F$). We believe that would be the proper time to consider this change in our general assumptions (cf. Note \ref{general assumptions}), since all those results we wanted to cover about the static GE has been stated by then.\\
Let us start with a lemma that provides a chain rule for strict graphical derivatives.

\begin{lem} [\textbf{Matrix Operation over Paratingent Cone Relations}] \citep[Lemma 1, p. 4]{cibulkaRoubal} \label{Lemma for SGD sum}
Let $ E \in \R^{k \times d} $ be any matrix, let $G \in \R^{l \times d} $ be injective, and let $\Gamma$ be a subset of rge $E$. Put $ \Xi := E^{-1}(\Gamma)$  and $ \Lambda := G(\Xi)$. For $ \bar{x} \in \Lambda $ denote by $\bar{y}$ the (unique) point in $\Xi$ with $G(\bar{y}) = \bar{x}$. Then
\begin{equation*}
\widetilde{T}(\bar{x}; \Lambda ) = \{ u \in \R^l~:~ \exists \, w \in \R^d \mathrm{~~such~that~} u = Gw \mathrm{~and~} Ew \in \widetilde{T}(E \bar{y}; \Gamma ) \, \}.
\end{equation*}
\end{lem}

\begin{proof}
The proof would be very similar to the proof of Lemma \ref{Lemma 4.1}. One only needs to consider the small necessary changes when replacing 
$T $ with $\widetilde{T}$.\\ 
\textbf{First}, we claim that $ \widetilde{T}(\bar{y}; \Xi ) = \{ w \in \R^d ~:~ Ew \in \widetilde{T}( E \bar{y}; \Gamma ) \} $.\\
Take any $ w \in  \widetilde{T}(\bar{y}; \Xi ) $, Find $ (t_n)_{n \in \N} $ in $ (0, \infty) $, $(y_n)_{n \in \N} $ in $ \Xi $, and $ (w_n)_{n \in \N} $ in $ \R^d $ converging to
$ 0, \bar{y} $, and $w$, respectively, such that $ y_n + t_n w_n \in \Xi $ whenever $ {n \in \N} $. Then we have that 
$$ E y_n + t_n E w_n = E (y_n + t_n w_n ) \in \Gamma \mathrm{~~for~each~} {n \in \N}. $$
Hence $ Ew \in \widetilde{T}( E \bar{y}; \Gamma ) $.\\
On the other hand, let $ w\in \R^d $ be such that $ Ew \in \widetilde{T}( E \bar{y}; \Gamma ) $. By definition, pick $ (t_n)_{n \in \N} $ in $ (0, \infty) $, and 
$ (u_n)_{n \in \N} $ in $\Gamma$, and $ (v_n)_{n \in \N} $ in $ \R^k $ converging to $ 0, E \bar{y} $, and $ E w $, respectively, such that $ u_n + t_n v_n \in \Gamma $ whenever $ {n \in \N} $ (see explanations in Footnote \ref{OMT}). 
$\Gamma \subset \mathrm{\, rge \,} E$, and since $\mathrm{\, rge \,} E$ is a closed subspace of $ \R^k $; one infers that $ v_n \in \mathrm{rge \,} E $ for each $ n \in \N $.  Therefore, by Banach open mapping theorem there are sequences $ (y_n)_{n \in \N} $ converging to $ \bar{y} $ and $ (w_n)_{n \in \N} $ converging to $w$, both in $ \R^d $, such that $ E y_n = u_n $ and $ E w_n = v_n $ for each $n \in \N$.\\
Thus, for an arbitrary index n, we have
\begin{equation*}
\begin{matrix}
u_n = E y_n \in \Gamma ~~ ~ & \mathrm{~so~} & y_n \in E^{-1} (\Gamma) = \Xi \\
u_n + t_n v_n = E ( y_n + t_n w_n ) \in \Gamma & \mathrm{~so~} &  y_n + t_n w_n \in \Xi ~~~~
\end{matrix}
\end{equation*}
Thus, $ w \in  \widetilde{T}(\bar{y}; \Xi ) $, and the claim is proved.\\
\textbf{Second}, we show that $ \widetilde{T}(\bar{x}; \Lambda ) = \{ G w ~:~ w \in \widetilde{T}(\bar{y}; \Xi ) \} $.\\
To prove that $ G \big( \widetilde{T}(\bar{y}; \Xi ) \big) \subset  \widetilde{T}(\bar{x}; \Lambda ) $, pick any $ w \in G \big( \widetilde{T}(\bar{y}; \Xi ) \big) $, and find 
$ v \in G \big( \widetilde{T}(\bar{y}; \Xi ) \big) $ with $ G v = w $. 
Thus there is 
$ (t_n)_{n \in \N} $ in $ (0, \infty) $ converging to $0$, and $ (y_n)_{n \in \N} $ in $\Xi$ converging to $\bar{y}$ ,  and $ (v_n)_{n \in \N} $ in $ \R^d $ converging to $ v $, such that $ y_n + t_n v_n \in \Xi $ whenever $n \in \N$. 
For each $n \in \N$, put $ u_n = G y_n $ and $ w_n = G v_n $. Clearly, $ (u_n)_{n \in \N} $ converges to $ G \bar{y} = \bar{x} $ and $ (w_n)_{n \in \N} $ converges to $ w = G v $. Moreover, \\
$$ u_n + t_n w_n = G (y_n + t_n v_n) \in G( \Xi ) = \Lambda \mathrm{~whenever~} n \in \N. $$
Thus, $ w \in \widetilde{T}(\bar{x}; \Lambda ) $.\\
To see the opposite inclusion, pick any $ w \in \widetilde{T}(\bar{x}; \Lambda ) $. Find ($ (t_n)_{n \in \N} $ in $ (0, \infty) $ converging to $0$, $ (x_n)_{n \in \N} $ in $\Lambda $ converging to $\bar{x}$, and $ (w_n)_{n \in \N} $ in $ \R^l $ converging to $ w $, such that $ x_n + t_n w_n \in \Lambda = G( \Xi ) $ for each $n \in \N$.\\
For each $n \in \N$, find $y_n \in \Xi $ such that $ x_n = G y_n $. Then, $ (y_n)_{n \in \N }$ is bounded.\\
Indeed, if this is not the case, find a cluster point $\bar{h}$ of $ \big( y_n / \| y_n \| \big)_{n \in \N } $. Let $N$ be an infinite subset of $\N$ such that 
$ \displaystyle \lim_{N \ni n \to \infty } \frac{y_n}{ \| y_n \| } =  \bar{h} $. Then
$$ 0 = \lim_{N \ni n \to \infty } \frac{x_n}{ \| y_n \| } = \lim_{N \ni n \to \infty } G \Bigg( \frac{y_n}{ \| y_n \| } \Bigg) = G(\bar{h}). $$
This contradicts the injectivity of $G$ because $ \norm{\bar{h} } = 1 $. Therefore there is an infinite subset $N$ of $\N$ such that $ (y_n)_{ n \in N } $ converges to $ \tilde{y} \in \R^d $, say. Then 
$$ G \bar{y} = \bar{x} = \lim_{N \ni n \to \infty } x_n = \lim_{N \ni n \to \infty } G y_n = G \tilde{y}. $$
Employing, the injectivity once more, we get $ \bar{y} = \tilde{y} $.\\
For each $n \in N$, find $ v_n $ in $\Xi$ such that $ w_n = G \big( \frac{v_n - y_n}{t_n} \big) $, and put $u_n = \frac{v_n - y_n}{t_n}$. Similar argument as in the case of $ (y_n) $ shows that $ (u_n)_{n \in N} $ is bounded. Therefore there is an infinite subset $N'$, of $N$ such that $(u_n)_{n \in N'} $ converges to some $ u \in \R^d $. For each $ n \in N' $, we have $ y_n + t_n u_n = v_n \in \Xi $, therefore $ u \in \widetilde{T}(\bar{y}; \Xi ) $; Moreover, $ w = G(u) $ and thus, 
$ w \in G \big( \widetilde{T}(\bar{y}; \Xi ) \big) $. The proof of second part is finished. \\
Combining the two parts yields the assertion.
\end{proof}

\begin{prop} [\textbf{Strict Graphical Derivative of Sum for GE}] \citep[Proposition 1, p. 5]{cibulkaRoubal} \label{SMR-Proposition 1}
Under the assumptions (A1) - (A2), for any $b \in \R^n$ one has
\begin{equation*}
\widetilde{D} \Phi \rfpp{z}{p} (b)  = \nabla f(\bar{z}) b + B \, \widetilde{D} F_C (C\bar{z}| \, \bar{v}) (Cb).
\end{equation*}
\end{prop}

\begin{proof}
From Proposition \ref{Sum Rule for Strict Graphical Derivatives}, one obtains
\begin{equation*}
\widetilde{D} \Phi \rfpp{z}{p} (b)  = \nabla f(\bar{z}) b + B \, \widetilde{D} Q \big( \bar{z}| \, \bar{p} - f(\bar{z}) \big) (b) ~~ \mathrm{~~for~each~} b \in \R^n .
\end{equation*}
Moreover, observe that
\begin{equation*}
\begin{split}
\gph{Q} =   
\Bigg \{ ~&
 \begin{pmatrix}
 u \\ v
 \end{pmatrix}
 \in \R^{2n} ~:~ \mathrm{there~exists }
 \begin{pmatrix}
 b \\ c
 \end{pmatrix} 
 \in \R^m \times \R^n \mathrm{\, such~that~}\\
 & \, 
 \begin{pmatrix}
 u \\ v
 \end{pmatrix}
 = G 
 \begin{pmatrix}
 b \\ c
 \end{pmatrix} 
 \mathrm{~and~} E
 \begin{pmatrix}
 b \\ c
 \end{pmatrix} 
 \in \gph{F_C} ~~~~ ~~~~~~ ~~~ ~~~~~~~~~
\Bigg \},
\end{split}
\end{equation*}
with 
\begin{equation*}
G := 
\begin{pmatrix}
I_n & 0 \\
0 & B
\end{pmatrix},~~ 
E := 
\begin{pmatrix}
C & 0 \\
0 & I_m
\end{pmatrix}.
\end{equation*}
As $B$ is injective, so is $G$. Using Lemma \ref{Lemma for SGD sum} with $ k := 2m, l := 2n, d := n + m $, $ \Gamma := \gph{F_C}, \bar{x} := \big( \bar{z}, \bar{p} - f (\bar{z}) \big) ^T $, and $ \bar{y} := \rfp{z}{v} ^T $ reveals that
\begin{equation*}
\widetilde{T} \Big( \big( \bar{z}, \bar{p} - f (\bar{z}) \big); \, \gph{Q} \Big) =
\begin{Bmatrix}
\begin{pmatrix}
b \\
Bc
\end{pmatrix}
:
\begin{pmatrix}
Cb \\
c
\end{pmatrix}
\in \widetilde{T} \big(  (C \bar{z}, \bar{v} ) ; \, \gph{F_C} \big)
\end{Bmatrix}.
\end{equation*}
This means that $ \widetilde{D} Q \big( \bar{z}, \bar{p} - f (\bar{z})\big) (b) = B \widetilde{D} F_C ( C\bar{z}, \bar{v})(Cb) $. The assertion is proved.
\end{proof}

The following proposition is the strict graphical version of Theorem \ref{CDCMR theorem}, and relates the strong metric regularity of a set-valued map at a reference point to a property of its strict graphical derivative. Remember that the necessary condition was already obtained in Theorem \ref{SGDCSMR}.

\begin{prop} [\textbf{Strong Metric Regularity in terms of Strict Graphical Derivative}] \citep[Proposition 2, p.5]{cibulkaRoubal} \label{SMR-Proposition 2} \\
Consider a set-valued mapping $\mmap{H}{n}{n}$ and a point $ \rfp{x}{y} \in \gph{~H} $. Then $H$ is strongly metrically regular at $\bar{x}$ for $\bar{y}$ if and only if it satisfies the following three conditions:
\begin{enumerate}[topsep=-1ex,itemsep=-1ex,partopsep=1ex,parsep=1ex]
\item[(a)] for each neighborhood $U$ of $\bar{x}$ there is a neighborhood $V$ of $\bar{y}$ such that \\
$H^{-1} (y) \cap U \neq \emptyset$ whenever $y \in V$;
\item[(b)] the set $\gph{~H} \cap \left ( \B_r (\bar{x}) \times \B_r(\bar{y}) \right )$ is closed for some $r > 0$;
\item[(c)] $0 \in \widetilde{D} H(\bar{x} | \, \bar{y})(u) \Longrightarrow \, u = 0 $.
\end{enumerate}
\end{prop}

\begin{proof}
Suppose that $H$ is strongly metrically regular \at{x}{y}. Then (c) holds by Theorem \ref{SGDCSMR}. Observe also that $H$ has necessarily locally closed graph at the reference point. Finally, (a) is satisfied since $H$ is open at $ \rfp{x}{y} $, that is, for any neighborhood $ U $ of $\bar{x}$, the set $ V := H(U) $ is a neighborhood of $\bar{y}$. \\
The converse implication is proved in Theorem \ref{SGDCSMR}.
\end{proof}

\begin{thm} [\textbf{Strong Metric Regularity Criterion for GE}] \label{SMR-Theorem 1} \hfill \\ \citep[Theorem 1, p. 6]{cibulkaRoubal} 
Assume that (A1) - (A3) hold true. Then $\Phi$ is strongly metrically regular at $\bar{z}$ for $\bar{p}$ if and only if 
\begin{enumerate} [topsep=-1ex,itemsep=-1ex,partopsep=1ex,parsep=1ex]
\item[(a)] for each neighborhood $U$ of $\bar{z}$ there is a neighborhood $V$ of $\bar{p}$ such that\\
 $\Phi^{-1} (p) \cap U \neq \emptyset$ whenever $p \in V$;
\item[(b)] $0 \in \nabla f(\bar{z}) b + B \, \widetilde{D} F_C (C\bar{z}| \, \bar{v}) (Cb) ~~ \Longrightarrow ~~ b = 0  $.
\end{enumerate}
Moreover, its regularity modulus is given by
\begin{equation*}
\mathrm{reg~} (\Phi; \bar{z}| \bar{p}) = \sup \left \{ \, \left \| b \right \| \, : \left (   \nabla f(\bar{z}) b + B \, \widetilde{D} F_C (C\bar{z}| \, \bar{v}) (Cb) \, \right ) \cap \B \neq \, \emptyset ~ \right \}.
\end{equation*}
\end{thm}

\begin{proof}
Having in mind that $\Phi$ has closed graph, one only needs to combine Proposition \ref{SMR-Proposition 1} and Proposition \ref{SMR-Proposition 2} to get the desired equivalence. The formula for the regularity modulus follows from Theorem \ref{SGDCSMR}.
\end{proof}

\begin{eg} \textbf{(A Simple Circuit with DIAC)} \label{A Simple Circuit with DIAC-chap3-SMR}\\
Consider the circuit in Figure \ref{fig: DIAC-circuit-chap03} with a DIAC. In previous sections we have studied the Aubin property, isolated calmness, and calmness of the solution mapping at different reference points. Let us now, investigate the strong metric regularity of the map $ \Phi $ in view of Theorem \ref{SMR-Theorem 1}. 
Note that $ \mmap{\Phi}{}{} $ is given by $ \Phi(z) = f(z) + F(z) $, where 
\begin{eqnarray*}
f(z) :=
\left \{ \begin{matrix}
R z - \dfrac{V}{\sqrt{1 - \frac{2az}{V}}} + V & & z < 0,\\
R z + \dfrac{V}{\sqrt{1 + \frac{2az}{V}}} - V & & z \geq 0,
\end{matrix} \right.
\mathrm{~~~~and~~~~}
F(z) :=
\left \{ \begin{matrix}
-V ~~~ & &  z < 0,\\
[- V, V] & & z = 0,\\
V ~~~ & & z > 0,
\end{matrix} \right.
\end{eqnarray*} 
in which $ R $, $a$, and $V$ are positive constants (ref. Example \ref{DIAC-Chap2}). Assumptions (A1) - (A3) hold true.  We will consider two reference points. \\
For $ \rfp{z}{p} = (0, 0) $, we would have $ f'(0) = R - a $, $ \bar{v} = 0 $. A close look at graph of $\Phi$ (see Figure \ref{fig: DIAC-circuit02-chap02}) reveals that condition $(a)$ is satisfied. Condition $(b)$ reads as
\begin{equation*}
0 \in (R - a ) \beta + \widetilde{D} F ( \,0\,|\, 0\,) (\beta) ~ \Longrightarrow ~ \beta = 0.
\end{equation*}
Since 
$ {\scriptsize
\widetilde{T} \Big( (0,0); \gph{F} \Big) = \R 
\begin{pmatrix}
0 \\ 1
\end{pmatrix}
}
$, also this condition is satisfied, and $\Phi$ is strongly metrically regular at $0$ for $0$. \\
For $ \rfp{z}{p} = (0, -V) $, parameters $R$, and $a$ play an important role, as we might have expected. If $ a > R $, the point $ (0, -V) $ is a local minimum of $\gph{\Phi} $, and hence condition $(a)$ is not satisfied. When $ a \leq R $, there would be no more problems for this condition (see Figure  \ref{fig: DIAC-circuit02-chap02}). 
To check condition $(b)$, notice that
\begin{equation*}
\begin{split}
\widetilde{D} F ( \,0\,|\, -V\,) (\beta) & = \begin{Bmatrix}
v \in \R ~|~ \begin{scriptsize}
\begin{pmatrix}
\beta \\ v
\end{pmatrix}
\end{scriptsize} 
\in \widetilde{T} \Big( (0,\, -V); \gph{F} \Big)
\end{Bmatrix} \\
& = 
\begin{Bmatrix}
v ~|~ \begin{scriptsize}
\begin{pmatrix}
\beta \\ v
\end{pmatrix} 
\in \R \begin{pmatrix}
0 \\ 1
\end{pmatrix} \, \bigcup \, \R 
\begin{pmatrix}
-1 \\ 0
\end{pmatrix}
\end{scriptsize} 
\end{Bmatrix}.
\end{split}
\end{equation*}
Hence, $ \big( \beta,\, (a - R) \, \beta \big)^T $ should belong to the latter union. If $ a < R $, the only possibility would be $ \beta = 0 $, and thus, the condition holds true. But in case of  $ a = R $, it is possible to have $ \beta \neq 0 $. \\
Therefore, $\Phi$ will not have the strong metric regularity property at $0$ for $-V$ if $ a \geq R$.
\end{eg}

\begin{thm} [\textbf{Summary of Metric Regularity Criteria for GE}] \label{SMR-Theorem 2} \hfill \\ \citep[Theorem 2, p. 6]{cibulkaRoubal}
Suppose that the assumptions (A1) - (A4) hold true. Then
\begin{enumerate} [topsep=-1ex,itemsep=-1ex,partopsep=1ex,parsep=1ex]
\item[(i)] $\Phi$ is metrically regular at $\bar{z}$ for $\bar{p}$ if and only if
 \begin{equation*}
\left.\begin{matrix}
((C C^T)^{-1} C \nabla f(\bar{z})^T \xi, B^T \xi ) \in - N ((C \bar{z}, \bar{v}); \gph{F}) ~~\\
~~~~~~ \nabla f(\bar{z})^T \xi \in \mathrm{~rge~} C^T 
\end{matrix}\right\}
\Rightarrow \xi = 0;
\end{equation*}

\item[(ii)] $\Phi$ is strongly metrically sub-regular at $\bar{z}$ for $\bar{p}$ if and only if
 \begin{equation*}
\left.\begin{matrix}
(Cb, -(B^T B)^{-1} B^T \nabla f(\bar{z}) b ) \in T ((C \bar{z}, \bar{v}); \gph{F}) ~~\\
~~~~~~~~~~ \nabla f(\bar{z}) b \in \mathrm{~rge~} B 
\end{matrix}\right\}
 \Rightarrow b = 0;
 \end{equation*}

\item[(iii)] $\Phi$ is strongly metrically regular at $\bar{z}$ for $\bar{p}$ if and only if
\begin{itemize}[topsep=-1ex,itemsep=-1ex,partopsep=1ex,parsep=1ex]
\item[(a)] for each neighborhood $U$ af $\bar{z}$ there is a neighborhood $V$ of $ \bar{p} $ such that \\
 $ \Phi^{-1}(p) \cap U \not = \emptyset $ whenever $p \in V$;
\item[(b)]
 \begin{equation*}
\left.\begin{matrix}
(Cb, -(B^T B)^{-1} B^T \nabla f(\bar{z}) b ) \in \widetilde{T} ((C \bar{z}, \bar{v}); \gph{F}) ~~\\
~~~~~~~~~~ \nabla f(\bar{z}) b \in \mathrm{~rge~} B 
\end{matrix}\right\}
 \Rightarrow b = 0.
 \end{equation*}
\end{itemize}
\end{enumerate}
\end{thm}

\begin{proof}
The statement (i) is Corollary \ref{Corollary 3.1}, the proof of which was based on Mordukhovich coderivative criterion. Whereas (ii) is Corollary \ref{Corollary 4.1}. To see
the last one, note that if $C$ is surjective, then $ F_C = F $. Moreover, (A1) ensures that $ B^T B \in \R^{m \times m} $ is non-singular. It suffices to show that (b) is equivalent to (b) in Theorem \ref{SMR-Theorem 1}.\\
First, let $ b \in \R^n $ be such that $ 0 \in \nabla f(\bar{z}) b + B \, \widetilde{D} F (C\bar{z}| \, \bar{v}) (Cb) $. Find a point 
$ w \in \widetilde{D} F (C\bar{z}| \, \bar{v}) (Cb) $ with $ \nabla f(\bar{z}) b + B w = 0 $. Thus, $  -(B^T B)^{-1} B^T \nabla f(\bar{z}) b  $ is in 
$ \widetilde{D} F (C\bar{z}| \, \bar{v}) (Cb) $. Clearly, we have $ \nabla f(\bar{z}) b \in \mathrm{rge} \, B $ and the definition of the strict graphical derivative of $F $ yields the rest.\\
On the other hand, pick any $ b \in \R^n $ with $ \big( Cb, -(B^T B)^{-1} B^T \nabla f(\bar{z}) b \big) $ in $ \widetilde{T} ((C \bar{z}, \bar{v}); \gph{F}) $ and 
$ \nabla f(\bar{z}) b \in \mathrm{rge} \, B $. The definition of the strict graphical derivative says that 
$$ w :=  -(B^T B)^{-1} B^T \nabla f(\bar{z}) b \in   \widetilde{D} F (C\bar{z}| \, \bar{v}) (Cb). $$ 
Thus we have $ B^T B w = - B^T \nabla f(\bar{z}) b $. So, $ B w + \nabla f(\bar{z}) \in \mathrm{ker} \, B^T \cap \mathrm{rge} \, B = \{ 0 \} $. 
Therefore, $ 0 \in \nabla f(\bar{z}) b + B \, \widetilde{D} F (C\bar{z}| \, \bar{v}) (Cb) $.
\end{proof}

\begin{cor} \citep[Corollary 1, p. 7]{cibulkaRoubal} \label{Corollary-SMR} \\
In addition to (A1) - (A5), assume that $ n = m$, $B = C = I_n$, $ \nabla f(\bar{z}) $ is a P-matrix, and for each $ i \in \{ 1, 2, \cdots ,n \} $, the mapping
$ \mmap{F_i}{}{} $ is maximal monotone. Then $\Phi$ is strongly metrically regular \at{z}{p}.
\end{cor}

\begin{proof}
For any $ x  = (x_1, \cdots, x_n), y = (y_1, \cdots, y_n) \in \R^n $, define the function $\varphi : \R^n \times \R^n \to (\R^2)^n$ as 
$ \varphi (x, y) := \Big( (x_1, y_1), \cdots, (x_n, y_n) \Big) $.\\
Clearly, $\varphi$ is linear and one-to-one. Also, by (A5), we have $ \displaystyle \prod_{i=1} ^n \gph{F_i} = \varphi( \gph{F} ) $. The definition of the paratingent cone, Proposition \ref{product rule}, and
Lemma \ref{Lemma for SGD sum} (with $E = I$, and $G$ be the representative matrix of $\varphi$) imply that
\begin{equation*}
\begin{split}
\varphi \Big( \widetilde{T} \big( \rfp{z}{v} ; \gph{F} \big) \Big) & = \widetilde{T}  \Big(  \varphi \rfp{z}{v} ; \varphi ( \gph{F} ) \Big) \\
& = \widetilde{T}  \Big(  \varphi \rfp{z}{v} ; \prod_{i=1} ^n \gph{F_i} \Big) \subset \prod_{i=1} ^n \widetilde{T} \big( \rfp{z_i}{v_i} ; \gph{F_i} \big).
\end{split} 
\end{equation*}
Also, it is well-known that
\begin{equation*}
\prod_{i=1} ^n N \big( \rfp{z_i}{v_i} ; \gph{F_i} \big) = N \Big(  \varphi \rfp{z}{v} ; \prod_{i=1} ^n \gph{F_i} \Big)
= \varphi \Big( N \big( \rfp{z}{v} ; \gph{F} \big) \Big).
\end{equation*}
As all $ F_i $'s are maximal monotone, we have
\begin{eqnarray*}
\begin{split}
N \big( \rfp{z_i}{v_i} ; \gph{F_i} \big) & \subset \{ (a, b) \in \R^2 ~:~ a b \leq 0 \}, \\
\widetilde{T} \big( \rfp{z_i}{v_i} ; \gph{F_i} \big) & \subset \{ (a, b) \in \R^2 ~:~ a b \geq 0 \}, \mathrm{~~for~each~~} i \in \{1, \cdots, n \}.
\end{split}
\end{eqnarray*}
Fix any non-zero $ \eta \in \R^n $. Since $ \nabla f(\bar{z}) $ is a P-matrix, so is $ \nabla f(\bar{z})^T $. There are $ k, l \in \{1, \cdots, n \} $ such that 
$ \eta_k ( \nabla f(\bar{z}) \eta )_k > 0 $ and $ \eta_l ( \nabla f(\bar{z})^T \eta )_l  > 0 $, which means that
$  \Big( \eta_k, - ( \nabla f(\bar{z}) \eta )_k \Big) \not \in \widetilde{T} \big( \rfp{z_k}{v_k} ; \gph{F_k} \big) $ and
$ \Big( ( \nabla f(\bar{z})^T \eta )_l, \eta_l \Big) \not \in - N \big( \rfp{z_l}{v_l} ; \gph{F_l} \big) $.\\
The above relations for the normal and paratingent cone and the fact that $\varphi$ is one-to-one, imply that conditions $(i)$ and $(iii-b)$ in Theorem \ref{SMR-Theorem 2} hold. Then, $(iii-a)$ is a result of $(i)$\footnote{
This fact is indeed, the very first lines of the proof of Theorem \ref{MR-AP}, (equivalence of metric regularity and inverse Aubin property). Here is a brief review.\\
Let $U = \B_a (\bar{z}) $, $V = \B_b (\bar{p}) $, $ \kappa > \mathrm{reg} \, (\Phi; \bar{z} \, | \, \bar{p} ) $, where $a$ and $b$ are positive constants such that $ b < \frac{a}{\kappa} $.\\
Choose an arbitrary $ p \in \B_b (\bar{p}) $. If $\Phi ^{-1} (p) \cap \B_a (\bar{z})  = \emptyset $, then $ d \big( \bar{z}, \Phi ^{-1} (p) \big) \geq a $. Thus, by definition of metric regularity one obtains
$$ a \leq  d \big( \bar{z}, \Phi ^{-1} (p) \big) \leq \kappa d \big( p, \Phi(\bar{z}) \big) \leq \kappa \norm{ p - \bar{p} } \leq \kappa b < a,  $$
which is a contradiction. So, $(iii-a)$ holds.
}, 
and the proof is complete.
\end{proof}

\begin{eg}
Let us have a closer look at Example \ref{A Simple Circuit with DIAC-chap3-SMR}. We have $ n = m = 1$, $ B = C = I_1$, and not only assumptions (A1)-(A3), but also assumptions (A4) and (A5) hold true. It is also easy to verify that $F$ is maximal monotone. Hence, in view of Corollary \ref{Corollary-SMR}, the strong metric regularity of $\Phi$ could be concluded if $f'(\bar{z}) $ is a P-matrix. \\
From Definition \ref{P-matrix}, and calculations in Example \ref{DIAC-Chap2}, we obtain that everything reduces to check whether 
$$ f'(\bar{z}) = R - \dfrac{a}{\paren{ \dfrac{2 a |\bar{z}|}{V} + 1 } ^{3/2}} $$
is positive or not. For $\bar{z} = 0$, $f'(0) = R - a > 0 $ whenever $ a < R $, which confirms the results of Example \ref{A Simple Circuit with DIAC-chap3-SMR}.\\
Note that the corollary is a sufficient condition and thus, does not say anything about the strong metric regularity when $ f'(\bar{z}) \, \leq \, 0 $, as we observed that for example, at the reference point $ (0, 0)$, the strong metric regularity always holds no matter what the relation of $R$ and $a$ is.
\end{eg}

\subsection{Strong Metric (Sub-) Regularity with a Nonsmooth Single-valued Part} \label{Strong Metric (Sub-) Regularity with a Nonsmooth Single-valued Part}

The rest of this section discusses the case of a non-smooth $f$ in the setting (\ref{inclusion}), that is when $f$ is not continuously differentiable on $\R^n$ as it was assumed by condition (A2) in Note \ref{general assumptions}. Since there is no differentiability of $f$ any more to use, one can think of generalized Jacobians as a possible alternative. \\
We shall recall that from the circuit interpretation point of view, in general, $F$ represents the $ i - v $ characteristic of diodes, transistors and such components and $f$ indicates all the other components and how they are connected to each other. However, as we have explained before in Chapter \ref{Chapter2}, one can use a simplification technique in order to simplify the graph of set-valued part (which has a computational importance) and thus, add some points to $f$. Hence, it is very probable to obtain a single-valued part which is not smooth enough (see Example \ref{A simple circuit with SCR}). \\
We will discuss the results obtained by Ismailov \citep{izmailov2014} for strong metric regularity and the similar results considering his method for strong metric 
sub-regularity \citep{cibulka2016strong}. \\ 
It is worth mentioning that as Ismailov claimed, the setting 
\begin{equation} \label{robinson-extended}
f (x) + F (x) \ni 0, 
\end{equation}
with a non-differentiable function $ \smap{f}{n}{n}$ and a set-valued map $\mmap{F}{m}{m}$, could be considered as a way of unifying two classical theorems of variational analysis. \\
The first theorem is due to Robinson \citep{robinson} and is based on the assumptions we were working with till now (see also Subsection \ref{Review on an existence theorem}). Thus, Theorem \ref{Izmailov's theorem} extends Theorem \ref{SMR-Theorem 2} to the case of a non-smooth function $f$.\\
The second result is Clarke’s inverse function theorem \citep{clarke1976}, which is concerned with the case of a usual nonlinear equation
$$ f(x) = 0, $$
corresponds to (\ref{robinson-extended}) with $ F (\cdot) \equiv \{ 0 \} $; but assumes local Lipschitz continuity of $f$ rather than smoothness. Thus, Theorem \ref{Izmailov's theorem} at the same time, extends Clarke’s theorem (see also \citep[Theorem 7.1.1, p. 253]{clarke1990}) from usual equations to GEs.\\
Thus, we assume the following condition instead of (A2):
\begin{itemize}[nolistsep]
\item[]$(\widetilde{\mathrm{A2}})$ $f$ is locally Lipschitz continuous on $\R^n$.
\end{itemize}
Let us first provide a perturbation results like Theorem \ref{theorem 5G.3}, with a strongly metrically sub-regular set-valued map and a calm single-valued term.

\begin{lem} [\textbf{Stability of Strong Metric Sub-regularity Under Single-valued Calm Perturbations}] \citep[Lemma 2, p. 8]{cibulkaRoubal} \label{Stability of Strong Metric Sub-regularity Under Single-valued Calm Perturbations}\\
Let $ \rfp{x}{y} \in \R^n \times \R^m $ and $ \mmap{G}{n}{m} $ be such that $ \bar{y} \in G(\bar{x}) $. Suppose that $G$ is strongly metrically sub-regular \at{x}{y} (that is, there is $ \kappa > 0 $ along with a neighborhood $U$ of $\bar{x}$ such that
$ \norm{ x - \bar{x} } \, \leq \, \kappa \, d \big( \bar{y}, G(x) \big) $ whenever $ x \in U $).\\
Then, for any function $ \smap{g}{n}{m} $ which is calm at $\bar{x}$ relative to $U \subset \mathrm{dom} \, g $ with the constant $ \mu < \frac{1}{\kappa} $, one has
\begin{equation*}
\norm{ x - \bar{x} } \, \leq \, \dfrac{\kappa}{1 - \kappa \mu} d \big( \bar{y} + g(\bar{x}),  g(x) + G(x) \big) \mathrm{~~for~each~~} x \in U.
\end{equation*}
That is, $ g + G $ is strongly metrically sub-regular at $\bar{x}$ for $\bar{y} + g(\bar{x})$.
\end{lem}

\begin{proof}
Fix any $ x \in U $. The calmness of $g$ means that $ \norm{ g(x) - g(\bar{x}) } \, \leq \, \mu \norm{ x - \bar{x} } $. Therefore, 
\begin{eqnarray*}
\begin{split}
\norm{ x - \bar{x} } \, & \leq \, \kappa \, d \big( \bar{y}, G(x) \big) \, \leq \, \kappa \, \Big [ d \big( \bar{y}, \, \bar{y} - g(x) + g(\bar{x}) \big)  + d \big( \bar{y} - g(x) + g(\bar{x}), G(x) \big)  \Big ] \\
& \leq \kappa \, \norm{ g(x) - g(\bar{x}) } + \, \kappa \, d \big( \bar{y} +  g(\bar{x}), \, g(x) + G(x) \big) \\
& \leq \kappa \mu \, \norm{ x - \bar{x} } + \, \kappa \, d \big( \bar{y} + g(\bar{x}), \, g(x) + G(x) \big).
\end{split}
\end{eqnarray*}
Performing a small rearrangement and dividing by $ 1 - \kappa \mu > 0 $, we obtain the desired inequality.
\end{proof}

Now, we go one step ahead with dropping the calmness condition on $g$.

\begin{thm} [\textbf{Stability of  SMSR Under Single-valued Perturbations}] \label{SMSR-Theorem 3} \hfill \\ \citep[Theorem 3.5, page 11]{cibulka2016strong} 
Let $ \rfp{x}{y} \in \R^n \times \R^l $, $ \smap{g}{n}{l} $, and $ \mmap{G}{n}{l} $ be such that $  \bar{y} \in g(\bar{x}) + G(\bar{x}) $. Suppose that there exists a compact subset $ \mathcal{A} $ of $ \R^{l \times n} $ such that
\begin{enumerate}[topsep=-1ex,itemsep=-1ex,partopsep=1ex,parsep=1ex]
\item[(a)] there are $ \varepsilon > 0 $ and $ r > 0 $ such that for each $ u \in \mathrm{int}\, \B_{r} (\bar{x}) $, one can find $ A \in \mathcal{A} $ such that 
$ \norm{ g(u) - g(\bar{x}) - A ( u - \bar{x}) } \leq \varepsilon \norm{ u - \bar{x} } $;
\item[(b)] for every $ A \in \mathcal{A} $ the mapping $ G_A : \R^n \ni x \mapsto g(\bar{x}) + A (x - \bar{x}) + G(x) \subset \R^l $ is strongly metrically sub-regular \at{x}{y} and let $ \displaystyle m := \sup_{ A \in \mathcal{A} } \mathrm{subreg} \, (G_A; \bar{x} \, | \, \bar{y} ) < \dfrac{1}{2 \varepsilon} $.
\end{enumerate}
Then $ g + G $ is strongly metrically sub-regular \at{x}{y}; and $$ \mathrm{subreg} \, (g + G; \bar{x} \, | \, \bar{y} ) \leq \dfrac{m}{1 - 2 m \varepsilon}. $$
\end{thm}

\begin{proof}
Without loss of generality assume that $ \bar{y} = 0 $ and note that $ m < \infty $. Fix any $ \kappa > m $ such that $ 2 \kappa \varepsilon < 1 $. Let $r > 0$ be as in $(a)$. First, we show that there exists $ a \in (0, r] $ such that
\begin{equation} \label{SMSR-1}
\norm{ x - \bar{x} } \, \leq \, \dfrac{\kappa}{1 - \kappa \varepsilon} d \big(0, G_A(x) \big) \mathrm{~~whenever~~} x \in \mathrm{int} \, \B_a (\bar{x}), \mathrm{~and~} A \in \mathcal{A}.
\end{equation}
As $\mathcal{A}$ is compact, there is a finite set $ \mathcal{A}_F \subset \mathcal{A} $ such that
\begin{equation} \label{SMSR-2}
\mathcal{A} \subset \mathcal{A}_F + \varepsilon \B.
\end{equation}
Pick any $ \widetilde{A} \in \mathcal{A}_F $. Then, by assumption (b), there exists $ \alpha_{\widetilde{A}} > 0 $ such that
\begin{equation*}
\norm{ x - \bar{x} } \, \leq \, \kappa \, d \big(0, G_{\widetilde{A}}(x) \big) \mathrm{~~whenever~~} x \in \mathrm{int} \, \B_{\alpha_{\widetilde{A}}} (\bar{x}).
\end{equation*}
Fix any $ A' \in \varepsilon \B $. As $ G_{\widetilde{A} + A' }  = G_{\widetilde{A} } + A' ( x - \bar{x}) $, Lemma \ref{Stability of Strong Metric Sub-regularity Under Single-valued Calm Perturbations} with $g(\cdot) = A' ( \cdot - \bar{x})  $ reveals that
\begin{equation*}
\norm{ x - \bar{x} } \, \leq \, \dfrac{\kappa}{1 - \kappa \varepsilon} d \big(0, G_{\widetilde{A} + A' }(x) \big) \mathrm{~~for~any~~} x \in \mathrm{int} \, \B_{\alpha_{\widetilde{A}}} (\bar{x}).
\end{equation*}
Thus for any $ \widetilde{A} \in \mathcal{A}_F $, there is $ \alpha_{\widetilde{A}} > 0 $ such that for each $ A' \in \varepsilon \B $ the above inequality holds. 
Let $ \displaystyle a = \min \Big \{ r , \min_{\widetilde{A} \, \in \, \mathcal{A}_F} \alpha_{\widetilde{A}} \Big \} $. Taking into account \eqref{SMSR-2}, we obtain \eqref{SMSR-1}.\\
Fix any $ x \in \mathrm{int} \, \B_a (\bar{x} ) $. Use $(a)$ to find $ A \in \mathcal{A} $ such that $ \norm{ g(x) - g(\bar{x}) - A ( x - \bar{x}) } \leq \varepsilon \norm{ x - \bar{x} } $.
This and \eqref{SMSR-1} imply that
\begin{eqnarray*}
\begin{split}
\norm{ x - \bar{x} } \, & \leq \, \dfrac{\kappa}{1 - \kappa \varepsilon} d \big(0, G_{A}(x) \big) = 
 \dfrac{\kappa}{1 - \kappa \varepsilon} d \big(- g(\bar{x}) - A ( x - \bar{x}) , G(x) \big) \\
& \leq \, \dfrac{\kappa}{1 - \kappa \varepsilon} \Big [ d \big(- g(x), G(x) \big) + \norm{ g(x) - g(\bar{x}) - A ( x - \bar{x}) }  \Big]\\
& \leq \, \dfrac{\kappa}{1 - \kappa \varepsilon} d \big( 0, g(x) + G(x) \big) + \dfrac{\kappa}{1 - \kappa \varepsilon} \varepsilon \norm{ x - \bar{x} }
\end{split}
\end{eqnarray*}
Since $ 2 \kappa \varepsilon < 1 $, and thus, $ \dfrac{\kappa \varepsilon}{1 - \kappa \varepsilon} < 1 $, we get that 
$$ \norm{ x - \bar{x} } \, \leq \, \dfrac{\kappa}{1 - 2 \kappa \varepsilon} d \big(0, \, g(x) + G(x) \big). $$
Thus, $ g + G $ is strongly metrically sub-regular at $ \bar{x} $ for $0$. As $ \kappa \in (m, \frac{1}{2 \varepsilon} ) $ was arbitrary, we get the desired estimate on the sub-regularity modulus.
\end{proof}

Unfortunately, Theorem \ref{SMSR-Theorem 3} does not say anything about the possible choices of subsets $\mathcal{A}$. The following corollary suggests a family of subsets $\mathcal{A}$, for which the condition (a) in Theorem \ref{SMSR-Theorem 3} is satisfied.\\

\begin{cor} \citep[Corollary 2, p. 10]{cibulkaRoubal} \label{Corollary2}
Let $ \rfp{x}{y} \in \R^n \times \R^l $, $ \smap{g}{n}{l} $, and $ \mmap{G}{n}{l} $ be such that $  \bar{y} \in g(\bar{x}) + G(\bar{x}) $. Suppose that $g$ is locally Lipschitz continuous at $\bar{x}$ and that for every $ \varepsilon > 0 $ there exists $ r > 0 $ along with a selection $h$ for $ \partial_B g $ such that
\begin{equation} \label{Corollary2-01}
\norm{ g(u) - g(\bar{x}) - h(u) ( u - \bar{x}) } \, \leq \, \varepsilon \, \norm{ u - \bar{x} } \mathrm{~~whenever~~} u \in \mathrm{int} \, \B_r (\bar{x}).
\end{equation}
Assume that the assumption (b) in Theorem \ref{SMSR-Theorem 3} is satisfied with $ \mathcal{A} := \partial_B g(\bar{x}) $. 
Then $ g + G $ is strongly metrically sub-regular \at{x}{y}; and $ \mathrm{subreg} \, (g + G; \bar{x} \, | \, \bar{y} ) \, \leq  m $.
\end{cor}

\begin{proof}
Let $ \gamma \in (0,1) $ be such that $ \gamma ( m + \gamma ) < 1 $. Set $ \mathcal{A} = \partial_B g(\bar{x}) + \gamma \B $. 
Then $\mathcal{A}$ is compact (as the sum of two compact sets). By Lemma \ref{Stability of Strong Metric Sub-regularity Under Single-valued Calm Perturbations}, for any
$ \widetilde{A} \in \partial_B g(\bar{x}) $ and any $ A' \in \gamma \B $, the mapping $ G_{ \widetilde{A} + A' } = G_{ \widetilde{A}} + A' (x  - \bar{x}) $ is strongly metrically sub-regular \at{x}{y} with the modulus at most $\frac{m + \gamma}{1 - (m + \gamma) \gamma } $. Thus, for every $ A \in \mathcal{A} $ the mapping $ G_A $ is strongly metrically sub-regular \at{x}{y}; and
$$ m' = \sup_{A \in \mathcal{A}} \mathrm{subreg} \, (G_A; \bar{x} \, | \, \bar{y} ) \leq \dfrac{m + \gamma}{1 - (m + \gamma) \gamma } $$
Let $ \varepsilon \in (0, \gamma ) $ be arbitrary. By the outer semi-continuity of $\partial_B g$ and (\ref{Corollary2-01}), there is $ r > 0 $ and a selection $h$ for 
$\partial_B g$ such that, for each $ u \in \mathrm{int} \, \B_r (\bar{x}) $, one has
\begin{equation*}
h(u) \in \partial_B g(u) \subset \partial_B g(\bar{x}) + \varepsilon \B \subset \mathcal{A} \mathrm{~~and~~}
\norm{ g(u) - g(\bar{x}) - h(u) ( u - \bar{x}) } \, \leq \, \varepsilon \, \norm{ u - \bar{x} }.
\end{equation*}
Thus, for each $ u \in \mathrm{int} \, \B_r (\bar{x}) $ one can find $ A \in \mathcal{A} $ such that $\norm{ g(u) - g(\bar{x}) - A ( u - \bar{x}) } \, \leq \, \varepsilon \, \norm{ u - \bar{x} }$. Theorem \ref{SMSR-Theorem 3} implies that the mapping $g + G$ is strongly metrically sub-regular \at{x}{y}; and 
$ \mathrm{subreg} \, (g + G; \bar{x} \, | \, \bar{y} ) \, \leq  m' $. As $ \gamma > 0 $ can be arbitrarily small, the proof is finished.
\end{proof}

Now, we will express the Izmailov's theorem \citep{izmailov2014} about the strong metric regularity in case of a non-smooth single-valued part, and then summarize the results of this subsection in one theorem.\\

\begin{thm} [\textbf{SMR of Perturbed GE with a Nonsmooth Single-valued Part}] \citep[Theorem 3, p. 583]{izmailov2014} \label{Izmailov's theorem}
Let $\smap{f}{n}{n}$ be a Lipschitz continuous function in a neighborhood of $ \bar{x} \in \R^n$, and $\mmap{F}{n}{n}$ be a set-valued map such that $ 0 \in f(\bar{x}) + F( \bar{x})$, and for each $ A \in \partial f(\bar{z}) $ the mapping
$$ \R^n \ni x \mapsto J_A(x) := f(\bar{x}) + A (x - \bar{x}) + F(x)  $$
is strongly metrically regular at $\bar{x}$ for $0$.\\
Then there exist neighborhoods $U$ of $\bar{x}$ and $V$ of $0$ such that for every $y \in V$, there exists a unique $x(y) \in U$ satisfying the perturbed generalized equation $ y \in f(x) + F( x ) $  , and the mapping $ y \mapsto x(y) $ is Lipschitz continuous on $V$.\\
\end{thm}

\begin{thm} [\textbf{Summary of MR Criteria for GE in Case of a Nonsmooth Single-valued Part}] \citep[Theorem 4, p. 10]{cibulkaRoubal} \label{Summary of MR Criteria for GE in Case of a Nonsmooth Single-valued Part}\\
Under the assumptions (A1), $\widetilde{(\mathrm{A2})} $, (A3), and (A4), for any $ A \in \partial f(\bar{z}) $, define the mapping 
$$ J_A : \R^n \ni z \mapsto f(\bar{z}) + A ( z - \bar{z}) + BF(Cz). $$
\begin{enumerate}[topsep=0pt,itemsep=1ex,partopsep=1ex,parsep=1ex]
\item[(a)] The mapping $\Phi$ is strongly metrically sub-regular \at{z}{p} provided that for each $ A \in \partial f(\bar{z}) $, one has that
 \begin{equation*}
\left.\begin{matrix}
(Cb, -(B^T B)^{-1} B^T A b ) \in T ((C \bar{z}, \bar{v}); \gph{F}) ~~\\
~~~~~~~~~~ ~ A b \in \mathrm{~rge~} B ~~
\end{matrix}\right\}
 \Rightarrow b = 0;
 \end{equation*}
\item[(b)] The mapping $\Phi$ is strongly metrically regular at $\bar{z}$ for $\bar{p}$ provided that for each $ A \in \partial f(\bar{z}) $, one has that
\begin{itemize}[topsep=-1ex,itemsep=-1ex,partopsep=1ex,parsep=1ex]
\item[(i)] for each neighborhood $U$ af $\bar{z}$ there is a neighborhood $V$ of $ \bar{p} $ such that \\
 $ J_A ^{-1}(p) \cap U \not = \emptyset $ whenever $p \in V$;
\item[(ii)]
 \begin{equation*}
\left.\begin{matrix}
(Cb, -(B^T B)^{-1} B^T Ab ) \in \widetilde{T} ((C \bar{z}, \bar{v}); \gph{F}) ~~\\
~~~~~~~~~~ Ab \in \mathrm{~rge~} B 
\end{matrix}\right\}
 \Rightarrow b = 0.
 \end{equation*}
\end{itemize}
\end{enumerate}
\end{thm}

\begin{proof}
(a) For each $ A \in \partial f(\bar{z}) $, the function $f^{*} (\cdot) : = f(\bar{z}) + A ( \cdot - \bar{z} ) $ is continuously differentiable on $\R^n$; so, condition (A2) is satisfied and the mapping $J_A$ is strongly metrically sub-regular \at{z}{p} by Theorem \ref{SMR-Theorem 2} (ii) with $ \Phi := J_A = f^{*} + BF(C \, \cdot)$. Now, apply Theorem \ref{SMSR-Theorem 3} with $ g= f$, $G= F$, and $ \mathcal{A} = \partial f(\bar{z}) $, to get the conclusion.\\
(b) The conditions (i) and (ii) guarantee that, for each $ A \in \partial f(\bar{z}) $, the mapping $J_A$ is strongly metrically regular \at{z}{p} (using Theorem \ref{SMR-Theorem 2} (iii) ). By Theorem \ref{Izmailov's theorem}, $\Phi$ is strongly metrically regular \at{z}{p}.
\end{proof}

\begin{eg} \textbf{(A simple circuit with SCR)} \citep{cibulkaRoubal} \label{A simple circuit with SCR} \\
Consider the circuit in Figure \ref{fig: SCR-circuit-chap3}, with an SCR, whose $ i - v $ characteristic is given by 

\begin{figure}[ht]
  \begin{minipage}{.50\textwidth}
    \begin{eqnarray*}
		G(z) :=
		\left \{ \begin{array}{lcl}
		a z + V_1  & & z < 0, \\
		$[$ V_{1}, \, \varphi(0) $]$  & &  z = 0, \\
		\varphi(z) & & z \in [0, \alpha], \\
		a (z - \alpha) + \varphi (\alpha) & & z > \alpha,
		\end{array} \right.
	\end{eqnarray*} 
  \end{minipage}
  \begin{minipage}{.45\textwidth}
    \centering
		\includegraphics[width=0.87\textwidth]{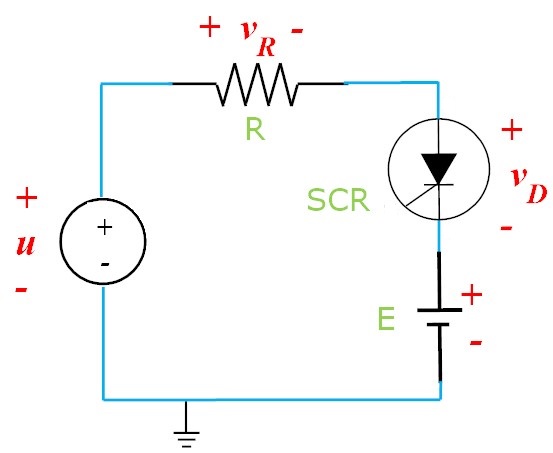}
  \end{minipage}
\caption{A simple circuit with SCR}
\label{fig: SCR-circuit-chap3}
\end{figure}

where $ -V_1, \, \alpha $, and $ a $ are positive constants, and $ \varphi : \R \longrightarrow (0, \infty) $ is a continuously differentiable function with 
$ V_1 < \varphi(0) $, $ \varphi(\alpha) < \varphi(0)$, $ \varphi'(0) > 0 $, and $ a > \varphi'(\alpha) > 0 $.
Note that $ G = g+ F $ with
\begin{eqnarray*}
		~F(z) :=
		\left \{ \begin{array}{ll}
		V_1  & z < 0, \\
		$[$ V_{1}, \, \varphi(0) $]$ ~ &  z = 0, \\
		\varphi (0)  & z >0,
		\end{array} \right.
		\mathrm{~~and~~~}
		g(z) :=
		\left \{ \begin{array}{ll}
		a \, z  &  z < 0, \\
		\varphi(z) - \varphi(0)  &   z = [0, \, \alpha ], \\
		a (z - \alpha) + \varphi (\alpha) - \varphi (0) &  z > \alpha.
		\end{array} \right.
\end{eqnarray*} 
A review of Example \ref{DIAC-Chap2} reveals that by setting $ p = u - E $ and $ z = i $, we get \eqref{inclusion} with $ m = n = 1 $, $ B = C = I_1 $, and 
$ f(z) = Rz + g(z) $, for $ z \in \R $. \\
Then $f$ is locally Lipschitz continuous on $ \R $ with
\begin{equation*}
\partial f(z) = \left \{ \begin{array}{lll}
R +  a & & z < 0, \\
$[$ R + a, \, \varphi'(0) $]$ & & z= 0,\\
R + \varphi'(z) & & z \in (0 , \alpha),\\
$[$ R + \varphi'(\alpha), R + a $]$  & & z = \alpha, \\
R +  a & & z > \alpha.
\end{array} \right.
\end{equation*}
Suppose that $  \varphi'(z) > - R $ for each $ z\in (0, \alpha) $. Then, for any $ z \in \R $, all the elements of $ \partial f(z) $ are positive. Assumptions (A1), $(\widetilde{\mathrm{A2}})$, (A3), and (A4), hold true. 
Given $ \rfp{z}{p} \in \gph{\Phi} $, we get that the assumptions of Theorem \ref{Summary of MR Criteria for GE in Case of a Nonsmooth Single-valued Part} $(b)$ are satisfied. Thus, $\Phi$ is strongly metrically regular at any reference point. \\
Let us do the computations for two reference points. For $ \bar{z} = 0, \, \bar{p} = V_1 $, we have $ f(\bar{z}) = 0 $, $\bar{v} = V_1$, $ J_A (z) = A \, z + F(z) $ where $ A \in [ R + a, \, \varphi'(0) ] $, and condition $(i)$ is satisfied (see Figure \ref{fig: SCR-circuit02-chap3}). Condition 
$ (ii) $ reads as
$$ (b, \, - Ab) \in \widetilde{T} \Big( (0, V_1); \, \gph{F} \Big) \Longrightarrow \, b= 0,  \mfa A \in [ R + a, \, \varphi'(0) ]. $$
Since $ \widetilde{T} \Big( (0, V_1); \, \gph{F} \Big) = \R 
{\scriptsize
\begin{pmatrix}
0 \\ 1
\end{pmatrix} \, \bigcup \, \R \begin{pmatrix}
-1 \\ 0
\end{pmatrix}}
$, 
and $ A > 0 $, the condition is satisfied. 

\begin{figure}[ht]
\centering
		\includegraphics[width=0.55\textwidth]{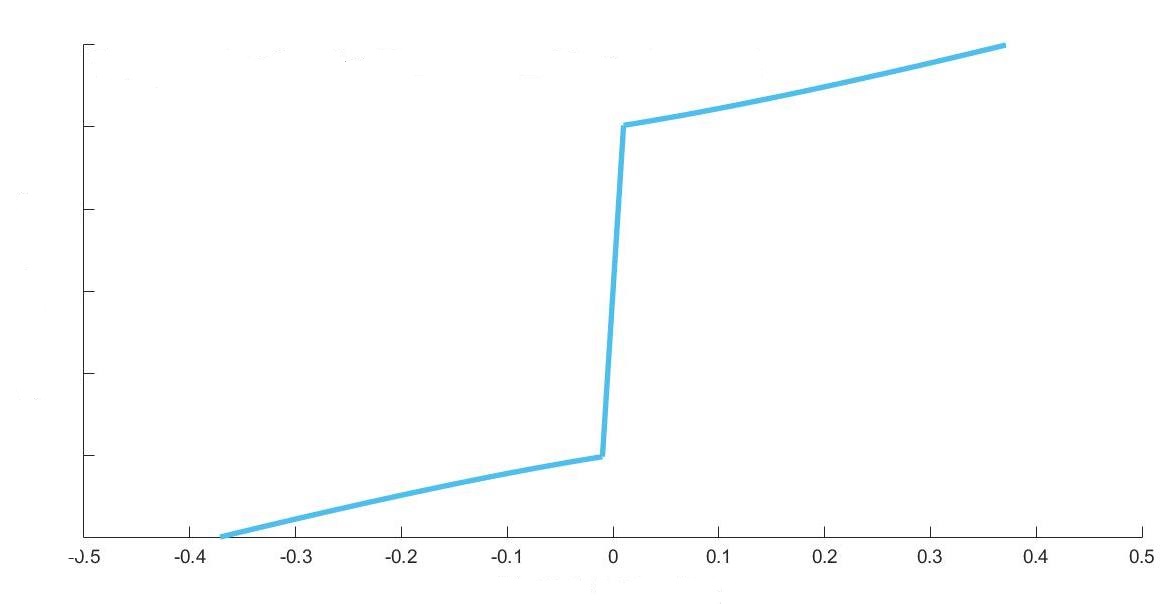}
\caption{General form of the map $ J_A $ in Example \ref{A simple circuit with SCR} }
\label{fig: SCR-circuit02-chap3}
\end{figure}

For $ \bar{z} = \alpha, \, \bar{p} = R \alpha + \varphi(\alpha) $, we have $\bar{v} = \varphi(0) $, and 
$ f(\alpha) = R \alpha + \varphi(\alpha) - \varphi(0) > 0 $. Then $ J_A (z) = f(\alpha) +  A \, (z - \alpha) + F(z) $ where $ A \in [ R + \varphi'(\alpha), R + a ]$, which is again of the form shown in Figure \ref{fig: SCR-circuit02-chap3}. Hence, condition $(i)$ is satisfied. \\
Moreover, 
$ \widetilde{T} \Big( \big( \alpha, \varphi(0) \big); \, \gph{F} \Big) = \R 
{\scriptsize
\begin{pmatrix}
1 \\ 0
\end{pmatrix} }
$, 
and thus, condition $(ii)$ also holds. 
\end{eg}

%% file: Chapters/Chapter04.tex
\chapter{Metric Regularity in the Case of Time-Varying Sources} 
\label{Chapter4} 
\lhead{Chapter 4. \emph{Metric Regularity in Presence of Time-Varying Sources}} 

\setlength{\epigraphwidth}{.49\textwidth}
\epigraph{\emph{Some mathematician,
I believe, has said that true pleasure lies not in the discovery of truth, but in the search for it.}}{\textit{Leo Tolstoy}}

In this chapter, we try to answer the question of perturbation effect when the vector $p$ becomes a function of time, which corresponds to the electronic circuits working with AC sources. We start with obtaining a proper model for the electrical circuit in this case, and explain the shift of our goal from studying the solution mapping to solution trajectories in Section \ref{Obtaining a Proper Model}. 
In Subsection \ref{Review on an existence theorem}, we express an important theorem (cf. Theorem \ref{theorem 6G.1}) to build our structure. \\
Section \ref{Results} starts with a proposition that has a simple proof, but contains a novel idea which considerably eases the study of the problem in the time varying case. In Subsection \ref{Continuity of Solution Trajectories}, we provide results that highlight smoothness dependence of trajectories on the input signal. \\
In Subsection \ref{Uniform Strong Metric Regularity-subsection}, we deviate from the study about solution trajectories for a while to present a uniform strong metric regularity result along a trajectory, with assumptions weaker than Theorem \ref{theorem 6G.1}.  It is necessary for obtaining the results about perturbation effect in Subsection \ref{Perturbations of the Input Signal}, where we provide a theorem that guarantees the existence of a solution trajectory for the perturbed problem which is continuous, and whose distance from the solution trajectory of the (non-perturbed) generalized equation is controlled by the distance between the input signal and its perturbed function.

\section{Obtaining a Proper Model} \label{Obtaining a Proper Model}
In this section, we would reconsider the situation described in Example \ref{Effect of AC Sources}, in order to introduce and study in depth the problem in the case of time-varying sources. We would see that there is a need to change the formulation and thus, the point of view towards the perturbation study. A primary result that guarantees the existence of continuous solution trajectories will be given in Subsection \ref{Review on an existence theorem}, Theorem \ref{theorem 6G.1}. This theorem will be a starting point for the detailed study of next sections.

\begin{figure}[ht]
	\centering
		\includegraphics[width=12.45cm]{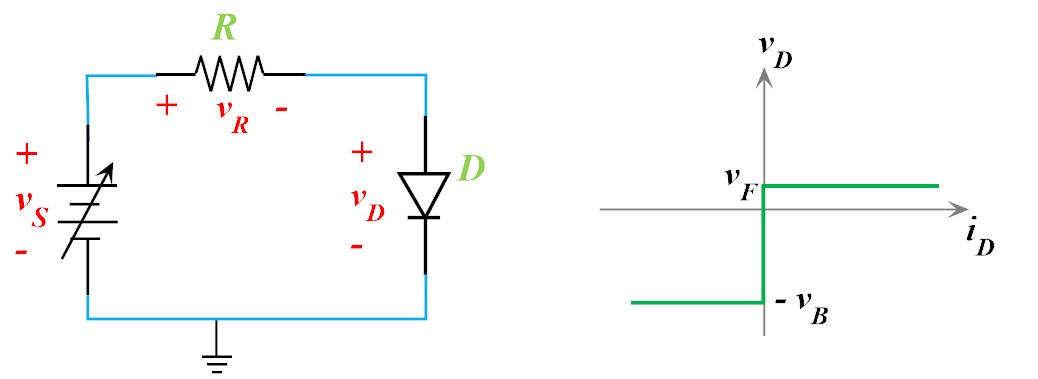}
	\caption{A regulator circuit with multiple DC sources}
	\label{fig: A regulator circuit with multiple DC sources}
\end{figure}

Let us start with an example. In Figure \ref{fig: A regulator circuit with multiple DC sources}, a simple regulator circuit with a practical model for the diode is shown. 
The voltage source is made of $n$ batteries connected to each other in a serial scheme, that provides $ n + 1 $ different levels between $ 0 \, V_s $, and $ 1 \, V_s$.

\begin{figure}[ht]
  \begin{minipage}{.40\textwidth}
  		\includegraphics[width=5.85cm]{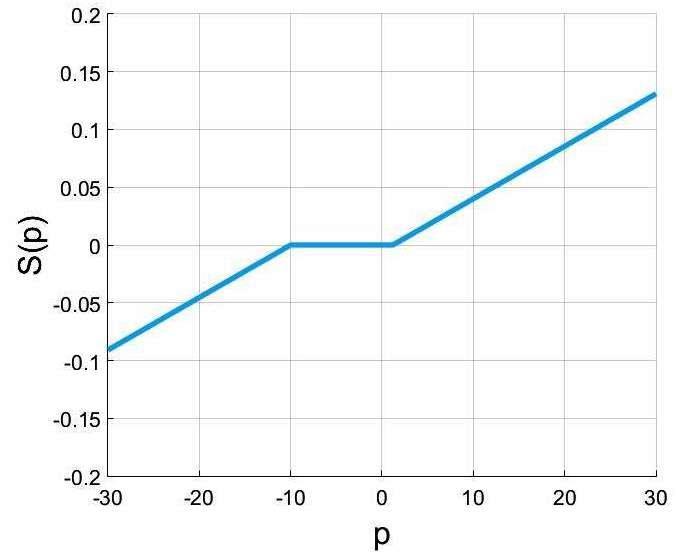}
  \end{minipage}
  \begin{minipage}{.60\textwidth}
    	\centering
		\includegraphics[width=8.45cm]{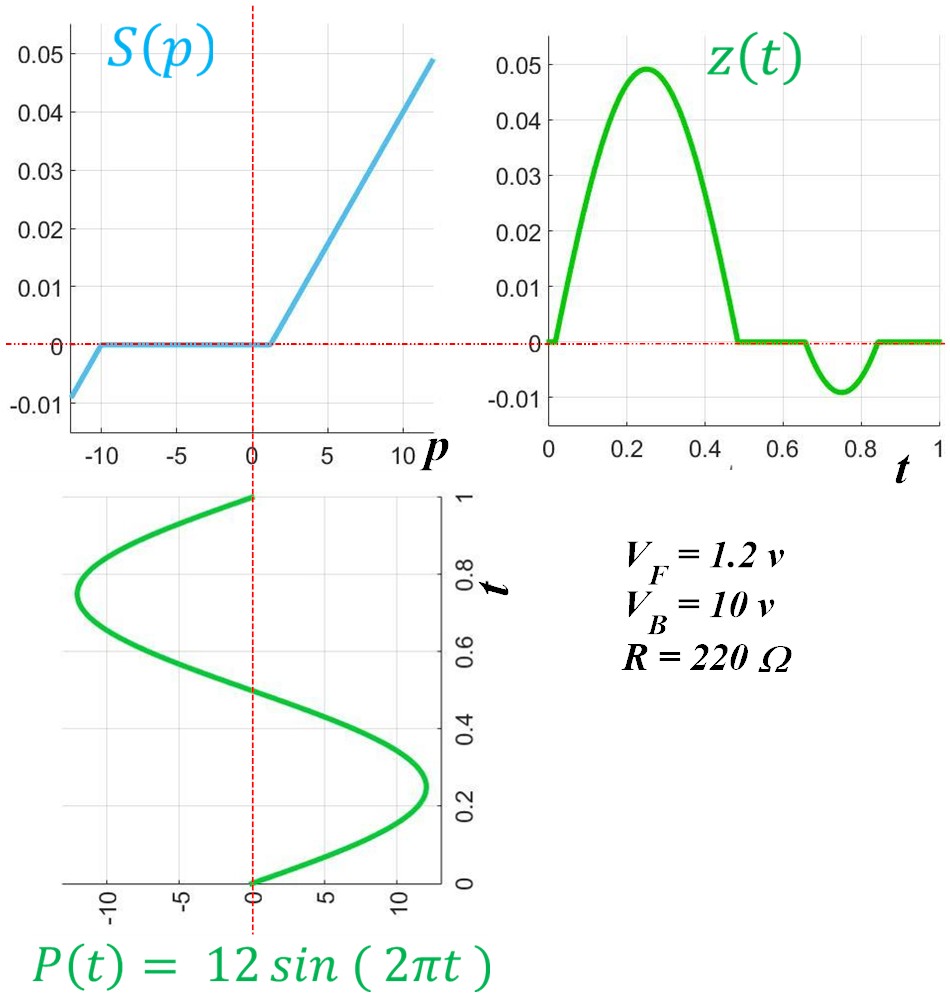}
  \end{minipage}
	\caption{Solution mapping for the circuit in Figure \ref{fig: A regulator circuit with multiple DC sources} (left), graphical method to find the output for a typical input (right)}
	\label{fig: Solution mapping for the circuit}
\end{figure}

Using Kirchhoff's laws and $ i- v $ characteristics of diode and resistor, we obtain that:
\begin{equation}
\left.\begin{matrix}
\mathrm{KVL: }  - \frac{m}{n} V_S + V_R+ V_D = 0~ \\ 
\mathrm{KCL: }~~~ I_{D} = I _R =: I~~~~~~~~~~~\\
~~~~~~\, V_R = R \, I_R ~~~~~~~~~~ \\
~~~~~~\, V_{D} \in F( I_{D}) ~~~~~~ ~~\\
\end{matrix}\right\}~~\Rightarrow
 \quad 0 \in - p + R z + F (z),
\end{equation}
where $ p = \frac{m}{n} V_S $, $ z = I $, and $ m = 0, 1, \, \dots , n $, indicates the number of turned-on batteries in the circuit. Then, the solution mapping would be 
$$ S(p) = \{ \, z \in \R \, | \, - p + R z + F (z) \ni 0 \, \}. $$
In order to find $S(p)$, we can consider the three parts of $\gph{F}$ separately to solve the generalized equation, fortunately, analytically this time.
\begin{enumerate}[topsep=-1ex,itemsep=-1ex,partopsep=1ex,parsep=1ex, leftmargin = 5ex]
\item[1.] For $ z > 0 $, with $ F(z) = \{ v_F \} $. \\
Then, we would have an equation, $ - p + R z + v_F = 0 $. Thus, $ z = \dfrac{p - v_F}{R}$ which is only valid for $ z > 0 $, that is, when $ p > v_F $.

\item[2.] For $ z = 0 $, with $ F(z) = [ - v_B, v_F ] $. \\
Then, $ - p + 0 + [ - v_B, v_F ] \ni 0 $. That is, $ z = 0 $ for $ p \in [ - v_B, v_F ] $.

\item[3.] For $ z < 0 $, with $ F(z) = \{ - v_B \} $. \\
Then, again we would have an equation, $ - p + R z - v_B = 0 $. Thus, $ z = \dfrac{p + v_B}{R}$, as long as $ p < - v_B $.
\end{enumerate} 
Therefore, $S$ is a single-valued map in this problem, with the graph shown in Figure \ref{fig: Solution mapping for the circuit} (left), and the rule given as:
\begin{equation}
S(p) = \left \{ 
\begin{matrix}
~ \Big \{ \dfrac{p + v_B}{R} \Big \} & & p < - v_B \\
 & & \\
~ \Big \{ 0 \Big \} & & ~~~~~~~p \in [ - v_B, v_F ] \\
 & & \\
~ \Big \{ \dfrac{p - v_F}{R} \Big \} & & p > v_F
\end{matrix} \right.
\end{equation}
When we deal with an AC voltage source, theoretically we can follow the same procedure. For any $ t \in [0, 1]$, use the specific value $p(t)$ and the transformation graph to find the value of $z$ at that time, that is $z(t)$. Then, we can obtain the graph of $z(\cdot)$ with respect to time, similar to the one shown in Figure \ref{fig: Solution mapping for the circuit} (right) for a sinusoid signal.\\
There are two interesting facts to highlight here:
\begin{enumerate}[topsep=-1ex,itemsep=-1ex,partopsep=1ex,parsep=1ex, leftmargin = 5ex]
\item[\textbf{(a)}] Very naturally, instead of asking for the graph of the solution mapping with respect to the input signal, we focused on the graph of the solution mapping with respect to the time. Of course, when the solution mapping is not a function like this problem, the latter expression needs a clarification. 

\item[\textbf{(b)}] Dealing with a function as the input signal, we searched for a function as the output signal. In order to keep the notations consistent, yet without ambiguity, we will refer to these functions as $p(\cdot)$, $\smap{p}{}{n}$, $z(\cdot)$, and so on. 
\end{enumerate}

Therefore, in the case of time varying sources, we can assume that $ t $ as a parameter, belongs to a set like $[0, 1]$, and it would be more appropriate to consider the solution mapping as the (generally set-valued) map that associates to every $ t \in [0, 1]$, the set of all possible vectors $ z$ in $\R^n$ that fits the generalized equation
\begin{equation}\label{pge}
f(z) - p(t) + F(z) \ni 0.
\end{equation}
The \emph{solution mapping} is therefore given by
\begin{equation}
S : t \mapsto S(t) = \{ z \in \R^n ~|~ f(z) - p(t) + F(z) \ni 0 \},
\end{equation}
and a \emph{\textbf{solution trajectory}}\index{solution trajectory} over $[0,1]$ is, 
in this case, a function $\bar{z}(\cdot)$ such that $ \bar{z}(t) \in S(t)$  for all $ t \in [0, ~ 1] $, that is, $\bar{z}(\cdot)$ is a selection for $S$ over $[0,1]$. 
Every solution trajectory defined in this way, would fit in Equation \eqref{time-varying-case-first-formulation} and fulfils our intuition.\\
Defining the solution mapping in terms of the parameter $t$, and not directly of the input signal $p$, will cause some difficulties to study the perturbation problem. \\
Comparing to the situation in Chapter \ref{Chapter3}, although for each $t$ one needs to solve a generalized equation of the type discussed in depth in Chapter \ref{Chapter3}, the Aubin property\footnote{
or any other local stability property like calmness or isolated calmness of $S$ or equivalently, metric regularities (all four different definitions) of $S^{-1}$.
} 
of the solution mapping at a certain point is not sufficient any more to guarantee the stability of the output with respect to the perturbations of the input signal. In other words, the relation between $S$ and $p(\cdot)$ is not explicitly expressed now.\\
So the first question one needs to answer is whether some pointwise local stability property of the new solution mapping can be related to a general statement about it or not. We will tackle this problem in Subsections \ref{Review on an existence theorem} and \ref{Uniform Strong Metric Regularity-subsection}. \\
Next, we will focus on the solution trajectories. In particular, we ask questions related to the existence of selections which are smooth functions, their relationship to the input signal, and their reaction to the small perturbations of the input signal. Section \ref{Results} is devoted to provide answers to these questions.

Let us note that under the general assumptions:
\begin{enumerate}[topsep=-1ex, itemsep=-1ex, partopsep=1ex, parsep=1ex, leftmargin = 9ex]
\item[(A2)] $f$ is continuously differentiable in $\R^n$; 
\item[(A3)] $F$ has closed graph;
\end{enumerate}
when $p$ is a continuous function, the map $S$ has closed graph.\\
Indeed, consider $ \{ (t_n), (z_n) \}_{n \in \N} $ such that $(t_n, z_n) \in \gph{\, S}~~\mathrm{for~every~} n$; and $ (t_n, z_n) \longrightarrow (t_0, z_0) $.
We show that $(t_0, z_0) \in \gph{\, S} $.\\
For each $n$, consider $ w_n := - f(z_n) + p(t_n) $. 
Since, by assumption, $z_n \in S(t_n) $, we get that $w_n \in F(z_n)$. Now, by the continuity of $p$ and $f$, we obtain that 
$ w_n \longrightarrow w_0 =  - f(z_0) + p(t_0)$.\\
Since $F$ has closed graph we get that $( w_0, z_0 ) \in \gph{\, F}$, i.e. $w_0 \in F(z_0)$. Therefore, $ F(z_0)  + f(z_0) - p(t_0) \ni 0 $, that is, $z_0 \in S(t_0)$.\\

Consider a function $ h: \R \times \R^n \longrightarrow \R^n $, defined as $h(t, v) = - p (t) + f(v) $. Then, the generalized equation (\ref{pge}) can be written as 
\begin{equation} \label{parametric generalized equation}
h(t,v) + F(v) \ni 0,
\end{equation}  
which is called a \emph{parametric generalized equation}, since the function $h$ now depends on a scalar parameter $t \in [0, 1] $.
For any given $(t,z) \in \gph{S}$, define the mapping
\begin{equation} \label{robinson formulation}
v \mapsto G_{t,z}(v) := h(t,z) + \nabla_z h(t,z) \, (v-z) + F(v).
\end{equation}
A point $ (t,z) \in \R^{1+n}$ is said to be a \emph{strongly regular}\index{strongly regular point} point\footnote{
\textbf{Historical note}.
This term was first introduced at 1980 by Robinson \citep{robinson}, and then became popular in the literature. In his work, the set-valued mapping $F$ is always a normal cone to a non-empty closed convex subset $C$ of a normed linear space, and he found this setting a convenient tool for formulating some problems in complementarity and in mathematical programming, as well as variational inequalities.\\
As he mentions, the idea behind this condition is that it is analogous to the non-singularity condition imposed in the usual implicit function theorem for nonlinear equations, and indeed, it reduces to that condition (in his setting) when $C$ is the whole space (so that the generalized equation reduces to the equation $ f(x) = 0 $).\\
Note that applying a condition on the \lq\lq linearization\rq\rq \,of \eqref{parametric generalized equation} will allow the numerical methods to enter the scene more effectively.
} 
for the generalized equation (\ref{pge}, or equivalently, \ref{parametric generalized equation}) when $ (t,z) \in \gph {S}$ and the mapping $G_{t,z}$ is strongly metrically regular at $z$ for $0$. That is, there exist constants $a_t, ~ b_t, ~ \lambda_t > 0$ such that the mapping
\begin{equation} \label{smr-constants}
\B_{b_t} (0) \ni y \longmapsto G_{t,z} ^{-1} (y) \cap \B_{a_t} (z(t) )
\end{equation}
is a Lipschitz continuous function with a Lipschitz constant $\lambda_t$.

The following theorem will show that the \emph{strong regularity condition} ensures a good behaviour of the solutions of the nonlinear problem. Although it is not explicitly mentioned, a closer look at the definition of estimators (Definition \ref{Estimators}) and the note after it (Note \ref{use of linearization as an estimator}), reveals its role behind the scene.

\begin{thm} [\textbf{Implicit Function Theorem for Generalized Equations}] \label{theorem 2B.7} \hfill \\  \cite[Theorem 2B.7, p. 89]{implicit} 
Consider a function $ f : \R^d \times \R^n \longrightarrow \R^n $ and a mapping $ \mmap{F}{n}{n} $ with $ \rfp{p}{x} \in \mathrm{int~dom}\, f $ and 
$ f\rfp{p}{x} + F(\bar{x}) \ni 0 $, and suppose that $ \widehat{\mathrm{lip}}_p (f; \rfp{p}{x}) \leq \gamma < \infty $.
Let $h$ be a strict estimator of $f$ with respect to $x$ uniformly in $p$ at $\rfp{p}{x}$ with constant $\mu$. 
Suppose that $ (h+F)^{-1} $ has a Lipschitz continuous single-valued localization $\sigma$ around $0$ for $\bar{x}$ with $ \mathrm{lip}(\sigma;0) \leq \kappa $ for a constant $\kappa$ such that $ \kappa \mu < 1$. Then the solution mapping
\begin{equation*}
S : p \longmapsto  \{ x \in \R^n ~|~ f(p,x) + F(x) \ni 0 \} ~~~~~\mathrm{for}~ p \in \R^d
\end{equation*}
has a Lipschitz continuous single-valued localization $s$ around $\bar{p}$ for $\bar{x}$ with
$$ \mathrm{lip}(s;\bar{p}) \leq \dfrac{\kappa \gamma}{ 1 - \kappa \mu}.$$
\end{thm}

\subsection{ Review on an Existence Theorem } \label{Review on an existence theorem}

From Theorem \ref{theorem 2B.7}, one obtains that when $(\bar{t},\bar{z})$ is a strongly regular point for (\ref{pge}), there are open neighborhoods $T$ of $\bar{t}$ and $U$ of $\bar{z}$ such that the mapping
\begin{equation}
T \cap [0, 1] \ni \tau \mapsto S(\tau) \cap U 
\end{equation}
is single-valued and Lipschitz continuous on $T \cap [0, ~1]$.\\
The theorem which follows shows that if each point in $ \gph{S} $ is strongly regular, then there are finitely many Lipschitz continuous solution trajectories defined on $ [0, 1] $ whose graphs never intersect each other. In addition, along any such trajectory $ \bar{u}(\cdot) $ the mapping $ G_{t, \bar{u}(t)} $ is strongly regular uniformly in $ t \in [0, 1] $, meaning that the neighborhoods and the constants involved in the definition do not depend on $t$.\\
Since we will work on the idea and assumptions of this theorem, and its proof needs some clarification, we add the proof here for the sake of completeness. We preferred to insert the added details as footnotes as long as it was possible.\footnote{
Check also \href{https://sites.google.com/site/adontchev/}{this web page} for Errata and Addenda.
}

\begin{thm} [\textbf{Uniform Strong Metric Regularity}] \index{metric regularity ! uniform strong metric regularity}\label{theorem 6G.1} \cite[Theorem 6G.1, p. 426]{implicit}\\
Suppose that there exists a bounded set $ C \subset \R^n$ such that, for each $ t \in [0,1]$, the set $S(t)$ is non-empty and contained in $C$ for all $t \in [0,1]$. Also, suppose that each point in $ \gph{S}$ is strongly regular. Then there are finitely many Lipschitz continuous functions $ \bar{u}_j : [0,1] \rightarrow \R^n, ~ j=1, 2, ... ,M $ such that, for each $t \in [0,1]$, one has 
$ \displaystyle S(t)= \bigcup_{1 \leq j \leq M} \{ \bar{u}_j(t) \}$. \\
Moreover, the graphs of the functions $ \bar{u}_j $ are isolated from each other, in the sense that there exists $ \delta > 0 $ such that
\begin{equation*}
\| \, \bar{u}_{j'}(t) - \bar{u}_j(t) \, \| \geq \delta \mathrm{~~for~every~} j' \not = j \mathrm{~~and~every~} t \in [0, 1].
\end{equation*}
Furthermore, there exist positive constants $ a, b$ and $ \lambda $ such that, for each such function $\bar{u}_i$, and for each $ t \in [0,1]$ the mapping
\begin{equation*}
\B_{b} (0) \ni w \longmapsto G_{t, \bar{u}_i (t)} ^{-1} (w) \cap \B_{a}  (\bar{u}_i (t) )
\end{equation*}
is a Lipschitz continuous function with a Lipschitz constant $\lambda$.
\end{thm}

\begin{proof}
From the assumed uniform boundedness of the solution mapping $S$ and the continuity of $f$ and its derivatives, we get the existence of a constant $ K > 0 $ such that
\begin{equation}
\sup_{ t \in [0,1], v \in C} \Big( \norm{ \nabla_t f(t,v) } + \norm{ \nabla_u f(t,v) }+ \norm{ \nabla^2 _{uu} f(t,v)} +|\norm{ \nabla^2 _{ut} f(t,v)}  \Big) \, \leq \, K.
\end{equation}
Let $ (t,v) \in \gph{S} $. Then, according to Theorem \ref{theorem 2B.7} there exists a neighborhood $ T_{t,v} $ of $t$ which is open relative to $ [0,1] $ and an open neighborhood $U_{t,v}$ of $v$ such that the mapping $T_{t,v} \ni \tau \mapsto S(\tau) \cap U_{t,v} $ is a function, denoted by $u_{t,v} (\cdot)$, which is Lipschitz continuous on $T_{t,v}$ with Lipschitz constant $L_{t,v}$. From the open covering $ \big\{ T_{t,v} \times U_{t,v} \big\}_{ (t,v) \, \in \, \gph{S}}$ of the graph of $S$, which is a compact set in $ \R^{1+n} $ (due to the boundedness assumption on $C$), we can extract a finite subcovering $ \Big\{ T_{t_j,v_j} \times U_{t_j, v_j} \Big\} _{j=1} ^M $. Let $ L = \displaystyle \max_{1 \leq j \leq M} L_{t_j, v_j} $.

We will prove the theorem in three steps corresponding to the following claims:\\
1. existence of finitely many Lipschitz continuous trajectories; \\
2. isolation of the graphs of the trajectories;\\
3. existence of uniform bounds.

\textit{\textbf{STEP 1.}} Let $ \tau \in [0,1] $ and choose any $ \bar{u} \in S(\tau) $. Now we will prove that there exists a Lipschitz continuous function $ \bar{u}(\cdot) $ with Lipschitz constant $L$ such that $ \bar{u}(t) \in S(t) $ for all $ t \in [0,1] $ and $ \bar{u}(\tau) = \bar{u} $.\\
Assume $ \tau < 1 $. Then there exists $ j \in \{ 1, \cdots ,M \} $ such that $ ( \tau, \bar{u}) \in T_{t_j,v_j} \times U_{t_j, v_j} $. Define $ \bar{u}(t) = u_{t_j, v_j} (t) $ for all $ t \in (t'_j, t''_j) := T_{t_j, v_j} $. Then $ \bar{u}(\tau) = u_{t_j, v_j } (\tau) = \bar{u} $ [first part of the claim] and $\bar{u}(\cdot)$ is Lipschitz continuous on $[t'_j, t''_j]$.
\footnote{
Let $ t \in (t'_j, t''_j) $, and consider a sequence $ (t^n _j) \to t''_j $ such that $ t^n _j < t''_j $ for each $n$. Then, by the Lipschitz continuity of $\bar{u}$ inside the interval, and the continuity of the norm we get
\begin{eqnarray*}
\begin{split}
\norm{ \bar{u}(t''_j) - \bar{u} (t) }  & = \norm{ \bar{u}(\lim_{n \to \infty} t^n _j ) - \bar{u}(t) } = \norm{  \lim_{n \to \infty} \bar{u} (t^n _j ) - \bar{u}(t) } = \lim_{n \to \infty} \norm{\bar{u}( t^n _j ) - \bar{u}(t)} \\
& \leq \, \lim_{n \to \infty} L \norm{t^n _j - t } = L \norm { t''_j - t }.
\end{split}
\end{eqnarray*}
So $\bar{u}(\cdot)$ is Lipschitz continuous at $t''_j$ (and also at $t'_j$ with a similar reasoning). It remains to show that $ \big( t''_j , \bar{u}(t''_j) \big) \in \gph{S} $. 
Since $ \gph{S}$ is a compact set and for each $n$, we have $ \big( t^n _j , \bar{u}(t^n _j) \big) \in \gph{S} $;
$$ t^n _j \rightarrow t''_j, \mathrm{~and~} \bar{u}(t^n _j) \rightarrow \bar{u}(t''_j) \mathrm{~as~} n \to \infty $$
we get the assertion.\\}
If $ t''_j < 1 $ then there exists some $ i \in \{1, \cdots ,M \} $ such that $ (t''_j, \bar{u}(t''_j)) \in T_{t_i, v_i} \times U_{t_i, v_i} := (t'_i, t''_i) \times U_{t_i, v_i} $. Then of course $ u_{t_i, v_i}(t''_j) = \bar{u}(t''_j) $ [single-valuedness of $ S \cap U_{t_i, v_i}$ ]; and we 
can extend $ \bar{u}(\cdot) $ to $ [ t_j, t''_i] $ as $ \bar{u}(t) = u_{t_i, v_i} (t) $ for $t \in [t_{j'}, t_{i''} ] $. After at most $M$ such steps we extend $ \bar{u}(\cdot) $ to $ [t_{j'}, 1] $. By repeating the same argument on the interval $ [0, \tau] $ we extend $ \bar{u}(\cdot) $ on the entire interval $[0, 1]$ thus obtaining a Lipschitz continuous selection for $S$. If $ \tau = 1 $ then we repeat the same argument on $ [0, 1] $ starting from $1$ and going to the left.\\
To finish the proof of the first claim, we assume that $ ( \tau, \bar{u}) $ and $ (\theta, \widetilde{u}) $ are two points in $ \gph{S} $ and let $ \bar{u}( \cdot) $ and $ u(\cdot) $ be the functions determined by the above procedure such that $ \bar{u}(\tau) = \bar{u} $ and $ \widetilde{u} (\theta) = \widetilde{u} $. \\
Assume that $ \bar{u}(0) \neq \widetilde{u}(0) $ and the set $ \Delta  := \{ t \in [0, 1] \, | \, \bar{u}(t) = \widetilde{u}(t) \} $ is non-empty. 
Since $ \Delta $ is closed\footnote{
In fact, $\Delta $ is the zero level set of the continuous function $\bar{u} - \widetilde{u}$.\\
}, 
$ \inf \Delta := v > 0 $ is attained and then we have that $\bar{u}(v) = \widetilde{u}(v)$ and $ \bar{u}(t) \not = \widetilde{u}(t) $ for $ t \in [0, v) $. But then $ \big(v, \bar{u}(v) \big) \in \gph{S} $ cannot be a strongly regular point of $S$, a contradiction\footnote{
To be more clear, in order to have $ ( v, \bar{u}(v) ) $ as a strongly regular point of $\gph{S}$, we should have that the mapping $ \tau \mapsto S(\tau) \cap U $ is single-valued in a vicinity of $(v, \bar{u}(v) )$, which is not.\\
}. 
Thus, the number of different Lipschitz continuous functions $ \bar{u}(\cdot) $ constructed from points $ (\tau, \bar{u}) \in \gph{S} $ is not more than the number of points in $ S(0) $\footnote{
Suppose $ ( 0 , u_i) \in \gph{S} $. Since every point in $\gph{S}$ is a strongly regular point, one has that the mapping $ T_{0 , u_i} \ni \tau \mapsto S(\tau) \cap U_{0 , u_i} $
is single-valued, especially, $ S(0) \cap U_{0 , u_i} $ is a singleton. Thus, $ S(0) = \displaystyle \bigcup_{i \in I} \big( S(0) \cap U_{0 , u_i} \big) $, for an arbitrary index set $I$. But since graph of $S$ is a compact subset of $\R^{n+1}$, there exists 
a finite set $ J \subset I $ such that $ \displaystyle S(0) \subset  \bigcup_{i \in J } \big( S(0) \cap U_{0 , u_i} \big) = \bigcup_{i \in J } \{ u_i \} $.
}. 
Hence there are finitely many Lipschitz continuous functions $ \bar{u}(\cdot) $ such that for every $ t \in [0, 1] $ one has $S(t)= \displaystyle \bigcup_{1 \leq j \leq M} \{ \bar{u}_j(t) \}$. This proves the first part of the theorem.

\textit{\textbf{STEP 2.}}
The fact that the solutions are isolated from each other, is implicitly shown in the contradiction above and follows from the fact that there are ``finitely many'' $\bar{u}_j(\cdot)$.\\
The argument is that if there are two different functions $ \bar{u}(\cdot) $, and $ \widetilde{u}(\cdot)$ with an intersection point, say $ \bar{u}(\tau) = \widetilde{u}(\tau) = \bar{u} $ for a point $ (\tau, \bar{u}) \in \gph{S} $, then $\Delta \not = \emptyset $ and we must have $  \bar{u}(0) = \widetilde{u}(0)  $, and so on for every $t$. Thus, the two functions have either one point of intersection, and are totally the same, or they do not have any intersection at all. 

\textit{\textbf{STEP 3.}}
Choose a Lipschitz continuous function $ \bar{u}(\cdot) $ whose graph is in the graph of $S$, that is, $ \bar{u}(\cdot) $ is one of the functions $ \bar{u}_j(\cdot) $ and its Lipschitz constant is $L$. Let $ t \in (0,1) $ and let $ G_t= G_{t, \bar{u}(t)} $\footnote{
One should be careful not to confuse this change of notation, which is only valid for the rest of this proof, with the auxiliary map $G_t$ that will be defined by \eqref{my auxiliary map} in Section \ref{Results} and will be used till the end of this chapter.\\
}, 
for simplicity. Let $ a_t, b_t $ and $ \lambda_t $ be positive constants such that the mapping
\begin{equation}
\B_{b_t} (0) \ni w \longmapsto G_{t} ^{-1} (w) \cap \B_{a_t}  (\bar{u}(t) )
\end{equation}
is a Lipschitz continuous function with Lipschitz constant $ \lambda_t $. Make $b_t > 0 $ smaller if necessary so that
\begin{equation} \label{6G.1-0}
2 b_t \lambda_t < a_t.
\end{equation}
Let $ \rho_t > 0 $ be such that $ L \rho_t < \dfrac{a_t }{2} $. Then, from the Lipschitz continuity of $\bar{u}$ around $t$ we have that 
$ \displaystyle \B_{\frac{a_t}{2}} (\bar{u}(\tau)) \subset \B_{a_t} (\bar{u}(t)) $ for all $ \tau \in (t - \rho_t, t + \rho_t) $\footnote{
Let $ z \in \B_{\frac{a_t}{2}} \bar{u}(\tau)) $ be an arbitrary point. Then
$$ \norm{ z - \bar{u}(t) } \, \leq \, \norm{ z - \bar{u}(\tau) } + \norm{ \bar{u}(\tau) - \bar{u}(t)} \, \leq \, \frac{a_t}{2} + L \norm{ \tau - t } \, < \, \frac{a_t}{2} + L \rho_t \, < \, \frac{a_t}{2} + \frac{a_t}{2}  $$
Thus, $ z \in \B_{a_t} ( \bar{u}(t)) $.\\
}. 
Make $ \rho_t > 0 $ smaller if necessary so that
\begin{equation} \label{6G.1-1}
K (L+1) \rho_t < \dfrac{1}{\lambda_t}, \mathrm{~~~and~~~}  \lambda_t \rho_t < \dfrac{a_t}{16 K}
\end{equation}
Our aim is to apply Theorem \ref{theorem 5G.3} (the strong regularity part) in order to show that there exist a neighborhood $O_t$ of $t$ and positive constants $ \alpha_t $ and $\beta_t$  such that for each $ \tau \in O_t \cap [0, 1] $ the mapping 
\begin{equation} \label{6G.1-4}
\B_{\beta_t} (0) \ni w \longmapsto G_{\tau} ^{-1} (w) \cap \B_{\alpha_t}  (\bar{u}(t) )
\end{equation}
is a Lipschitz continuous function. \\
Consider the function $ g_{t, \tau} : \R^n \longrightarrow \R^n $ defined as 
\begin{eqnarray}
\begin{split}
g_{t, \tau} (v) = & ~ f( \tau , \bar{u}(\tau)) - f(t, \bar{u}(t)) \\ 
& + \Big( \nabla _u f(\tau , \bar{u}(\tau)) - \nabla _u f(t, \bar{u}(t)) \Big) v \\
& + \nabla _u f(t, \bar{u}(t)) \bar{u}(t) - \nabla _u f(\tau , \bar{u}(\tau)) \bar{u}(\tau).
\end{split}
\end{eqnarray}
For each $v$ we have $G_{\tau} (v) = G_t (v) + g_{t, \tau} (v) $. Since $ g_{t, \tau} $ is an affine map, it is Lipschitz continuous. Let us show that the Lipschitz constant is bounded by the expression on the left of (\ref{6G.1-1}) \footnote{
Here, we have used the fact that $f$ is a twice continuously differentiable function. The following proof requires only Lipschitz continuity of $\nabla_u f$.
\begin{eqnarray*}
\begin{split}
\norm{g_{t, \tau} (v) - g_{t, \tau} (v') }  & \, \leq \,  \norm{ \nabla _u f(\tau , \bar{u}(\tau)) - \nabla _u f(t, \bar{u}(t)) } \norm{ v - v' }\\
& \, \leq \, L_{ \nabla_u f} \big( \norm{\tau - t} + \norm{ \bar{u}(\tau) - \bar{u}(t) } \big) \norm{ v - v' } \\
& \, < \, K' \big( \rho_t + L \rho_t \big) \norm{ v - v' } ~ = ~  K' (L + 1 ) \rho_t \norm{ v - v' }
\end{split}
\end{eqnarray*}
In which $L_{ \nabla_u f}$ is the Lipschitz constant of $\nabla_u f$, and $K'$ could be chosen in accordance with $\rho_t$, in a way that $ K' (L+1) \rho_t < \frac{1}{\lambda_t} $ holds true. \\
},
 by using a mean-value theorem in $\R^n$ \footnote{ \label{mean-value theorem}
For more details on mean value theorems in $\R^n$, refer to \cite[Section 3.2]{ortega1970}, specially Theorem 3.2.3, page 69. We used in the above proof the main idea of that theorem (instead of using the end result of it, to give the reader a scheme of the proof). To be more precise, for a G-differentiable function $\smap{g}{n}{m}$ defined on a convex subset $D_0 $ of $\R^n$, one has
$$ \norm{g(x) - g(y) } \leq \sup_{ 0 \leq s \leq 1} \norm{g'( x + s (y - x) )} \, \norm{x - y} \mfa \, x,y \in D_0. $$
}. 
For simplicity, let $p_{\tau} := (\tau, \bar{u}(\tau) )$, and $p_t := (t, \bar{u}(t) ) $.
\begin{align*}
 \norm{g_{t, \tau} (v) - g_{t, \tau} (v') }  \, \leq \, & \norm{ \nabla _u f(p_{\tau}) - \nabla _u f(p_{t}) } \norm{ v - v' } \\
 \leq \, & \norm{  \int_{0}^{1} \dfrac{d}{ds} \Big[ \nabla_u f \big( p_{t} + s (p_{\tau} - p_t) \big) \Big] ds \,} \norm{ v - v' } \\
 \leq \, & \int_{0}^{1} \norm{ \dfrac{d}{ds} \Big[ \nabla_u f \big( t + s (\tau - t), \, \bar{u}(t) + s ( \bar{u}(\tau) - \bar{u}(t) ) \big) \Big] } ds ~ \norm{ v - v' } \\
 \begin{split}
    \leq \, & \int_{0}^{1} \Big \| \nabla^2 _{ut} f (p_{t} + s (p_{\tau} - p_t))  (\tau - t) \\
              & \qquad \quad +  \nabla^2 _{uu}  f (p_{t} + s (p_{\tau} - p_t)) (\bar{u}(\tau) - \bar{u}(t)) \Big \| \, ds ~ \norm{ v - v' }
  \end{split}\\
  \begin{split}
	\leq \, & \Big( \sup_{ 0 \leq s \leq 1} \norm{ \nabla^2 _{ut} f \Big( p_{t} + s (p_{\tau} - p_t) \Big) } \norm{\tau - t} \\
	& \qquad \quad + \sup_{ 0 \leq s \leq 1} \norm{ \nabla^2 _{uu} f \Big( p_{t} + s (p_{\tau} - p_t) \Big) } \norm{ \bar{u}(\tau) - \bar{u}(t) } \Big) \norm{ v - v' } 
	\end{split}\\
\leq \, & K \big( \norm{\tau - t} + \norm{ \bar{u}(\tau) - \bar{u}(t) } \big) \norm{ v - v' } \\
\leq \, & K \big( \rho_t + L \rho_t \big) \norm{ v - v' } \\
\leq \, &  K (L + 1 ) \rho_t \norm{ v - v' }.
\end{align*}
To show the boundedness of $ \norm{ g_{t, \tau} (\bar{u}(t)) }$, one can use elementary calculus and the mean value theorem mentioned in Footnote (\ref{mean-value theorem}), again, to argue as follows
\begin{align*}
\norm{ g_{t, \tau} (\bar{u}(t)) } = \, & \norm{ f( \tau , \bar{u}(\tau)) - f(t, \bar{u}(t)) + \nabla _u f(\tau , \bar{u}(\tau)) \big( \bar{u}(t)-  \bar{u}(\tau) \big) } \\
= \, & \Big \| \int_{0}^{1} \dfrac{d}{ds} f \big( t + s  (\tau - t), \, \bar{u}(t) + s ( \bar{u}(\tau) - \bar{u}(t) ) \big) \,  ds 
+ \nabla _u f(p_\tau) \big( \bar{u}(t)-  \bar{u}(\tau) \big) \Big \| \\
	\begin{split}
	  = \, & \Big \| \int_{0}^{1}  (\tau - t) \nabla_t f \big( t + s  (\tau - t), \, \bar{u}(t) + s ( \bar{u}(\tau) - \bar{u}(t) ) \big) \, ds \\
	  & \qquad \quad + \int_{0}^{1} \big( \bar{u}(\tau)-  \bar{u}(t) \big) \nabla_u f \big( t + s  (\tau - t), \, \bar{u}(t) + s ( \bar{u}(\tau) - \bar{u}(t) ) \big) \, ds \\
	  & \qquad \quad - \int_{0}^{1} \nabla _u f(\tau, \bar{u}(\tau)) \big( \bar{u}(\tau) - \bar{u}(t)  \big) \, ds ~ \Big \|
	\end{split}\\
	\begin{split}
		\leq \, & \norm{\tau - t} \sup_{ 0 \leq s \leq 1} \norm{\nabla_t f \big( p_t + s ( p_{\tau} - p_t) \big)} \\
		& \qquad \quad + \norm{\bar{u}(t) - \bar{u}(\tau)} \norm{ \int_{0}^{1} \nabla_u f \big( p_t + s ( p_{\tau} - p_t) \big) - \nabla_u f(p_{\tau}) ds }
	\end{split}\\
\leq \, & K \rho_t + L \rho_t ~ \int_{0}^{1} K \norm{(s-1) (p_{\tau} - p_t)} ds \\
\leq \, & K \rho_t \big( 1 + \dfrac{1}{2} L \rho_t + \dfrac{1}{2} L^2 \rho_t \big).
\end{align*}
We apply Theorem \ref{theorem 5G.3} with $ F = G_t, ( \bar{x}, \bar{y} ) = (\bar{u}(t), 0), g = g_{t, \tau}, a = a_t, b = b_t, \kappa = \lambda_t $, and defining 
$ \mu = \mu_t := K (L+1) \rho_t $, and 
\begin{equation}
\kappa' = \lambda'_t : = \dfrac{3 \lambda_t}{2(1-K(L+1) \rho_t \lambda_t)} > \dfrac{\lambda_t}{1 - \mu_t \lambda_t }.
\end{equation}
For that purpose we need to show that there exist constants $\alpha_t$ and $\beta_t$ that satisfy the inequalities
\begin{equation} \label{6G.1-2}
\alpha_t \, \leq \, \dfrac{a_t}{2}, ~~~~ 2 \mu_t \alpha_t + 2 \beta_t \, \leq \, b_t, ~~~~ 2 \lambda'_t  \beta_t \, \leq \, \alpha_t, ~~~~ \norm{ g_{t, \tau} (\bar{u}(t)) } \leq \beta_t.
\end{equation}
Choose $ \rho_t $ smaller if necessary such that $  \dfrac{1}{2} L \rho_t + \dfrac{1}{2} L^2 \rho_t < 1 $, then the above calculations reveal that
 $ \norm{ g_{t, \tau} (\bar{u}(t)) } \leq 2 K \rho_t $.\\
 Denoting $ A := K (1+L) $ and $ B := 2K $ we have
 \begin{equation*}
 \mu_t = A \rho_t \mathrm{~~and~~}  \norm{ g_{t, \tau} (\bar{u}(t)) } \leq B \rho_t.
 \end{equation*}
Set $\beta_t := B \rho_t $. We will now show that there exists a positive $ \alpha_t $ which satisfies all inequalities in (\ref{6G.1-2}). \\
Substituting the already chosen $ \mu_t $ and $ \beta_t $ in (\ref{6G.1-2}), we obtain that $\alpha_t$ should satisfy 
\begin{equation} \label{6G.1-3}
\left \{ \begin{matrix}
\alpha_t \leq \dfrac{a_t}{2} ~~~~~~~~~~~~~~~~~~~~~  \\
~2 A \rho_t  \alpha_t + 2 B \rho_t \, \leq \, b_t ~~~~~~ \\
~3 B \lambda_t \rho_t \leq \alpha_t (1 - A \rho_t \lambda_t ).
\end{matrix} \right.
\end{equation}
The system (\ref{6G.1-3}) has a solution $ \alpha_t > 0 $ provided that
\begin{equation*}
\dfrac{3 B \lambda_t \rho_t}{1 - A \lambda_t \rho_t} \, \leq \, \dfrac{b_t - 2 B \rho_t}{2 A \rho_t} \mathrm{~~~and~~~} \dfrac{b_t - 2 B \rho_t}{2 A \rho_t} \, \leq \, \dfrac{a_t}{2}.
\end{equation*}
Thus, everything comes down to checking whether this system of inequalities is consistent. But this system is consistent whenever
\begin{equation*}
\rho_t (A \lambda_t b_t + 2B + 4 A B \lambda_t \rho_t) \, \leq \, b_t \, \leq \, \rho_t ( 2B + A a_t ),
\end{equation*}
which holds when $ A \lambda_t b_t + 4 A B \lambda_t \rho_t \leq  A a_t $; which in turn always holds because of the assumptions (\ref{6G.1-0}) and (\ref{6G.1-1}).\\
We are now ready to apply Theorem \ref{theorem 5G.3} from which we conclude that the mapping in (\ref{6G.1-4}), which is 
$$\B_{\beta _t} (0) \ni w \longmapsto G^{-1}_{\tau} (w) \cap \B_{\alpha _t } (\bar{u}(t))$$
is a Lipschitz continuous function with Lipschitz constant $ \lambda'_t$. The next step is to move the center of the second ball in the above intersection to $\bar{u}(\tau)$.
We claim that the mapping 
\begin{equation} \label{6G.1-5}
\B_{\frac{\beta _t}{2}} (0) \ni w \longmapsto G^{-1}_{\tau} (w) \cap \B_{\frac{\alpha _t }{2}} (\bar{u}(\tau))
\end{equation}
is a Lipschitz continuous function with Lipschitz constant $\lambda' _t $ (or maybe less).\\
First we observe that $ G^{-1}_{\tau} (w) \not = \emptyset $, for $ w \in \B_{\frac{\beta _t}{2}} (0) $.\\
Let $ u_1, u_2 \in G^{-1}_{\tau} (w) \cap \B_{\frac{\alpha _t }{2}} (\bar{u}(\tau)) $. From the Lipschitz continuity of $\bar{u}(\cdot)$ around $t$ we have that 
$\B_{\frac{\alpha _t }{2}} (\bar{u}(\tau)) \subset \B_{\alpha _t } (\bar{u}(t)) $ (choosing $\rho_t $ smaller if necessary such that $ L \rho_t < \alpha _t / 2 $). 
Hence, we get $ u_1, u_2 \in G^{-1}_{\tau} (w) \cap \B_{\alpha _t } (\bar{u}(t)) $, which is a single-valued map and thus, we get $u_1 = u_2$.\\
The last part is to show $ G^{-1}_{\tau} (w) \cap \B_{\frac{\alpha _t }{2}} (\bar{u}(\tau)) \not = \emptyset $.\\
For any $ w_0 \in  \B_{\frac{\beta _t}{2}} (0) $, we know that there exists a unique $u_0 := G^{-1}_{\tau} (w_0) \cap \B_{\alpha _t } (\bar{u}(t))$.\\
We also know that, by definition, $ \bar{u}(\tau) \in G^{-1}_{\tau} (0) $. Since $ \norm{\bar{u}(\tau) - \bar{u}(t)} < L \rho_t < \dfrac{ \alpha _t }{2}  $, we get $\bar{u}(\tau) \in \B_{\alpha _t } (\bar{u}(t)) $. So, $ \bar{u}(\tau) = G^{-1}_{\tau} (0) \cap \B_{\alpha _t } (\bar{u}(t)) $.\\
\begin{equation*}
\begin{split}
\norm{ u_0 - \bar{u}(\tau) } & = \norm{ [ G^{-1}_{\tau} (w_0) \cap \B_{\alpha _t } (\bar{u}(t)) ] - [  G^{-1}_{\tau} (0) \cap \B_{\alpha _t } (\bar{u}(t)) ]  }\\
& \leq \lambda'_t \norm{ w_0 - 0 } \\
& \leq \lambda'_t \, \frac{\beta_t}{2}  \, \leq \, \frac{\alpha_t}{2}.
\end{split}
\end{equation*}
Thus, $u_0 \in \B_{\frac{\alpha _t }{2}} (\bar{u}(\tau)) $. In fact, we get $u_0 = G^{-1}_{\tau} (w_0) \cap \B_{\frac{\alpha _t }{2}} (\bar{u}(\tau)) $. \\
A similar reasoning provides the Lipschitz property of the map in that neighborhood\footnote{
In fact, one may argue that in both (\ref{6G.1-4}) and (\ref{6G.1-5}), we deal with the graph of $ G^{-1}_{\tau} $, only the localization (i.e. the cutting area) changes. Since 
$ \Big( \B_{\frac{\beta _t}{2}} (0) \times  \B_{\frac{\alpha _t }{2}} (\bar{u}(\tau)) \Big) \subset \Big( \B_{\beta _t} (0) \times \B_{\alpha_t}  (\bar{u}(t)) \Big) $, if there exists a point in (\ref{6G.1-5}), it is the same point of (\ref{6G.1-4}); and the only thing that needs to be checked is whether the new localization is empty-valued or not.
}. \\
Now we are ready to unify the constants. From the open covering $ \displaystyle \bigcup_{t \, \in \, [0,1] } (t - \rho_t, t + \rho_t) $ of $[ \, 0, 1] $ choose a finite subcovering of open intervals $ ( t_i - \rho_{t_i}, t_i + \rho_{t_i}), ~ i = 1, 2, ..., m $. \\
Let 
\begin{equation*}
\begin{split}
& a := \min \big \{ \frac{\alpha_{t_i} }{2} ~|~ i = 1, ..., m \big \} ~~~~~\\
& \lambda := \max \{ \lambda'_{t_i} ~|~  i = 1, ..., m \} ~~~~~~\\
\mathrm{and~} & b := \min \Big \{~ \dfrac{a}{\lambda} ,~ \min \big \{ \frac{\beta_{t_i}}{2} ~|~  i = 1, ..., m \big \} \Big \}.
\end{split}
\end{equation*}
Since $ b \leq a / \lambda $, the observation in Lemma \ref{ratio law} (possibility to reduce \lq\lq proportionally\rq\rq \,the radii of balls in the definition of a strongly metrically regular map) applies; hence, for each $ \tau \in  ( t_i - \rho_{t_i}, t_i + \rho_{t_i}) \cap [0,1]$ the mapping
$$ \B_{b} (0) \ni w \longmapsto G^{-1}_{\tau} (w) \cap \B_{a} (\bar{u}(\tau)) $$
 is a Lipschitz continuous function with Lipschitz constant $ \lambda$.
Let $ t \in [0, 1] $; then $ t \in (t_i - \rho_{t_i}, t_i + \rho_{t_i} ) $ for some $ i \in \{1, \cdots ,m \} $. Hence the mapping 
$$ \B_{b} (0) \ni w \longmapsto G^{-1}_{t} (w) \cap \B_{a} (\bar{u}(t)) $$
is a Lipschitz continuous function with Lipschitz constant $\lambda$. The proof is complete.
\end{proof}

Although Theorem \ref{theorem 6G.1} is an important result in our study, there is an unpleasant assumption in that theorem which is: $S(t)$  is uniformly bounded. Even if this condition is fulfilled, it is hard to be checked since it requires the whole set $S(t)$ to be clarified and available for any $t \in [0, 1]$. \\
In the next section it will be shown that the uniform bound could be obtained without this extra assumption. Moreover, we will provide some results about the nature of the solution trajectories related to specific classes of input functions and their perturbations. 

\section{A Study about the Solution Trajectories} \label{Results}

As a completion to the shift of interest from studying the solution mapping of a generalized equation to searching for certain properties of the solution trajectories, we described in the previous section; in this section, we will present some results about the behaviour of the solution trajectories, first for the generalized equation (\ref{pge}), and then for the perturbed generalized equation. \\
Though we are inspired by Robinson's idea of strongly regular points in defining the auxiliary map (\ref{robinson formulation}), and the techniques in Theorems \ref{theorem 2B.7} and \ref{theorem 6G.1}, we find it more convenient to do some modifications in the setting in order to adapt it to our problem.\\
Since our aim in this section is the study of the solution trajectories with respect to variations of the input function, $p (t)$, and since working with the function $f(z)$ or with its first order approximation does not play an important role in our case (the proof of this statement will follow soon), we assume to deal with $f(z)$ itself and so to consider the auxiliary mapping
\begin{equation} \label{my auxiliary map}
v \mapsto G_{t}(v) := f(v) - p(t) + F(v).
\end{equation}
For more details on different possible choices of auxiliary maps and how the strong metric regularity would be affected, we state the following proposition:\\

\begin{prop} [\textbf{Different Auxiliary Maps}] \label{different auxiliary maps} \hfill \\
Given the generalized equation $ f(t,u) + F(u) \ni 0 $, where $ f: \R \times \R^n \longrightarrow \R^n $ is a function, and $ \mmap{F}{n}{n}$ is a set-valued map with closed graph, let us denote by $h$ a strict estimator of $f$ with respect to $u$ uniformly in $t$, at $\rfp{t}{u}$ with a constant $ \mu $ (cf. Definition \ref{Estimators}). Consider the following auxiliary maps:
\begin{eqnarray}
G_{\bar{t},\bar{u}}(v) =h(v) + F(v), \, \\
G_{\bar{t}}(v) = f(\bar{t}, v) + F(v). 
\end{eqnarray}
Then, $G_{\bar{t},\bar{u}}$ is SMR at $\bar{u}$ for $0$, if and only if $G_{\bar{t}}$ is SMR at $\bar{u}$ for $0$, provided that the regularity modulus of each map times 
$ \mu $ is less than $1$.
\end{prop}

\begin{proof}
First observe that, by definition of a strict estimator, $h(\bar{u}) = f(\bar{t},\bar{u})$ and so, $ 0 \in G_{\bar{t},\bar{u}}(\bar{u}) $ is equivalent to $0 \in G_{\bar{t}}(\bar{u}) $. Now taking into account the pointwise relation
$$ G_{\bar{t}}(v) =  f(\bar{t}, v) + F(v) = G_{\bar{t},\bar{u}}(v)  + f(\bar{t}, v) - h(v), $$
one can define a map $ \smap{g_{\bar{t},\bar{u}}}{n}{n} $ with $g_{\bar{t},\bar{u}}(v) = f(\bar{t}, v) - h(v) $. 
For any $ v_1, v_2 \in U$ (a neighborhood of $\bar{u}$), we get
\begin{equation*}
\begin{split}
\norm{ g_{\bar{t},\bar{u}}(v_1) - g_{\bar{t},\bar{u}}(v_2) } & = \norm{ f(\bar{t}, v_1) - h(v_1) - f(\bar{t}, v_2) + h(v_2) } \\
& = \norm{ e(\bar{t}, v_1) - e(\bar{t}, v_2) } \\
& \leq \mu \norm{v_1 - v_2 }
\end{split}
\end{equation*}
where the last inequality is obtained by definition of strict estimator. Thus, $g_{\bar{t},\bar{u}}$ is Lipschitz continuous around $\bar{u}$.\\
Now one can use Theorem \ref{Inverse Function Theorem for Set-Valued Mappings} with $G = G_{\bar{t},\bar{u}} $ and $g = g_{\bar{t},\bar{u}} $ and, by assuming that 
$ ~ \mu \, \mathrm{reg} \, (G_{\bar{t},\bar{u}} \,;\, \bar{u}) < 1 $, to conclude that $ g + G = G_{\bar{t}} $ has a Lipschitz continuous single-valued localization around $ 0 + g(\bar{u})$ for $\bar{u} $. 
Since $g_{\bar{t},\bar{u}} (\bar{u}) = h(\bar{u}) - f(\bar{t},\bar{u}) = 0 $, the latter could be expressed as the SMR of $G_{\bar{t}}$ at $\bar{u}$ for $0$.\\
The converse implication is satisfied in a similar way by letting $G = G_{\bar{t}} $, $g = - g_{\bar{t},\bar{u}} $ and assuming
 $ \mu \, \mathrm{reg} \, (G_{\bar{t}} \,;\, \bar{u})  < 1  $.
\end{proof}

\rem
\textbf{(a)} A closer look at the proof reveals that if $h$ is a strict estimator, then the regularity modulus of $G_{\bar{t}}$ and $ G_{\bar{t},\bar{u}} $ are related to each other with $ \kappa' = \dfrac{\kappa}{1 - \kappa \mu} $.\\
Considering a partially first order approximation of $f$ like $h(v) = f(\bar{t},\bar{u}) + \nabla f_u (\bar{t},\bar{u}) (v - \bar{u}) $, will result in the same modulus for auxiliary maps (since $\mu = 0$ in this case).

\textbf{(b)} 
One should note that in general, $h_1(\cdot) = f(\bar{t},\cdot)$ is not a strict estimator of $f$ at the reference point. To guarantee this, one needs an extra assumption like the following:
\begin{center}
$f(t,\cdot)$ is Lipschitz continuous, for any $ t $ in a neighborhood of $\bar{t}$.
\end{center}
However, this is not a necessary condition. For example, in the specific case we are interested in, that is $ h(t,u) = f(u) - p(t) $, $h_1$ is automatically a strict estimator with $\mu = 0$ (in fact, a partial first order approximation), without requiring a Lipschitz continuity assumption on $p(\cdot)$.\\

\subsection{Continuity of Solution Trajectories} \label{Continuity of Solution Trajectories}
In this subsection we will discuss the smoothness relation between the input signal and solution trajectories under the strong metric regularity assumption of the auxiliary map \eqref{my auxiliary map}. 
For future reference we remind that the \emph{solution mapping}, in our setting, is defined as
\begin{equation} \label{solution mapping for time-varying case}
S : t \mapsto S(t) = \{ z \in \R^n ~|~ G_t (z) \ni 0 \},
\end{equation}
and a function $ z : [0, 1] \to \R^n $ is called a \emph{solution trajectory} if
\begin{equation} \label{solution trajectory for time-varying case}
z(t) \in S(t),  \mfa t \in [0, 1].
\end{equation}
Throughout the whole subsection we will assume that, given a function $p(\cdot)$, a solution trajectory $z(\cdot)$ exists. The reason we mentioned the existence result in Subsection \ref{Review on an existence theorem} is to guarantee that this assumption is not nonsense. \\
Let us start with a simple observation that will be used several times in this chapter. The following lemma will provide a rule for moving from one auxiliary map to another. This simple yet handy result is a consequence of our choice of auxiliary map and our setting.

\begin{lem} \label{G-relations}
Consider the generalized equation \eqref{pge}, and the auxiliary map (\ref{my auxiliary map}). For arbitrary points $t, t' \in [0, 1] $, the following equalities hold
\begin{eqnarray} \label{set-relations} 
\boxed{ G_{t} (v) = G_{t'} (v) + p(t') - p(t) } \, ~~\\
\boxed{ G_{t} ^{-1} (w) = G_{t'} ^{-1} \big( w + p(t) - p(t') \big) }
\end{eqnarray}
\end{lem}

\begin{proof}
The first equality is trivial. One only needs to write down the definition of the auxiliary map $G_t$. For the second one, consider $ v \in G_{t} ^{-1} (w) $. Then, 
$ w \in G_t(v)  = G_{t'} (v) + p(t') - p(t) $, implies $ v \in G_{t'} ^{-1} \big( w + p(t) - p(t') \big) $. \\
Since there is nothing special about $ t $ and $t'$, the inverse inclusion also holds true, and hence the equality is proved.
\end{proof}

In the following proposition, we will prove a continuity result for a given solution trajectory under suitable assumptions. One of the assumptions is that \lq\lq different\rq\rq \,trajectories, that is, trajectories without intersections, may not get arbitrary close to each other. In Theorem \ref{theorem 6G.1}, we have already seen this assumption as a result, which comes after continuity of the trajectories. Here, as a somehow inverse statement, we start from \lq\lq isolated trajectories\rq\rq \,and prove their continuity. \\
Since continuity of the trajectories is not assumed any more, one should be careful about how to formulate this property. We use an expression based on the graphs of trajectories (see \citep{MRDGE2016}). Example \ref{Isolated trajectories} after the proposition will provide some ideas for the difficulties that may arise by \lq\lq bad\rq\rq \,formulations.

\begin{prop} [\textbf{Smoothness Dependence of Trajectories on Input Signal}] \label{claim1}
For the generalized equation \eqref{pge}, and the solution mapping \eqref{solution mapping for time-varying case}, assume that 
\begin{itemize} [topsep=-1ex,itemsep=0ex,partopsep=1ex,parsep=1ex, leftmargin = 7ex]
\item[(i)] $ z(\cdot) $ is a given solution trajectory which is \textbf{isolated}\index{isolated trajectories} from other trajectories; that is, 
there is an open set $ \mathcal{O} \in \R^{n+1} $ such that
\begin{equation} \label{isolated-solution-real}
\{ (t, v) ~|~  t \in [0, 1] \mathrm{~and~} 0 \in G_t (v) \} \cap \mathcal{O} = \gph{ z }.
\end{equation}
\item[(ii)] $ p (\cdot) $ is a continuous function;
\item[(iii)] $G_t$ is pointwise strongly metrically regular; i.e. for any $ t \in [0, 1]$, $G_t$ is strongly metrically regular at $z(t) $ for $0$, with constants $a_t, b_t$, and $\kappa_t > 0 $ defined as \eqref{smr-constants}.
\end{itemize}
Then $ z(\cdot) $ is a continuous function.
\end{prop}

\begin{proof}
Fix $t \in [0, 1]$. We know that $ (t, z(t) ) \in \gph{S} $, so $ 0 \in G_t (z(t)) $ or $ z (t) \in G_t ^{-1} (0)$. For any $\epsilon > 0$, let 
 $ \epsilon_1 := \min \, \{ \frac{\epsilon}{\kappa_t}, b_t \}$, in which $b_t$ is the radius of the neighborhood around $0$ in the assumption $(iii)$.
 By the uniform continuity of $p(\cdot)$, there exists $\delta_p > 0 $ such that 
\begin{equation*}
\norm{ p(t) - p(\tau) } < \epsilon_1  \mathrm{~~whenever~~} \| \, \tau - t \, \| < \delta_p .
\end{equation*}
Let $ \delta < \delta_p $ and consider $\tau \in [0, 1] $ such that $ \| \, \tau - t \, \| < \delta $. By definition, $z(\tau) \in G_{\tau} ^{-1} (0)$. Using Lemma \ref{G-relations} we obtain 
$ z(\tau) \in G_t ^{-1} ( p(\tau) - p(t) ) $. \\
On the other hand, by assumption $ (i) $ we also know that $z(\tau) \in \B_{a_t} (z(t)) $. \\
Indeed, assuming $z(\tau) \not \in \B_{a_t} (z(t)) $, allows us to define a Lipschitz continuous function $\widetilde{z}$ as 
$$ \widetilde{z} (\tau) := G_t ^{-1} \big( p(\tau) - p(t)  \big) \cap \B_{a_t} (z(t)) $$
on $\B_{\delta} (t) $. By Lemma \ref{G-relations}, $ \widetilde{z} (\tau) \in G_{\tau} ^{-1} (0) $ and thus, is (part of) a solution trajectory. 
Now, consider a sequence $ (t_n) $ in $\B_{\delta} (t) $ converging to $t$, and recall that, by definition, 
$ z(t) = G_t ^{-1} ( 0 ) \cap \B_{a_t} (z(t))$. Thus, $ \widetilde{z} (t_n) \longrightarrow z(t) $.\\
This means that $\widetilde{z}$ is a solution trajectory that could get arbitrarily close to $z(\cdot)$ at $ \big(t, z(t) \big)$, which contradicts assumption $(i)$.\\
So $ z(\tau) \in G_t ^{-1} ( y ) \cap \B_{a_t} (z(t))$ where $ y \in \B_{b_t} (0)$.\\
Now, by assumption $(iii)$, the mapping $  G_t ^{-1} ( \cdot ) \cap \B_{a_t} (z(t)) $ is single-valued and Lipschitz continuous on $  \B_{b_t} (0) $ with Lipschitz constant $\kappa_t$. So
\begin{equation*}
\left \| z(t) - z(\tau) \right \| \, \leq \, \kappa_t \left \| p(t) - p(\tau) \right \| \, < \, \epsilon.
\end{equation*} 
Since $t$ was an arbitrary point in $[0, 1]$, the proof is complete.
\end{proof}

\rem
If we assume that $p(\cdot)$ is a Lipschitz continuous function, then following the
previous proof by considering $ \tau_1, \tau_2 \in \B_{\delta} (t) $, we can obtain
\begin{equation*}
\left \| z(\tau_1) - z(\tau_2) \right \| \, \leq \, \kappa_t \norm { \, [p(\tau_1) - p(t) ] - [ p(\tau_2) - p(t) ] \, } \, \leq \, \kappa_t L_p  \, \norm{ \tau_1 - \tau_2 }.
\end{equation*}
This means that $ z(\cdot) $ is locally Lipschitz on $ [0, 1] $ which is a compact set; so it is globally Lipschitz and we can restate the proposition as the following corollary.

\begin{cor} \label{claim1-corollary}
Assume that
\begin{enumerate}[topsep=-1ex,itemsep=-1ex,partopsep=1ex,parsep=1ex, leftmargin = 7ex]
\item [(i)] $z(\cdot) $ is a given continuous solution trajectory; 
\item [(ii)] $p (\cdot) $ is a Lipschitz continuous function;
\item [(iii)] $G_t$ is pointwise strongly metrically regular at $z(t) $ for $0$.
\end{enumerate}
Then $ z(\cdot) $ is a Lipschitz continuous function.
\end{cor}

\begin{eg} \label{Isolated trajectories} 
Having a look at Theorem \ref{theorem 6G.1}, one can think of defining isolation of trajectories in this way:
$\bar{u}_{j}$ is isolated from other trajectories, in the sense that there exists $ \delta > 0 $ such that
\begin{equation} \label{isolated-solution-fake}
\| \, \bar{u}_{j'}(t) - \bar{u}_j(t) \, \| \geq \delta \mathrm{~~for~every~} j' \not = j \mathrm{~~and~every~} t \in [0, 1].
\end{equation}
Since we have not proved the continuity of trajectories and a selection could be made as a function with \lq\lq jump\rq\rq, this condition could not be valid for any trajectory. Consider the functions $z_1, z_2$, and $\bar{z}$ in Figure \ref{fig: isolated-trajectories}.
\begin{equation*}
\left \{\begin{matrix}
z_1 : [0, 1] \to \R \\
z_1(t) = 1~~~~~~
\end{matrix} \right.,
~~~~~~~
\left \{\begin{matrix}
z_2 : [0, 1] \to \R \\
z_2(t) = 0 ~~~~~~
\end{matrix} \right.,
~~~~~~~
\bar{z}(t) = 
\left \{\begin{matrix}
z_1 (t) & t \in [0, t_0) \\
z_2 (t) & t \in [t_0, 1] 
\end{matrix} \right.
\end{equation*}

\begin{figure}[ht]
	\centering
		\includegraphics[width=0.43\textwidth]{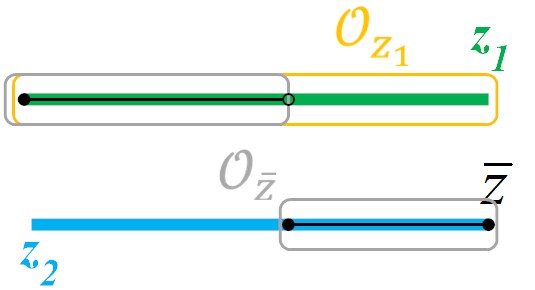}
	\caption{Difficulty in defining isolated trajectories}
	\label{fig: isolated-trajectories}
\end{figure}
\end{eg}
In this example, solution $\bar{z}$ is not isolated from $ z_1 $, nor $ z_2 $ in the sense of \eqref{isolated-solution-fake}. In fact, there is no solution trajectory which is isolated from all other trajectories (infinitely many others!). Thus, criteria \eqref{isolated-solution-fake} is nonsense. \\
However, using the formulation of \eqref{isolated-solution-real}, one can see that $ z_1$ is isolated from $z_2$, and $ \bar{z} $; the same holds for $z_2$. 
While for the trajectory $ \bar{z} $, the situation is a bit different. Any open set $ \mathcal{O} \in \R^{2} $ around the $\gph{\bar{z}}$ (like the grey rectangles in the figure), may include points of $ \gph{z_1} $, or $ \gph{z_2}$ which does not belong to $\gph{\bar{z}}$. Hence, $ \bar{z} $ (and any other solution trajectory that could be made in this way) is not isolated from $ z_1 $, nor $z_2$, and automatically drops out of the discussion of Proposition \ref{claim1}.  

\subsection{Uniform Strong Metric Regularity} \label{Uniform Strong Metric Regularity-subsection}
In this subsection we focus our attention on the uniform strong metric regularity of $ G_t $. 
One can consider this subsection as a quite independent section, but since we need the results we obtain here in the following subsection where some properties of the solution trajectories will be discussed, we prefer to keep logical hierarchy, and put this subsection under the general name of \lq\lq a study about the solution trajectories\rq\rq. \\
We have already seen Theorem \ref{theorem 6G.1} about uniform strong metric regularity. Our aim is to provide statements under simpler conditions, adapted to our particular setting. \\
In order to clarify the next statement, we remind that pointwise strong metric regularity of $ G_t $ for all $ t \in [0, 1]$, guarantees for each $t \in [0, 1]$ the existence of constants $a_t, \, b_t, \, \kappa_t > 0$ such that the mapping
\begin{equation*}
\B_{b_t} (0) \ni y \longmapsto G_{t} ^{-1} (y) \cap \B_{a_t} (z(t) )
\end{equation*}
is single valued and Lipschitz continuous with a Lipschitz constant $\kappa_t$.\\
It is worthwhile noting that the radii $a_t, \, b_t $ can be decreased provided that a suitable proportion is kept. The details are expressed in following lemma.

\begin{lem} [\textbf{Proportional Reduction of Radii}] \label{ratio law}
Let $H$ be a strongly metrically regular map at $\bar{x}$ for $\bar{y}$ with a Lipschitz constant $\kappa > 0$ and neighborhoods $ \B_a (\bar{x}) $ and  $ \B_b (\bar{y})$. Then for every positive constants
\begin{center}
 $ a' \leq a $ and $ b' \leq b $ such that $ \kappa b' \leq a' $, 
\end{center}the mapping $H$ is strongly metrically regular with the corresponding Lipschitz constant $ \kappa $ and neighborhoods $ \B_{a'} (\bar{x}) $ and  $ \B_{b'} (\bar{y}) $.
\end{lem}

\begin{proof}
Since $ B_{b'} (\bar{y}) \subset B_{b} (\bar{y})$ by assumption, $ H^{-1} (y) \cap \B_a (\bar{x}) = : x $ for every $y \in B_{b'} (\bar{y})$. 
Taking into account that $ H^{-1} (\cdot) \cap \B_a (\bar{x}) $ is a Lipschitz continuous function on $ B_{b} (\bar{y}) $, and by definition, $ \bar{x} : = H^{-1} (\bar{y}) \cap \B_a (\bar{x}) $, we get:
\begin{equation*}
\norm{ x - \bar{x} } = \norm{ \big( H^{-1} (y) \cap \B_a (\bar{x}) \big) - \big( H^{-1} (\bar{y}) \cap \B_a (\bar{x}) \big) } \, \leq \, \kappa \norm{ y - \bar{y} } \, \leq \, \kappa b'.
\end{equation*}
So, under the condition $ \kappa b' \leq a' $, we get $ x \in \B_{a'} (\bar{x}) $. \\
Indeed, in this case any $ y \in B_{b'} (\bar{y}) $ will be in the domain of $ H ^{-1} (\cdot) \cap \B_{a'} (\bar{x}) $. Then, the single-valuedness and Lipschitz continuity are the consequences of dealing with the same map (that is, $\gph{H}$).
\end{proof}

\rem \label{radii effect}
Fixing a point $t$, and the corresponding radii  $ a_t $, and $b_t$, we want to study the effect that increasing or decreasing these radii can have on the properties of the graphical localization.\\ Choosing a bigger $a$ (letting $ a > a_t $) may cause the graphical localization loose single valuedness by adding extra points to the localized map; while decreasing $a$ may omit some points from the localized map and therefore, cause emptiness. \\
Increasing $b$, enlarges the domain of the localized map and may cause both unpleasant situations; multi-valuedness and emptiness. But even if the new points remain properly (i.e. in a single-valued manner) inside the second neighbourhood, the function may act not enough smoothly and cause loosing the Lipschitz property of the localized map. On the other hand, decreasing $b$ may not harm anything as it only drops some well-behaved points out of our scope. 

\begin{figure}[ht]
	\centering
		\includegraphics[width=0.85\textwidth]{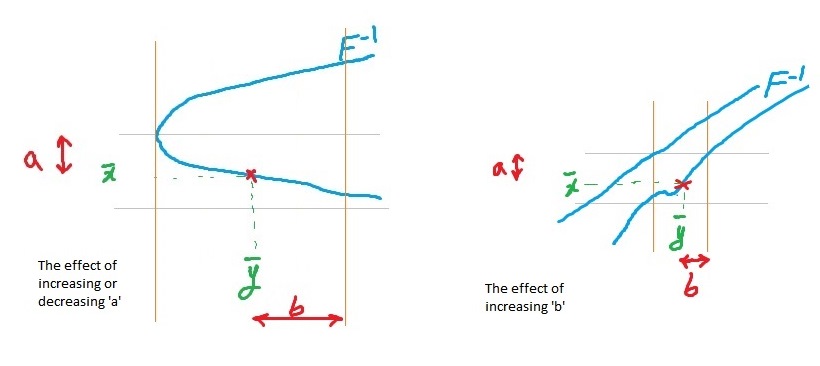}
	\caption{The side effects of changing radii independently}
	\label{fig: radii-effect}
\end{figure}

The previous lemma suggests that working with both radii $a$, and $b$ at the same time is allowed (with a careful control, of course) in order to keep the desired properties of the localized map. \\

\begin{thm} [\textbf{Uniform Strong Metric Regularity}] \label{claim2} \hfill \\
For the generalized equation \eqref{pge}, and the solution mapping \eqref{solution mapping for time-varying case}, assume that 
\begin{enumerate}[topsep=-1ex, itemsep=-1ex, partopsep=1ex, parsep=1ex, leftmargin = 7ex]
\item [(i)] $z(\cdot) $ is a given continuous solution trajectory;
\item [(ii)] $p (\cdot) $ is a continuous function;
\item [(iii)] $G_t$ is pointwise strongly metrically regular at $z(t) $ for $0$.
\end{enumerate}
Then there exist constants $a, \, b, \, \kappa > 0$ such that for any $ t \in [0, \,1] $, the mapping
\begin{equation*}
\B_{b} (0) \ni y \longmapsto G_{t} ^{-1} (y) \cap \B_{a} (z(t) )
\end{equation*}
is single valued and Lipschitz continuous with a Lipschitz constant $\kappa$.
\end{thm}

\begin{proof}
We prove the statement in two steps. First, by showing the mentioned map must be single-valued without caring about the Lipschitz regularity, and then by proving it is a Lipschitz continuous function. 
\par
\textit{\textbf{STEP 1. Single-valuedness: }}\\
We show that there exist $ a, b > 0 $ such that for any $ t \in [0, 1]$, the map 
\begin{equation} \label{usmr-01}
\B_b (0) \ni y \longmapsto G_{t} ^{-1} (y) \cap \B_{a} (z(t) ) 
\end{equation}
is single-valued.
We argue by contradiction, by assuming that for any $ a, \, b > 0 $, there exists $t_{a,b} \in [0,\, 1] $ such that \eqref{usmr-01} does not hold. In particular, take 
$ a_n =\frac{1}{n},  b_n = \frac{1}{n^3} $; then, for every $n \in \N$, there exists $t_n (:= t_{a_n, b_n} ) \in [0, \, 1] $ such that
\begin{equation} 
\B_{b_n} (0) \ni y \longmapsto G_{t_n} ^{-1} (y) \cap \B_{a_n} (z(t_n) ) 
\end{equation}
is not single-valued, which is equivalent to 
\begin{enumerate}[topsep=-1ex, itemsep=-1ex, partopsep=1ex, parsep=1ex, leftmargin = 10ex]
\item[\textit{Case 1.}] there exists $ y_n \in \B_{b_n} (0) $ such that the cardinality of the set $  G_{t_n} ^{-1} (y_n) \cap \B_{a_n} (z(t_n) ) $ is grater than one, or 
\item[\textit{Case 2.}] there exists $ y'_n \in \B_{b_n} (0) $ such that the set $  G_{t_n} ^{-1} (y'_n) \cap \B_{a_n} (z(t_n) ) $ is empty\footnote{
In other words, the mapping
\begin{equation*}
y \longmapsto G_{t_n } ^{-1} (y) \cap \B_{a_n} (z(t_n ) )
\end{equation*}
for at least a point $ y \in \B_{b_n} (0) $, is empty, or multivalued, that is, it has at least two values.
}. 
\end{enumerate}
By replacing $ (t_n) $ with a subsequence (if necessary), from the compactness of $[0, 1] $ in $\R$, we can assume that $ t_n \longrightarrow t_0 \in [0, \, 1]$. We will try to reach a contradiction in each case.

\textit{Case 1. Multi-valuedness}\\
Suppose that, for any $ n \in \N$, there exist $t_n \in [0, \, 1] $ and at least a $y_n \in \B_{b_n } (0)$ such that 
 $z_n ^1, z_n ^2 \in G_{t_n } ^{-1} (y_n) \cap \B_{a_n} (z(t_n ) ) $ with $z_n ^1 \not = z_n ^2$. \\
By assumption $(iii)$, there exist constants $a_{t_0}, b_{t_0}, \kappa_{t_0} > 0 $ such that the mapping 
\begin{equation*}
\B_{b_{t_0}} (0) \ni w \longmapsto G_{t_0} ^{-1} (w) \cap \B_{a_{t_0}} (z(t_0) )
\end{equation*}
is single valued and Lipschitz continuous with Lipschitz constant $ \kappa_{t_0} $.\\
Make $b_{t_0} > 0$ smaller if necessary so that
\begin{equation}\label{relation1}
 \kappa_{t_0} b_{t_0} \leq a_{t_0}.
\end{equation}
For $n$ large enough (i.e. $ n > N_0 \in \N$), one can have the following:
\begin{equation}\label{relation2}
b_n < \dfrac{b_{t_0}}{2},~ ~ \norm{ p(t_n) - p(t_0) } < \dfrac{b_{t_0}}{2},~ ~ \norm{ z(t_n) - z(t_0) } < \dfrac{a_{t_0}}{2}, ~~ a_n < \dfrac{a_{t_0}}{2}, ~~ \kappa_{t_0} < n^2,
\end{equation} 
in which the second and third inequalities are the results of continuity assumptions of $p(\cdot)$ and $z(\cdot)$, respectively. Then,
\begin{equation*}
\norm{ z_n ^1 - z(t_0) } \leq \norm{z_n ^1 - z(t_n)} + \norm{z(t_n) - z(t_0)} \leq  a_n + \frac{a_{t_0}}{2} < a_{t_0}.
\end{equation*}
The same holds for $z_n ^2 $; thus, $z_n ^1, z_n ^2  \in \B_{a_{t_0}} (z(t_0 ) ) $. On the other hand, $ z_n ^i \in G_{t_n } ^{-1} (y_n)$ for $ i = 1, 2$, and Lemma \ref{G-relations} implies that $ z_n ^i \in G_{t_0 } ^{-1} \big(y_n + p(t_n) - p(t_0) \big) $.
But 
$$ \norm{ y_n + p(t_n) - p(t_0) } \, \leq \, \norm{y_n - 0} + \norm{p(t_n) - p(t_0)} \, \leq \, b_n + \dfrac{ b_{t_0} } {2} \, < \,  b_{t_0}. $$
Thus, 
 $ (y_n + p(t_n) - p(t_0)) \in \B_{b_{t_0}}(0)$, which is a contradiction since, in that neighborhood, 
$ G_{t_0 } ^{-1} (\cdot) \cap  \B_{a_{t_0}}(z(t_0)) $ is single-valued. \\

\textit{Case 2. Emptiness}\\
Let us now suppose that, for any $n \in \N$, there exist $t_n \in [0, \, 1] $ and at least a point $y'_n \in \B_{b_n } (0)$ such that $ G_{t_n } ^{-1} (y'_n) \cap \B_{a_n}  (z(t_n ) )$ is empty. \\
For $n$ large enough, the inequalities in (\ref{relation1}) and (\ref{relation2}) hold, and we have already proved that $ y \in \B_{b_n} (0) $ implies 
$ y + p(t_n) - p(t_0)) \in \B_{b_{t_0}} (0) $. Therefore, since $y'_n \in \B_{b_n} (0)$, the mapping $  G_{t_0} ^{-1} (y'_n + p(t_n) - p(t_0)) \cap \B_{a_{t_0}} (z(t_0) ) $ is single-valued. In particular, it implies that $G_{t_0} ^{-1} (y'_n + p(t_n) - p(t_0)) \not = \emptyset $.
Let $z$ be a point in $ G_{t_0} ^{-1} (y'_n + p(t_n) - p(t_0))$. Then, by using Lemma \ref{G-relations}, we obtain 
$ y'_n \in G_{t_n} (z) $, in particular, $G_{t_n} ^{-1} (y'_n) $ is not empty. The contradiction assumption implies that
\begin{equation}\label{contradiction point}
\norm{ z - z(t_n) } > a_n.
\end{equation}
We will show the inconsistency between the contradiction assumption and the assumptions of the theorem with this inequality. In order to proceed, let us first prove that the mapping 
\begin{equation} \label{usmr-02}
\B_{\frac{b_{t_0}}{2}} (0) \ni y \longmapsto G_{t_n} ^{-1} (y) \cap \B_{\frac{a_{t_0}}{2}} (z(t_n) )
\end{equation}
is single-valued and Lipschitz continuous with Lipschitz constant $  \kappa_{t_0} $. \\
As a matter of fact, we have already seen that $ G_{t_n} ^{-1} (y) \neq \emptyset $ for every $ y \in \B_{b_n} (0) $, and $ G_{t_n} ^{-1} (y) \cap \B_{a_{t_0} / 2} (z(t_n) ) $ is not multi-valued. Thus, it only remains to show that 
$ G_{t_n} ^{-1} (y) \cap \B_{a_{t_0} / 2} (z(t_n) ) \not = \emptyset $ for every $ y \in \B_{b_{t_0} / 2} (0) $.\\
Denote by $z_y$ the point $ z_y := G_{t_0} ^{-1} \big( y + p(t_n) - p(t_0) \big) \cap \B_{a_{t_0}} (z(t_0) ) $. First observe that, by Lemma \ref{G-relations}, $ z_y \in G_{t_n} ^{-1} (y) $.\\
On the other hand, by definition, $ z(t_n) \in G_{t_n} ^{-1} (0) $ and by using Lemma \ref{G-relations}, we get $ z(t_n) \in G_{t_0} ^{-1} (p(t_n) - p(t_0)) $. We also know that $ z ( t_n) \in \B_{a_{t_0}} (z(t_0) ) $ (from the inequalities in \eqref{relation2}). The single-valuedness of $ G_{t_0} ^{-1} (.) \cap \B_{a_{t_0}} (z(t_0) ) $ over $ \B_{b_{t_0}} (0) $ allows us to write 
$ z(t_n) = G_{t_0} ^{-1} \big( p(t_n) - p(t_0) \big) \cap \B_{a_{t_0}} (z(t_0) ) $ without ambiguity. Thus, we have:
\begin{equation*}
\begin{split}
\norm{ z_y - z(t_n) } & = \scalebox{0.9}{ $ \norm{ [ \, G_{t_0} ^{-1}\paren{ y + p(t_n) - p(t_0) } \cap \B_{a_{t_0}} (z(t_0) ) \, ] -  [ \, G_{t_0} ^{-1} (p(t_n) - p(t_0)) \cap \B_{a_{t_0}} (z(t_0) ) \, ] } $ } \\
& \leq \kappa_{t_0} \norm{  y + p(t_n) - p(t_0) - ( p(t_n) - p(t_0) ) } \\
& \leq \kappa_{t_0} \norm{y - 0 } \\
& \leq \, \kappa_{t_0} \frac{b_{t_0}}{2} \, \leq \, \frac{1}{2} a_{t_0} .
\end{split}
\end{equation*}
Which means $ z_y \in G_{t_n} ^{-1} (y'_n) \cap  \B_{\frac{a_{t_0}}{2}} (z(t_n) )$.\\
A similar reasoning reveals the Lipschitz continuity of the map $ G_{t_n} ^{-1} (\cdot) \cap \B_{\frac{a_{t_0}}{2}} (z(t_n) ) $.\\ 
Indeed, taking any $y_i \in \B_{b_{t_0}/2} (0) $, one can define $z_i:= G_{t_n} ^{-1} (y_i) \cap \B_{a_{t_0}/2} (z(t_n) ) $ for $i =1, 2$ without ambiguity. Using the second and third inequalities in (\ref{relation2}), we have
\begin{equation*}
\begin{split}
& z_i \in G_{t_n} ^{-1} (y_i) = G_{t_0} ^{-1}\paren{ y_i + p(t_n) - p(t_0) } \mathrm{~~~~and~~~~}  y_i + p(t_n) - p(t_0)  \in \B_{b_{t_0}} (0) \\
& z_i \in  \B_{a_{t_0}/2} (z(t_n) ), \mathrm{~~~and~~~} \norm{ z(t_n) - z(t_0) } < a_{t_0}/2, \mathrm{~~so~~} z_i \in  \B_{a_{t_0}} (z(t_0) ).
\end{split}
\end{equation*}
Thus we are allowed to use the single-valuedness and Lipschitz property of $G_{t_0} ^{-1}$ to obtain
\begin{equation*}
\begin{split}
\norm{ z_1 - z_2 } & = \norm{ [ \, G_{t_n} ^{-1} (y_1) \cap \B_{a_{t_0}/2} (z(t_n) ) \, ] -  [ \, G_{t_n} ^{-1} (y_2) \cap \B_{a_{t_0}/2} (z(t_n) ) \, ] } \\
& = \scalebox{0.88}{ $\norm{ [ \, G_{t_0} ^{-1}\paren{ y_1 + p(t_n) - p(t_0) } \cap \B_{a_{t_0}} (z(t_0) ) \, ] -  [ \, G_{t_0} ^{-1} (y_2 + p(t_n) - p(t_0)) \cap \B_{a_{t_0}} (z(t_0) ) \, ] } $ } \\
& \leq \kappa_{t_0} \norm{y_1 - y_2 },
\end{split}
\end{equation*}
which reveals the Lipschitz property of the map in \eqref{usmr-02}. 

Now, having the strong metric regularity of $ G_{t_n}(\cdot)$ with constants $ \frac{a_{t_0}}{2}, \frac{b_{t_0}}{2}, \kappa_{t_0} $ in hand, we use Lemma \ref{ratio law} with $ a' = a_n = \frac{1}{n} \leq \frac{a_{t_0}}{2} $, $ b' = b_n = \frac{1}{n^3} \leq \frac{b_{t_0}}{2} $, to obtain the strong metric regularity of 
$ G_{t_n}(\cdot)$ with constants $ a_n, b_n, \kappa_{t_0} $ (reminding that the last inequality of (\ref{relation2}) guarantees $  \kappa_{t_0} b' \leq a' $).
Now for the specific $y'_n \in \B_{b_n } (0) $ claimed before, there exists $ z \in G_{t_n} ^{-1} (y'_n) \cap \B_{a_n} (z(t_n) )$ which 
contradicts (\ref{contradiction point}).\\
Therefore, till now we have proved that there exist $ a^*,b^*  > 0 $ such that the mapping
\begin{equation*}
\B_{b^*} (0) \ni y \longmapsto G_{t } ^{-1} (y) \cap \B_{a^*} (z(t) )
\end{equation*}
is single-valued for any $ t \in [0, 1] $. \\

\textit{\textbf{STEP 2. Lipschitz Continuity:}}\\
Being sure that we deal with a function, now we proceed by claiming that there exist $ b \leq b^* $, and  $ \kappa > 0 $ such that the mapping
\begin{equation*}
\B_{b} (0) \ni y \longmapsto G_{t } ^{-1} (y) \cap \B_{a^*} (z(t) )
\end{equation*}
is Lipschitz continuous with Lipschitz constant $\kappa $ for all $ t \in [0, 1] $.\\

We will prove the assertion by contradiction. Suppose the claim is false; then, for any $ b \leq b^* $, and any $ \kappa > 0 $, there exists $t_{b, \kappa}  \in [0, \, 1] $ such that the mapping 
\begin{equation*}
\B_{b} (0) \ni y \longmapsto G_{t_{b, \kappa} } ^{-1} (y) \cap \B_{a^*} (z(t_{b, \kappa}) )
\end{equation*}
is not Lipschitz with constant $\kappa$. Since we already know that this map is single-valued, due to the explanations in Remark \ref{radii effect}, the contradiction assumption yields that for every $ \kappa > 0 $, there exist $y_1, y_2 \in \B_b (0) $, with $ y_1 \not = y_2 $ such that 
\begin{equation*}
\norm{ \, \left [  G_{t_{b, \kappa} } ^{-1} (y_1) \cap \B_{a^*} (z(t_{b, \kappa}) ) \right ] - \left [  G_{t_{b, \kappa} } ^{-1} (y_2) \cap \B_{a^*} (z(t_{b, \kappa}) ) \right ] \, } \, > \, \kappa \norm{ y_1 - y_2 }.
\end{equation*}
For any $ n \in \N $, let $ b_n : = \min \{\, \frac{1}{n^3}, b^* \, \}, \kappa_n := n $ and set $ t_n := t_{b_n, \kappa _n } $. Then there exist at least
$y_n ^1, y_n ^2 \in \B_{b_n} (0) $, with $ y_n ^1 \not = y_n ^2 $ such that
\begin{equation*}
\left \| y_n ^1 -  y_n ^2 \right \| \, n \, < \, \norm{ \, \left [  G_{t_n } ^{-1} (y_n ^1) \cap \B_{a^*} (z(t_n) ) \right ] - \left [  G_{t_n } ^{-1} (y_n ^2) \cap \B_{a^*} (z(t_n) ) \right ] \, }
\end{equation*}
Let $g_n ^i := G_{t_n } ^{-1} (y_n ^i) \cap \B_{a^*} (z(t_n) ) $ for $ i = 1, 2 $, and assume that $t_n$ converges to a point, say $t_0$. \\
For $n$ large enough, one has the following:
\begin{equation*}\label{relation4}
b_n < \dfrac{b_{t_0}}{2},~ ~ \norm{ p(t_n) - p(t_0) } < \dfrac{b_{t_0}}{2},~ ~ \norm{ z(t_n) - z(t_0) } < \dfrac{a^*}{2},~~ \kappa_{t_0} < n^2.
\end{equation*} 
On the one hand, $g_n ^i \in G_{t_n } ^{-1} (y_n ^i) = G_{t_0 } ^{-1}\paren{ \, y_n ^i + p(t_n) - p(t_0) \, } $ and the above conditions imply that 
 $ \paren{ \, y_n ^i + p(t_n) - p(t_0) \, } \in \B_{b_{t_0}} (0) $.\\
On the other hand, $ g_n ^i \in B_{a^*} (z(t_n) ) $. We will show that $ g_n ^i \in B_{a_{t_0}} (z(t_0) ) $. \\
Indeed, let $w_n ^i := G_{t_0 } ^{-1}  \paren{\, y_n ^i + p(t_n) - p(t_0)} \cap B_{a_{t_0}} (z(t_0) ) $ for $ i = 1, 2 $. \\
Since $  \paren{\, y_n ^i + p(t_n) - p(t_0)} \longrightarrow 0 $, by the continuity of 
$G_{t_0 } ^{-1}  (\cdot) \cap B_{a_{t_0}} (z(t_0) ) $ around $0$, we get 
$$ w_n ^i \longrightarrow z(t_0) =  G_{t_0 } ^{-1}  (0) \cap B_{a_{t_0}} (z(t_0) ). $$
Thus, for any $\epsilon > 0 $, there exists $N_{\epsilon} \in \N$ such that, for $ n > N_{\epsilon}$, one has 
$ \norm{ w_n ^i - z(t_0) } < \epsilon $. Let $ \epsilon = a^* / 2 $. Then, 
\begin{equation*}
\begin{split}
\norm{ w_n ^i - z(t_n) } & \leq \norm{ w_n ^i - z(t_0) } + \norm{ z(t_0) - z(t_n) } \\
& <  a^* / 2 +  a^* / 2 \\
& <  a^*,
\end{split}
\end{equation*}
which means that $ w_n ^i \in \B_{a^*} (z(t_n)) $. Combining with $ w_n ^i \in  G_{t_n } ^{-1} (y_n ^i) $ (obtained by using Lemma \ref{G-relations}), 
we get that
$  w_n ^i \in  G_{t_n } ^{-1} (y_n ^i) \cap \B_{a^*} (z(t_n)) $.
Hence, by the single-valuedness of $ G_{t_n } ^{-1} (\cdot) \cap \B_{a^*} (z(t_n)) $, we can conclude that $  w_n ^i =  g_n ^i \in B_{a_{t_0}} (z(t_0) ). $ 
Then, the assumption (iii) of the theorem results in $ \norm {g_n ^1 -  g_n ^2} \, \leq \, \kappa_{t_0} \norm{y_n ^1 -  y_n ^2} $. So
\begin{equation*}
\norm{y_n ^1 -  y_n ^2} \, n \, < \, \norm {g_n ^1 -  g_n ^2} \, \leq \, \kappa_{t_0} \norm{y_n ^1 -  y_n ^2},
\end{equation*}
which is a contradiction, since it implies boundedness of the sequence $ (\kappa_n ) := ( n )$. Combining the two steps ends the proof.
\end{proof}

\rem
A slightly different version of this Theorem could be stated and proved as follows. Under stronger assumptions on $p(\cdot)$ and $z(\cdot)$, a simpler and more direct proof can be provided. The proof is in the direction of Theorem \ref{theorem 6G.1}, except that the special structure of the single-valued part here allows us to bypass the use of Theorem \ref{theorem 5G.3}. \\

\begin{thm} [\textbf{Uniform Strong Metric Regularity}] \label{claim2-second version} \hfill \\
For the generalized equation \eqref{pge}, and the solution mapping \eqref{solution mapping for time-varying case}, assume that 
\begin{enumerate} [topsep=-1ex,itemsep=-1ex,partopsep=1ex,parsep=1ex, leftmargin = 7ex]
\item [(i)] $z(\cdot) $ is a Lipschitz continuous solution trajectory with Lipschitz constant $L_z$;
\item [(ii)] $p (\cdot) $ is a Lipschitz continuous function with Lipschitz constant $L_p$;
\item [(iii)] $G_t$ is pointwise strongly metrically regular at $z(t) $ for $0$.
\end{enumerate}
Then there exist constants $a, \, b, \, \kappa > 0$ such that  the mapping
\begin{equation*}
\B_{b} (0) \ni y \longmapsto G_{t} ^{-1} (y) \cap \B_{a} (z(t) )
\end{equation*}
is single valued and Lipschitz continuous with a Lipschitz constant $\kappa$ for any $ t \in [0, 1] $.
\end{thm}

\begin{proof}
For an arbitrary point $t_0 \in [0, 1]$, by assumption $(iii)$, there exist constants $a_{t_0}, b_{t_0}, \kappa_{t_0} > 0$ such that the mapping
$$ \B_{b_{t_0}} (0) \ni y \longmapsto G_{t_0} ^{-1} (y) \cap \B_{a_{t_0}} (z(t_0) ) $$
is single-valued and Lipschitz continuous with Lipschitz constant $\kappa_{t_0}$. Choose $b_{t_0}$ smaller if necessary, so that 
 $  b_{t_0} \kappa_{t_0} < a_{t_0}  $.\\
Fix $\rho_{t_0} >0$ such that the following conditions are satisfied:
\begin{equation} \label{thm-condition} 
L_z \rho_{t_0} < a_{t_0} /2,~~~ L_p \rho_{t_0} < b_{t_0} /2
\end{equation}
Now, for any $ \tau \in ( t_0 - \rho_{t_0} , t_0 + \rho_{t_0} ) \cap [0, 1] $, we claim that the mapping
\begin{equation} \label{usmr-03}
\B_{b_{t_0}/2} (0) \ni y \longmapsto G_{\tau} ^{-1} (y) \cap \B_{a_{t_0}/2} (z(\tau) )
\end{equation} 
is single-valued and Lipschitz continuous with Lipschitz constant $\kappa_{t_0}$.\\
The proof will include the following simple steps:
\begin{enumerate} [topsep=-1ex,itemsep=-1ex,partopsep=1ex,parsep=1ex, leftmargin = 7ex]
\item[1.] the set $  G_{\tau} ^{-1} ( y ) $ is not empty for any $ y \in \B_{b_{t_0}/2} (0) $;
\item[2.] the mapping \eqref{usmr-03} is not multivalued;
\item[3.] the set $ G_{\tau} ^{-1} ( y ) \cap \B_{a_{t_0}/2} (z(\tau)) $ is a singleton for any $ y \in \B_{b_{t_0}/2} (0) $;
\item[4.] the mapping \eqref{usmr-03} is a Lipschitz continuous function (with constant $\kappa_{t_0}$).
\end{enumerate}

\textbf{1. }Since $ \norm{ p(\tau) - p(t_0)} \leq L_p \norm{\tau - t_0 } < L_p \rho_{t_0} < b_{t_0}/2 $, for any $ y \in \B_{b_{t_0}/2} (0) $, one gets 
$ y + p(\tau) - p(t_0)  \in \B_{b_{t_0}} (0)$; so $G_{t_0} ^{-1} ( y + p(\tau) - p(t_0) ) \not = \emptyset $. Thus, from Lemma \ref{G-relations} one concludes 
that $ G_{\tau} ^{-1} ( y ) \not = \emptyset $.

\textbf{2. }Let $u_1, u_2 \in G_{\tau} ^{-1} ( y ) \cap \B_{a_{t_0}/2} (z(\tau) ) $. 
Observe that $u_i \in G_{t_0} ^{-1} ( y + p(\tau) - p(t_0) )$  for $i =1, 2$, by Lemma \ref{G-relations}. 
Since Lipschitz continuity of $z(\cdot)$ implies that
$ \B_{a_{t_0}/2} (z(\tau) ) $ is a subset of $ \B_{a_{t_0}} (z(t_0) ) $, one obtains that $u_i \in G_{t_0} ^{-1} ( y + p(\tau) - p(t_0) ) \cap \B_{a_{t_0}/2} (z(t_0) ) $ for $i =1, 2$. 
Thus, assumption $(iii)$ implies $u_1 = u_2$.

\textbf{3. }For any $w_0 \in \B_{b_{t_0}/2} (0)$, let $u_0 := G_{t_0} ^{-1} ( w_0 + p(\tau) - p(t_0) ) \cap \B_{a_{t_0}} (z(t_0))$. \\
By definition, $ z(\tau) \in G_{\tau} ^{-1} (0) $, so from \eqref{thm-condition} and Lemma \ref{G-relations} one can obtain that 
$ z(\tau) \in G_{t_0} ^{-1} (0 + p(\tau) - p(t_0) )) \cap \B_{a_{t_0}} (z(t_0)) $. Then,
\begin{equation*}
\begin{split}
\norm{ u_0 - z(\tau) } & = \scalebox{0.9}{
 $ \norm{ [G_{t_0} ^{-1} ( w_0 + p(\tau) - p(t_0) ) \cap \B_{a_{t_0}} (z(t_0))] - [G_{t_0} ^{-1} (0 + p(\tau) - p(t_0) )) \cap \B_{a_{t_0}} (z(t_0))]  }$ } \\
& \leq \kappa_{t_0} \norm{ w_0 + p(\tau) - p(t_0) - (0 + p(\tau) - p(t_0)) } \\
& < \kappa_{t_0} \frac{b_{t_0}}{2} < \frac{a_{t_0}}{2}.
\end{split}
\end{equation*}
Thus, $u_0 \in G_{\tau} ^{-1} ( w_0 ) \cap \B_{a_{t_0}/2} (z(\tau))$. In fact, by considering the result of the previous part of the proof we obtain 
$u_0 = G_{\tau} ^{-1} ( w_0 ) \cap \B_{a_{t_0}/2} (z(\tau))$.

\textbf{4. }Consider any $ w_1, w_2 \in \B_{b_{t_0}/2} (0) $, and let $u_i = G_{\tau} ^{-1} ( w_i ) \cap \B_{a_{t_0}/2} (z(\tau))$ for $i = 1, 2$. Then, by using Lemma \ref{G-relations} we can conclude that
\begin{equation*}
\begin{split}
\norm{u_1 - u_2} & = \norm{ [G_{\tau} ^{-1} ( w_1 ) \cap \B_{a_{t_0}/2} (z(\tau))] - [G_{\tau} ^{-1} ( w_2 ) \cap \B_{a_{t_0}/2} (z(\tau))] } \\
& = \scalebox{0.9}{$ \norm{ [G_{t_0} ^{-1} ( w_1 + p(\tau) - p(t_0) ) \cap \B_{a_{t_0}/2} (z(\tau))] - [G_{t_0} ^{-1} ( w_2 + p(\tau) - p(t_0) ) \cap \B_{a_{t_0}/2} (z(\tau))] } $ } \\
& \leq \kappa_{t_0} \norm{  w_1 + p(\tau) - p(t_0) - ( w_2 + p(\tau) - p(t_0) ) } \\
& \leq \kappa_{t_0} \norm{ w_1 - w_2 }.
\end{split}
\end{equation*}

Now that we have local uniformity, choose a finite subcover $  \bigcup_{i=1} ^{m} (t_i - \rho_{t_i}, \,t_i +  \rho_{t_i}) $ from the open covering of the compact interval $[0, 1]$. Let
\begin{equation} \label{constants} 
\begin{split}
& a := \min \{ \frac{a_{t_i}}{2} ~|~ i = 1, ..., m \}, ~~~~~\\
& \kappa := \max \{ \kappa_{t_i} ~|~  i = 1, ..., m \}, ~~~~~~\\
\mathrm{and~} & b := \min \Big \{ \dfrac{a}{\kappa} ,~ \min \{ \frac{b_{t_i}}{2} ~|~  i = 1, ..., m \} \Big \}.
\end{split}
\end{equation}
For any interval $ (t_i - \rho_{t_i}, \,t_i +  \rho_{t_i}) $, we have the strong metric regularity with constants $ \dfrac{a_{t_i}}{2}, \dfrac{b_{t_i}}{2} $, and $\kappa $.  From \eqref{constants} and Lemma \ref{ratio law}, we obtain the strong metric regularity with constants $a, b$, and $ \kappa $, which does not depend on $t_i$ any more. Thus, the proof is complete.
\end{proof}

\subsection{Perturbations of the Input Signal} \label{Perturbations of the Input Signal}
In this subsection we try to take into account the small variations of the function $ p (\cdot) $. More precisely, we consider a continuous function $\widetilde{p} (\cdot) $ such that $ \norm { \widetilde{p} (t) - p(t) } < \epsilon $ for any $ t \in [0, \, 1] $, and for a suitably small $ \epsilon > 0 $. 
The reason of such perturbations and the importance of this study in the case of electronic circuits was already discussed in Chapter \ref{Chapter2} (cf. Section \ref{Formulating the problem}).\\
Unlike the static case, where $p $ was a fixed vector in $\R^n$, since we changed the notation of solution mapping, the problem does not reduce to the study of the stability properties of $S(t)$; however, we can take advantage of those results by considering the fact that in the newly defined set-valued map $G_t$ the variable $t$ reduces to a parameter.\\
In this subsection we deal with the perturbed form of problem (\ref{pge}). To be more specific, we consider the generalized equation 
\begin{equation} \label{perturbed GE}
f(z) - \widetilde{p} (t) + F(z) \ni 0,
\end{equation}
denote the corresponding solution mapping with $\widetilde{S}$, 
\begin{equation} \label{perturbed S}
\widetilde{S} : t \mapsto \widetilde{S}(t) = \{ z \in \R^n ~|~ f(z) - \widetilde{p} (t) + F(z) \ni 0 \},
\end{equation}
and define the auxiliary map $ \mmap{ \widetilde{G_t} }{n}{n}$ as $  \widetilde{G_t} (v) = f(v) + F(v) - \widetilde{p}(t) $. \\
The easy-to-check equalities
\begin{eqnarray} 
\boxed{ \widetilde{G_t}(v) = G_t(v)  + p(t) - \widetilde{p}(t) }~\,~~ \label{perturbation-relation1} \\
\boxed{ \widetilde{G_t} ^{-1} (w) = G_t ^{-1}  \big( w + \widetilde{p}(t) - p(t) \big)  } \label{perturbation-relation2}
\end{eqnarray}
for each $ t \in [0, 1] $, will be useful for connecting the strong metric regularity properties of $\widetilde{G_t} $ to those of $G_t$ as described in the following lemma. \\
Once more, we want to indicate that the straightforward equalities \eqref{perturbation-relation1} and \eqref{perturbation-relation2} are a consequence of our choice of the auxiliary maps and the special form of the single-valued part of the generalized equation \eqref{perturbed GE}. 

\begin{lem} [\textbf{Perturbation Effect on the Auxiliary Map}] \label{Perturbation Effect on the Auxiliary Map} \hfill \\
Assume that $ p(\cdot) $ and $ \widetilde{p}(\cdot) $ are continuous functions from $ [0, 1]$ to $ \R^n$ with $ \norm { \widetilde{p} (t) - p(t) } < \epsilon $ for any 
$ t \in [0, 1] $. If $G_t$ is strongly metrically regular at $\bar{u}$ for $0$ 
\big(i.e. $ (\bar{u},0) \in \gph \, G_t $ and there exist constants $a_t, \, b_t, \, \kappa_t > 0$ such that the mapping
\begin{equation*}
\B_{b_t} (0) \ni y \longmapsto G_{t} ^{-1} (y) \cap \B_{a_t} (\bar{u} )
\end{equation*}
is single valued and Lipschitz continuous with Lipschitz constant $ \kappa_t $\big), 
then for any positive $ \epsilon < b_t $ the mapping 
\begin{equation} \label{perturbed G} 
w \longmapsto \widetilde{G_t} ^{-1} (w) \cap \B_{a_t} (\bar{u})
\end{equation}
is a Lipschitz continuous function on 
$\B_{b_t - \epsilon} (0)$ with Lipschitz constant $\kappa_t$. 
\end{lem}

\begin{proof}
The process of proof is very similar to those steps we provided for the mapping \eqref{usmr-03} in the proof of the previous theorem, except that here, the intersecting ball $ \B_{a_t} (\bar{u} )$ is the same for both maps $ G_{t} ^{-1} $ and $ \widetilde{G_t} ^{-1} $ and this makes the proof much easier. The proof should include the following steps:
\begin{enumerate} [topsep=-1ex,itemsep=-1ex,partopsep=1ex,parsep=1ex, leftmargin = 7ex]
\item[1.] the sets $  \widetilde{G_t} ^{-1} (w)  $ and also $ \widetilde{G_t} ^{-1} (w) \cap \B_{a_t} (\bar{u}) $ are not empty for any $ w \in \B_{b_t - \epsilon} (0) $;
\item[2.] the mapping \eqref{perturbed G} is not multivalued;
\item[3.] the mapping \eqref{perturbed G} is a Lipschitz continuous function (with constant $\kappa_t$).
\end{enumerate}

Choose any $w_1, w_2 \in \B_{b_t - \epsilon} (0)$. From assumption we get $ \big( w_i + \widetilde{p}(t) - p(t) \big) \in \B_{b_t} (0) $ for $ i = 1, 2 $. Then, the pointwise strong metric regularity of $G_t$, lets us define
$ u_i := G_t ^{-1}  ( w_i + \widetilde{p}(t) - p(t) ) \cap \B_{a_t} (\bar{u}) $ for $ i = 1, 2 $.
By using \eqref{perturbation-relation2}, one obtains $ u_i \in  \widetilde{G_t} ^{-1} (w_i) $. In fact, $ u_i =  \widetilde{G_t} ^{-1} (w_i) \cap \B_{a_t} (\bar{u}) $.
Thus, steps 1. and 2. are proved. \\
But pointwise strong metric regularity of $G_t$ provides more information, that is
$$ \norm{u_1 - u_2} \leq \kappa_t \norm{w_1 - w_2}. $$
Therefore, step 3. is also proved.
\end{proof}

\begin{rem} \label{Perturbation Effect on the Auxiliary Map-remark}
\textbf{(a)} A careful look at the proof reveals that the lemma could be also expressed in the following way:
\begin{center}
If $G_t$ is SMR at $\bar{u}$ for $0$, then $\widetilde{G_t}$ is SMR at $\bar{u}$ for $  p(t) - \widetilde{p}(t) $.
\end{center}
In this case, $\epsilon$ could be as big as $b_t$.\\
In fact, in this case one can consider Theorem \ref{Inverse Function Theorem with Strong Metric Regularity} with $ F = G_t,~ \rfp{x}{y} =(\bar{u}, 0 ), ~ \kappa = \kappa_t $, and $g(\cdot) = p(t) - \widetilde{p}(t) $ which is a constant function with respect to $u$, so is obviously Lipschitz with any $ \mu < \kappa^{-1} $, and immediately obtain the SMR at $\bar{u}$ for 
$  p(t) - \widetilde{p}(t) $ of the map $ ( G_t + p(t) - \widetilde{p}(t) ) $ which is exactly $\widetilde{G_t}$.\\
\textbf{(b)} Under the assumptions of Theorem \ref{claim2}, we would have uniform strong metric regularity for $G_t$ at $z(t) $ for $0$ and the proof shows that we obtain uniform strong metric regularity for
 $\widetilde{G_t}$ at $z(t) $ for $ p(t) - \widetilde{p}(t) $, too.
\end{rem}

Finally, we have provided enough instruments to declare the main result of this section, that is the existence of a solution trajectory $\widetilde{z} (\cdot) $ close to $z(\cdot) $ that inherits its continuity properties.\\
We may recall that, under the assumptions of Theorem \ref{claim2}, existence of positive constants $a, b $, and $\kappa$ not depending on $t$ is guaranteed for uniform strong metric regularity. Since the following theorem satisfies those assumptions, we will use the uniform constants without ambiguity.

\begin{thm} [\textbf{Existence of a Continuous Trajectory for the Perturbed Problem}] \label{Existence of a Continuous Trajectory for the Perturbed Problem}
For the generalized equations \eqref{pge}, and \eqref{perturbed GE} and the corresponding solution mappings \eqref{solution mapping for time-varying case}, and \eqref{perturbed S},  assume that 
\begin{enumerate} [topsep=-1ex, itemsep=-0.5ex, partopsep=1ex, parsep=1ex, leftmargin = 7ex]
\item [(i)] $z(\cdot) $ is a given continuous solution trajectory (for $S$);
\item [(ii)] $p (\cdot)$ and $ \widetilde{p}(\cdot) $ are continuous functions such that for every $t \in [ 0, \, 1 ] $, $ \norm { \widetilde{p} (t) - p(t) } < \epsilon $ (with 
$ \epsilon < b/4$);
\item [(iii)] $G_t$ is pointwise strongly metrically regular at $z(t) $ for $0$.
\end{enumerate} 
Then there exists a continuous solution trajectory $\widetilde{z} (\cdot) $ for $\widetilde{S} $ such that, for every $t \in [ 0, \, 1 ] $, we have $ \norm { \widetilde{z} (t) - z(t) } <  \frac{4 a \epsilon}{b} $.
\end{thm}

\begin{proof}
We will present two proofs for this theorem, both are constructional methods, yet with different approaches. Remark \ref{another proof} after the proof, will provide a comparison between the methods.

\emph{Method 1. Pointwise construction:}\\
Consider an arbitrary $t_0 \in [0, \, 1]$. Since $ \big( t_0, z(t_0) \big) \in \gph{\, S}$ and $G_{t_0}$ is strongly metrically regular at $z(t_0)$ for $0$, by using Theorem \ref{claim2}, we obtain that the mapping
\begin{equation*}
\B_{b} (0) \ni y \longmapsto G_{t_0} ^{-1} (y) \cap \B_{a} (z(t_0) )
\end{equation*}  
is single-valued and Lipschitz continuous with constant $\kappa$. 
Let $y_0 =\widetilde{p}(t_0) - p(t_0) $. For $\epsilon$ small enough (i.e. $\epsilon < b / 4$), we have $y_0 \in \B_{b} (0)$. Let 
\begin{equation} \label{Existence of a Continuous Trajectory for the Perturbed Problem-01}
\widetilde{z} (t_0) := G_{t_0} ^{-1} (y_0) \cap \B_{a} (z(t_0) ).
\end{equation}
Note that the right-hand side of this expression is a singleton and so 
$ \widetilde{z} (t_0) $ is exactly determined without ambiguity. Let us check if $ (t_0, \widetilde{z} (t_0) ) \in \gph{\, \widetilde{S} } $ or, equivalently, $ 0 \in \widetilde{G_{t_0}} (\widetilde{z} (t_0))$.\\
From the definition of $ \widetilde{z} (t_0) $ we have $ y_0 \in G_{t_0} ( \widetilde{z} (t_0) ) = f(\widetilde{z} (t_0) ) + F(\widetilde{z} (t_0) ) - p(t_0)$. Then, from
 $y_0 =\widetilde{p}(t_0) - p(t_0) $, one gets  $\widetilde{p}(t_0) \in f(\widetilde{z} (t_0) ) + F(\widetilde{z} (t_0) ) $ or 
$ 0 \in \widetilde{G_{t_0}} (\widetilde{z} (t_0)) $.
Since $t_0$ is an arbitrary point in $[0, 1]$, 
$$ t \longmapsto G_{t} ^{-1} \big( \widetilde{p}(t) - p(t) \big) \cap \B_{a} (z(t) )$$
defines a single-valued map $\widetilde{z}(\cdot)$.\\
To prove the continuity, consider a sequence $(t_n) \in [0, 1]$ converging to $t_0$. By continuity of $ \widetilde{p}(\cdot)$ and $p(\cdot)$ we know that  
$ y_n := \widetilde{p}(t_n) - p(t_n) ~ \longrightarrow ~ y_0 = \widetilde{p}(t_0) - p(t_0) $.\\
By definition, $ \widetilde{z} (t_n) := G_{t_n} ^{-1} (y_n) \cap \B_a (z(t_n) )$. Lemma \ref{G-relations} yields that 
$$ \widetilde{z} (t_n) \in G_{t_0} ^{-1} ( \, y_n + p(t_n) - p(t_0) \, ).$$ 
On the other hand, $  \widetilde{z} (t_n) \in \B_a (z(t_n) ) $. We claim that $ \widetilde{z} (t_n) \in \B_a (z(t_0) ) $.\\
Indeed, for $n$ large enough, one can have the following:
\begin{equation} \label{construction-M1}
\norm{p(t_n) - p(t_0)} < \frac{b}{4},~ \norm{z(t_n) - z(t_0)} < \frac{a}{2}.
\end{equation}
Considering Lemma \ref{ratio law}, with $b' = b/2$, and $a' = a/2$, we obtain that the mapping
$$ \B_{b/2} (0) \ni y \longmapsto G_{t_0} ^{-1} (y) \cap \B_{a/2} (z(t_0) ) $$
is single-valued and Lipschitz continuous with Lipschitz constant $\kappa$. 
Now observing that 
\begin{equation*}
\begin{split}
\norm{ y_n + p(t_n) - p(t_0) } & \leq ~ \norm{ \widetilde{p}(t_n) - p(t_n) } + \norm{ p(t_n) - p(t_0) } \\
 & < \epsilon + \frac{b}{4} \, < \, \frac{b}{2},
\end{split}
\end{equation*}
we can define $ w_n := G_{t_0} ^{-1} \big( y_n + p(t_n) - p(t_0) \big) \cap \B_{a/2} (z(t_0) ) $ without ambiguity.\\
On the one hand, Lemma \ref{G-relations} implies that $w_n \in  G_{t_n} ^{-1} ( y_n)$.\\
On the other hand, $w_n \in \B_{a/2} (z(t_0) ) \subset \B_{a} (z(t_n) ) $. Thus, $w_n \in G_{t_n} ^{-1} ( y_n) \cap \B_{a} (z(t_n) )$.\\
Since $G_{t_n} ^{-1}$ is single-valued and Lipschitz continuous when restricted to $ \B_{b} ( 0 ) \times \B_{a} (z(t_n) ) $, we obtain the equality 
$ w_n = G_{t_n} ^{-1} ( y_n) \cap \B_{a} (z(t_n) ) $, and thus, $w_n = \widetilde{z} (t_n) $. 
Therefore, $ \widetilde{z} (t_n) \in \B_{a/2} (z(t_0) ) \subset \B_{a} (z(t_0) ). $\\
The strong metric regularity of $G_{t_0} $ implies that
\begin{equation*}
\norm{ \widetilde{z} (t_n) - \widetilde{z} (t_0)} \, \leq \, \kappa \norm{ \, \paren{y_n + p(t_n) - p(t_0)} - y_0 \, }.
\end{equation*}
Hence, $ \norm{ \widetilde{z} (t_n) - \widetilde{z} (t_0)} $ converges to zero as $ n \rightarrow \infty $. \\
It only remains to remind that the estimate for the difference $\norm { \widetilde{z} (t) - z(t) }$ is a straightforward consequence of the way we constructed $\widetilde{z}$. 
Indeed, let $ \displaystyle r( \epsilon ) := \sup_{ t \, \in \, [0, 1] } \norm{ \widetilde{z} (t) - z(t) } $. Starting from $ \epsilon < \frac{b}{4} $, we obtained $ r(\epsilon) < a $. If we let $ \epsilon < \frac{b}{8} $, a deeper look into the proof reveals that we can obtain $ r( \epsilon ) < \frac{a}{2} $, and so on. Thus, the distance $ \norm{ \widetilde{z} (t) - z(t) } $ (for every $t \, \in \, [0, 1]$) is controlled linearly by $ \epsilon $ and the proof is complete. 

\emph{Method 2. Construction over an interval:}\\
Fix $t \in [0, 1]$, and let $b$ smaller if necessary such that $ \kappa b < a $. By Remark \ref{radii effect}, this will not affect the uniform strong metric regularity of $ G_t $ guaranteed by the assumptions of this theorem and Theorem \ref{claim2}. The uniform continuity of $ \widetilde{p} (\cdot)$, and $ z(\cdot) $ allows us to choose $ \rho >0 $ sufficiently small and independent of $t$, such that for any $ \tau \in (t - \rho, \, t + \rho ) $, the following hold:
\begin{equation} \label{construction-01}
\norm{ \widetilde{p}(\tau) - \widetilde{p} (t) } < \dfrac{b}{4}, ~~ \norm{ z(\tau) - z(t) } < \dfrac{a}{2}.
\end{equation}
Then, for any $ \tau \in (t - \rho, \, t + \rho ) $ the continuity of $ \widetilde{p}(\cdot) $ and its closeness to $p(\cdot) $ implies that 
\begin{equation*}
\norm{\widetilde{p}(\tau) - p(t)} \, \leq \, \norm{\widetilde{p}(\tau) -\widetilde{p}(t) } + \norm{\widetilde{p}(t) - p(t)} < \dfrac{b}{2},
\end{equation*}
and therefore, by using Lemma \ref{ratio law} with $ b' = \frac{b}{2}$, and $ a' = \frac{a}{2}$, we obtain that the set $  G_t ^{-1} \big(\, \widetilde{p}(\tau) - p(t) \,\big)  \cap \B_{\frac{a}{2}} (z(t) ) $ is a singleton. Thus, we can define 
\begin{equation} \label{Existence of a Continuous Trajectory for the Perturbed Problem-02}
\widetilde{z} (\tau) := G_t ^{-1} \big(\, \widetilde{p}(\tau) - p(t) \,\big)  \cap \B_{\frac{a}{2}} ( z(t) ) \mathrm{~~~for~ any~~ } \tau \in (t - \rho, \, t + \rho ), 
\end{equation}
without ambiguity. In order to prove the continuity of this function, consider a sequence $(\tau_n)$ in $(t - \rho, \, t + \rho )$ such that $ \tau_n \longrightarrow \tau $. 
Then, from the Lipschitz continuity of $ G_t ^{-1} ( \cdot ) \cap \B_a ( z(t) ) $ over $ \B_b ( 0 ) $ we obtain that
\begin{equation*}
\begin{split}
\norm{ \widetilde{z}(\tau_n) - \widetilde{z}(\tau) } & =
 \norm{ \Big[ G_t ^{-1} \big(\, \widetilde{p}(\tau_n) - p(t) \,\big)  \cap \B_a ( z(t) )  \Big] - \Big[  G_t ^{-1} \big(\, \widetilde{p}(\tau) - p(t) \,\big)  \cap \B_a ( z(t) ) \Big] } \\
& \leq \, \kappa \norm{ ~ \widetilde{p}(\tau_n) - p(t) - \paren{\widetilde{p}(\tau) - p(t) }\, } \\
& \leq \, \kappa \norm{ ~ \widetilde{p}(\tau_n) - \widetilde{p}(\tau) \, }.
\end{split}	
\end{equation*}
The continuity of $ \widetilde{p} $ implies that $ \norm{ \widetilde{z}(\tau_n) - \widetilde{z}(\tau) }  \longrightarrow 0 $ as $ \tau_n \longrightarrow \tau $. \\
It remains to show that $ \widetilde{z} $ is (part of) a solution trajectory, that is, $ \big( \tau, \widetilde{z}(\tau) \big) \in \gph{\, \widetilde{S}} $.\\
Since $ \widetilde{z} (\tau) \in G_t ^{-1} \big(\, \widetilde{p}(\tau) - p(t) \, \big) $, from Lemma \ref{G-relations} we get 
$ \widetilde{z} (\tau) \in G_{\tau} ^{-1} \big(\, \widetilde{p}(\tau) - p(\tau) \, \big) $. Then, Equation \eqref{perturbation-relation2} implies that 
$ \widetilde{z} (\tau) \in \widetilde{G_{\tau}} ^{-1} (0)$  for any $ \tau \in (t - \rho, \, t + \rho ) $, or equivalently, $ \widetilde{z}(\tau) \in \widetilde{S}(\tau) $.\\
Up to now, we have proved that for each $ t \in [0, 1]$, we can find a solution trajectory in the interval $ (t - \rho, \, t + \rho ) $. It remains to show that this construction over different intervals remains consistent. 
To be more clear, let us consider two points $t_1$, and $t_2$, with corresponding trajectory pieces $\widetilde{z_1}$, and 
$\widetilde{z_2}$. Suppose that $t_1 < t_2 $ and let us consider the situation where $ t_2 - \rho < \tau < t_1 + \rho $. We should prove that 
$ \widetilde{z_1} (\tau) = \widetilde{z_2} (\tau) $.

By definition, $ \widetilde{z_i} (\tau) =  G_{t_i} ^{-1} \paren{\, \widetilde{p}(\tau) - p(t_i) \,} \cap \B_a (z(t_i) ) $ for $i = 1, 2 $, and as already shown, Lemma \ref{G-relations}, and Equality \eqref{perturbation-relation2} imply that $ \widetilde{z_i} (\tau) \in \widetilde{G_{\tau}} ^{-1} (0) $ for $i = 1, 2 $.\\ 
On the other hand, the continuity of $ z (\cdot) $, and inequalities in \eqref{construction-01} reveal that 
\begin{equation*}
\norm{ z(\tau) - \widetilde{z_i}(\tau) } \, \leq \,  \norm{ z(\tau) - z(t_i) } + \norm{ z(t_i) - \widetilde{z_i}(\tau) } \\
< \, \dfrac{a}{2} \, + \, \dfrac{a}{2}.
\end{equation*}
Thus, $ \widetilde{z_i}(\tau) \in \B_a (z(\tau) ) $ for $i = 1, 2 $. Using Lemma \ref{Perturbation Effect on the Auxiliary Map} and Remark \ref{Perturbation Effect on the Auxiliary Map-remark} for $ \widetilde{G_{\tau}}$, we obtain that the mapping
$ \widetilde{G_{\tau}} ^{-1} (\cdot) \cap \B_a (z(\tau) ) $ is single-valued and Lipschitz continuous over $\B_{\frac{3b}{4}} (0)$. So, $ \widetilde{z_1} (\tau) = \widetilde{z_2} (\tau) $, and the proof is complete.
\end{proof}

\rem \label{another proof}
\textbf{(a)} A through observation reveals that in fact, the two methods produce the same function mainly because of the single-valuedness of the mapping 
\begin{equation*}
\B_b (0) \ni y \longmapsto G_t ^{-1} (y) \cap \B_a (z(t) ).
\end{equation*}  
To be more precise, let us denote the trajectory obtained from \emph{Method 1.} by $\widetilde{z}_{M1}$, and the other one by $ \widetilde{z}_{M2} $. 
Consider an arbitrary point $ t \in [0, 1] $ and a neighborhood $ (t - \rho, \, t + \rho ) $ with $ \rho > 0 $ defined in such a way that \eqref{construction-01} holds.
First observe that from Equations \eqref{Existence of a Continuous Trajectory for the Perturbed Problem-01}, and \eqref{Existence of a Continuous Trajectory for the Perturbed Problem-02} we obtain immediately that $\widetilde{z}_{M1} (t) = \widetilde{z}_{M2} (t)$. \\
Now for any $ \tau \in (t - \rho, \, t + \rho )$, we have 
$ \widetilde{z}_{M1} (\tau) = G_{\tau} ^{-1} \big( \widetilde{p}(\tau) - p(\tau) \big) \cap \B_a (z(\tau) ) $.\\
We have already seen in proof \emph{Method 1.} of the previous theorem that when inequalities in \eqref{construction-M1} are satisfied (which is the case, by Condition \ref{construction-01} on 
$ \rho$), it is possible to conclude that $ \widetilde{z}_{M1} (\tau) \in B_a (z(t) )$. \\
On the other hand, Lemma \ref{G-relations} implies that $ \widetilde{z}_{M1} (\tau) \in G_t ^{-1} (\widetilde{p}(\tau) - p(t))$. Thus, 
$ \widetilde{z}_{M1} (\tau) \in G_t ^{-1} (\widetilde{p}(\tau) - p(t))\cap B_a (z(t)) $, and by strong metric regularity of $G_t$ we can obtain the desired equality 
$ \widetilde{z}_{M1} (\tau) = \widetilde{z}_{M2} (\tau) $.\\

\textbf{(b)} It is worth mentioning that the method of construction over intervals shows explicitly that Lipschitz continuity of $ \widetilde{z} (\cdot) $ could be easily obtained from Lipschitz continuity of $ \widetilde{p}(\cdot) $. But this is not something new or more than what we can obtain from the method of pointwise construction, as it was implicitly mentioned there, too. Indeed, in view of Lemma \ref{Perturbation Effect on the Auxiliary Map}, Proposition \ref{claim1}, and Corollary \ref{claim1-corollary}, we get the same result.